\documentclass[11pt]{amsart}

\usepackage{mathtools}        
\usepackage{mathrsfs}         
\usepackage{tikz}
\usetikzlibrary{graphs, positioning}
\usepackage{tikz-network}
\usepackage{booktabs}
\usepackage{multirow}
\usepackage{float}
\usepackage{bbm}
\usepackage{amsfonts}
\usepackage{amsmath}
\usepackage{amssymb}
\usepackage{graphicx}

\usepackage[left=1.1in,right=1.1in,top=1.25in,bottom=1.25in,footskip=.25in]{geometry}
\usepackage{mathtools}
\usepackage{amsthm}
\usepackage{latexsym}
\usepackage{fancyhdr}
\usepackage{array}
\usepackage{amscd}

\usepackage{lscape}
\usepackage{tikz}
\usepackage{tikz-cd}
\usepackage{float}
\usepackage{bm}
\usepackage{subfigure}
\usepackage{grffile}
\usepackage{array}
\usepackage{longtable}
\usepackage{booktabs}
\usepackage{adjustbox}
\usepackage{diagbox}
\usepackage{enumitem}

\usepackage{algorithm}
 \usepackage{algpseudocode}
 \usepackage{algorithmicx}
 \algdef{SE}[DOWHILE]{Do}{doWhile}{\algorithmicdo}[1]{\algorithmicwhile\ #1}

\algrenewcommand\algorithmicrequire{\textbf{Input:}}
\algrenewcommand\algorithmicensure{\textbf{Output:}}

\usepackage{ dsfont }
\usepackage{lipsum}
\newcolumntype{C}[1]{>{\centering\arraybackslash}p{#1}}

\definecolor{navy}{HTML}{2F729C} 
\usepackage[hyperfootnotes=false, colorlinks, linkcolor={blue}, citecolor={magenta}, filecolor={blue}, urlcolor={blue}, plainpages=false, pdfpagelabels]{hyperref}

\usepackage{cleveref}
\usepackage[tableposition=above]{caption}

\usepackage{xparse}
\usepackage{etoolbox}

\DeclareMathOperator{\dist}{dist}
\newcommand{\Caterpillar}[1]{
\begin{tikzpicture}[scale=0.8, transform shape]
    \pgfmathtruncatemacro{\k}{#1}
    \pgfmathtruncatemacro{\kmone}{\k-1}
    \def\step{1.5}
    \def\leafy{1.5}

    \foreach \i in {0,...,\k} {
        \node[shape=circle, draw, fill=white, minimum size=5mm, inner sep=1pt, line width=1.2pt]
            (v\i) at ({\i*\step},0) {};
    }

    \foreach \i in {1,...,\kmone} {
        \node[shape=circle, draw, fill=white, minimum size=5mm, inner sep=1pt, line width=1.2pt]
            (u\i) at ({\i*\step},\leafy) {};
    }

    \foreach \i [evaluate=\i as \j using int(\i+1)] in {0,...,\kmone} {
        \draw[line width=1.2pt] (v\i) -- (v\j);
    }

    \foreach \i in {1,...,\kmone} {
        \draw[line width=1.2pt] (v\i) -- (u\i);
    }

    \node at ({0.5*\k*\step},-1) {\Large $\mathcal H^1_{2^{\k}}$};
\end{tikzpicture}
}

\tikzset{
  Style:circle/.style={
    draw,
    circle,
    fill=white,
    minimum size=5mm,
    inner sep=1pt,
    line width=1.2pt
  },
  Style:edge/.style={
    line width=1.2pt
  }
}

\tikzset{
  StyleEllipse/.style={
    draw,
    ellipse,
    fill=white,
    inner sep=2pt,
    outer sep=0pt,
    line width=1.2pt
  },
  StyleEdge/.style={
    line width=1.2pt
  }
}

\newcommand{\Q}{{\mathbb Q}}
\newcommand{\Z}{{\mathbb Z}}

\newcommand{\GL}{{\rm GL}}

\newcommand\sm[1]{\begin{psmallmatrix}#1\end{psmallmatrix}}

\newcommand\m[1]{\begin{pmatrix}#1\end{pmatrix}}

\newcommand{\lmfdb}[2]{\href{#1}{\texttt{#2}}}

\newtheorem{mainthm}{Theorem}
\newtheorem{maincor}[mainthm]{Corollary}

\theoremstyle{definition}
\newtheorem{defn}{Definition}[section]
\newtheorem{definition}[defn]{Definition}

\newtheorem{notation}[defn]{Notation}
\newtheorem{example}[defn]{Example}

\theoremstyle{plain}

\newtheorem{corollary}[defn]{Corollary}
\newtheorem{lemma}[defn]{Lemma}

\newtheorem{theorem}[defn]{Theorem}

\newtheorem{proposition}[defn]{Proposition}

\theoremstyle{remark}

\newtheorem{remark}[defn]{Remark}

\numberwithin{equation}{section}

\title[Isogeny graphs of elliptic curves in characteristic zero]{Isogeny graphs of elliptic curves \\ in characteristic zero}

\author{Alexander J. Barrios}
\address{University of St. Thomas, Department of Mathematics, St. Paul, MN 55105 USA}
\email{abarrios@stthomas.edu}

\author{Enrique Gonz\'{a}lez-Jim\'{e}nez}
\address{Universidad Aut\'{o}noma de Madrid, Departamento de Matem\'{a}ticas, Madrid, Spain}
\email{enrique.gonzalez.jimenez@uam.es}

\author{Ivan Novak}
\address{University of Zagreb, Bijeni\v{c}ka Cesta 30, 10000 Zagreb, Croatia}
\email{ivan.novak@math.hr}



\begin{document}

\begin{abstract}
For an elliptic curve $E$ defined over a field $K$ of characteristic $0$ with $\operatorname{End}_K \! E \cong \mathbb{Z}$, we classify which isogeny graphs $\mathcal{G}(E/K)$ can occur. We first show that $\mathcal{G}(E/K)$ decomposes as a weak Cartesian product of its $p$-primary isogeny graphs, one for each prime $p$, thereby reducing the problem to classifying $p$-primary isogeny graphs. We then show that each such graph is isomorphic, as an edge-weighted graph, to a member of an explicit family of edge-weighted graphs $\mathcal{H}_{p^k}^r$ and $\mathcal{H}_{p^{\infty,+}}^r$, every member of which occurs as a $p$-primary isogeny graph except for $\mathcal{H}_{2^k}^0$ for $k\ge 2$. The proof relies on a detailed study of the $p$-adic Galois representation attached to $E$, through which we identify each graph with a subgroup of $\operatorname*{GL}\nolimits_{2}(\mathbb{Z}_{p})$. 

More generally, we identify subgroups of $\operatorname*{GL}\nolimits_{2}(\widehat{\mathbb{Z}})$ for each possible isogeny graph and describe their corresponding modular curves, completing, in the genus $0$ case, the explicit parameterization of isogeny graphs via parameterized isogenous families of elliptic curves. We also introduce the $p$-blooming invariant $\mathfrak{I}_p(E/K)$, an isogeny class invariant determining the value of $r$ in the $p$-primary isogeny graph, and show that elliptic curves over fields with a real embedding attain the smallest possible value. As applications, we characterize the isogeny graphs of elliptic curves with potential complex multiplication; give an algorithm for determining the isogeny graph from the adelic Galois representation; recover the classification of rational isogeny graphs; and, under GRH, classify the isogeny graphs occurring over certain number fields.
\end{abstract}

\maketitle

\tableofcontents

\section{Introduction}
For an elliptic curve $E$ over a field $K$, its isogeny graph $\mathcal{G}(E/K)$ is the undirected edge-weighted graph whose vertices are the $K$-isomorphism classes $[E^{\prime}]_{K}$ of elliptic curves $E^{\prime}$ that are isogenous to $E$, and two vertices are joined by an edge of weight $p$ whenever there is an isogeny of prime degree $p$ between representatives of the corresponding classes. When the characteristic of $K$ is finite, these graphs have been studied extensively. Ordinary elliptic curves give rise to volcanoes \cite{Clarkvolcanoes,Clarkvolcanoes2,MireilleMorain,Kohel,DrewVolcanoes}, while supersingular elliptic curves give rise to highly structured graphs governed by the quaternionic nature of their endomorphism rings \cite{Voight,Waterhouse}. These graphs have long been of interest for their intrinsic arithmetic properties and have also found applications in cryptography \cite{GalCrypto,GalVer,KoblitzCrypto}, with supersingular isogeny graphs recently attracting particular attention in the search for quantum-safe cryptosystems \cite{SS2, DelfsGal,SS1}. In contrast, the characteristic zero case has not been studied to the same depth.

In this article, we investigate the structure of $\mathcal{G}(E/K)$ over fields $K$ of characteristic zero. We fix a prime $p$ and consider the $p$-primary (isogeny) graph $\mathcal{G}_{p}(E/K)$, defined as the full~subgraph of $\mathcal{G}(E/K)$ consisting of those vertices $[E^{\prime}]_{K}$ such that $E$ and $E^{\prime}$ are $p^k$-isogenous over $K$ for some positive integer $k$. Over an algebraically closed field, the structure of $\mathcal{G}_{p}(E/K)$ is especially transparent and is naturally governed by Bruhat-Tits theory \cite{BruhatTits,defeo,Trees}; in particular, the algebraically closed case provides the local combinatorial model that motivates the global study.

The situation becomes substantially subtler once $K$ is not algebraically closed. In the simplest case, isogeny graphs over $\mathbb{Q}$ seem to have first appeared in the so-called Antwerp tables \cite{MR0376533}, where the authors remark that the tables ``illustrate almost all the known ways in which isogenies can occur.'' For decades, it was folklore among experts that these configurations indeed exhausted all possibilities, although no general proof was available at the time.

The classification of rational isogeny graphs depends critically on the possible degrees of rational cyclic isogenies of elliptic curves over $\mathbb{Q}$. Mazur’s celebrated 1978 theorem \cite{MR482230} determined the possible prime degrees. Building on ideas and computations due to numerous mathematicians, including Fricke \cite{MR3221641}, Fricke and Klein \cite{MR3838339}, Kenku \cite{MR510395,MR588271,MR549292,MR616546,MR675184}, Klein \cite{MR1509988}, Kubert~\cite{MR0434947}, Ligozat \cite{MR0422158}, Mazur \cite{MR482230}, Mazur and V\'{e}lu \cite{MR320010}, and Ogg \cite{MR0337974}, among others, this line of work culminated in 1982 with Kenku’s result \cite{MR675184}, which completed the classification of all rational cyclic isogeny degrees. Finally, in 2022, Chiloyan and Lozano‑Robledo \cite[\S 6]{MR4203041} provided the first complete proof of the classification of rational isogeny graphs $\mathcal{G}(E/\mathbb{Q})$, as part of their broader study of rational isogeny–torsion graphs. Their work confirmed the long‑standing folklore expectation. Subsequently, Barrios \cite{Bariso} gave an explicit classification of rational isogeny graphs via parameterized isogenous families of elliptic curves. Motivated by these developments, we ask:
\begin{center}
\textit{Given a field $K$ of characteristic $0$, which isogeny graphs of elliptic curves occur over $K$, and can they be explicitly classified?}
\end{center}

In full generality, this problem appears to lie far beyond the reach of current techniques. Even isolating a single basic invariant of these graphs, namely the possible prime degrees of $K$-rational isogenies, leads to difficult and largely unresolved questions. More precisely, for a number field $K$, let $\operatorname{IsogPrimeDeg}(K)$ denote the set of primes $p$ such that some elliptic curve over $K$ admits a $K$-rational $p$-isogeny. When $K=\mathbb{Q}$, Mazur's theorem~\cite{MR482230} yields \[
\operatorname{IsogPrimeDeg}(\mathbb{Q})=\left\{  2,3,5,7,11,13,17,19,37,43,67,163\right\}.
\]
For $K\neq \mathbb{Q}$, the best general results are largely conditional on the generalized Riemann hypothesis (GRH). Larson and Vaintrob \cite{LarsonVaintrob} proved, assuming GRH, that $\operatorname{IsogPrimeDeg}(K)$ is finite if and only if $K$ does not contain the Hilbert class field of an imaginary quadratic field. The forward implication is unconditional, following from the theory of complex multiplication (CM), whereas the converse is deduced from results of Merel~\cite{Merel} and Momose~\cite{Momose} under GRH.

Even under GRH, making these results explicit is highly nontrivial. Combining the work of Larson and Vaintrob with that of David \cite{David} and Momose \cite{Momose}, yields an algorithm \cite[Algorithm~1.6]{BanwaitAlg} that produces explicit upper bounds for the largest prime in $\operatorname{IsogPrimeDeg}(K)$ for certain quadratic number fields $K$. However, these bounds can be astronomical: for instance, the bound obtained for $\mathbb{Q}(\sqrt{-10})$ is $3.20\times10^{316}$ \cite[Table~1]{BanwaitAlg}. Banwait \cite[Theorem~1.8]{BanwaitAlg} subsequently gave a refined GRH-dependent algorithm that dramatically improves these bounds; in the same example, the largest prime is shown to be at most $73$. With this refinement, he proves that $\operatorname{IsogPrimeDeg}(\mathbb{Q}(\sqrt{-10}))=\operatorname{IsogPrimeDeg}(\mathbb{Q})$ if GRH holds. More recently, Banwait, Najman, and Padurariu \cite{Banwait} determined the degrees of cyclic isogenies over several fixed quadratic fields, and Banwait and Derickx \cite{BanwaitDerickx} considered the problem for arbitrary number fields, giving a GRH-dependent algorithm for computing an explicit upper bound for the largest prime in $\operatorname{IsogPrimeDeg}(K)$ under suitable hypotheses.

These results illustrate the depth of the problem and the current limitations of our techniques, and suggest that a complete classification of isogeny graphs over number fields is presently out of reach. Accordingly, the main objective of this paper is to answer the following variant:
\begin{center}
\textit{Which isogeny graphs are realizable for elliptic curves $E$ over fields $K$ of characteristic $0$?}
\end{center}
The case of CM in characteristic zero has been studied in the literature \cite{Clarkvolcanoes,Clarkvolcanoes2,DrewVolcanoes}. We therefore focus on elliptic curves $E$ that satisfying $\operatorname{End}_K E\cong \mathbb{Z}$. The core of the article fully answers this question, obtaining a classification of the isogeny graphs that occur, expressed in terms of the image of the associated Galois representation. In particular, from our work, we recover the classification of rational isogeny graphs, as well as a classification, conditional on GRH, of the possible isogeny graphs over certain number fields. This is the context of Section~\ref{sec:nffields}, where we revisit the above discussion about number fields for which all possible degrees of cyclic isogenies have been determined under GRH, and place them within the context of our results.

The first step towards answering our motivating question is to reduce the study of isogeny graphs $\mathcal{G}(E/K)$ to an understanding of the $p$-primary graphs $\mathcal{G}_{p}(E/K)$ for each prime $p$. This is the context of our first result, which plays for isogeny graphs a role analogous to unique factorization in the Fundamental Theorem of Arithmetic.

\begin{mainthm}\label{mainthmCarPro}
Let $E$ be an elliptic curve over a field $K$ of characteristic zero with $\operatorname*{End}_{K}\! E\cong\mathbb{Z}$. For each prime $p$, let $(\mathcal{G}_{p}(E/K),[E]_K)$ denote the $p$-primary graph as a pointed graph. Then $\mathcal{G}(E/K)$ is graph-isomorphic to the following weak Cartesian product, taken over all primes~$p$:
\[
\underset{p}{\square}(\mathcal{G}_{p}(E/K),[E]_K)
\]
\end{mainthm}
Theorem~\ref{mainthmCarPro} appears as Theorem~\ref{Thm:CartesianPro}, and its proof relies on a study of the full torsion subgroup. Indeed, the study of isogenies of an elliptic curve $E$ is intertwined with the lattice of cyclic $G_K$-submodules of the full torsion subgroup $E(\overline{K})_{\text{tors}}$, where $G_K=\operatorname{Gal}(\overline{K}/K)$ is the absolute Galois group of $K$. Motivated by this connection, we begin in Section~\ref{Sec:StructureGraphs} by analyzing the $G_K$-invariant cyclic sublattice of $E(\overline{K})_{\text{tors}}$. In the finite setting, Suzuki~\cite{Suzuki} proved that the subgroup lattice of a finite group decomposes as a Cartesian product under the appropriate primary decomposition hypotheses, and the same viewpoint applies to the cyclic subgroup lattice of a finite group, where the subgraph lattice is the familiar divisibility lattice~\cite{Tarnauceanu}. In our setting, we are interested in the corresponding abelian torsion case, together with the additional $G_{K}$-module structure. Accordingly, we establish in Proposition~\ref{Lem:IsoGraPro} that the lattice of cyclic $G_K$-invariant submodules of an abelian torsion group decomposes as a weak Cartesian product of its primary components. 

This fits naturally into the classical graph-theoretic literature on Cartesian products of graphs: Sabidussi~\cite{Sabidussi} and Vizing~\cite{Vizing} independently proved unique prime factorization for connected finite graphs with respect to the Cartesian product. We note that a graph is said to be prime if it cannot be represented as the Cartesian product of two non-trivial graphs. In the infinite setting, Miller~\cite{Miller} and Imrich~\cite{Imrich1} independently proved an analogous unique prime factorization via the weak Cartesian product. These ideas were further extended by Imrich and Peterin~\cite{Imrich} to graphs with loops. Our contribution is to place the cyclic-submodule lattice in this framework and show that, for an abelian torsion group with primary decomposition, it admits a compatible decomposition in terms of the weak Cartesian product. In particular, while this article concerns itself with elliptic curves, our results can be adapted to study the possible subgraphs of cyclic isogenies in the isogeny graph of abelian varieties with no additional endomorphisms.

With Theorem~\ref{mainthmCarPro} on hand, a classification of the possible isogeny graphs reduces to classifying the $p$-primary graphs. To state this classification, we first recall that the isogeny class degree of $E/K$, which we denote by $\deg\mathcal{G}(E/K)$, is the least common multiple of the entries in the isogeny matrix $M_{E/K}$ of  $E$, where each entry is the minimal degree of a cyclic isogeny between two elliptic curves in the isogeny class. When $\operatorname*{End}_{K}\! E\cong\mathbb{Z}$ and $\mathcal{G}(E/K)$ is finite, we prove that the isogeny matrix $M_{E/K}$ is permutation-similar to the Kronecker product of the isogeny matrices associated to the $p$-primary graphs (see Proposition~\ref{prop:maxmat}). Consequently, we obtain that $\deg\mathcal{G}(E/K)$ is equal to the largest entry in the isogeny matrix.  We note that the hypothesis $\operatorname*{End}_{K}\! E\cong\mathbb{Z}$ is essential here, and also in Theorem~\ref{mainthmCarPro}, as demonstrated by Examples~\ref{Ex:CMnotCarPro} and~\ref{ex:CMvolc}.

For each prime $p$, let $\deg\mathcal{G}_{p}(E/K)$ denote the largest $p$-power degree of a cyclic isogeny occurring in $\mathcal{G}_p(E/K)$. Thus, $\deg\mathcal{G}_{p}(E/K)=p^{k_{p}}$ for some $k_{p}\in\mathbb{Z}_{\geq0}\cup\{\infty\}$. If $\mathcal{G}(E/K)$ is finite, then each $k_{p}\in\mathbb{Z}_{\geq0}$ and only finitely many $k_{p}$ are nonzero, and from Proposition~\ref{prop:maxmat}, we obtain:
\begin{equation}
\deg\mathcal{G}(E/K)=\prod_{p}\deg\mathcal{G}_{p}(E/K).\label{eq:degisoclass}
\end{equation}
Motivated by this, we extend the definition of isogeny class degree to the infinite setting by defining $\deg\mathcal{G}(E/K)$ as the supernatural number corresponding to the product in~\eqref{eq:degisoclass}.

To describe the graphs that occur in our classification, we introduce in Section~\ref{sec:graphs} a family of finite graphs $\mathcal{H}_{p^{k}}^{r}$ for $k,r\in\mathbb{Z}_{\geq0}$ with $r\leq\frac{k}{2}$, and in Section~\ref{sec:infgraphs}, a family of infinite graphs $\mathcal{H}_{p^{\infty}}^{r}$ and $\mathcal{H}_{p^{\infty,+}}^{r}$ with $r\in \mathbb{Z}_{\geq0}\cup\{\infty\}$. We note that $\mathcal{H}_{p^{0}}^{0}$ is the trivial graph with one vertex and no edges. The graph $\mathcal{H}_{p^{\infty}}^{\infty}$ is the $p$-Bruhat--Tits graph, and the graphs $\mathcal{H}_{p^{k}}^{0}$ are path graphs with $k+1$ vertices. Our classification shows that the $p$-primary graph $\mathcal{G}_{p}(E/K)$ is isomorphic to a member of one of these families. The quantity $r$ encodes information about the $p$-adic Galois image. More precisely, in Section~\ref{sec:bloominv}, we introduce and study the \textit{$p$-blooming invariant} $\mathfrak{I}_p(E/K)$ (see Definition~\ref{def:p-bloominginv}) of an elliptic curve $E$, and show that it is an isogeny class invariant. 

Before discussing our main result, we note that the only member of the family that is not realizable as the $p$-isogeny graph of some elliptic curve is the subfamily of path graphs $\mathcal{H}_{2^{k}}^{0}$ with $k\geq2$. This is because any elliptic curve that admits two distinct $2$-isogenies automatically admits three distinct $2$-isogenies. Consequently, for $p=2$ with $k\geq2$, the realizable graphs begin with $\mathcal{H}_{2^{k}}^{1}$. Over $\mathbb{Q}$, these graphs occur with $2\leq k\leq4$, and are given below.
\[
\adjustbox{scale=0.9}{$\begin{array}
[c]{ccc}
\Caterpillar{2}
& \qquad
\Caterpillar{3}
&\qquad
\Caterpillar{4}
\end{array}$}
\]
More generally, for $2\le k<\infty$, $\mathcal{H}_{2^{k}}^{1}$ is the caterpillar graph whose spine is a path on $k+1$ vertices, with each internal vertex having degree $3$. For instance, $\mathcal{H}_{2^{10}}^{1}$ is shown below:
\[
\adjustbox{scale=0.9}{\Caterpillar{10}}
\]

For $p=2$,  Arg\'{a}ez-Garc\'{\i}a and Cremona \cite{BlackBox} used the $2$-adic Galois representation to determine when a nontrivial $2$-primary graph is small or large. Using our terminology, a nontrivial $2$-primary is small precisely when $\mathcal{G}_{p}(E/K)\cong\mathcal{H}_{2}^{0}$, and is large otherwise. Sanna \cite{Sanna} obtained analogous results for $3$-primary graphs. 

With terminology established, we now state our classification of $p$-primary isogeny graphs:

\begin{mainthm}\label{mathmclass}
Let $E$ be an elliptic curve over a field $K$ of characteristic $0$ with $\operatorname*{End}_{K}\! E\cong\mathbb{Z}$.
\begin{enumerate}
    \item If  $\deg\mathcal{G}_{p}(E/K)=p^k$ with $k<\infty$, then $\mathcal{G}_p(E/K)$ is graph-isomorphic to $\mathcal{H}_{p^{k}}^{r}$ with $r = \min\{\mathfrak{I}_p(E/K),\lfloor k/2 \rfloor\}$.
    \item If  $\deg\mathcal{G}_{p}(E/K)=p^\infty$, then $\mathcal{G}_p(E/K)$ is graph-isomorphic to $\mathcal{H}_{p^{\infty}}^{\mathfrak{I}_p(E/K)}$ or $\mathcal{H}_{p^{\infty,+}}^{\mathfrak{I}_p(E/K)}$.
\end{enumerate}

Conversely, with the exception of $\mathcal{H}_{2^{k}}^{0}$ with $k\ge 2$, if $\mathcal{G}$ is one of the graphs $\mathcal{H}_{p^{k}}^{r}$ or $\mathcal{H}_{p^{\infty,+}}^{r}$ with $(k,r)$ satisfying one of the conclusions above, then there exists a field $K$ of characteristic $0$ and an elliptic curve $E$ defined over $K$ such that $\mathcal{G}$ is graph-isomorphic to $\mathcal{G}_p(E/K)$.

Further, for a prime $p\geq5$, if $r\geq1$, then $\sqrt{p^{\ast}}\in K$, where $p^*=(  -1)  ^{(p-1)/2}p$. In addition,
\begin{enumerate}
\item[(a)] if $p=2$, then

\begin{enumerate}
\item[(i)] $r\leq1$ if and only if either $i\not \in K$ or $k\leq3$;

\item[(ii)] $r=2$ if and only if $i\in K$ with $k\geq4$ and either $\zeta_{8}\not \in K$ or $k\in\{4,5\}$;

\item[(iii)] $r\geq3$ if and only if $\zeta_{8}\in K$ and $k\geq6$; 
\end{enumerate}

\item [(b)]if $p=3$, then

\begin{enumerate}
\item[(i)] $r=0$ if and only if either $\zeta_{3}\not \in K$ or $k=1$;

\item[(ii)] $r\geq1$ if and only if $\zeta_{3}\in K$ and $k\geq2$; 
\end{enumerate}
\end{enumerate}
\end{mainthm}
The finite and infinite cases of the theorem are treated separately in Sections \ref{sec:graphs} and \ref{sec:infgraphs}, respectively. In particular, parts (1) and (2) of Theorem \ref{mathmclass}, and their respective converse, are Theorems~\ref{classificationGpk} and \ref{thm:infinitegpks}, respectively. In the finite case, the determination of $r$ from the $p$-blooming invariant is the statement of Corollary~\ref{ex:CMvolc}. The further portion of the theorem is the focus of Corollary \ref{cor:mainthmfieldofdef}. The proof of Theorem \ref{mathmclass} relies on a combinatorial study of the families of the graphs $\mathcal{H}_{p^{k}}^{r}$ and $\mathcal{H}_{p^{\infty,+}}^{r}$, together with a study of the $p$-adic Galois representation attached to an elliptic curve. More precisely, in Sections \ref{subsec:modpGalois} and \ref{subsec:padicGalois} we investigate, respectively, the $\operatorname{mod}p^{k}$ and $p$-adic Galois representations attached to an elliptic curve. Recall that after a choice of basis, the $\operatorname{mod}p^{k}$ and $p$-adic Galois representation attached to $E$ are
\[
\rho_{E,p^{k}}:G_{K}\rightarrow\operatorname*{GL}\nolimits_{2}(\mathbb{Z}/p^{k}\mathbb{Z})\qquad\text{and}\qquad\rho_{E,p^{\infty}}:G_{K}\rightarrow\operatorname*{GL}\nolimits_{2}(\mathbb{Z}_{p}),
\]
respectively. Key ingredients in the proof of Theorem~\ref{mathmclass} are Propositions~\ref{prop:imagechangeuponisogeny} and~\ref{prop:infiniteimagechange}, which allow us to investigate the Galois images of elliptic curves isogenous to $E$ from a fixed basis for $E[p^{k}]$ or $T_{p}(E)$, where $T_{p}(E)$ denotes the Tate module of $E$. Moreover, the proof of Theorem~\ref{mathmclass} also relies on exploiting combinatorial information first explored by Novak~\cite{ivan-numberofiso}. In particular, we exploit the fact that for each positive integer $i$, the $p$-adic Galois representation uniquely determines the number of $K$-rational $p^{i}$-isogenies admitted by $E$. In fact, the proof of Theorem~\ref{mathmclass} proceeds through a careful study of the Galois images to establish that $\mathcal{G}_{p}(E/K)$ must contain a member of the family of graphs as a subgraph. We then adhere to the established combinatorial properties of the graphs to deduce our desired isomorphism. As a consequence, we obtain the exact number of elliptic curves in a given isogeny class. In particular, we obtain the following result, which is an immediate consequence of Theorems~\ref{mainthmCarPro} and~\ref{mathmclass}, together with Corollary~\ref{number_of_vertices}.
\begin{maincor}\label{Cor:isoclassize}
Let $E$ be an elliptic curve over a field $K$ of characteristic $0$ with $\operatorname*{End}_{K}E\cong\mathbb{Z}$. If the isogeny graph of $E/K$ is finite and satisfies
\[
\mathcal{G}(E/K)\cong \underset{p}{\square}\mathcal{H}_{p^{k_{p}}}^{r_{p}},
\]
then the size of the isogeny class of $E/K$ is
\[
\left\vert V(\mathcal{G}(E/K))\right\vert = \prod_{p}\left(  \left(k_{p}+1-2r_{p}\right)  p^{r_{p}}+2\cdot\frac{p^{r_{p}}-1}{p-1}\right)  .
\]
\end{maincor}

In Section \ref{sec:bloominv}, we introduce the $p$-blooming invariant $\mathfrak{I}_{p}(E/K)$ of an elliptic curve and study its properties. We note that $\mathfrak{I}_{p}(E/K)$ is the minimum of the $p$-adic valuations of the difference of the eigenvalues of $\rho_{E,p^{\infty}}(\sigma)$, as $\sigma$ ranges over $G_{K}$ (see Definition \ref{def:p-bloominginv}). We conclude the section by demonstrating that elliptic curves defined over fields admitting a real embedding have the smallest possible $\mathfrak{I}_{p}(E/K)$ (see Proposition \ref{Prop:BloomInvReal}). Consequently, the isogeny class degree uniquely determines the isogeny graph over such fields (see Corollary~\ref{Cor:BloomInvReal}).

This observation has particularly strong consequences for elliptic curves with potential CM. In this direction, let $E$ be a complex elliptic curve with CM by an imaginary quadratic field $K$, and let $F=\mathbb{Q}(j(E))$. Since $F$ admits a real embedding \cite[Remark 5.2]{SilverbergFoD}, Proposition~\ref{Prop:BloomInvReal} applies. In Section~\ref{sec:pCM}, we exploit this fact to analyze the isogeny graphs of potential CM elliptic curves. We begin by proving that, for each $n\geq3$, the CM field $K$ is contained in the splitting field of the level~$n$ modular polynomial~$\Phi_{n}(X,j(E))$, generalizing a result of Bourdon, Clark, and Stankewicz \cite[Lemma 3.15]{bourdonclarkstan}, which states that $K\subseteq F(E[n])$. We then use Kwon's characterization of $n$-isogenies~\cite{Kwon} to classify the $p$-primary isogeny graphs of potential CM elliptic curves over their field of definition (see Lemma~\ref{lem:potCM}). Building on this result, we determine the distribution of the $p$-adic valuations of the conductors of the endomorphism rings along maximal paths in these graphs (see Theorem~\ref{thm:CMdiscclass}). These results culminate in the following complete classification of isogeny graphs associated to potential CM elliptic curves:

\begin{mainthm}\label{mainthm:genkwon}
Let $E$ be a complex elliptic curve with $\operatorname*{End}E\otimes_{\mathbb{Z}}\mathbb{Q}\cong K$ an imaginary quadratic field. Suppose further that $E$ is defined over a field $F$ such that $\operatorname*{End}_{F}E\cong\mathbb{Z}$. For each $E^{\prime}\in V(\mathcal{G}(E/F))$, let $\mathfrak{f}_{E^{\prime}}$ denote the conductor of the order $\operatorname*{End}E^{\prime}$. For each prime $p$, let
\[
\mathfrak{F}_{p}=\sup\left\{  v_{p}(\mathfrak{f}_{E^{\prime}})\mid E^{\prime
}\in V(\mathcal{G}_{p}(E/F))\right\}  .
\]
Then, the isogeny graph $\mathcal{G}(E/F)\cong\square_{p}(\mathcal{G}_{p}(E/F),[E]_{F})$, where $\mathcal{G}_{p}(E/F)$ is determined as follows:

\begin{enumerate}
\item if $\mathfrak{F}_{p}=\infty$, then
\[
\mathcal{G}_{p}(E/F)\cong\left\{
\begin{array}
[c]{cl}
\mathcal{H}_{2^{\infty}}^{1} & \text{if }p=2,\\
\mathcal{H}_{p^{\infty}}^{0} & \text{if }p\geq3;
\end{array}
\right.
\]

\item if $\mathfrak{F}_{p}<\infty$, then there exists an elliptic curve $E_{0}\in V(\mathcal{G}_{p}(E/F))$ such that $\operatorname*{End}E_{0}$ has discriminant $\Delta=\mathfrak{f}_{E_{0}}^{2}\Delta_{K}$, where $\Delta_{K}$ is the discriminant of $K$ and $v_{p}(\mathfrak{f}_{E_{0}})=\mathfrak{F}_{p}$. Further, let $\mathcal{G}_{p}(E_{0}/\mathbb{Q}(j(E_{0}))$  be as given by Lemma~\ref{lem:specialdiscpotcm} or Theorem~\ref{thm:CMdiscclass}. Then,
\[
\mathcal{G}_{p}(E/F)\cong\mathcal{G}_{p}(E_{0}/\mathbb{Q}(j(E_{0}))).
\]
\end{enumerate}
\end{mainthm}
This is Theorem~\ref{thm:genkwon} in the text. As an application, we obtain necessary and sufficient conditions for a potential CM elliptic curve to admit a cyclic $n$-isogeny (see Corollary \ref{genKwonBC}), thereby extending theorems of Kwon~\cite{Kwon}, as well as Bourdon and Clark~\cite{BourdonClark1}. While our work focuses on elliptic curves satisfying $\operatorname*{End}_{F}\!E\cong\mathbb{Z}$, the case of CM has been extensively studied in the literature, including recent work of Clark \cite{Clarkvolcanoes}, as well as Clark and Saia~\cite{Clarkvolcanoes2}, and, for the characterization of $n$-isogenies for CM elliptic curves, by Bourdon and Clark \cite{BourdonClark2}. Aside from Examples~\ref{Ex:CMnotCarPro} and~\ref{ex:CMvolc}, we do not investigate the CM case any further.

Beyond classifying isogeny graphs, the proof of Theorem~\ref{mathmclass} naturally identifies a family of subgroups of $\operatorname*{GL}\nolimits_{2}(\mathbb{Z}/p^{k}\mathbb{Z})$ and $\operatorname*{GL}\nolimits_{2}(\mathbb{Z}_p)$ that govern the Galois representation attached to elliptic curves realizing each graph. These groups provide a bridge between the combinatorial classification of isogeny graphs and modular curves, forming the basis for the results discussed in Section~\ref{sec:ModularCurves}. More precisely, for each finite graph $\mathcal{H}_{p^{k}}^{r}$, we identify a subgroup $H_{p^{k}}^{r}\le \operatorname*{GL}\nolimits_{2}(\mathbb{Z}/p^{k}\mathbb{Z})$ with the property that whenever $\mathcal{G}_{p}(E/K)\cong\mathcal{H}_{p^{k}}^{r}$, the elliptic curve $E$ is isogenous to an elliptic curve $E'$ whose~$\operatorname{mod} p^k$ Galois image is conjugate to a subgroup of $H_{p^{k}}^{r}$. Moreover, if $(E_0,E_1,\ldots E_k)$ is a maximal path in $\mathcal{G}_{p}(E/K)$, we defined subgroups $H_{p^{k}}^{r}(j)$ describing the Galois image attached to each vertex along the path (see Corollary~\ref{cor:subgroup-subgraph-corr}). An analogous construction in the infinite case yields subgroups $H_{p^{\infty}}^{r}$ and $H_{p^{\infty,+}}^{r}$ of $\operatorname*{GL}\nolimits_{2}(\mathbb{Z}_p)$, which characterize the $p$-adic Galois image of elliptic curves realizing the infinite graphs (see Corollary~\ref{cor:infiniteimageconjugate}).

In sum, we show that each $p$-primary graph $\mathcal{G}_{p}(E/K)$ corresponds to a subgroup of $\operatorname*{GL}\nolimits_{2}(\mathbb{Z}_{p})$, where we identify those subgroups of $\operatorname*{GL}\nolimits_{2}(\mathbb{Z}/p^{k}\mathbb{Z})$ with their corresponding lift. Next, after choosing a basis, the adelic Galois representation attached to $E$ is
\[
\rho_{E}:G_{K}\longrightarrow\operatorname*{GL}\nolimits_{2}(\widehat{\mathbb{Z}}).
\]
The above discussion, together with Theorem~\ref{mainthmCarPro}, then shows that each isogeny graph $\mathcal{G}(E/K)$ can be associated to a subgroup of $\operatorname*{GL}\nolimits_{2}(\widehat{\mathbb{Z}})$. In particular, if
\[
\mathcal{G}(E/K)\cong\underset{p}{\square}\mathcal{H}_{p^{k_{p}}}^{r_{p}},
\]
then the image of the adelic Galois representation of some elliptic curve in the isogeny class of~$E$ is conjugate to a subgroup of $\prod_{p}H_{p^{k_{p}}}^{r_{p}}$. In Section \ref{sec:algorithm image->graph}, we discuss how to determine $\mathcal{G}(E/K)$ from the adelic Galois image. More precisely, we prove Theorem~\ref{thm:algorithm}, whose statement is given below:
\begin{mainthm}\label{mainthm:algorithm}
Let $E$ be an elliptic curve over a field $K$ of characteristic $0$ with $\operatorname*{End}_{K}\! E\cong~\mathbb{Z}$. There is an algorithm which outputs the pointed graph $\left(  \mathcal{G}(E/K),[E]_K\right)  $ from the following inputs:
\begin{itemize}
\item a positive integer $n$, the level of $\rho_{E}$;

\item a subset $S\subseteq\operatorname*{GL}\nolimits_{2}(\mathbb{Z}/n\mathbb{Z})$ such that $\rho_{E,n}(G_{K})\cong \langle S\rangle$ with respect to some basis of~$E[n]$.
\end{itemize}
\end{mainthm}

Next, we observe that the discussion preceding Theorem \ref{mainthm:algorithm} highlighted that each possible isogeny graph $\mathcal{G}(E/K)$ with $\operatorname*{End}_{K}E\cong\mathbb{Z}$ is associated to some subgroup of $\operatorname*{GL}\nolimits_{2}(\widehat{\mathbb{Z}})$. In Section~\ref{sec:ModularCurves}, we restrict to the open subgroups $H\leq\operatorname*{GL}\nolimits_{2}(\widehat{\mathbb{Z}})$ associated to isogeny graphs, and consider their associated modular curve $X_{H}$. In particular, we determine their field of definition (see Lemma~\ref{lem:fieldofdefmodcurve}), from which we deduce Corollary~\ref{cor:mainthmfieldofdef}, which establishes the further portion of Theorem~\ref{mathmclass}, thus completing its proof. We also classify the genus $0$ and $1$ modular curves associated to isogeny graphs (see Lemma~\ref{lem:gen01classmodcurves}). In particular, we prove that $X_{H}$ has genus~$0$ if and only if $X_{H}$ is isomorphic to either a genus $0$ modular curve $X_{0}(n)$ or $H$ is conjugate to one of $H_{9}^{1}$, $H_{16}^{2}$, $H_{2}^{0}\times H_{9}^{1}$, or $H_{25}^{1}$. 

By Lemma~\ref{lem:gen01classmodcurves}, the genus $0$ modular curves associated to isogeny graphs fall into two families: those corresponding to the classical modular curves $X_{0}(n)$, and those corresponding to the four groups $H_{9}^{1},H_{16}^{2},H_{2}^{0}\times H_{9}^{1}$, and $H_{25}^{1}$. We now recall the explicit parameterizations associated with the genus $0$ modular curves $X_{0}(n)$ with $n>1$. Let $\Phi_{n}(X,Y)$ denote the level $n$ modular polynomial. When $X_{0}(n)$ has genus $0$, the classical Fricke parameterizations \cite{MR3838339} provides rational functions $j_{n,1}(t)$ and $j_{n,2}(t)$ satisfying
\[
\Phi_{n}(j_{n,1}(t),j_{n,2}(t))=0,
\]
thereby yielding a rational parameterization of the non-cuspidal points on $X_{0}(n)$. Barrios~\cite{Bariso} proved that these parameterizations were realized explicitly by isogenous families of elliptic curves $\mathcal{C}_{n,j}(t,d)$ (see Table 5 of loc. cit. for the corresponding Weierstrass models). More precisely, if $E_{1}$ and $E_{2}$ are $n$-isogenous over a field $K$ of characteristic $0$ or relatively prime to $6n$, and their $j$-invariants are not both identically $0$ or $1728$, then $E_{1}$ and $E_{2}$ occur as specializations of the families $\mathcal{C}_{n,j}(t,d)$. By~\cite[Lemma 3.3]{Bariso}, the restriction on the $j$-invariants is only necessary when $n$ is prime, and the remaining prime cases were subsequently treated in~\cite[Lemma~4.5]{disctwins}. From the perspective of Galois representations, these results admit a natural reinterpretation. They provide an explicit description of the $K$-rational non-cuspidal points on the modular curve $X_{0}(n)$.

In Section \ref{sec:expcalssgenus0}, we establish analogous results for the four remaining genus~$0$ modular curves~$X_{H}$ associated to isogeny graphs, thereby concluding the explicit classification of genus~$0$ modular curves associated to isogeny graphs. We achieve this by extending the families $\mathcal{C}_{n,j}(t,d)$ to larger families $\mathcal{C}_{n,m}^{1}(t,d)$ (see Table \ref{ta:curves}) for $n\in\left\{  9,16,18,25\right\}$. For each such $n$, the extended family contains the previously constructed family $\mathcal{C}_{n,j}(t,d)$ while having the additional elliptic curves needed to parameterize each vertex on the associated isogeny graph. In particular, we prove:

\begin{mainthm}\label{mainthm:explicitclassK}
Let $n\in\left\{  9,16,18,25\right\}$ and let $H_n$ be defined as:
\[
{\renewcommand{\arraystretch}{1.2}\begin{array}
[c]{c|cccc}
n & 9 & 16 & 18 & 25\\\hline
H_n & H_{9}^{1} & H_{16}^{2} & H_{2}^{0}\times H_{9}^{1} & H_{25}^{1}
\end{array}}
\]
Let $E$ be an elliptic curve over a field $K$ of characteristic $0$ or relatively prime to $6n$, and suppose that $\rho_{E,n}(G_{K})$ is conjugate to $H_{n}\leq\operatorname*{GL}\nolimits_{2}(\mathbb{Z}/n\mathbb{Z})$. Then there exist $t\in K$ and $d\in K^{\times}/(K^{\times})^{2}$ such that $E$ is $K$-isomorphic to $\mathcal{C}_{n,1}^{1}(t,d)$. In addition, if $\operatorname*{End}_{K}E\cong\mathbb{Z}$, then the partial isogeny graph associated to $\{  [  \mathcal{C}_{n,m}^{1}(t,d)]_{K}\}  _{m}$, given in Table \ref{isographs}, is a subgraph of $\mathcal{G}(E/K)$.
\end{mainthm}
Theorem \ref{mainthm:explicitclassK} is Theorem \ref{thm:explicitclassK}. As a consequence, we explicitly determine the $n$-division field for $n\in\{3,4,5\}$ of elliptic curves $E/K$ such that $\rho_{E,n}(G_{K})$ is diagonalizable (see Corollary~\ref{cor:ndivisionclass}). 

Our earlier discussion highlighted the difficulty of determining all possible isogeny graphs over a given characteristic $0$ field $K$. That said, certain fields have the property that the isogeny class degree uniquely determines the isogeny graph. More precisely, we have the following result, which is immediate from Theorem~\ref{mathmclass}.
\begin{maincor}\label{maincor1}
Let $K$ be a field of characteristic $0$ such that $i\not \in K$ and $\sqrt{p^{\ast}}\not \in K$ for each prime $p$, where $p^{\ast}=(-1)^{(p-1)/2}p$. If $E$ is an elliptic curve over a field $K$ with $\operatorname*{End}_{K}E\cong\mathbb{Z}$, then the isogeny class degree of $E$ uniquely determines the isogeny graph $\mathcal{G}(E/K)$. That is, if $\deg\mathcal{G}(E/K)=\prod_{p}p^{k_{p}}$, then
\[
\mathcal{G}(E/K)\cong\left\{
\begin{array}
[c]{cl}
\underset{p}{\square}\mathcal{H}_{p^{k_{p}}}^{0} & \text{if }k_{2}\leq1,\\
\underset{p\neq2}{\square}\mathcal{H}_{p^{k_{p}}}^{0}\square\mathcal{H}
_{2^{k_{2}}}^{1} & \text{if }k_{2}\geq2.
\end{array}
\right.
\]
\end{maincor}
In particular, if $K$ is a field containing no quadratic subfields, then every elliptic curve $E$ defined over $K$ satisfies $\operatorname*{End}_{K}E\cong\mathbb{Z}$. We note that the same conclusion holds if $K$ admits a real embedding (see Corollary~\ref{Cor:BloomInvReal}). Consequently, for such fields, the isogeny class degree uniquely determines the isogeny graph. This observation allows us to recover the classification of rational isogeny graphs, which we prove below.

\begin{mainthm}[{\cite[Theorem~1.2]{MR4203041}}]\label{ThmQ}
The Isogeny graph of an elliptic curve defined over $\Q$ is uniquely determined by its isogeny class degree. In particular, up to weight-preserving graph isomorphism, there are $26$ possible isogeny graphs over $\mathbb{Q}$. Table~\ref{isographsQ} lists the $26$ possible graphs and their Chiloyan and Lozano-Robledo (CLR) label \cite{MR4203041}.
\end{mainthm}
{\begingroup
\renewcommand{\arraystretch}{1.2}
 \begin{longtable}{ccc|ccc}
 	\caption{Isogeny graphs over $\mathbb{Q}$}\label{isographsQ}\\
	\hline
	$\deg\mathcal{G}(E/\mathbb{Q})$ & $\mathcal{G}(E/\mathbb{Q})$ & $\text{CLR label}$ & $\deg\mathcal{G}(E/\mathbb{Q})$ & $\mathcal{G}(E/\mathbb{Q})$ & $\text{CLR label}$ \\
	\hline

	\endfirsthead
	\hline
	$\deg\mathcal{G}(E/\mathbb{Q})$ & $\mathcal{G}(E/\mathbb{Q})$ & $\text{CLR label}$ & $\deg\mathcal{G}(E/\mathbb{Q})$ & $\mathcal{G}(E/\mathbb{Q})$ & $\text{CLR label}$ \\
	\hline
	\endhead
	\hline

	\multicolumn{4}{r}{\emph{continued on next page}}
	\endfoot
	\hline
	\endlastfoot

$1$ & $\mathcal{H}_{1}^{0}$ & $L_{1}$ & 4 & $\mathcal{H}_{4}^{1}$ & $T_{4}$\\\hline
$p\in\operatorname*{IsoPrimeDeg}(\mathbb{Q})$ & $\mathcal{H}_{p}^{0}$ & $L_{2}(p)$ & 8 & $\mathcal{H}_{8}^{1}$ & $T_{6}$\\\hline
$p^{2}$ for $p\in\{3,5\}$ & $\mathcal{H}_{p^{2}}^{0}$ & $L_{3}(p^2)$ & 12 &
$\mathcal{H}_{4}^{1}\square\mathcal{H}_{3}^{0}$ & $S$\\\hline
$27$ & $\mathcal{H}_{27}^{0}$ & $L_{4}$ & 16 & $\mathcal{H}_{16}^{1}$ &
$T_{8}$\\\hline
$
\begin{array}
[c]{c}
pq\text{\ for\ }p,q\text{ primes such that}\\
pq\in\{6,10,14,15,21\}
\end{array}
$ & $\mathcal{H}_{p}^{0}\square\mathcal{H}_{q}^{0}$ & $R_{4}(pq)$ & 18 &
$\mathcal{H}_{2}^{0}\square\mathcal{H}_{9}^{0}$ & $R_{6}$

    \end{longtable}
    \endgroup}

\begin{proof}
By the classification of rational cyclic isogeny degrees, there are exactly $26$ positive integers $n$ for which $X_0(n)$ has a non-cuspidal $\mathbb{Q}$-rational point. Thus, there are $26$ distinct positive integers occurring as isogeny class degrees over $\mathbb{Q}$. The result now follows by Corollary~\ref{maincor1}.
\end{proof}
In Section~\ref{sec:nffields}, we consider results of a similar flavor by considering those number fields for which the possible degrees of cyclic isogenies have been classified under the assumption of GRH.

We conclude by noting that the same Galois-representation viewpoint also points toward a natural base-change problem. In the torsion setting, the growth of torsion subgroups under base change has been studied extensively in the literature, including work of Bruin and Najman~\cite{BruinNajman}, as well as Gonz\'{a}lez-Jim\'{e}nez and Najman~\cite{GonzalezJimenezNajman1,GonzalezJimenezNajman2}, among others. Analogous questions have been considered in the cyclic isogeny setting, with work of Furio and Lombardo~\cite{FurioLombardo}, Najman~\cite{NajmanIsoj}, Novak~\cite{Novakisoprimedeg}, Vukorepa~\cite{Vukorepa}, and others addressing how cyclic isogenies can arise after base change. Motivated by these developments, our forthcoming article studies how the isogeny graph $\mathcal{G}(E/K)$ changes under base change to an extension $F/K$, and to what extent the resulting graphs can be classified when the isogeny graphs over the original field are already understood.

\textbf{Outline of paper.} Section \ref{sec:prelim} reviews some standard facts that are assumed throughout the paper. Section \ref{Sec:StructureGraphs} studies the structure of $\mathcal{G}(E/K)$, and contains the proof of Theorem~\ref{mainthmCarPro}. In Section~\ref{sec:graphs}, we classify the finite $p$-primary graphs. Section~\ref{subsec:finitegraphs} introduces the family of finite graphs $\mathcal{H}_{p^{k}}^{r}$, and Section~\ref{subsec:modpGalois} studies the $\operatorname{mod}p^{k}$ Galois representation, resulting in the proof of Theorem \ref{mathmclass} (1), and its corresponding converse. Section~\ref{sec:bloominv} introduces the $p$-blooming invariant of an elliptic curve and its properties. In particular, we show that it is an isogeny class invariant and that elliptic curves defined over fields admitting a real embedding have smallest possible $p$-blooming invariant. Section~\ref{sec:infgraphs} mirrors Section~\ref{sec:graphs}: after introducing the infinite families of graphs $\mathcal{H}_{p^{\infty}}^{r}$ and $\mathcal{H}_{p^{\infty,+}}^{r}$, we study the associated $p$-adic Galois representations, proving Theorem~\ref{mathmclass} (2), and its corresponding converse. Section~\ref{sec:pCM} applies the preceding theory to elliptic curves with potential CM and culminates with the proof of Theorem~\ref{mainthm:genkwon}. Section~\ref{sec:algorithm image->graph} presents an algorithm for determining the isogeny graph of an elliptic curve from a subset generating its adelic Galois image. Section~\ref{sec:ModularCurves} introduces modular curves associated to isogeny graphs via the families of subgroups $H_{p^k}^r$ arising from our classification. Section~\ref{sec:expcalssgenus0} completes the explicit classification of isogeny graphs associated to genus $0$ modular curves. We conclude with Section \ref{sec:nffields}, which provides classification results, conditional on GRH, for the possible isogeny graphs over certain number fields. Finally, we include two appendices. Appendix~\ref{LMFDBzoo} contains a table with one elliptic curve for each isogeny graph appearing in the LMFDB. Appendix~\ref{appendix_para} contains parametrizations of the elliptic curves whose isogeny graphs correspond to the genus $0$ modular curves considered in Section~\ref{sec:expcalssgenus0}.

\textbf{Computational resources.} The explicit computations in this paper were carried out using the
computer algebra systems \texttt{Magma}~\cite{magma} and
\textsc{SageMath}~\cite{sagemath}. The code and data required to
reproduce and verify these computations are publicly available at the GitHub repository
\begin{center}
\url{https://github.com/enrique-gonzalez-jimenez/isogeny-graphs}.
\end{center}
Throughout the paper, we use the LMFDB \cite{lmfdb} labeling conventions for elliptic curves and number fields.

\textbf{Acknowledgments.} The authors would like to thank Abbey Bourdon, John Cullinan, Maarten Derickx, and Filip Najman for helpful comments in the preparation of this article. Barrios and Novak also thank the Universidad Aut\'{o}noma de Madrid for their hospitality during visits to work on this project. During the carrying out of this work, Barrios was supported through research grants from the University of St. Thomas and an AMS Simons Research Enhancement Grant; Gonz\'{a}lez-Jim\'{e}nez was supported by Grant PID2022-138916NB-I00 funded by MCIN/AEI/10.13039/501100011033 and by ERDF A way of making Europe; and Novak was financed by the Croatian
Science Foundation under the project no. IP-2022-10-5008 and by the project “Implementation of cutting-edge research and its application as part of the Scientific Center of Excellence for Quantum and Complex Systems, and Representations of Lie Algebras,'' Grant No. PK.1.1.10.0004, co-financed by the European Union through the European Regional Development Fund - Competitiveness and Cohesion Programme 2021-2027.


\section{Preliminaries}\label{sec:prelim}

We begin with a brief review of the basic notions from the theory of elliptic curves, Galois representations, and graphs that will be used throughout this article. Unless stated otherwise, the background material can be found in \cite{DelRap, graphthy,GraphHandbook,ladiclass,Trees,ladics,MR2514094}. Henceforth, $K$ denotes a field of characteristic $0$, and $G_{K}=\operatorname*{Gal}(\overline{K}/K)$ is the absolute Galois group of $K$.

\subsection{Elliptic curves and isogenies}
Let $E$ be an elliptic curve over $K$, and denote its \textit{$j$-invariant} by $j(E)$. For an elliptic curve $E'$, an \textit{isogeny} between $E$ and $E'$ is a surjective morphism of elliptic curves $\pi:E\rightarrow E^{\prime}$. Since $K$ has characteristic $0$, every isogeny is separable, and so its degree is the order of its kernel. If $\ker \pi\cong\mathbb{Z}/n\mathbb{Z}$ for some positive integer $n$, then the isogeny is said to be cyclic, and we call it an $n$-isogeny. If $\pi$ is defined over $K$, or equivalently, if $\ker\pi$ is $G_{K}$-invariant, then $\pi$ is a $K$-rational isogeny. In that case, we also say that $E$ and $E'$ are $n$-isogenous over $K$.

Now let $n$ be a positive integer and let $E[n]=\left\{  P\in E(\overline{K})\mid\left[  n\right]  P=\mathcal{O}\right\}$ be the $n$-torsion subgroup of $E$. This is the kernel of the multiplication-by-$n$ map $\left[  n\right]  :E\rightarrow E$. The \textit{$n$-th division field} of $E/K$ is the extension $K(E[n])/K$, and the Weil pairing shows that $K(\mu_{n})\subseteq K(E[n])$, where $\mu_{n}$ denotes a primitive $n$-th root of unity.

We write $\operatorname*{Hom}(E,E^{\prime})$ for the group of all isogenies $\varphi:E\rightarrow E^{\prime}$, and $\operatorname*{Hom}_{K}(E,E^{\prime})$ for the subgroup consisting of isogenies defined over $K$. Next, for an integer $n$, let $\left[  n\right]  _{E}:E\rightarrow E$ and $\left[  n\right]  _{E^{\prime}}:E^{\prime}\rightarrow E^{\prime}$ denote the corresponding multiplication-by-$n$ maps on $E$ and $E^{\prime}$, respectively. Then, for $\varphi\in\operatorname*{Hom}_{K}(E,E^{\prime})$, we have that for each $P\in E$,
\[
\left(  \varphi\circ\left[  n\right]  _{E}\right)  \left(  P\right)
=\varphi\left(  nP\right)  =n\varphi(P)=\left(  \left[  n\right]  _{E^{\prime
}}\circ\varphi\right)  \left(  P\right)  .
\]
Consequently, for each $n\in\mathbb{Z}$, it is the case that $\varphi \circ\left[  n\right]  _{E}=\left[  n\right]  _{E^{\prime}}\circ\varphi\in\operatorname*{Hom}_{K}(E,E^{\prime})$.

When $E'=E$, an isogeny is an endomorphism of $E$. We write $\operatorname*{End}E$ for the endomorphism ring over $\overline{K}$, and $\operatorname*{End}\nolimits_{K}(E)$ for the subring consisting of endomorphisms defined over $K$. The ring $\operatorname*{End}E$ is isomorphic to either $\mathbb{Z}$ or an order in an imaginary quadratic field; in the latter case, we say that $E$ has complex multiplication (CM). In our setting, the assumption that $\operatorname*{End}_{K}E\cong\mathbb{Z}$ restricts $\operatorname*{Hom}_{K}(E,E^{\prime})$, as illustrated by our next result.

\begin{lemma}
\label{Lem:IsogenyGpRk1}Let $E$ and $E'$ be isogenous elliptic curves over a field $K$ of characteristic $0$. If $\operatorname*{End}_{K}E\cong\mathbb{Z}$, then $\operatorname*{Hom}_{K}(E,E^{\prime})$ is a torsion-free $\mathbb{Z}$-module of rank~$1$. 

In particular, if $\gamma\in\operatorname*{Hom}_{K}(E,E^{\prime})$ is a cyclic isogeny, then, for each $\varphi\in\operatorname*{Hom}_{K}(E,E^{\prime})$, there exist $n\in\mathbb{Z}$ such that $\varphi=\left[  n\right]  _{E^{\prime}}\circ\gamma=\gamma\circ\left[  n\right]  _{E}$.
\end{lemma}

\begin{proof}
By \cite[Proposition III.4.2]{MR2514094}, $\operatorname*{Hom}(E,E^{\prime})$ is a torsion-free $\mathbb{Z}$-module. Since $\operatorname*{Hom}_{K}(E,E^{\prime})\leq\operatorname*{Hom}(E,E^{\prime})$ is a nonzero submodule, it follows that $\operatorname*{Hom}_{K}(E,E^{\prime})$ is also a torsion-free $\mathbb{Z}$-module. It thus suffices to show that $\operatorname*{Hom}_{K}(E,E^{\prime})\otimes_{\mathbb{Z}}\mathbb{Q}$ is a $1$-dimensional vector space. To this end, let $\gamma:E\rightarrow E^{\prime}$ be a cyclic $K$-rational isogeny, and let $\widehat{\gamma}:E^{\prime}\rightarrow E$ denote the dual isogeny. In particular, $\gamma\circ\widehat{\gamma}=\left[  \deg\gamma\right]  _{E}$ and $\widehat{\gamma}\circ\gamma=\left[  \deg\gamma\right]  _{E^{\prime}}$. Now consider the $\mathbb{Z}$-module homomorphisms $f:\operatorname*{Hom}_{K}(E,E^{\prime})\rightarrow \operatorname*{End}_{K}\! E$ and $g:\operatorname*{End}_{K}\! E\rightarrow \operatorname*{Hom}_{K}(E,E^{\prime})$ defined by $f(\varphi)=\widehat{\gamma}\circ\varphi$ and $g(\psi)=\gamma\circ\psi$. After tensoring with $\mathbb{Q}$, we obtain the following $\mathbb{Q}$-linear maps:
\[
\begin{array}
[c]{rclcrcl}
\operatorname*{Hom}_{K}(E,E^{\prime})\otimes_{
\mathbb{Z}}\mathbb{Q}& \overset{f_{\mathbb{Q}}}{\longrightarrow} & \operatorname*{End}_{K}\! E\otimes_{\mathbb{Z}}\mathbb{Q}& \quad & \operatorname*{End}_{K}\! E\otimes_{\mathbb{Z}}\mathbb{Q}& \overset{g_{\mathbb{Q}}}{\longrightarrow} & \operatorname*{Hom}_{K}(E,E^{\prime})\otimes_{\mathbb{Z}}\mathbb{Q}\\
\varphi\otimes q & \longmapsto & \left(  \widehat{\gamma}\circ\varphi\right)\otimes q &  & \psi\otimes q & \longmapsto & \left(  \gamma\circ\psi\right)\otimes q
\end{array}
\]
Hence, for $\varphi\in\operatorname*{Hom}_{K}(E,E^{\prime})$ and $\psi\in\operatorname*{End}_{K}\! E$, we obtain
\begin{align*}
g_{\mathbb{Q}}\!\left(  f_{\mathbb{Q}}\!\left(  \varphi\otimes1\right)  \right)  
& = 
g_{\mathbb{Q}}\!\left(  \left(  \widehat{\gamma}\circ\varphi\right)  \otimes1\right) =
\left(  \gamma\circ\left(  \widehat{\gamma}\circ\varphi\right)  \right)\otimes1 
 =
\left(  \left[  \deg\gamma\right]_{E}\circ\varphi\right)\otimes1 
=
\deg\gamma\left(  \varphi\otimes1\right), \\
f_{\mathbb{Q}}\!\left(  g_{\mathbb{Q}}\!\left(  \psi\otimes1\right)  \right)  
& =
f_{\mathbb{Q}}\!\left(  \left(  \gamma\circ\psi\right)  \otimes1\right)  
 =
\left(\widehat{\gamma}\circ\left(  \gamma\circ\psi\right)  \right)  \otimes1 
=
\left(\left[  \deg\gamma\right]  _{E^{\prime}}\circ\psi\right)  \otimes1 
=
\deg\gamma\left(  \psi\otimes1\right)  .
\end{align*}
It thus follows that
\begin{align*}
g_{\mathbb{Q}}\circ f_{\mathbb{Q}}  & =\left[  \deg\gamma\right]  :\operatorname*{Hom}\nolimits_{K}(E,E^{\prime})\otimes_{\mathbb{Z}}\mathbb{Q}\longrightarrow\operatorname*{Hom}\nolimits_{K}(E,E^{\prime})\otimes_{\mathbb{Z}}\mathbb{Q},\\
f_{\mathbb{Q}}\circ g_{\mathbb{Q}}  & =\left[  \deg\gamma\right]  :\operatorname*{End}\nolimits_{K}(E)\otimes_{\mathbb{Z}}\mathbb{Q}\longrightarrow\operatorname*{End}\nolimits_{K}(E)\otimes_{\mathbb{Z}}\mathbb{Q}.
\end{align*}
In particular, $g_{\mathbb{Q}}\circ f_{\mathbb{Q}}$ and $f_{\mathbb{Q}}\circ g_{\mathbb{Q}}$ are endomorphisms of the corresponding $\mathbb{Q}$-vector spaces. Since $\deg\gamma\neq0$, we conclude that $f_{\mathbb{Q}}$ and $g_{\mathbb{Q}}$ are isomorphisms, as their corresponding inverses are obtained by dividing by $\deg\gamma$. Thus,
\[
\operatorname*{Hom}\nolimits_{K}(E,E^{\prime})\otimes_{\mathbb{Z}}\mathbb{Q}\cong\operatorname*{End}\nolimits_{K}(E)\otimes_{\mathbb{Z}}\mathbb{Q}
\]
as $\mathbb{Q}$-vector spaces. By assumption, $\operatorname*{End}\nolimits_{K}(E)\otimes_{\mathbb{Z}}\mathbb{Q}\cong\mathbb{Q}$, which shows that $\operatorname*{Hom}\nolimits_{K}(E,E^{\prime})$ is torsion-free $\mathbb{Z}$-module of rank $1$.

We now show that a cyclic isogeny $\gamma\in\operatorname*{Hom}_{K}(E,E^{\prime})$ can be taken as a generator for $\operatorname*{Hom}_{K}(E,E^{\prime})$. To this end, let $\alpha\in\operatorname*{Hom}_{K}(E,E^{\prime})$ be a generator, so that $\operatorname*{Hom}_{K}(E,E^{\prime})\cong\alpha\mathbb{Z}$ as $\mathbb{Z}$-modules. In particular, $\gamma=\alpha\circ\left[  n\right]  _{E}$ for some nonzero integer $n\in\mathbb{Z}$. Observe that if $n\geq2$, then $P\in E[n]$ satisfies $\gamma(P)=\alpha(nP)=\mathcal{O}$. Thus, $\left(\mathbb{Z}/n\mathbb{Z}\right)  ^{2}\cong E[n]\leq\ker\gamma$, which contradicts our assumption that $\gamma$ is cyclic. Therefore, $\left\vert n\right\vert =1$, which shows that without loss of generality, we may take $\alpha=\gamma$. The result now follows.
\end{proof}

Next, we recall that an elliptic curve over $K$ admits either $0,1,$ or $3$ distinct $K$-rational $2$-isogenies. For an odd prime $p$, the number of distinct $K$-rational $p$-isogenies is either $0,1,2,$ or $p+1$ distinct $p$-isogenies. Novak \cite{ivan-numberofiso} has recently extended this result to determine the number of $K$-rational $p^{k}$-isogenies for any positive integer $k$. We record this result below, and our classification of primary isogeny graphs builds on ideas developed in the proof of said result.

\begin{proposition}
[{\cite[Propositions 1.1 and 1.2]{ivan-numberofiso}}]Let $E$ be an elliptic curve over a field $K$ of characteristic $0$. For a prime $p$ and positive integer $k$, let $m$ be the number of distinct $p^{k}$-isogenies admitted by $E$. Then,
\begin{enumerate}
\item if $p=2$ and $k\ge 2$, then $m\in\left\{  0\right\}  \cup\left\{  2^{j}\mid j\in\left[  1,k\right]  \cap\mathbb{Z}\right\}  \cup\left\{  2^{k}+2^{k-1}\right\}  $;

\item if $p\geq3$, then $m\in\left\{  0\right\}  \cup\left\{  p^{j}\mid j\in\lbrack0,k)\cap\mathbb{Z}\right\}  \cup\left\{  2p^{j}\mid j\in\lbrack0,k)\cap\mathbb{Z}\right\}  \cup\left\{  p^{k}+p^{k-1}\right\}  $.
\end{enumerate}
\end{proposition}

Now let $\left[  E\right]_{K}$ denote the $K$-isomorphism class of $E$. The \textit{isogeny graph} of $E/K$, denoted $\mathcal{G}(E/K)$, is the undirected edge-weighted graph whose vertices are the $K$-isomorphism classes $\left[  E^{\prime}\right]_{K}$ of elliptic curves $E'$ isogenous to $E$, and whose edges are the prime-degree isogenies between representatives of the corresponding classes. An edge is assigned weight $p$ when it arises from an isogeny of prime degree $p$. The \textit{isogeny class} of $E/K$ is thus the vertex set $V(\mathcal{G}(E/K))$. For a prime $p$, we define the \textit{$p$-primary (isogeny) graph}, denoted $\mathcal{G}_{p}(E/K)$, as the full subgraph of $\mathcal{G}(E/K)$ consisting of those vertices $[E^{\prime}]_{K}$ such that $E$ and $E^{\prime}$ are $p^k$-isogenous for some positive integer $k$. Further note that $\mathcal{G}(E/K)$ and $\mathcal{G}_{p}(E/K)$ are connected graphs since all isogenies considered are separable. 

Next, we suppose $\mathcal{G}(E/K)$ is finite, and write 
\[
V(\mathcal{G}(E/K))=\left\{  \left[  E_{1}\right]_{K}, \left[  E_{2}\right]_{K}, \ldots, \left[  E_{k}\right]_{K}  \right\}.
\]
The \textit{isogeny matrix} of $E/K$ is the symmetric $k\times k$ matrix $M_{E/K}=\left(  a_{i,j}\right)  $, where $a_{i,j}$ denotes the smallest degree of a cyclic isogeny between $E_{i}$ and $E_{j}$. The \textit{isogeny class degree} of $E/K$, denoted $\deg\mathcal{G}(E/K)$, is the least common multiple of the entries in $M_{E/K}$. Note that the isogeny matrix depends on the chosen ordering of the vertices of $\mathcal{G}(E/K)$; different orderings produce permutation-similar matrices. In particular, the isogeny class degree is independent of this choice.

We note that by Faltings \cite{FaltingsThm}, the isogeny class $V(\mathcal{G}(E/K))$ of an elliptic curve $E$ over a number field $K$ is finite. However, if $E$ has CM over $K$, the set of edges may be infinite (see Examples~\ref{Ex:CMnotCarPro} and~\ref{ex:CMvolc}). However, in this setting, the isogeny matrix is still well-defined, as it records the smallest degree of a cyclic isogeny between two vertices. In particular, the isogeny class degree is well-defined so long as the isogeny class of $E$ is finite.

Next, for a prime $p$, we recall that the $p$-Bruhat-Tits tree $\mathcal{T}_{p}$ associated to $\operatorname*{PGL}_{2}(\mathbb{Q}_{p})$. It is a $p+1$-regular infinite tree (see Figure~\ref{fig:BruhatTits} for an illustration of the $2$-Bruhat-Tits tree), whose vertices are homothety classes of $\mathbb{Z}_{p}$-lattices in $\mathbb{Q}_{p}^{2}$, and two homothety classes $[L]$ and $[L^{\prime}]$ are connected by an edge if their representatives may be chosen so that $pL\subsetneq L^{\prime}\subsetneq L$. Equivalently, $[L:L^{\prime}]=p$.

Now suppose that $E$ does not have CM. Then, over an algebraic closure $\overline{K}$ of $K$, the $p$-primary graph $\mathcal{G}_{p}(E/\overline{K})$ is isomorphic to $\mathcal{T}_{p}$. More precisely, this identification comes from the $p$-adic Tate module $V_{p}(E)=T_{p}(E)\otimes_{\mathbb{Z}_{p}}\mathbb{Q}_{p}$, whose $\mathbb{Z}_{p}$-lattices correspond to vertices of $\mathcal{T}_{p}$, and whose $p$-isogenies correspond to edges \cite{ssgraphs,defeo}. In particular, the $K$-rational $p$-primary graph $\mathcal{G}_{p}(E/K)$ embeds as a connected subgraph of $\mathcal{T}_{p}$, and hence is a tree. Moreover, the hypothesis $\operatorname*{End}_{K}E\cong\mathbb{Z}$ ensures that there are no loops coming from nontrivial $K$-endomorphisms~\cite[Lemma in IV-10]{ladics}, and so the resulting $p$-primary graph $\mathcal{G}_{p}(E/K)$ also embeds as a connected subgraph of $\mathcal{T}_{p}$, and thus is also a tree. 

Now suppose that $\mathcal{G}(E/K)$ is finite, and suppose that $v_{p}(\deg\mathcal{G}(E/K))=k$. By definition of the isogeny class degree and the fact that $\mathcal{G}_{p}(E/K)$ is a tree, we obtain that the longest path in $\mathcal{G}_{p}(E/K)$ has length $k$. In particular, there is some elliptic curve $E^{\prime}\in V(\mathcal{G}_{p}(E/K))$ that admits a $K$-rational $p^{k}$-isogeny, and no elliptic curve in $V(\mathcal{G}_{p}(E/K))$ admits a $K$-rational $p^{k+1}$-isogeny. For later reference, we record the following lemma, which summarizes this discussion.

\begin{lemma}\label{Lem:BruhatTitsTree}
Let $E$ be an elliptic curve over a field $K$ of characteristic $0$ with $\operatorname*{End}_{K}E\cong~\mathbb{Z}$. For each prime $p$, the $p$-primary graph $\mathcal{G}_{p}(E/K)$ is a tree. Further, if $v_{p}(\deg\mathcal{G}(E/K))=k$, then there exists an elliptic curve $E^{\prime}\in V(\mathcal{G}_{p}(E/K))$ such that $E^{\prime}$ admits a $K$-rational $p^{k}$-isogeny, and no elliptic curve in $V(\mathcal{G}_{p}(E/K))$ admits a $K$-rational $p^{k+1}$-isogeny.
\end{lemma}

\subsection{Galois representations and modular curves}
For a positive integer $n$, the group of $n$-torsion points $E[n]$ is a free module over $\mathbb{Z}/n\mathbb{Z}$ of rank $2$. Fixing a basis $\left(  P,Q\right)  $ identifies $E[n]$ with $\left(\mathbb{Z}/n\mathbb{Z}\right)^2$ by sending $P\longmapsto\left(  1,0\right)  $ and $Q\longmapsto\left(0,1\right)  $. Under this identification, we obtain an isomorphism $\iota:\operatorname*{Aut}(E[n])\rightarrow\operatorname*{GL}\nolimits_{2}(\mathbb{Z}/n\mathbb{Z})$, where an automorphism $\varphi \in \operatorname*{Aut}(E[n])$ is represented by the matrix
\[
\iota(\varphi)=\left(
\begin{array}
[c]{cc}
a & b\\
c & d
\end{array}
\right)
\]
with $\varphi(P)=aP+cQ$ and $\varphi(Q)=bP+dQ$. We thus obtain the \textit{$\operatorname{mod}n$ Galois representation attached to $E$},
\[
\rho_{E,n}:G_{K}\longrightarrow\operatorname*{Aut}(E[n])\cong\operatorname*{GL}\nolimits_{2}(\mathbb{Z}/n\mathbb{Z}),
\]
defined by $\rho_{E,n}(\sigma)=\iota(\sigma_{n})$, where $\sigma_{n}$ is the automorphism of $E[n]$ induced by $\sigma$. We note that our choice of $\iota$ corresponds to the left action of $\operatorname*{GL}\nolimits_{2}(\mathbb{Z}/n\mathbb{Z})$ on $\left(\mathbb{Z}/n\mathbb{Z}\right)^2$.

To define the adelic Galois representation, we choose a compatible system of bases for the torsion groups $E[n]$ as $n$ varies. Concretely, for each positive integer $n$, we choose a basis of $E[n]$ in such a way that whenever $n=n_1 n_2$, the basis of $E[n_{1}]$ (resp. $E[n_{2}]$) is obtained from the basis of $E[n]$ by multiplication-by-$n_2$ (resp. $n_1$). In this way, the chosen bases are compatible under the natural transition maps among the torsion subgroups. Let $\iota_{E}:\operatorname*{Aut}\!\left(  E(\overline{K})_{\text{tors}}\right)  \rightarrow\operatorname*{GL}\nolimits_{2}(\widehat{\mathbb{Z}})$ be the isomorphism determined by this system of bases, where $\widehat{\mathbb{Z}} = \underleftarrow{\lim}_n \mathbb{Z}/n\mathbb{Z}$. This gives rise to
the \textit{adelic Galois representation attached to $E$}
\[
\rho_{E}:G_{K}\longrightarrow\underset{n}{\underleftarrow{\lim}}\operatorname*{GL}\nolimits_{2}(\mathbb{Z}/n\mathbb{Z})\cong\operatorname*{GL}\nolimits_{2}(\widehat{\mathbb{Z}}),
\]
defined by $\rho_{E}(\sigma)=\iota_{E}(\sigma)$. 

The \textit{level} of an open subgroup $H\leq\operatorname*{GL}\nolimits_{2}(\widehat{\mathbb{Z}})$ is the least positive integer $n$ such that $H$ contains the kernel of the projection $\pi_{n}:\operatorname*{GL}\nolimits_{2}(\widehat{\mathbb{Z}})\rightarrow\operatorname*{GL}\nolimits_{2}(\mathbb{Z}/n\mathbb{Z})$. Equivalently, $H=\pi_{n}^{-1}(\pi_{n}(H))$, and so the subgroup $H$ is determined by its image in $\operatorname*{GL}\nolimits_{2}(\mathbb{Z}/n\mathbb{Z})$. Moreover,
\[
\left[  \operatorname*{GL}\nolimits_{2}(\widehat{\mathbb{Z}}):H\right]  =\left[  \operatorname*{GL}\nolimits_{2}(\mathbb{Z}/n\mathbb{Z}):\pi_{n}(H)\right]. 
\]
Thus, for open subgroups $H\leq\operatorname*{GL}\nolimits_{2}(\widehat{\mathbb{Z}})$, it suffices to consider them as subgroups of $\operatorname*{GL}\nolimits_{2}(\mathbb{Z}/n\mathbb{Z})$ for some positive integer $n$. For a prime $p$, we also consider the projection $\pi_{p^{\infty}}:\operatorname*{GL}\nolimits_{2}(\widehat{\mathbb{Z}})\rightarrow\operatorname*{GL}\nolimits_{2}(\mathbb{Z}_{p})$. This gives rise to the \textit{$p$-adic Galois representation attached to $E$},
\[
\rho_{E,p^{\infty}}=\pi_{p^{\infty}}\circ\rho_{E}:G_{K}\longrightarrow
\underset{k}{\underleftarrow{\lim}}\operatorname*{GL}\nolimits_{2}(\mathbb{Z}/p^{k}\mathbb{Z})\cong\operatorname*{GL}\nolimits_{2}(\mathbb{Z}_{p}).
\]

For a positive integer $n$, let $B_0(n)\leq \operatorname*{GL}\nolimits_{2}(\mathbb{Z}/n\mathbb{Z})$ denote the Borel subgroup of upper triangular matrices. When $p$ is prime, we let $B_0(p^\infty)$ denote the Borel subgroup of upper triangular matrices in $\operatorname*{GL}\nolimits_{2}(\mathbb{Z}_p)$.  In our setting, we are motivated by the following basic lemma regarding the $\operatorname{mod}n$ Galois representation.

\begin{lemma}\label{lem:Borel}
Let $E$ be an elliptic curve defined over a field of characteristic $0$. If $E$ admits an $n$-isogeny over $K$, then the image of $\rho_{E,n}$ is conjugate to a subgroup of $B_0(n)$.
\end{lemma}

\begin{proof}
Let $\pi:E\rightarrow E^{\prime}$ be an $n$-isogeny, and let $C=\ker\pi$. Since $\pi$ is defined over $K$, the subgroup $C \subseteq E[n]$ is $G_{K}$-stable and cyclic of order $n$. Choose a basis $\left(  P,Q\right)$ of $E[n]$ such that $C=\left\langle P\right\rangle $. Then, for every $\sigma\in G_{K}$, we have that $\sigma(P)\in C$, and thus the image of $P$ under $\rho_{E,n}(\sigma)$ is a multiple of $P$. Thus, with respect to the basis $\left(  P,Q\right)  $, every matrix in $\rho_{E,n}(G_{K})$ is upper triangular. Hence, $\rho_{E,n}(G_{K})$ is contained in a conjugate of $B_0(n)$.
\end{proof}

In this article, we will work with subgroups of $\operatorname*{GL}\nolimits_{2}(\widehat{\mathbb{Z}})$, up to conjugacy. Thus, we view each subgroup $H$ as a representative of its conjugacy class. In particular, all inclusions $H_{1}\leq H_{2}$ of subgroups of $\operatorname*{GL}\nolimits_{2}(\widehat{\mathbb{Z}})$ should be understood to mean that $H_{1}$ is conjugate to a subgroup of $H_{2}$.

Now let $H\leq\operatorname*{GL}\nolimits_{2}(\widehat{\mathbb{Z}})$ be an open subgroup of level $n$. An \textit{$H$-level structure} on an elliptic curve $E$ is an $H$-orbit of isomorphisms $\iota:E[n]\rightarrow\left(\mathbb{Z}/n\mathbb{Z}\right)  ^{2}$. In particular, if $\iota'$ is another such isomorphism, we write $\iota\sim_{H}\iota^{\prime}$ if and only if $\iota=h\circ\iota^{\prime}$ for some $h\in H$, and we denote the corresponding equivalence class by $[\iota]_{H}$. 

The modular curve $Y_{H}(\overline{K})$ consists of equivalence classes of pairs $\left(  E,[\iota]_{H}\right)  $, where $\left(E,[\iota]_{H}\right)  \sim\left(  E^{\prime},[\iota^{\prime}]_{H}\right)  $ if there is an isomorphism $\varphi:E\rightarrow E^{\prime}$ such that the induced isomorphism $\varphi_{n}:E[n]\rightarrow E^{\prime}[n]$ satisfies
$\iota\sim\iota^{\prime}\circ\varphi_{n}$. More broadly, the modular curve $Y_{H}$ (resp. $X_{H}$) is defined to be the coarse moduli space of the stack $\mathcal{M}_{H}^{0}$ (resp. $\mathcal{M}_{H}$) over $\operatorname*{Spec}\mathbb{Z}\!\left[1/n\right]  $, that parameterizes elliptic curves (resp. generalized elliptic curves) with $H$-level structure. By \cite[IV-3.1]{DelRap}, this is equivalent to $X_{H}$ being isomorphic to the coarse space of the stack quotient $X(n)/H$, where $X(n)$ is the classical modular curve parameterizing full level $n$ structure. We record the following result regarding the $K$-rational points of $X_H$.

\begin{lemma}[{\cite[Lemma 2.1]{JRDZB}}]\label{lem:modularcurveimage}
Let $E$ be an elliptic curve over a field $K$ of characteristic~$0$, and let $H\leq\operatorname*{GL}\nolimits_{2}(\widehat{\mathbb{Z}})$ be an open subgroup of level $n$. Then there exists an $\iota$ such that $\left(  E,[\iota]_{H}\right)  \in X_{H}(K)$ if and only if the image of $\rho_{E,n}$ is contained in a subgroup conjugate to~$H$. 
\end{lemma}

Now suppose that $H\leq\operatorname*{GL}\nolimits_{2}(\widehat{\mathbb{Z}})$ is an open subgroup of level $n$ such that $-I\in H$, and let $\zeta_{n}$ be a primitive $n$-th root of unity. From \cite[\S 4]{Deligne-Rapoport} and \cite[\S 2]{MR482230}, we obtain that the field of definition of the modular curve $X_{H}$ is $\mathbb{Q}(\zeta_{n})^{\det H}$. Since $\operatorname*{Gal}(\mathbb{Q}(\zeta_{n})/\mathbb{Q})\cong\left(\mathbb{Z}/n\mathbb{Z}\right)  ^{\times}$, we observe that the field of definition of $X_{H}$ is $\mathbb{Q}$ if and only if $\det H=\left(\mathbb{Z}/n\mathbb{Z}\right)  ^{\times}$.

Now suppose that $H_{1},H_{2}\leq\operatorname*{GL}\nolimits_{2}(\widehat{\mathbb{Z}})$ are open subgroups containing $-I$, of coprime levels $n_{1}$ and $n_{2}$, respectively. Then $H_{1}\cap H_{2}$ is an open subgroup of level $n_{1}n_{2}$, and the associated modular curve satisfies $X_{H_{1}\cap H_{2}}\cong X_{H_{1}}\times_{X(1)}X_{H_{2}}$, where the fiber product is taken over the $j$-line $X(1)$. Next, let $\pi_{n}:\operatorname*{GL}\nolimits_{2}(\widehat{\mathbb{Z}})\rightarrow\operatorname*{GL}\nolimits_{2}(\mathbb{Z}/n\mathbb{Z})$ be the projection modulo $n$, and let $\overline{H}_{i}=\pi_{n_{i}}(H_{i})$ for $i\in\{1,2\}$. We can then identify $\overline{H}_{1}\times\overline {H}_{2}$ with a subgroup of $\operatorname*{GL}\nolimits_{2}(\mathbb{Z}/n_{1}n_{2}\mathbb{Z})$ via the Chinese remainder theorem. Moreover, $\pi_{n_{1}n_{2}}^{-1}(\overline{H}_{1}\times\overline{H}_{2})=H_{1}\cap H_{2}$. Therefore, the modular curve associated to $H_{1}\cap H_{2}$ may equivalently be viewed as the modular curve associated to the subgroup $\overline{H}_{1}\times \overline{H}_{2}$ of $\operatorname*{GL}\nolimits_{2}(\mathbb{Z}/n_{1}n_{2}\mathbb{Z})$. In what follows, we set $X_{\overline{H}_{1}\times\overline{H}_{2}}:=X_{\overline{H}_{1}}\times_{X(1)}X_{\overline{H}_{2}}$, identifying the finite-level subgroup with its full inverse image in $\operatorname*{GL}\nolimits_{2}(\widehat{\mathbb{Z}})$.

\subsection{Graph theory}
We now recall the graph-theoretic notions used throughout the paper. By a
\textit{weighted graph} $\mathcal{G}$ we mean a set of vertices $V(\mathcal{G})$, a set of edges $E(\mathcal{G})$, and a weight function $w:E(\mathcal{G})\rightarrow\mathbb{N}$. For a vertex $v$ of $\mathcal{G}$, we write $\deg(v)$ for the \textit{degree} of $v$, that is, the number of edges incident to $v$. For vertices $v,w\in V(\mathcal{G})$, we write $\operatorname{dist}(v,w)$ for the \textit{distance} between $v$ and $w$, namely the minimum number of edges in a path joining $v$ and $w$. Unless stated otherwise, all distances considered throughout this article are unweighted; that is, edge weights are ignored when computing distances. A graph is said to be \textit{weight-isomorphic} to another if there is a graph isomorphism between the graphs that preserves edge weights. Since the graphs in this article are weighted graphs, we reserve the notation $\mathcal{G}\cong\mathcal{H}$ to mean that the graphs $\mathcal{G}$ and $\mathcal{H}$ are weight-isomorphic. 

We now fix some terminology. A \textit{path graph} is a graph with vertex set $\{v_0,v_1,\ldots,v_k\}$ and edges $\{v_{i-1},v_{i} \}$ for $1\le i \le k$. A \textit{line graph} is a bi-infinite path graph with vertex set $\left\{v_{i}\right\}  _{i\in\mathbb{Z}}$ and edges $\left\{  v_{i},v_{i+1}\right\}  $ for each $i\in \mathbb{Z}$. A \textit{ray graph} is a one-sided infinite path graph with vertex set $\left\{  v_{i}\right\}  _{i\ge 0}$ and edges $\left\{  v_{i},v_{i+1}\right\}  $ for each $i\ge 0$.

For later use, we also introduce the notion of a \textit{$k$-spine} for a tree graph $\mathcal{G}$, which is a path on $k+1$ vertices in $\mathcal{G}$. This terminology will be convenient when we discuss our family of finite graphs~$\mathcal{H}_{p^{k}}^{r}$. We note that if $\mathcal{G}$ is an infinite tree, then it must contain a ray graph as a subgraph, as illustrated by the next result.
\begin{lemma}[{\cite[ VI.\S2, Kőnig's infinity lemma]{konig1936theorie}}]\label{koenig}
    Let $G$ be an infinite, connected graph in which every vertex has finite degree. Then $G$ contains a ray graph as a subgraph.
\end{lemma}

Now let $\left\{  \mathcal{G}_{i}\right\}_{i\in I}$ be a family of weighted graphs indexed by a set $I$, and write $\square_{i\in I}\mathcal{G}_{i}$ for their \textit{Cartesian product}. Its vertex set is $\prod_{i\in I}V(\mathcal{G}_{i})$, and two vertices $x=\left(  x_{i}\right)_{i\in I}$ and $y=\left(  y_{i}\right)_{i\in I}$ are adjacent if and only if there exists a unique $j\in I$ such that $x_{j}$ and $y_{j}$ are adjacent in $\mathcal{G}_{j}$, while $x_{i}=y_{i}$ for each $i\neq j$. In that case, the edge weight of $\left\{  x,y\right\}  $ is the weight of the edge $\left\{x_{j},y_{j}\right\}  $ in $\mathcal{G}_{j}$. In particular, graphs with a single vertex do not contribute to the product. Accordingly, if $J\subseteq I$ denotes the set of indices such that $\mathcal{G}_{j}$ has at least two vertices, then
\[
\underset{i\in I}{\square}\mathcal{G}_{i}\cong\underset{j\in J}{\square}\mathcal{G}_{j}.
\]
Thus, if $J$ is finite and each $\mathcal{G}_{j}$ is finite, then the Cartesian product is finite, even if $I$ is infinite.

Now suppose that $x=(x_i)_{i \in I}$ and $y=(y_i)_{i\in I}$ are vertices in $\square_{i\in I}\mathcal{G}_{i}$ such that there exists $j\in I$ with $x_j=y_j$ and $x_i \neq y_i$ for all $i\neq j$. Then any path from $x_j$ to $y_j$ in $\mathcal{G}_j$ lifts to a path of the same length from $x$ to $y$ in the Cartesian product. In particular, the distance between $x$ and $y$ in the Cartesian product is the same as the distance between $x_j$ and $y_j$ in $\mathcal{G}_j$. 

When $I$ is infinite, however, the full Cartesian product may be much larger than the graphs we will encounter in this paper. Indeed, if each $  \mathcal{G}_{i}$ is infinite, the vertices of $\square_{i\in I}\mathcal{G}_{i}$ will be uncountable. Since isogeny graphs are countable, the full Cartesian product is not the correct object for our study of infinite isogeny graphs. Instead, we will use the weak Cartesian product relative to a chosen basepoint.

In this direction, we define a \textit{pointed weighted graph} to be the pair $(\mathcal{G},v)$, where $\mathcal{G}$ is a weighted graph and $v\in V(\mathcal{G})$ is a distinguished vertex, called the \textit{basepoint}. Suppose now that $\{(\mathcal{G}_i,v_i)\}_{i \in I}$ is a family of pointed weighted graphs indexed by $I$, and let $v=(v_i)_{i \in I} \in \prod_{i\in I}V(\mathcal{G}_i) $. The \textit{weak Cartesian product} of the pointed graphs $(\mathcal{G}_i,v_i)$, denoted $\square_{i\in I}(\mathcal{G}_{i},v_i)$, is the induced subgraph of $\square_{i\in I}\mathcal{G}_{i}$ whose vertices are the tuples $x=(x_i)_{i\in I}$ such that $x_i = v_i$ for all but finitely many $i \in I$. Equivalently, it consists of those vertices that differ from the basepoint $v$ in only finitely many coordinates. The edges, and their corresponding weight, are exactly those inherited from the ambient Cartesian product. Note that whenever the weak Cartesian product is used, it is with respect to pointed graphs. Whereas the notation $\square_{i\in I} \mathcal{G}_i$ is understood to denote the ordinary Cartesian product of graphs.

The weak Cartesian product is the graph-theoretic analogue of a direct sum. Further, if each $\mathcal{G}_i$ is connected, then the weak Cartesian product is exactly the connected component of the basepoint $v$ in the full Cartesian product. When $I$ is finite, the weak Cartesian product agrees with the Cartesian product. Thus, in the finite case, it suffices to consider the Cartesian product of graphs. We next record the following standard result for later reference.
\begin{lemma}\label{Lem:IsoCartGraphs}
Let $(\mathcal{G}_{i},v_i)$ and $(\mathcal{H}_{i},w_i)$ denote
pointed weighted graphs for each $i \in I$. If $f_{i}:\mathcal{G}_{i}\rightarrow\mathcal{H}_{i}$ is a weight-preserving isomorphism of pointed weighted graphs with $f(v_i)=w_i$ for each $i\in I$,
then there is a weight-preserving isomorphism
\[
f:\underset{i\in I}{\square}(\mathcal{G}_{i},v_i)\longrightarrow\underset{i \in I}{\square}(\mathcal{H}_{i},w_i).
\]
\end{lemma}

The next example illustrates the Cartesian product of graphs.

\begin{example}
Let $\mathcal{G}_{1},\mathcal{G}_{2},$ and $\mathcal{G}_{3}$ be path graphs with two vertices. Set $V(\mathcal{G}_{1})=\{x_1,x_2\}$, $V(\mathcal{G}_{2})=\{y_1,y_2\}$, and $V(\mathcal{G}_{3})=\{z_1,z_2\}$. Let $w_i$ be the weight of the unique edge in $\mathcal{G}_i$. Thus, our graphs are as given below:
\[
\begin{tikzcd}[column sep=large]
|[Style:circle, draw=blue, text=blue]| x_1
  \arrow[r, "w_1", no head, color=blue, Style:edge]
& |[Style:circle, draw=blue, text=blue]| x_2
&
|[Style:circle, draw=red, text=red]| y_1
  \arrow[r, "w_2", no head, color=red, Style:edge]
& |[Style:circle, draw=red, text=red]| y_2
&
|[Style:circle, draw=teal!70!black, text=teal!70!black]| z_1
  \arrow[r, "w_3", no head, color=teal!70!black, Style:edge]
& |[Style:circle, draw=teal!70!black, text=teal!70!black]| z_2
\end{tikzcd}
\]
Next, we consider the graph $\mathcal{G}=\square_{i=1}^{3}\mathcal{G}_{i}$. Its vertex set $\prod_i ^3 V(\mathcal{G}_i)$ consists of eight $3$-tuples, and \eqref{eq:square} gives the graph $\mathcal{G}$. 
\begin{equation}\label{eq:square}
\adjustbox{scale=0.7}{\begin{tikzcd}[column sep=huge, row sep=huge, transform shape]
& |[StyleEllipse]| {(x_1,y_2,z_2)}
  \arrow[rr, "w_1", no head, color=blue, StyleEdge]
  \arrow[dd, no head, dashed, color=teal!70!black, StyleEdge]
& &
|[StyleEllipse]| {(x_2,y_2,z_2)} \\
|[StyleEllipse]| {(x_1,y_1,z_2)}
  \arrow[rr, no head, color=blue, StyleEdge]
  \arrow[ru, "w_2", no head, color=red, StyleEdge]
& &
|[StyleEllipse]| {(x_2,y_1,z_2)}
  \arrow[ru, "w_2", no head, color=red, StyleEdge]
  \arrow[dd, no head, color=teal!70!black, StyleEdge]
& \\
& |[StyleEllipse]| {(x_1,y_2,z_1)}
  \arrow[rr, no head, dashed, color=blue, StyleEdge]
& &
|[StyleEllipse]| {(x_2,y_2,z_1)}
  \arrow[uu, "w_3"', no head, color=teal!70!black, StyleEdge] \\
|[StyleEllipse]| {(x_1,y_1,z_1)}
  \arrow[uu, "w_3", no head, color=teal!70!black, StyleEdge]
  \arrow[rr, "w_1", no head, color=blue, StyleEdge]
  \arrow[ru, "w_2", no head, dashed, color=red, StyleEdge]
& &
|[StyleEllipse]| {(x_2,y_1,z_1)}
  \arrow[ru, "w_2", no head, color=red, StyleEdge]
&
\end{tikzcd}}
\end{equation}
Lastly, let $\theta$ denote a root of $x^2-x+2$. In particular, $\theta \in K=\mathbb{Q}(\sqrt{-7})$. If $(w_1,w_2,w_3) = (2,3,5)$, then $\mathcal{G}$ is weight-isomorphic to the isogeny graph of the elliptic curve $E/K$ with LMFDB label \href{https://www.lmfdb.org/EllipticCurve/2.0.7.1/324.2/a/6}{324.2-a6}. A Weierstrass model for $E$ is
\[
y^2 + xy + \theta y = x^3 - x^2 +\theta x - \theta.
\]
\end{example}


\section{Primary decomposition of isogeny graphs}\label{Sec:StructureGraphs}

Let $E$ be an elliptic curve defined over a field of characteristic $0$ with $\operatorname*{End}_{K}E\cong\mathbb{Z}$. The goal of this section is to study the structure of the isogeny graph $\mathcal{G}(E/K)$ of $E$ and prove Theorem~\ref{mainthmCarPro}. In this direction, we begin by recalling that the full torsion subgroup of $E$ is
\[
E(\overline{K})_{\text{tors}}=\bigcup_{n\in\mathbb{N}}E[n].
\]
The full torsion subgroup of $E$ is naturally a $G_{K}$-module. An elliptic curve $E_{1}$ is $n$-isogenous to~$E$ if and only if there exists a cyclic subgroup $C_{1}\leq E(\overline{K})_{\text{tors}}$ of order $n$ such that $E_{1}\cong E/C_{1}$. If $C_{1}$ is $G_{K}$-invariant, then $E_{1}$ is defined over $K$, and the corresponding isogeny $\varphi_{1}:E\rightarrow E_{1}$ is $K$-rational. More generally, if $C_{2}\leq C_{1}$, then $\varphi$ factors through the intermediate elliptic curve $E_{2}\cong E/C_{2}$. In other words, there exists $K$-rational isogenies $\varphi_{1}:E\rightarrow E_{2}$ and $\varphi_{2}:E_{2}\rightarrow E_{1}$ such that $\varphi=\varphi_{2}\circ\varphi_{1}$. Further, the degrees of $\varphi_{1}$ and $\varphi_{2}$ are $\left\vert C_{2}\right\vert $ and $\left[  C_{1}:C_{2}\right]  $, respectively. Thus, an understanding of cyclic isogenies is equivalent to understanding the $G_K$-invariant cyclic submodules of $E(\overline{K})_\text{tors}$. With that, we now introduce the following definition:

\begin{definition}
Let $K$ be a field of characteristic $0$, and let $M$ be a $G_{K}$-module. Let $\mathcal{L}_{K}(M)$ be the undirected weighted graph whose vertices are $G_{K}$-invariant finite cyclic submodules of~$M$, and two vertices $C_1$ and $C_2$ are connected by an edge if and only if one of these is a submodule of the other of prime index. The weight of each edge corresponds to this index.
\end{definition}

In our setting, $M$ will be a $G_K$-submodule of $E(\overline{K})_{\text{tors}}$, as illustrated by the following example.

\begin{example}\label{latticeexample1}
Consider the elliptic curve $E:y^{2}=x^{3}-x$ over $\Q$, given by LMFDB label \href{https://www.lmfdb.org/EllipticCurve/Q/32/a/3}{32.a3}. The $4$-torsion subgroup $E[4]$ is generated by $P=\left(  \sqrt{2}+1,\sqrt{2}+2\right)  $ and $Q=(i,i-1)$. Consequently, the $4$-th division field of $E$ is $\mathbb{Q}(E[4])=\mathbb{Q}(\zeta_{8})$. In particular, if $K$ is a field such that $\zeta_{8}\in K$, then $G_K$ fixes each element of $E[4]$. Therefore, $\mathcal{L}_{K}(E\left[  4\right]  )=\mathcal{L}_{\mathbb{Q}(\zeta_{8})}(E\left[  4\right]  )$, which is given below:
\begin{equation}\label{ex:latQE4}
    \adjustbox{scale=0.7}{
\begin{tikzcd}
<P> & <P+2Q>                                                                             & <P+Q> &                                                                                                                               & <3P+Q> & <2P+Q>                                                                             & <Q> \\
    & {<(1,0)>} \arrow[lu, "2" description, no head] \arrow[u, "2" description, no head] &       & {<(-1,0)>} \arrow[lu, "2" description, no head] \arrow[ru, "2" description, no head]                                          &        & {<(0,0)>} \arrow[u, "2" description, no head] \arrow[ru, "2" description, no head] &     \\
    &                                                                                    &       & <\mathcal{O}> \arrow[llu, "2" description, no head] \arrow[u, "2" description, no head] \arrow[rru, "2" description, no head] &        &                                                                                    &    
\end{tikzcd}}
\end{equation}
Now suppose that $F$ is a field such that $\sqrt{2}\in F$ and $i\not \in F$. Then the subgroups generated by $2P+Q=\left(  -i,i+1\right)  $ and $Q$ are not $G_F$-invariant. Now consider the subgroups generated by $P+Q=\left(  \sqrt{2}-1,-\sqrt{2}i+2i\right)  $ and $3P+Q=\left(-\sqrt{2}-1,-\sqrt{2}i-2i\right)  $. Since the $x$-coordinate of both of these points are in $F$, we see that the subgroups generated by these points are $G_F$-invariant. Lastly, $P,P+2Q=\left(-\sqrt{2}+1,-\sqrt{2}+2\right)  \in F^2$, and so $\mathcal{L}_{F}(E\left[4\right]  )$ is given by the following graph:
\begin{equation}\label{ex:latQ2noti}
    \adjustbox{scale=0.7}{
\begin{tikzcd}
<P> & <P+2Q>                                                                             & <P+Q> &                                                                                                                               & <3P+Q> &           \\
    & {<(1,0)>} \arrow[lu, "2" description, no head] \arrow[u, "2" description, no head] &       & {<(-1,0)>} \arrow[lu, "2" description, no head] \arrow[ru, "2" description, no head]                                          &        & {<(0,0)>} \\
    &                                                                                    &       & <\mathcal{O}> \arrow[llu, "2" description, no head] \arrow[u, "2" description, no head] \arrow[rru, "2" description, no head] &        &          
\end{tikzcd}
}\end{equation}
Lastly, if $L$ is a field such that $\sqrt{2},i\not \in L$, then the only subgroups that are $G_L$-invariant are the subgroups generated by $\left(  1,0\right)  ,\left(  -1,0\right)  ,$ and $\left(  0,0\right)  $. Thus, $\mathcal{L}_{L}(E\left[  4\right]  )$ is as follows:
\begin{equation}\label{ex:latQ}
    \adjustbox{scale=0.9}{
\begin{tikzcd}
{<(1,0)>} &  & {<(-1,0)>}                                                                                                                    &  & {<(0,0)>} \\
          &  & <\mathcal{O}> \arrow[llu, "2" description, no head] \arrow[u, "2" description, no head] \arrow[rru, "2" description, no head] &  &          
\end{tikzcd}
}\end{equation}
\end{example}

We are now ready to establish our first result, which shows that $\mathcal{L}_{K}(E(\overline{K})_{\text{tors}})\cong\mathcal{G}(E/K)$.

\begin{proposition}\label{prop:lattice-graph}
Let $E$ be an elliptic curve over a field $K$ of characteristic $0$ with $\operatorname*{End}_{K}E\cong\mathbb{Z}$. Then, the mapping $f:\mathcal{L}_{K}(E(\overline{K})_{\text{tors}})\rightarrow\mathcal{G}(E/K)$ defined by $f(C)=[E/C]_{K}$ is a weight-preserving graph isomorphism.
\end{proposition}

\begin{proof}
Let $C_{1},C_{2}\in V(\mathcal{L}_{K}(E(\overline{K})_{\text{tors}}))$ such that $\left\{  C_{1},C_{2}\right\}  $ is an edge in $\mathcal{L}_{K}(E(\overline{K})_{\text{tors}})$. Then, $C_{1}$ and $C_{2}$ are $G_{K}$-invariant, and without loss of generality, we may assume that $\left[C_{1}:C_{2}\right]  =~p$ for some prime $p$. Since $\{\mathcal{O}\}\leq C_2 \leq C_1$, there is a $K$-rational isogeny $\varphi_{1}:E\rightarrow E/C_{1}$ that factors through a $p$-isogeny $\psi:E/C_{2}\rightarrow E/C_{1}$. Since $\psi$ is necessarily $K$-rational, we have that 
\[
\left\{  f(C_{1}),f(C_{2})\right\}  =\left\{  \left[E/C_{1}\right]  _{K},\left[  E/C_{2}\right]  _{K}\right\}
\]
is an edge in $\mathcal{G}(E/K)$. Thus, $f$ is a graph homomorphism. Note that this also shows that $f$ is weight-preserving.

Next, let $f_{V}$ be the corresponding function on the sets of vertices. Since each $\left[  E^{\prime}\right]  _{K}\in V(\mathcal{G}(E/K))$ corresponds to a cyclic subgroup $C\leq E(\overline{K})_{\text{tors}}$ such that $E^{\prime}\cong E/C$ and $C$ is $G_{K}$-stable, we have that $f_{V}$ is surjective. Now suppose that $C_{1},C_{2}\in V(\mathcal{L}_{K}(E(\overline{K})_{\text{tors}}))$ such that $f_V(C_{1})=f_V(C_{2})$. Then, there exists a $K$-isomorphism $\psi:E/C_{1}\rightarrow E/C_{2}$. We further have $K$-rational cyclic isogenies $\varphi_{1}:E\rightarrow E/C_{1}$ and $\varphi_{2}:E\rightarrow E/C_{2}$. Let $\widehat{\varphi_{1}}:E/C_{1}\rightarrow E$ and $\widehat{\varphi_{2}}:E/C_{2}\rightarrow E$ denote the dual isogeny of $\varphi_{1}$ and $\varphi_{2}$, respectively. From Lemma~\ref{Lem:IsogenyGpRk1}, we have that $\operatorname*{Hom}_{K}(E/C_{1},E)\cong\left\langle \widehat{\varphi_{1}}\right\rangle $ as $\mathbb{Z}$-modules. Since $\psi$ is an isomorphism, we have that $\widehat{\varphi_{2}}\circ\psi\in\operatorname*{End}_{K}(E/C_{1},E)$ is cyclic, and there exist a nonzero integer $n$ such that $\widehat{\varphi_{2}}\circ\psi=\widehat{\varphi_{1}}\circ\left[  n\right]  _{E/C_{1}}$. In particular, $\ker\left[  n\right]_{E/C_{1}}\leq\ker\!\left(  \widehat{\varphi_{2}}\circ\psi\right)  $. Since $\widehat{\varphi_{2}}\circ\psi$ is cyclic, it follows that $\left\vert n\right\vert =1$. By taking duals of $\widehat{\varphi_{2}}\circ\psi=\widehat{\varphi_{1}}\circ\left[  \pm 1\right]  _{E/C_{1}}$, we obtain
\[
\psi^{-1}\circ\varphi_{2}=\left[  \pm1\right]  _{E/C_{1}}\circ\varphi_{1}.
\]
Thus, $C_{1}=\ker\left(  \varphi_{1}\right)  =\ker\varphi_{2}=C_{2}$, which shows that $f_{V}$ is bijective. 

It remains to show that the map induced by $f$ on the sets of edges is surjective. To this end, let $\left\{  \left[  E_{1}\right]  _{K},\left[E_{2}\right]  _{K}\right\}  $ be an edge in $\mathcal{G}(E/K)$. In particular, $E_{1}$ and $E_{2}$ are non-isomorphic and $p$-isogenous over $K$. Since $f_{V}$ is bijective, there exist distinct $G_{K}$-invariant cyclic subgroups $C_{1},C_{2}\leq E(\overline{K})_{\text{tors}}$ such that $\left[E_{1}\right]  _{K}=\left[  E/C_{1}\right]  _{K}$ and $\left[  E_{2}\right]_{K}=\left[  E/C_{2}\right]  _{K}$. Next, let $\varphi_{1}:E\rightarrow E/C_{1}$ and $\varphi_{2}:E\rightarrow E/C_{2}$ be the corresponding cyclic isogenies. By assumption, there is a $K$-rational $p$-isogeny $\psi:E/C_{1}\rightarrow E/C_{2}$. Let $\widehat{\psi}:E/C_{2}\rightarrow E/C_{1}$ be its corresponding dual isogeny. By Lemma \ref{Lem:IsogenyGpRk1}, $\operatorname*{Hom}_{K}(E,E/C_{1})\cong\left\langle \varphi_{1}\right\rangle$ and $\operatorname*{Hom}_{K}(E,E/C_{2})\cong\left\langle \varphi_{2}\right\rangle $. Consequently, since $\widehat{\psi}\circ\varphi_{2}\in\operatorname*{Hom}_{K}(E,E/C_{1})$ and $\psi\circ\varphi_{1}\in\operatorname*{Hom}_{K}(E,E/C_{2})$, there exists nonzero integers $a$ and $b$ such that $\widehat{\psi}\circ\varphi_{2}=\varphi_{1}\circ\left[a\right]  _{E}$ and $\psi\circ\varphi_{1}=\varphi_{2}\circ\left[  b\right]_{E}$. In particular, we have the following commutative diagram.
\[
\begin{tikzcd}
E \arrow[r, "{[b]_E}"] \arrow[d, "\varphi_1"] & E \arrow[d, "\varphi_2"] \arrow[r, "{[a]_E}"] & E \arrow[d, "\varphi_1"] \\
E/C_1 \arrow[r, "\psi"]                     & E/C_2 \arrow[r, "\widehat{\psi}"]           & E/C_1                   
\end{tikzcd}
\]
Now observe that
\[
\left[  p\right]  _{E/C_{1}}\circ\varphi_{1}=\widehat{\psi}\circ\psi\circ\varphi_{1}=\widehat{\psi}\circ\varphi_{2}\circ\left[  b\right]_{E}=\varphi_{1}\circ\left[  a\right]  _{E}\circ\left[  b\right]  _{E}=\left[ab\right]  _{E/C_{1}}\circ\varphi_{1}.
\]
Thus, $ab=p$, and without loss of generality we may assume that $\left(|a|,|b|\right)  =\left(  p,1\right)  $. Then, $\psi\circ\varphi_{1}=[\pm 1]_{E/C_{2}}\circ\varphi_{2}$, which implies that $C_{1}=\ker\varphi_{1}\leq\ker\varphi_{2}=C_{2}$. Now observe that
\begin{align*}
\deg\varphi_{2}=\deg\psi\deg\varphi_{1}\qquad & \Longleftrightarrow
\qquad\left\vert C_{2}\right\vert =p\left\vert C_{1}\right\vert \\
& \Longleftrightarrow\qquad\left[  C_{2}:C_{1}\right]  =p.
\end{align*}
Consequently, $\left\{  C_{1},C_{2}\right\}  $ is an edge in $\mathcal{L}_{K}(E(\overline{K})_{\text{tors}})$, which is mapped to the edge $\left\{\left[  E_{1}\right]  _{K},\left[  E_{2}\right]  _{K}\right\}  $ by $f$. Thus, $f$ is a weight-preserving graph isomorphism.
\end{proof}

\begin{example}\label{example-32a}
For a field $K$ containing $\mathbb{Q}$, Example~\ref{latticeexample1} studied the structure of $\mathcal{L}_K(E[4])$, where $E:y^{2}=x^{3}-x$ is the elliptic curve given by LMFDB label
\href{https://www.lmfdb.org/EllipticCurve/Q/32/a/3}{32.a3}. To apply Proposition~\ref{prop:lattice-graph}, we need knowledge of the isogeny class degree of $E/K$, and $E$ must satisfy $\operatorname{End}_K(E)\cong \mathbb{Z}$. We start by noting that $j(E)=1728$, and thus $E$ has CM, with $\operatorname{End}(E)\cong\mathbb{Z}[i]$. So if $i \not \in K$, then $\operatorname{End}_K(E)\cong\mathbb{Z}$ is satisfied and Proposition~\ref{prop:lattice-graph} can be applied if $\deg\mathcal{I}(E/K)=4$.

We start by observing that $\deg\mathcal{G}(E/\mathbb{Q})=4$, and by Proposition~\ref{prop:lattice-graph}, we have a graph isomorphism between $\mathcal{G}(E/\mathbb{Q})$ and $\mathcal{L}_\mathbb{Q}(E[4])$, where the latter graph was given in \eqref{ex:latQ}. Below, we give  $\mathcal{G}(E/\mathbb{Q})$, where the vertices are given by their corresponding LMFDB labels and oriented to match its corresponding isomorphic graph  $\mathcal{L}_\mathbb{Q}(E[4])$ in \eqref{ex:latQ}.

\begin{equation}
\adjustbox{scale=0.8}{\begin{tikzpicture}[
    vertex/.style={inner sep=2pt},
    edge/.style={thick},
    every node/.style={font=\normalsize}
]

\node[vertex] (a1) at (0,0)
{{[\href{https://www.lmfdb.org/EllipticCurve/Q/32/a/1}{32.a1}]$_{\mathbb{Q}}$}};

\node[vertex] (a3) at (3,0)
{{[\href{https://www.lmfdb.org/EllipticCurve/Q/32/a/3}{32.a3}]$_{\mathbb{Q}}$}};

\node[vertex] (a4) at (6,0)
{{[\href{https://www.lmfdb.org/EllipticCurve/Q/32/a/4}{32.a4}]$_{\mathbb{Q}}$}};

\node[vertex] (a2) at (3,2)
{{[\href{https://www.lmfdb.org/EllipticCurve/Q/32/a/2}{32.a2}]$_{\mathbb{Q}}$}};

\draw[edge] (a1) -- node[above] {$2$} (a3);
\draw[edge] (a3) -- node[above] {$2$} (a4);
\draw[edge] (a3) -- node[right] {$2$} (a2);

\end{tikzpicture}}
\end{equation}

Next, let $F=\mathbb{Q}(\sqrt{2})$. Base changing $E$ to $F$ yields the elliptic curve given by LMFDB label \href{https://www.lmfdb.org/EllipticCurve/2.2.8.1/32.1/a/4}{2.2.8.1-32.1-a4}. We note that $\deg\mathcal{G}(E/F)=16$ and that $\mathcal{G}(E/F)$ is graph isomorphic to $\mathcal{L}_{F}(E[4])$ by Proposition~\ref{prop:lattice-graph}. The former graph is given in~\eqref{ex:latQ2noti}, and below, we give the graph of $\mathcal{G}(E/F)$, where the vertices are given by their LMFDB label with respect to the number field $F$, and the vertices are oriented to coincide with its isomorphic
graph in~\eqref{ex:latQ2noti}.

\begin{equation}\label{ex:gra32Qsq2}
\adjustbox{scale=0.85}{
\begin{tikzpicture}[
    vertex/.style={inner sep=2pt},
    edge/.style={thick},
    every node/.style={font=\normalsize}
]

\node[vertex] (a7) at (0,0) 
{{$[\href{https://www.lmfdb.org/EllipticCurve/2.2.8.1/32.1/a/7}{32.1\text{-}a7}]_{F}$}};

\node[vertex] (a5) at (3,0) 
{{$[\href{https://www.lmfdb.org/EllipticCurve/2.2.8.1/32.1/a/5}{32.1\text{-}a5}]_{F}$}};

\node[vertex] (a4) at (6,0) 
{{$[\href{https://www.lmfdb.org/EllipticCurve/2.2.8.1/32.1/a/4}{32.1\text{-}a4}]_{F}$}};

\node[vertex] (a6) at (9,0) 
{{$[\href{https://www.lmfdb.org/EllipticCurve/2.2.8.1/32.1/a/6}{32.1\text{-}a6}]_{F}$}};

\node[vertex] (a8) at (12,0) 
{{$[\href{https://www.lmfdb.org/EllipticCurve/2.2.8.1/32.1/a/8}{32.1\text{-}a8}]_{F}$}};

\node[vertex] (a1) at (3,2) 
{[$\href{https://www.lmfdb.org/EllipticCurve/2.2.8.1/32.1/a/1}{32.1\text{-}a1}]_{F}$};

\node[vertex] (a3) at (6,2) 
{{$[\href{https://www.lmfdb.org/EllipticCurve/2.2.8.1/32.1/a/3}{32.1\text{-}a3}]_{F}$}};

\node[vertex] (a2) at (9,2) 
{{$[\href{https://www.lmfdb.org/EllipticCurve/2.2.8.1/32.1/a/2}{32.1\text{-}a2}]_{F}$}};

 horizontales
\draw[edge] (a7) -- node[above] {$2$} (a5);
\draw[edge] (a5) -- node[above] {$2$} (a4);
\draw[edge] (a4) -- node[above] {$2$} (a6);
\draw[edge] (a6) -- node[above] {$2$} (a8);

 verticales
\draw[edge] (a5) -- node[right] {$2$} (a1);
\draw[edge] (a4) -- node[right] {$2$} (a3);
\draw[edge] (a6) -- node[right] {$2$} (a2);

\end{tikzpicture}}
\end{equation}

Lastly, let $L=\mathbb{Q}(i)$ so that $\operatorname*{End}_{L}E\cong\mathbb{Z}\lbrack i]$. In this case, $[\href{https://www.lmfdb.org/EllipticCurve/Q/32/a/1}{32.a1}]_{L}=[\href{https://www.lmfdb.org/EllipticCurve/Q/32/a/2}{32.a2}]_{L}$ and $[\href{https://www.lmfdb.org/EllipticCurve/Q/32/a/3}{32.a3}]_{L}=[\href{https://www.lmfdb.org/EllipticCurve/Q/32/a/4}{32.a4}]_{L}$. Consequently, the isogeny class degree of $E$ is $2$, and we
obtain the linear graph $\mathcal{G}(E/L)$ given below, which illustrates the
necessity of the condition $\operatorname*{End}_{K}E\cong\mathbb{Z}$ in Proposition~\ref{prop:lattice-graph}.
\begin{equation}\label{eq:isograQi}
\begin{tikzcd}
{[\href{https://www.lmfdb.org/EllipticCurve/Q/32/a/1}{32.a1}]_L} \arrow[r, "2", no head] & {[\href{https://www.lmfdb.org/EllipticCurve/Q/32/a/3}{32.a3}]_L}
\end{tikzcd}
\end{equation}
\end{example}

Next, let $p$ be a prime. The $p$-primary torsion subgroup is defined as $E[p^{\infty}]=\bigcup_{k\geq1}E[p^{k}]$. In particular, $E[p^{\infty}]\leq E(\overline{K})_{\text{tors}}$. Now observe that by definition, $\mathcal{L}_{K}(E[p^{\infty}])$ is a subgraph of $\mathcal{L}_K(E(\overline{K})_{\text{tors}})$. In fact, it is the full subgraph of $\mathcal{L}_{K}(E(\overline{K})_{\text{tors}})$ consisting of those vertices $C$ such that $\left\vert C\right\vert =p^{k}$ for some positive integer $k$. Further note that the isomorphism in Proposition~\ref{prop:lattice-graph} satisfies $f(\mathcal{O})=[E]_{K}$. This discussion, together with loc. cit., Lemma~\ref{Lem:BruhatTitsTree}, and the definition of $\mathcal{G}_{p}(E/K)$, results in the following immediate consequence.

\begin{corollary}\label{cor:fpdecom}
Let $E$ be an elliptic curve over a field $K$ of characteristic $0$ with $\operatorname*{End}_{K}E\cong\mathbb{Z}$. For each prime $p$, the primary graph $\mathcal{G}_{p}(E/K)$ is a tree. Further, the mapping $f_{p}:\left(  \mathcal{L}_{K}(E[p^{\infty}]),\{\mathcal{O}\}\right)  \rightarrow\left(  \mathcal{G}_{p}(E/K),[E]_{K}\right)  $ defined by $f_{p}(C)=[E/C]_{K}$ is a weight-preserving graph isomorphism of pointed graphs.
\end{corollary}

The $G_{K}$-module $E(\overline{K})_{\text{tors}}$ has a natural primary decomposition. More generally, the same discussion applies to any torsion abelian $G_{K}$-module, and it is in this level of generality that we state our first decomposition result. In this direction, let $A$ be a torsion abelian $G_{K}$-module. For a prime $p$ and positive integer $k$, let $A[p^{k}]=\left\{  a\in A\mid p^{k}a=0\right\}  $ and set $A[p^{\infty}]=\bigcup_{k\geq1}A[p^{k}]$. Then we have a natural primary decomposition of $G_{K}$-modules. Namely,
\[
A\cong\bigoplus_{p}A[p^{\infty}].
\]
Now let $C\leq A$ be a cyclic group, and for each prime $p$, the $p$-primary component of $C$ is $C[p^{\infty}]=C\cap A[p^{\infty}]$. Then, $C$ decomposes naturally as
\[
C\cong\bigoplus_{p}C[p^{\infty}].
\]
In fact, this is an isomorphism of $G_{K}$-modules, and thus, $C$ is $G_{K}$-invariant if and only if each $C[p^{\infty}]$ is $G_{K}$-invariant. We summarize this in the following lemma.

\begin{lemma}\label{lem:EKtorDec}
Let $A$ be a torsion abelian $G_{K}$-module. Then, there is a natural isomorphism of $G_{K}$-modules
\[
A\cong\bigoplus_{p}A[p^{\infty}].
\]
Further, if $C\leq A$ is cyclic, then $C$ is $G_{K}$-invariant if and only if each of its $p$-primary components $C[p^{\infty}]$ is $G_{K}$-invariant.
\end{lemma}

We now prove our first primary decomposition result.

\begin{proposition}\label{Lem:IsoGraPro}
Let $A$ be a torsion abelian $G_{K}$-module. Then, there is a weight-preserving isomorphism of pointed graphs
\[
f:\left(  \mathcal{L}_{K}\!\left(  A\right)  ,\{0\}\right)  \longrightarrow
\underset{p}{\square}\left(  \mathcal{L}_{K}\!\left(  A[p^{\infty}]\right)
,\{0\}\right).
\]
\end{proposition}

\begin{proof}
By definition of the weak Cartesian product, we have that $V\!\left(  \square_{p}\left(  \mathcal{L}_{K}\!\left(  A[p^{\infty}]\right)  ,\{0\}\right)  \right)  $ is
\[
\left\{  \left.  \left(  C_{p}\right)  _{p}\in\prod_{p}V\!\left(\mathcal{L}_{K}\!\left(  A[p^{\infty}]\right)  \right)  \right\vert C_{p}=\left\{  0\right\}  \text{ for all but finitely many }p\right\}  .
\]
In particular, for a vertex $\left(  C_{p}\right)  _{p}$ in $\square_{p}\left(  \mathcal{L}_{K}\!\left(  A[p^{\infty}]\right)  ,\{0\}\right)  $,
we have that $\prod_{p}C_{p}=\bigoplus_{p}C_{p}$ is a finite cyclic group.
Next, we observe that Lemma \ref{lem:EKtorDec} implies that there exists a $G_{K}$-equivariant isomorphism
\[
g:A\longrightarrow\bigoplus_{p}A[p^{\infty}]
\]
satisfying $g(C)=\bigoplus_{p}C[p^{\infty}]$ for each cyclic subgroup $C\leq A$. Further, $C$ is $G_{K}$-invariant if and only if $C[p^{\infty}]$ is $G_{K}$-invariant for each prime $p$. In particular, $g$ induces a function
\[
f:\left(  \mathcal{L}_{K}\!\left(  A\right)  ,\left\{  0\right\}  \right) \longrightarrow\underset{p}{\square}\left(  \mathcal{L}_{K}\!\left(A[p^{\infty}]\right)  ,\left\{  0\right\}  \right)
\]
defined by $f(C)=\left(  C_{p}\right)  _{p}$ where $C_{p}=C[p^{\infty}]$ for each prime $p$. By assumption, $C_{p}$ is trivial for all but finitely many primes $p$. Since $g$ is a $G_{K}$-equivariant isomorphism, we observe that the induced map on vertices $f_{V}$ is bijective. Indeed, if $\left(C_{p}\right)  _{p}$ is a vertex in $\square_{p}\left(  \mathcal{L}_{K}\!\left(  A[p^{\infty}]\right)  ,\left\{  0\right\}  \right)  $, then each $C_{p}$ is $G_{K}$-invariant and $\bigoplus_{p}C_{p}\leq\bigoplus_{p}A[p^{\infty}]$ is a $G_K$-invariant cyclic group. Thus, there exists $C\leq A$ such that $f(C)=\left(  C_{p}\right)_{p}$.

Now suppose that $\left\{  C,D\right\}  $ is an edge in $\mathcal{L}_{K}\!\left(  A\right)  $. Then $C$ and $D$ are $G_{K}$-invariant, and without loss of generality, we may assume that $D\leq C$ with $\left[  C:D\right]  =q$ for some prime $q$. In particular, $f(D)\leq f(C)$ with $\left[f(D):f(C)\right]  =q$. It follows that if $f(C)=\left(  C_{p}\right)  _{p}$ and $f(D)=\left(  D_{p}\right)  _{p}$, then $C_{p}=D_{p}$ for each prime $p\neq q$ and $\left[  D_{q}:C_{q}\right]  =q$. This shows that $\left\{f(C),f(D)\right\}  $ is an edge in $\square_{p}\left(  \mathcal{L}_{K}\!\left(  A[p^{\infty}]\right)  ,\left\{  0\right\}  \right)  $. Thus, $f$ is a weight-preserving graph homomorphism. It remains to show that $f$ is surjective on edges. To this end, let $\{\left(  C_{p}\right)  _{p},\left(D_{p}\right)  _{p}\}$ be an edge in $\square_{p}\left(  \mathcal{L}_{K}\!\left(  A[p^{\infty}]\right)  ,\left\{  0\right\}  \right)  $. By definition, we have that there exists a prime $q$ such that $C_{q}\neq D_{q}$ and $C_{p}=D_{p}$ for all primes $p\neq q$. Further, $\left\{  C_{q},D_{q}\right\}  $ is an edge in $\mathcal{L}_{K}\!\left(  A[p^{\infty}]\right)  $. So without loss of generality, we may assume that $D_{q}\leq C_{q}$ with $\left[  C_{q}:D_{q}\right]  =q$.\ Consequently, $C=\bigoplus_{p}C_{p}$ and $D=\bigoplus_{p}D_{p}$ are $G_{K}$-invariant with $D\leq C$ and $\left[C:D\right]  =q$. Then, $\left[  f_{V}^{-1}(C):f_{V}^{-1}(D)\right]  =q$ and thus, $\left\{  f_{V}^{-1}(C),f_{V}^{-1}(D)\right\}  $ is an edge in $\mathcal{L}_{K}\!\left(  A\right)  $, which concludes the proof.
\end{proof}

We now obtain Theorem \ref{mainthmCarPro}:

\begin{theorem}[Theorem~\ref{mainthmCarPro}]\label{Thm:CartesianPro}
Let $E$ be an elliptic curve over a field $K$ of characteristic $0$ with $\operatorname*{End}_{K}E\cong~\mathbb{Z}$. Then, there is a weight-preserving isomorphism of pointed graphs
\[
f:\left(  \mathcal{G}(E/K),[E]_{K}\right)  \longrightarrow\underset{p}{\square}\left(  \mathcal{G}_{p}(E/K),[E]_{K}\right)  .
\]
\end{theorem}

\begin{proof}
By Proposition \ref{prop:lattice-graph} and Lemma \ref{Lem:IsoGraPro}, we have the following weight-preserving graph isomorphisms of pointed graphs:
\begin{align*}
j  & :\left(  \mathcal{G}(E/K),[E]_{K}\right)  \longrightarrow\left(
\mathcal{L}_{K}(E(\overline{K})_{\text{tors}}),\{\mathcal{O\}}\right)  ,\\
h  & :\left(  \mathcal{L}_{K}(E(\overline{K})_{\text{tors}}),\{\mathcal{O\}}\right)  \longrightarrow\underset{p}{\square}\left(  \mathcal{L}_{K}\!\left(E[p^{\infty}]\right)  ,\left\{  \mathcal{O}\right\}  \right)  .
\end{align*}
Note that $h(C)=\left(  C_{p}\right)  _{p}$, and $C_{p}$ is trivial for all but finitely many $p$. By Corollary \ref{cor:fpdecom}, we have that for each prime $p$, there is a weight-preserving isomorphism of pointed graphs
\[
g_{p}:\left(  \mathcal{L}_{K}\!\left(  E[p^{\infty}]\right)  ,\{\mathcal{O}\}\right)  \longrightarrow\left(  \mathcal{G}_{p}(E/K),[E]_{K}\right)
\]
defined by $g_{p}(C)=[E/C]_{K}$. It then follows from Lemma~\ref{Lem:IsoCartGraphs} that there is a weight preserving graph isomorphism
\[
g:\underset{p}{\square}\left(  \mathcal{L}_{K}\!\left(  E[p^{\infty}]\right),\{\mathcal{O}\}\right)  \longrightarrow\underset{p}{\square}\left(\mathcal{G}_{p}(E/K),[E]_{K}\right)  .
\]
Our desired isomorphism is then
\[
f=g\circ h\circ j:\left(  \mathcal{G}(E/K),[E]_{K}\right)  \longrightarrow
\underset{p}{\square}\left(  \mathcal{G}_{p}(E/K),[E]_{K}\right). \qedhere
\]
\end{proof}

As a consequence of Theorem~\ref{mainthmCarPro} and Lemma~\ref{Lem:BruhatTitsTree}, we obtain that each elliptic curve in the isogeny class of $E$ has the same $p$-primary graph:

\begin{corollary}\label{cor:sameprimary}
Let $E$ be an elliptic curve over a field $K$ of characteristic $0$ with $\operatorname*{End}_{K}E\cong~\mathbb{Z}$. For each prime $p$, if $E^{\prime}$ is isogenous to $E$ over $K$, then $\mathcal{G}_{p}(E/K)\cong\mathcal{G}_{p}(E^{\prime}/K)$. In particular, if $v_{p}(\deg\mathcal{G}(E/K))=k$, then there exists an elliptic curve $E^{\prime}\in V(\mathcal{G}_{p}(E/K))$ such that $E^{\prime}$ admits a $K$-rational $p^{k}$-isogeny, and no elliptic curve in $V(\mathcal{G}_{p}(E/K))$ admits a $K$-rational $p^{k+1}$-isogeny.
\end{corollary}

\begin{example}
Consider the elliptic curve $E:y^{2}+xy+y=x^{3}-2731x-55146$, given by LMFDB label \href{https://www.lmfdb.org/EllipticCurve/Q/14/a/1}{14.a1}. We note that $E$ does not have CM, and that $\mathcal{G}_{p}(E/\mathbb{Q})$ is trivial for all $p\neq2,3$. The $2$- and $3$-isogeny graphs of $E/\mathbb{Q}$ are given below, where the vertices are given by their LMFDB labels.
\begin{equation*}
    \adjustbox{scale=0.7}{\begin{tikzcd}
{[\href{https://www.lmfdb.org/EllipticCurve/Q/14/a/1}{14.a1}]_\mathbb{Q}} \arrow[d, "2"', no head] & {[\href{https://www.lmfdb.org/EllipticCurve/Q/14/a/1}{14.a1}]_\mathbb{Q}} \arrow[r, "3", no head] & {[\href{https://www.lmfdb.org/EllipticCurve/Q/14/a/3}{14.a3}]_\mathbb{Q}} \arrow[r, "3", no head] & {[\href{https://www.lmfdb.org/EllipticCurve/Q/14/a/4}{14.a4}]_\mathbb{Q}} \\
{[\href{https://www.lmfdb.org/EllipticCurve/Q/14/a/2}{14.a2}]_\mathbb{Q}}                          &                                                                                                   &                                                                                                   &                                                                          
\end{tikzcd}}
\end{equation*}
By Theorem~\ref{Thm:CartesianPro}, $\mathcal{G}(E/\mathbb{Q})=\mathcal{G}_{2}(E/\mathbb{Q})\square\mathcal{G}_{3}(E/\mathbb{Q})$, which is depicted below.
\begin{equation*}
    \adjustbox{scale=0.7}{\begin{tikzcd}
{[\href{https://www.lmfdb.org/EllipticCurve/Q/14/a/1}{14.a1}]_\mathbb{Q}} \arrow[r, "3", no head] \arrow[d, "2", no head] & {[\href{https://www.lmfdb.org/EllipticCurve/Q/14/a/3}{14.a3}]_\mathbb{Q}} \arrow[r, "3", no head] \arrow[d, "2", no head] & {[\href{https://www.lmfdb.org/EllipticCurve/Q/14/a/4}{14.a4}]_\mathbb{Q}} \arrow[d, "2", no head] \\
{[\href{https://www.lmfdb.org/EllipticCurve/Q/14/a/2}{14.a2}]_\mathbb{Q}} \arrow[r, "3", no head]                         & {[\href{https://www.lmfdb.org/EllipticCurve/Q/14/a/6}{14.a6}]_\mathbb{Q}} \arrow[r, "3", no head]                         & {[\href{https://www.lmfdb.org/EllipticCurve/Q/14/a/5}{14.a5}]_\mathbb{Q}}                        
\end{tikzcd}}
\end{equation*}

Next, let $K=\mathbb{Q}(\sqrt{2})$. Then, $E/K$ is given by the LMFDB label \href{https://www.lmfdb.org/EllipticCurve/2.2.8.1/98.1/a/10}{2.2.8.1-98.1-a10}. Below, we omit reference to the label \href{https://www.lmfdb.org/NumberField/2.2.8.1}{2.2.8.1}, which corresponds to the number field $K$. Over $K$, we have that $\mathcal{G}_{p}(E/K)\cong\mathcal{G}_{p}(E/\mathbb{Q})$ for each odd prime $p$. The $2$-isogeny graph of $E/K$ includes two new elliptic curves, which do not arise from elliptic curves definer over $\Q$ under base change:
\begin{equation*}  \adjustbox{scale=0.75}{
\begin{tikzcd}
                                                                                                 & {[\href{https://www.lmfdb.org/EllipticCurve/2.2.8.1/98.1/a/7}{98.1\text{-}a7}]_K} \arrow[d, "2", no head] &                                                                                     \\
{[\href{https://www.lmfdb.org/EllipticCurve/2.2.8.1/98.1/a/1}{14.a2}]_K} \arrow[r, "2", no head] & {[\href{https://www.lmfdb.org/EllipticCurve/Q/14.a1/}{14.a1}]_K} \arrow[r, "2", no head]                  & {[\href{https://www.lmfdb.org/EllipticCurve/2.2.8.1/98.1/a/12}{98.1\text{-}a12}]_K}
\end{tikzcd}}
\end{equation*}
Consequently, $\mathcal{G}(E/K)=\mathcal{G}_{2}(E/K)\square\mathcal{G}_{3}(E/K)$, shown below, has isogeny class degree $36$, and consists of $12$ $K$-isomorphism classes of elliptic curves.
\begin{equation*}  \adjustbox{scale=0.75}{
\begin{tikzcd}
                                                                 &                                                                                          & {[\href{https://www.lmfdb.org/EllipticCurve/2.2.8.1/98.1/a/7}{98.1\text{-}a7}]_K} \arrow[r, "3", no head] \arrow[d, "2", no head]             & {[\href{https://www.lmfdb.org/EllipticCurve/2.2.8.1/98.1/a/6}{98.1\text{-}a6}]_K} \arrow[r, "3", no head] \arrow[d, "2", no head]             & {[\href{https://www.lmfdb.org/EllipticCurve/2.2.8.1/98.1/a/11}{98.1\text{-}a11}]_K} \arrow[d, "2", no head]          &                                                                                                            &                                                                                   \\
                                                                 &                                                                                          & {[\href{https://www.lmfdb.org/EllipticCurve/Q/14/a/1}{14.a1}]_K} \arrow[r, "3", no head] \arrow[lld, "2"', no head] \arrow[rrd, "2", no head] & {[\href{https://www.lmfdb.org/EllipticCurve/Q/14/a/3}{14.a3}]_K} \arrow[r, "3", no head] \arrow[rrd, "2", no head] \arrow[lld, "2"', no head] & {[\href{https://www.lmfdb.org/EllipticCurve/Q/14/a/4}{14.a4}]_K} \arrow[rrd, "2", no head] \arrow[lld, "2", no head] &                                                                                                            &                                                                                   \\
{[\href{https://www.lmfdb.org/EllipticCurve/Q/14/a/2}{14.a2}]_K} & {[\href{https://www.lmfdb.org/EllipticCurve/Q/14/a/6}{14.a6}]_K} \arrow[l, "3", no head] & {[\href{https://www.lmfdb.org/EllipticCurve/Q/14/a/5}{14.a5}]_K} \arrow[l, "3", no head]                                                      &                                                                                                                                               & {[\href{https://www.lmfdb.org/EllipticCurve/2.2.8.1/98.1/a/12}{98.1\text{-}a12}]_K} \arrow[r, "3"', no head]         & {[\href{https://www.lmfdb.org/EllipticCurve/2.2.8.1/98.1/a/9}{98.1\text{-}a9}]_K} \arrow[r, "3"', no head] & {[\href{https://www.lmfdb.org/EllipticCurve/2.2.8.1/98.1/a/4}{98.1\text{-}a4}]_K}
\end{tikzcd}}
\end{equation*}

Now let $F=\mathbb{Q}(\sqrt{-3})$. Then $E/F$ is given by the LMFDB label \href{https://www.lmfdb.org/EllipticCurve/2.0.3.1/196.2/a/10}{2.0.3.1-196.2-a10}. The label \href{https://www.lmfdb.org/NumberField/2.0.3.1}{2.0.3.1} corresponds to the number field $F$ and is omitted below. Over $F$, we have that $\mathcal{G}_{p}(E/F)\cong\mathcal{G}_{p}(E/\mathbb{Q})$ for each prime $p\neq3$. The $3$-isogeny graph of $E/F$, shown below, has two new elliptic curves that are not base changes from elliptic curves defined over $\Q$.
\begin{equation*}\adjustbox{scale=0.75}{
\begin{tikzcd}
                                                                 & {[\href{https://www.lmfdb.org/EllipticCurve/2.0.3.1/196.2/a/5}{196.2\text{-}a5}]_F} \arrow[d, "3", no head]                               &                                                                  \\
{[\href{https://www.lmfdb.org/EllipticCurve/Q/14.a1/}{14.a1}]_F} & {[\href{https://www.lmfdb.org/EllipticCurve/Q/14.a3/}{14.a3}]_F} \arrow[r, "3", no head] \arrow[l, "3"', no head] \arrow[d, "3", no head] & {[\href{https://www.lmfdb.org/EllipticCurve/Q/14.a4/}{14.a4}]_F} \\
                                                                 & {[\href{https://www.lmfdb.org/EllipticCurve/2.0.3.1/196.2/a/6}{196.2\text{-}a6}]_F}                                                       &                                                                 
\end{tikzcd}}
\end{equation*}
The isogeny graph $\mathcal{G}(E/F)$ is then given by $\mathcal{G}_{2}(E/F)\square\mathcal{G}_{3}(E/F)$, which is given below. In particular, the isogeny class degree of $E/F$ is $18$, and there are $10$ distinct $F$-isomorphism classes in the isogeny class of $E/F$.
\begin{equation*}\adjustbox{scale=0.75}{
   \begin{tikzcd}
                                                                                             &                                                                                                                                                                        &                                                                                              &                                                                                          & {[\href{https://www.lmfdb.org/EllipticCurve/2.0.3.1/196.2/a/8}{196.2\text{-}a8}]_F} \arrow[d, "3", no head]      &                                                                  \\
                                                                                             & {[\href{https://www.lmfdb.org/EllipticCurve/2.0.3.1/196.2/a/5}{196.2\text{-}a5}]_F} \arrow[d, "3"', no head] \arrow[rrru, "2"', no head]                               &                                                                                              & {[\href{https://www.lmfdb.org/EllipticCurve/Q/14.a2/}{14.a2}]_F} \arrow[r, "3", no head] & {[\href{https://www.lmfdb.org/EllipticCurve/Q/14.a6/}{14.a6}]_F} \arrow[d, "3", no head] \arrow[r, "3", no head] & {[\href{https://www.lmfdb.org/EllipticCurve/Q/14.a5/}{14.a5}]_F} \\
{[\href{https://www.lmfdb.org/EllipticCurve/Q/14.a1/}{14.a1}]_F} \arrow[rrru, "2"', no head] & {[\href{https://www.lmfdb.org/EllipticCurve/Q/14.a3/}{14.a3}]_F} \arrow[r, "3"', no head] \arrow[l, "3", no head] \arrow[d, "3"', no head] \arrow[rrru, "2"', no head] & {[\href{https://www.lmfdb.org/EllipticCurve/Q/14.a4/}{14.a4}]_F} \arrow[rrru, "2"', no head] &                                                                                          & {[\href{https://www.lmfdb.org/EllipticCurve/2.0.3.1/196.2/a/9}{196.2\text{-}a9}]_F}                              &                                                                  \\
                                                                                             & {[\href{https://www.lmfdb.org/EllipticCurve/2.0.3.1/196.2/a/6}{196.2\text{-}a6}]_F} \arrow[rrru, "2"', no head]                                                        &                                                                                              &                                                                                          &                                                                                                                  &                                                                 
\end{tikzcd}}
\end{equation*}

Lastly, we consider the field $KF=\mathbb{Q}(\sqrt{2},\sqrt{-3})$. This time around, $\mathcal{G}_{2}(E/KF)\cong\mathcal{G}_{2}(E/K)$, $\mathcal{G}_{3}(E/KF)\cong\mathcal{G}_{3}(E/F)$, and $\mathcal{G}_{p}(E/KF)$ remains trivial for each prime $p>3$. Thus, $\mathcal{G}(E/KF)=\mathcal{G}_{2}(E/KF)\square\mathcal{G}_{3}(E/KF)$. This isogeny graph consists of $20$ distinct $KF$-isomorphism classes of elliptic curves, four of which are not base changes from elliptic curves defined over $\Q,K,$ or $F$. We assign to these four elliptic curves the label 9604.ai for $1\leq i \leq 4$, where $9604$ references the norm of the conductor of $E/KF$. We note that this label is not an LMFDB label, as these elliptic curves are not currently in their database. With this notation, we give the graph of $\mathcal{G}(E/KF)$ below.
\begin{equation*}\adjustbox{scale=0.57}{
\begin{tikzcd}
                                                                                                                         &                                                                                                                                                   &                                                                                                 &                                                                                                                                                       &                                                                                                                                                                               &                                                                                                                              &                                                                                                                 & {[9604.1\text{-}a1]_{KF}} \arrow[d, "3", no head]                                                                                      &                                                                                        \\
                                                                                                                         &                                                                                                                                                   &                                                                                                 &                                                                                                                                                       & {[\href{https://www.lmfdb.org/EllipticCurve/2.0.3.1/196.2/a/5}{196.2\text{-}a5}]_{KF}} \arrow[d, "3"', no head] \arrow[rrru, "2", no head] \arrow[rrrdd, "2"', no head]       &                                                                                                                              & {[\href{https://www.lmfdb.org/EllipticCurve/2.2.8.1/98.1/a/7}{98.1\text{-}a7}]_{KF}} \arrow[r, "3", no head]    & {[\href{https://www.lmfdb.org/EllipticCurve/2.2.8.1/98.1/a/6}{98.1\text{-}a6}]_{KF}} \arrow[d, "3", no head] \arrow[r, "3", no head]   & {[\href{https://www.lmfdb.org/EllipticCurve/2.2.8.1/98.1/a/11}{98.1\text{-}a11}]_{KF}} \\
                                                                                                                         & {[\href{https://www.lmfdb.org/EllipticCurve/2.0.3.1/196.2/a/8}{196.2\text{-}a8}]_{KF}} \arrow[d, "3"', no head] \arrow[rrru, "2"', no head]       &                                                                                                 & {[\href{https://www.lmfdb.org/EllipticCurve/Q/14.a1/}{14.a1}]_{KF}} \arrow[r, "3"', no head] \arrow[rrru, "2"', no head] \arrow[rrrdd, "2"', no head] & {[\href{https://www.lmfdb.org/EllipticCurve/Q/14.a3/}{14.a3}]_{KF}} \arrow[r, "3"', no head] \arrow[rrru, "2"', no head] \arrow[d, "3", no head] \arrow[rrrdd, "2"', no head] & {[\href{https://www.lmfdb.org/EllipticCurve/Q/14.a4/}{14.a4}]_{KF}} \arrow[rrru, "2"', no head] \arrow[rrrdd, "2"', no head] &                                                                                                                 & {[9604.1\text{-}a2]_{KF}}                                                                                                              &                                                                                        \\
{[\href{https://www.lmfdb.org/EllipticCurve/Q/14.a2/}{14.a2}]_{KF}} \arrow[r, "3"', no head] \arrow[rrru, "2"', no head] & {[\href{https://www.lmfdb.org/EllipticCurve/Q/14.a6/}{14.a6}]_{KF}} \arrow[d, "3"', no head] \arrow[r, "3"', no head] \arrow[rrru, "2"', no head] & {[\href{https://www.lmfdb.org/EllipticCurve/Q/14.a5/}{14.a5}]_{KF}} \arrow[rrru, "2"', no head] &                                                                                                                                                       & {[\href{https://www.lmfdb.org/EllipticCurve/2.0.3.1/196.2/a/6}{196.2\text{-}a6}]_{KF}} \arrow[rrru, "2"', no head] \arrow[rrrdd, "2"', no head]                               &                                                                                                                              &                                                                                                                 & {[9604.1\text{-}a3]_{KF}} \arrow[d, "3"', no head]                                                                                     &                                                                                        \\
                                                                                                                         & {[\href{https://www.lmfdb.org/EllipticCurve/2.0.3.1/196.2/a/9}{196.2\text{-}a9}]_{KF}} \arrow[rrru, "2"', no head]                                &                                                                                                 &                                                                                                                                                       &                                                                                                                                                                               &                                                                                                                              & {[\href{https://www.lmfdb.org/EllipticCurve/2.2.8.1/98.1/a/12}{98.1\text{-}a12}]_{KF}} \arrow[r, "3"', no head] & {[\href{https://www.lmfdb.org/EllipticCurve/2.2.8.1/98.1/a/9}{98.1\text{-}a9}]_{KF}} \arrow[d, "3"', no head] \arrow[r, "3"', no head] & {[\href{https://www.lmfdb.org/EllipticCurve/2.2.8.1/98.1/a/4}{98.1\text{-}a4}]_{KF}}   \\
                                                                                                                         &                                                                                                                                                   &                                                                                                 &                                                                                                                                                       &                                                                                                                                                                               &                                                                                                                              &                                                                                                                 & {[9604.1\text{-}a4]_{KF}}                                                                                                              &                                                                                       
\end{tikzcd}
}\end{equation*}

\end{example}

For each of the isogeny graphs considered in Example~\ref{example-32a}, observe that $\deg\mathcal{G}(E/K)=\max M_{E/K}$. While this does not hold in general (see Example\ \ref{Ex:CMnotCarPro}), it does hold whenever $\operatorname{End}_K E \cong \mathbb{Z}$. This is demonstrated by the following proposition.

\begin{proposition}\label{prop:maxmat}
Let $E$ be an elliptic curve over a field $K$ of characteristic $0$ with $\operatorname*{End}_{K}E\cong\mathbb{Z}$, and suppose that $\mathcal{G}(E/K)$ is finite. For each prime $p$, let $M_{E/K}^{(p)}$ denote the isogeny matrix of $\mathcal{G}_{p}(E/K)$. Then, the isogeny matrix $M_{E/K}$ of $\mathcal{G}(E/K)$ is permutation-similar to the following Kronecker product:
\[
\bigotimes_{p}M_{E/K}^{(p)}.
\]
In particular, $\deg\mathcal{G}(E/K)=\max M_{E/K}$, and if $n=\deg\mathcal{G}(E/K)$, then there exists an elliptic curve $E^{\prime}$ in the isogeny class of $E$ such that $E^{\prime}$ admits a $K$-rational $n$-isogeny.
\end{proposition}

\begin{proof}
Let $P$ be the finite set of primes $p$ dividing $\deg\mathcal{G}(E/K)$ so that
\[
\underset{p}{\square}\mathcal{G}_{p}(E/K)\cong\underset{p\in P}{\square
}\mathcal{G}_{p}(E/K).
\]
By Theorem \ref{Thm:CartesianPro}, we may identify $\mathcal{G}(E/K)$ with $\square_{p\in P}\mathcal{G}_{p}(E/K)$. For each $p\in P$, fix an ordering
\[
V(\mathcal{G}_{p}(E/K))=\left\{  \left[  E_{1}^{(p)}\right]  _{K},\left[E_{2}^{(p)}\right]  _{K},\ldots,\left[  E_{n_{p}}^{(p)}\right]  _{K}\right\},
\]
and order the vertices of the Cartesian product lexicographically. Now let $$[F]_{K},[F^{\prime}]_{K}\in V(\square_{p\in P}\mathcal{G}_{p}(E/K)),$$
and with respect to the lexicographic ordering, we have tuples
\[
\lbrack F]_{K}=\left(  \left[  E_{i_{p}}^{(p)}\right]  _{K}\right)  _{p\in P}\qquad\text{and}\qquad\lbrack F^{\prime}]_{K}=\left(  \left[  E_{j_{p}}^{(p)}\right]  _{K}\right)  _{p\in P}.
\]

For each $p$, set $d_{p}=(M_{E/K}^{(p)})_{i_{p},j_{p}}$. By definition, $d_{p}$ is the smallest degree of a cyclic isogeny between $E_{i_{p}}^{(p)}$ and $E_{j_{p}}^{(p)}$, and $d_{p}$ is a power of $p$. We claim that the smallest degree of a cyclic isogeny between the elliptic curves $F$ and $F^{\prime}$ is $\prod_{p\in P}d_{p}$. In this direction, let $\phi :F\rightarrow F^{\prime}$ be any cyclic isogeny, and let $n=\deg\phi$. For each $p\in P$, there is a $G_{K}$-invariant cyclic subgroup $C_{p}$ of order $p^{v_{p}(n)}$ such that $C_{p}\leq\ker\phi$. Next, set $F_{p}=F/C_{p}$ so that $F$ and $F_{p}$ are $p^{v_{p}(n)}$-isogenous. Therefore, $F_{p}$ and $F^{\prime}$ are isogenous of degree coprime to $p$. In particular, the only contribution to the $p$-primary graph $\mathcal{G}_{p}(E/K)$ comes from the $p^{v_{p}(n)}$-isogeny $F\rightarrow F_{p}$, and so the vertices $[E_{i_{p}}^{(p)}]$ and $[E_{j_{p}}^{(p)}]$ in $\mathcal{G}_{p}(E/K)$ are connected by a path of length at most $v_{p}(n)$, and therefore $v_{p}(d_{p})\leq v_{p}(n)$ by minimality of $d_{p}$. Since $d_{p}$ is a power of $p$, we have that $d_{p}$ divides $p^{v_{p}(n)}$. Since this holds for each $p\in P$, we deduce that $\prod_{p\in P}d_{p}$ divides $n$.

We now show that equality must hold. To this end, set $P=\left\{  p_{1},p_{2},\ldots,p_{m}\right\}  $. Now let $\left[  F_{0}\right]  _{K}=\left[F\right]  _{K}$, and for $1\leq t\leq m$, let
\[
\left[  F_{t}\right]  _{K}=\left(  \left[  E_{j_{p_{1}}}^{(p_{1})}\right]_{K},\ldots,\left[  E_{j_{p_{t}}}^{(p_{t})}\right]  _{K},\left[E_{i_{p_{t+1}}}^{(p_{t+1})}\right]  _{K},\ldots,\left[  E_{i_{p_{m}}}^{(p_{m})}\right]_K  \right)  .
\]
In particular, $\left[  F_{m}\right]  _{K}=\left[  F^{\prime}\right]  _{K}$, and for $1\leq t\leq m$, the tuples $\left[  F_{t-1}\right]  _{K}$ and $\left[  F_{t}\right]  _{K}$ differ in exactly one coordinate, namely the coordinate corresponding to $p_{t}$. By definition of $d_{p_{t}}$, representatives corresponding to the vertices $[E_{i_{p_{t}}}^{(p_{t})}]_{K}$ and $[E_{j_{p_{t}}}^{(p_{t})}]_{K}$ in $\mathcal{G}_{p_{t}}(E/K)$ are joined by a cyclic isogeny of degree $d_{p_{t}}$. Since all other coordinates are unchanged, this lifts to a cyclic isogeny $F_{t-1}\rightarrow F_{t}$ of degree $d_{p_{t}}$ in the Cartesian product. Since $d_{p_{1}},d_{p_{2}},\ldots,d_{p_{m}}$ are powers of distinct primes, they are pairwise coprime. Therefore the composition
\[
\begin{tikzcd}
F=F_0 \arrow[r, "d_{p_1}"] & F_1 \arrow[r] & \cdots \arrow[r] & F_{m-1} \arrow[r, "d_{p_m}"] & F_m=F'
\end{tikzcd},
\]
and has degree $\prod_{p\in P}d_{p}$, and its kernel is cyclic since a product of cyclic groups of pairwise coprime order is cyclic. This shows that the smallest degree of a cyclic isogeny between $F$ and $F^{\prime}$ is $\prod_{p\in P}d_{p}$. 

Next, we observe that
\begin{equation}
\prod_{p\in P}(M_{E/K}^{(p)})_{i_{p},j_{p}}=\prod_{p\in P}d_{p}.\label{kronmatrix}
\end{equation}
By induction on \cite[Equation (13)]{Kroneckerlexi}, we have that with the lexicographic ordering, (\ref{kronmatrix}) is the $([F]_{K},[F^{\prime}]_{K})$ entry in $M_{E/K}$. Consequently,
\[
M_{E/K}=\bigotimes_{p\in P}M_{E/K}^{(p)}.
\]
Finally, since the entries of $M_{E/K}^{(p)}$ are powers of $p$, we have that $\operatorname{lcm}M_{E/K}^{(p)}=\max M_{E/K}^{(p)}$. Therefore,
\[
\deg\mathcal{G}(E/K)=\operatorname{lcm}M_{E/K}=\prod_{p\in P} \operatorname{lcm}M_{E/K}^{(p)}=\prod_{p\in P}\max M_{E/K}^{(p)}=\max M_{E/K}.
\]
Thus, if $n=\deg\mathcal{G}(E/K)$, then there exists an elliptic curve $E^{\prime}$ in the isogeny class of $E$ such that $E^{\prime}$ admits a $K$-rational $n$-isogeny, which concludes the proof.
\end{proof}

Motivated by Proposition \ref{prop:maxmat}, we now extend our definition of
isogeny class degree to infinite graphs:

\begin{definition}
Let $E$ be an elliptic curve over a field $K$ of characteristic $0$ with
$\operatorname*{End}_{K}E\cong\mathbb{Z}$. For a prime $p$, if $\mathcal{G}_{p}(E/K)$ is finite, let $M_{E/K}^{(p)}$ be the isogeny matrix associated to $\mathcal{G}_{p}(E/K)$. For each prime $p$, we define the \textit{$p$-isogeny class degree of $E/K$}, denoted $\deg\mathcal{G}_{p}(E/K)$, to be the supernatural number
\[
\deg\mathcal{G}_{p}(E/K)=\left\{
\begin{array}
[c]{cl}
\max M_{E/K}^{(p)} & \text{if }\mathcal{G}_{p}(E/K)\text{ is finite,}\\
p^{\infty} & \text{if }\mathcal{G}_{p}(E/K)\text{ is infinite.}
\end{array}
\right.
\]
The \textit{isogeny class degree of $E/K$}, denoted $\deg\mathcal{G}(E/K)$, is then defined as the
supernatural number
\[
\deg\mathcal{G}(E/K)=\prod_{p}\deg\mathcal{G}_{p}(E/K).
\]
\end{definition}

We conclude by noting that from Theorem~\ref{mainthmCarPro}, we have that if $\operatorname*{End}_{K}E\cong\mathbb{Z}$, then $\mathcal{G}(E/K)$ is graph isomorphic to the Cartesian product of its $p$-primary graphs. The following example demonstrates the necessity of the assumption. In particular, it illustrates that elliptic curves with complex multiplication may have isogeny graphs that are not expressible as a Cartesian product of $p$-primary graphs.

\begin{example}\label{Ex:CMnotCarPro}
Let $E:y^{2}=x^{3}-x$ be the elliptic curve with LMFDB label \href{https://www.lmfdb.org/EllipticCurve/Q/32/a/3}{32.a3}. Over the field $F=\mathbb{Q}(\sqrt{2})$, $E$ is given by LMFDB label \href{https://www.lmfdb.org/EllipticCurve/2.2.8.1/32.1/a/4}{2.2.8.1-32.1-a4}. Below, we omit reference to the label \href{https://www.lmfdb.org/NumberField/2.2.8.1}{2.2.8.1}, which corresponds to the number field $F$. The isogeny graph of $E/F$ is given in (\ref{ex:gra32Qsq2}). Since $\operatorname*{End}_{F}E\cong\mathbb{Z}$ and $\deg\mathcal{G}(E/F)=4$, we have that $\mathcal{G}(E/F)$ is graph isomorphic to $\mathcal{L}_{F}(E[4])$, whose graph is given in (\ref{ex:latQ2noti}). We now consider the field extension $K=\mathbb{Q}(i,\sqrt{2})=\mathbb{Q}(\zeta_{8})$. Then $\operatorname*{End}_{K}E\cong\mathbb{Z}\lbrack i]$, and we have the following equalities of $K$-isomorphism classes:
\[
\begin{array}
[c]{ccc}
\left[  \href{https://www.lmfdb.org/EllipticCurve/2.2.8.1/32.1/a/1}{32.1\text{-}a1}\right]  _{K}=\left[  \href{https://www.lmfdb.org/EllipticCurve/2.2.8.1/32.1/a/2}{32.1\text{-}a2}\right]  _{K}, &
\qquad & \left[  \href{https://www.lmfdb.org/EllipticCurve/2.2.8.1/32.1/a/3}{32.1\text{-}a3}\right]  _{K}=\left[   \href{https://www.lmfdb.org/EllipticCurve/2.2.8.1/32.1/a/4}{32.1\text{-}a4}\right]
_{K},\\
\left[  \href{https://www.lmfdb.org/EllipticCurve/2.2.8.1/32.1/a/5}{32.1\text{-}a5}\right]  _{K}=\left[  \href{https://www.lmfdb.org/EllipticCurve/2.2.8.1/32.1/a/6}{32.1\text{-}a6}\right]  _{K}, &  &
\left[  \href{https://www.lmfdb.org/EllipticCurve/2.2.8.1/32.1/a/7}{32.1\text{-}a7}\right]  _{K}=\left[  \href{https://www.lmfdb.org/EllipticCurve/2.2.8.1/32.1/a/8}{32.1\text{-}a8}\right]  _{K}.
\end{array}
\]
Consequently, the $2$-isogeny graph of $E/K$ is:
\begin{equation}\label{eq:CMgraphsec3}
\adjustbox{scale=0.8}{\begin{tikzcd}
                                                                                                          &                                                                                                                                    & {[\href{https://www.lmfdb.org/EllipticCurve/2.2.8.1/32.1/a/2}{32.1\text{-}a2}]_K} \\
{[\href{https://www.lmfdb.org/EllipticCurve/2.2.8.1/32.1/a/4}{32.1\text{-}a4}]_K} \arrow[r, "2", no head] & {[\href{https://www.lmfdb.org/EllipticCurve/2.2.8.1/32.1/a/6}{32.1\text{-}a6}]_K} \arrow[ru, "2", no head] \arrow[r, "2", no head] & {[\href{https://www.lmfdb.org/EllipticCurve/2.2.8.1/32.1/a/8}{32.1\text{-}a8}]_K}
\end{tikzcd}} 
\end{equation}
The determination of the full isogeny graph $\mathcal{G}(E/K)$ now falls under the theory of isogeny volcanoes \cite{Clarkvolcanoes,DrewVolcanoes}. Since $K$ contains the CM field, every elliptic curve in the isogeny class admits infinitely many $n$-isogenies. The discriminants of the corresponding endomorphism rings are given in the following table:
\[
\begin{tabular}
[c]{ccccc}
$E'$ & \href{https://www.lmfdb.org/EllipticCurve/2.2.8.1/32.1/a/4}{32.1\text{-}a4} & \href{https://www.lmfdb.org/EllipticCurve/2.2.8.1/32.1/a/6}{32.1\text{-}a6} & \href{https://www.lmfdb.org/EllipticCurve/2.2.8.1/32.1/a/2}{32.1\text{-}a2} & \href{https://www.lmfdb.org/EllipticCurve/2.2.8.1/32.1/a/8}{32.1\text{-}a8}\\\hline
$\Delta(\operatorname{End}E')$& $-4$ & $-16$ & $-64$ & $-64$
\end{tabular}
\]
Using the terminology of loc. cit., each of the $2$-isogenies occurring in~\eqref{eq:CMgraphsec3} is a vertical isogeny. Moreover, the discriminants $-4$ and $-16$ each occur only once in the isogeny class. Consequently, every $n$-isogeny from either \href{https://www.lmfdb.org/EllipticCurve/2.2.8.1/32.1/a/4}{32.1\text{-}a4} or \href{https://www.lmfdb.org/EllipticCurve/2.2.8.1/32.1/a/6}{32.1\text{-}a6} is necessarily an endomorphism. By contrast, the elliptic curves \href{https://www.lmfdb.org/EllipticCurve/2.2.8.1/32.1/a/2}{32.1\text{-}a2} and \href{https://www.lmfdb.org/EllipticCurve/2.2.8.1/32.1/a/8}{32.1\text{-}a8} have isomorphic endomorphism rings, so any isogeny between them is horizontal.

By \cite[Theorem 6.18]{BourdonClark2}, these two elliptic curves admit an $n$-isogeny over $K$ if and only if $-64$ is a square modulo $4n$. Some of these cyclic isogenies are endomorphisms, while the remaining ones are horizontal isogenies between the two elliptic curves. For example, among primes $p<100$, there is a horizontal $p$-isogeny precisely for $p\in\left\{  5,13,29,37,53,61\right\}  $. Since there are infinitely many such primes, there are infinitely many horizontal edges joining the vertices corresponding to \href{https://www.lmfdb.org/EllipticCurve/2.2.8.1/32.1/a/2}{32.1\text{-}a2} and \href{https://www.lmfdb.org/EllipticCurve/2.2.8.1/32.1/a/8}{32.1\text{-}a8}. Accordingly, we represent this infinite family of horizontal isogenies by a single dashed edge labeled by the smallest prime degree, namely $5$.
\begin{equation*}
\adjustbox{scale=0.8}{\begin{tikzcd}
                                                                                                          &                                                                                                                                    & {[\href{https://www.lmfdb.org/EllipticCurve/2.2.8.1/32.1/a/2}{32.1\text{-}a2}]_K} \arrow[d, "5", no head, dashed] \\
{[\href{https://www.lmfdb.org/EllipticCurve/2.2.8.1/32.1/a/4}{32.1\text{-}a4}]_K} \arrow[r, "2", no head] & {[\href{https://www.lmfdb.org/EllipticCurve/2.2.8.1/32.1/a/6}{32.1\text{-}a6}]_K} \arrow[ru, "2", no head] \arrow[r, "2", no head] & {[\href{https://www.lmfdb.org/EllipticCurve/2.2.8.1/32.1/a/8}{32.1\text{-}a8}]_K}                                
\end{tikzcd}}
\end{equation*}
Now let $E'$ be an elliptic curve in the isogeny class. If $E'$ is either \href{https://www.lmfdb.org/EllipticCurve/2.2.8.1/32.1/a/4}{32.1\text{-}a4} or \href{https://www.lmfdb.org/EllipticCurve/2.2.8.1/32.1/a/6}{32.1\text{-}a6}, then $\mathcal{G}_{5}(E^{\prime}/K)$ is trivial. On the other hand, if $E'$ is either \href{https://www.lmfdb.org/EllipticCurve/2.2.8.1/32.1/a/2}{32.1\text{-}a2} or \href{https://www.lmfdb.org/EllipticCurve/2.2.8.1/32.1/a/8}{32.1\text{-}a8}, then $\mathcal{G}_{5}(E^{\prime}/K)$ consists of two vertices joined by a horizontal $5$-isogeny. Consequently,  Corollary~\ref{cor:sameprimary} fails since $\operatorname*{End}_{K}E\not\cong\mathbb{Z}$. In particular, $\mathcal{G}(E/K)$ cannot be expressed as a Cartesian product of its $p$-primary graphs. 

Nevertheless, although $\mathcal{G}(E/K)$ has infinitely many edges, its isogeny class degree is $20=\operatorname{lcm}M_{E/K}$, since the isogeny matrix $M_{E/K}$ records only the least degree of a cyclic isogeny between any two vertices. Thus, the hypothesis $\operatorname*{End}_{K}E\cong \mathbb{Z}$ in Theorem~\ref{mainthmCarPro} is essential. For an additional example involving CM graphs, see Example~\ref{ex:CMvolc}.
\end{example}


\section{Finite \texorpdfstring{$p$}{p}-primary isogeny graphs} \label{sec:gps}\label{sec:graphs}

With Theorem \ref{mainthmCarPro} established, our classification is reduced to classifying the possible $p$-primary isogeny graphs that could occur over fields of characteristic $0$ for each prime $p$. In this section, we consider the case when the $p$-primary isogeny graphs are finite. We conclude by proving the main result of this section, Theorem~\ref{classificationGpk}. In particular, this result implies Theorem~\ref{mathmclass} (1), and the converse portion in the case when $\mathcal{H}_{p^{k}}^{r}$ is a finite graph. In this direction, we begin by introducing the family of graphs $\mathcal{H}_{p^{k}}^{r}$.

\subsection{The family of finite graphs \texorpdfstring{$\mathcal{H}_{p^{k}}^r$}{Hrpk}}\label{subsec:finitegraphs}
Before introducing the family of graphs $\mathcal{H}_{p^{k}}^{r}$, we first establish some necessary graph-theoretic terminology. Recall that by a $k$-spine, with $k$ a positive integer, we mean a path graph with vertex set $\left\{  v_{0},v_{1},\ldots,v_{k}\right\}  $ and edges $\left\{  v_{i},v_{i+1}\right\}  $ for $0\leq i\leq k-1$. 

Our next definition provides the local growth operation that will be used in the construction of the graphs $\mathcal{H}_{p^k}^r$.

\begin{definition}
Let $v$ be a vertex of a graph $\mathcal{G}$, and let $p$ be a prime. We say that $v$ is \textit{$p$-bloomed} if $\deg v = p+1$. If $\deg v<p+1$, then a \textit{$p$-blossoming at $v$} is the operation of adjoining exactly $p+1-\deg v$ new vertices, each adjacent only to $v$. In the resulting graph, $v$ is $p$-bloomed. A \textit{$p$-blossoming of $\mathcal G$} is the operation of $p$-blossoming at every vertex of $\mathcal G$ simultaneously. The resulting graph is called the $p$-blossom of $\mathcal G$.
\end{definition}

We now define the graphs $\mathcal{H}_{p^k}^r$ recursively.

\begin{definition}
    Let $p$ be a prime number and $k$ a positive integer, and $r$ an integer satisfying $r \leq \frac{k}{2}$. Let $\mathcal P_0$ be the path on $k+1-2r$ vertices and $k-2r$ edges. For $1 \leq i \leq r$, let $\mathcal{P}_i$ be the $p$-blossom of $\mathcal{P}_{i-1}$. We define $\mathcal H_{p^k}^r$ to be the graph $\mathcal P_r$, with the additional property that each edge is assigned a weight of $p$. 
    
    The path $\mathcal{P}_0$ is called the \textit{skeleton} of $\mathcal H_{p^k}^r$, and the $k+1-2r$ vertices of $\mathcal{P}_0$ are called the \textit{skeletal vertices} of $\mathcal H_{p^k}^r$.
\end{definition}

\begin{remark}
While the graph $\mathcal{H}_{p^k}^r$ is a weighted graph, with each edge assigned weight $p$, all distances and path lengths are measured in the underlying unweighted tree. The edge weights play no role in the combinatorial arguments below.
\end{remark}

\begin{example}\label{ex:blossoming}
We now illustrate the construction of $\mathcal{H}_{3^{6}}^{2}$. Since $k=6$ and $r=2$, we begin with a path $\mathcal{P}_{0}$ on $3$ vertices. To obtain $\mathcal{P}_{1}$, we perform a $3$-blossoming of each vertex of $\mathcal{P}_{0}$. The graphs $\mathcal{P}_{0}$ and $\mathcal{P}_{1}$ are given below.

\[
\adjustbox{scale=0.85}{\begin{tikzpicture}[scale=1, transform shape]

    \Vertex[shape=circle, size=0.3, color=white, style={draw=blue}, x=1.5]{1}
    \Vertex[shape=circle, size=0.3, color=white, style={draw=blue}, x=3]{4}
    \Vertex[shape=circle, size=0.3, color=white, style={draw=blue}, x=4.5]{6}

    \Edge[style={color=blue}](1)(4)
    \Edge[style={color=blue}](4)(6)

        \node at (3.1, -1.5) { $\mathcal P_0$};

\end{tikzpicture}
\qquad \qquad
\begin{tikzpicture}[scale=1, transform shape]

    \Vertex[shape=circle, size=0.3, color=white, style={draw=red}]{3}

    \Vertex[shape=circle, size=0.3, color=white, style={draw=blue}, x=1.5]{1}
    \Vertex[shape=circle, size=0.3, color=white, style={draw=blue}, x=3]{4}
    \Vertex[shape=circle, size=0.3, color=white, style={draw=blue}, x=4.5]{6}
    \Vertex[shape=circle, size=0.3, color=white, style={draw=red}, x=6]{8}
    \Vertex[shape=circle, size=0.3,  color=white, style={draw=red}, x=1.5, y=1]{2}
    
    \Vertex[shape=circle, size=0.3,  color=white, style={draw=red}, x=3, y=1]{5}
    \Vertex[shape=circle, size=0.3,  color=white, style={draw=red}, x=4.5, y=1]{7}

    \Vertex[shape=circle, size=0.3,  color=white, style={draw=red}, x=3, y=-1]{15}

    \Vertex[shape=circle, size=0.3,  color=white, style={draw=red}, x=1.5, y=-1]{12}

        \Vertex[shape=circle, size=0.3,  color=white, style={draw=red}, x=4.5, y=-1]{17}

    \Edge[style={color=red}](3)(1)
    \Edge[style={color=blue}](1)(4)
    \Edge[style={color=blue}](4)(6)
    \Edge[style={color=red}](6)(8)

    \Edge[style={color=red}](1)(2)
    \Edge[style={color=red}](4)(5)
    \Edge[style={color=red}](6)(7)

    \Edge[style={color=red}](2)(12)
    \Edge[style={color=red}](5)(15)
    \Edge[style={color=red}](7)(17)
  
        \node at (3.1, -1.5) { $\mathcal P_1$};

\end{tikzpicture}}
\]
Each of the new vertices in $\mathcal{P}_{1}$ is a leaf, whereas each skeletal vertex is $3$-bloomed. To obtain~$\mathcal{P}_{2}$, we then apply a $3$-blossoming at each leaf of $\mathcal{P}_{1}$, which gives the graph $\mathcal{H}_{3^{6}}^{2}$. In particular, each leaf of $\mathcal{P}_{1}$ is now $3$-bloomed, and the new vertices are leaves in $\mathcal{H}_{3^{6}}^{2}$.
\[
\adjustbox{scale=0.85}{\begin{tikzpicture}[scale=1, transform shape]

    \Vertex[shape=circle, size=0.3, color=white, style={draw=red}]{3}
    \Vertex[shape=circle, size=0.3,  color=white, x=-1.5]{3a}
    \Vertex[shape=circle, size=0.3,  color=white, y=1]{3b}
   
    \Vertex[shape=circle, size=0.3, color=white, style={draw=blue}, x=1.5]{1}
    \Vertex[shape=circle, size=0.3, color=white, style={draw=blue}, x=3]{4}
    \Vertex[shape=circle, size=0.3, color=white, style={draw=blue}, x=4.5]{6}
    \Vertex[shape=circle, size=0.3, color=white, style={draw=red}, x=6]{8}
    \Vertex[shape=circle, size=0.3,  color=white, style={draw=red}, x=1.5, y=1]{2}
    \Vertex[shape=circle, size=0.3,  color=white, x=1.2, y=2]{2a}
    \Vertex[shape=circle, size=0.3,  color=white, x=1.9, y=2]{2b}
    \Vertex[shape=circle, size=0.3,  color=white, x=1.55, y=2]{2c}
    
    \Vertex[shape=circle, size=0.3,  color=white, style={draw=red}, x=3, y=1]{5}
    \Vertex[shape=circle, size=0.3,  color=white, style={draw=red}, x=4.5, y=1]{7}

    \Vertex[shape=circle, size=0.3,  color=white, x=4.15, y=2]{7a}
    \Vertex[shape=circle, size=0.3,  color=white, x=4.5, y=2]{7b}
    \Vertex[shape=circle, size=0.3,  color=white, x=4.85, y=2]{7c}

    \Vertex[shape=circle, size=0.3,  color=white, x=7.5]{8a}
    \Vertex[shape=circle, size=0.3,  color=white, x=6, y=1]{8b}

    \Vertex[shape=circle, size=0.3,  color=white, x=2.6, y=2]{5a}
    \Vertex[shape=circle, size=0.3,  color=white, x=3.3, y=2]{5b}
    \Vertex[shape=circle, size=0.3,  color=white, x=2.95, y=2]{5c}

    \Vertex[shape=circle, size=0.3,  color=white, style={draw=red}, x=3, y=-1]{15}
    \Vertex[shape=circle, size=0.3,  color=white, x=2.6, y=-2]{15a}
    \Vertex[shape=circle, size=0.3,  color=white, x=3.3, y=-2]{15b}
    \Vertex[shape=circle, size=0.3,  color=white, x=2.95, y=-2]{15c}
    
    \Vertex[shape=circle, size=0.3,  color=white, x=4.5, y=-1]{17}

    \Vertex[shape=circle, size=0.3,  color=white, style={draw=red}, x=1.5, y=-1]{12}
    \Vertex[shape=circle, size=0.3,  color=white, x=1.2, y=-2]{12a}
    \Vertex[shape=circle, size=0.3,  color=white, x=1.9, y=-2]{12b}
    \Vertex[shape=circle, size=0.3,  color=white, x=1.55, y=-2]{12c}

    \Vertex[shape=circle, size=0.3,  color=white, y=-1]{13b}

        \Vertex[shape=circle, size=0.3,  color=white, style={draw=red}, x=4.5, y=-1]{17}

    \Vertex[shape=circle, size=0.3,  color=white, x=4.15, y=-2]{17a}
    \Vertex[shape=circle, size=0.3,  color=white, x=4.5, y=-2]{17b}
    \Vertex[shape=circle, size=0.3,  color=white, x=4.85, y=-2]{17c}

    \Vertex[shape=circle, size=0.3,  color=white, x=6, y=-1]{18b}

    \Edge[style={color=red}](3)(1)
    \Edge[style={color=blue}](1)(4)
    \Edge[style={color=blue}](4)(6)
    \Edge[style={color=red}](6)(8)

    \Edge[style={color=red}](1)(2)
    \Edge[style={color=red}](4)(5)
    \Edge[style={color=red}](6)(7)

    \Edge(8)(8a)
    \Edge(8)(8b)

    \Edge(5)(5a)
    \Edge(5)(5b)
    \Edge(5)(5c)

    \Edge(7)(7a)
    \Edge(7)(7b)
    \Edge(7)(7c)
    
    \Edge(3)(3a)
    \Edge(3)(3b)

    \Edge(2)(2a)
    \Edge(2)(2b)
    \Edge(2)(2c)

    \Edge[style={color=red}](2)(12)
    \Edge[style={color=red}](5)(15)
    \Edge[style={color=red}](7)(17)

    \Edge(12)(12a)
    \Edge(12)(12b)
    \Edge(12)(12c)

    \Edge(15)(15a)
    \Edge(15)(15b)
    \Edge(15)(15c)

    \Edge(17)(17a)
    \Edge(17)(17b)
    \Edge(17)(17c)

    \Edge(8)(18b)
    \Edge(3)(13b)
  
        \node at (3.1, -2.8) {\Large $\mathcal H^2_{3^6}$};

\end{tikzpicture}}
\]

\end{example}

\begin{lemma}\label{Lem:longpath}\label{Lem:skel}
The longest possible path in the graph $\mathcal{H}_{p^{k}}^{r}$ has $k$ edges, and any path with $k$ edges contains the skeleton.
\end{lemma}

\begin{proof}
Let $\mathcal{P}_{0},\mathcal{P}_{1},\ldots,\mathcal{P}_{r}$ be the sequence of graphs used in the construction of $\mathcal{H}_{p^{k}}^{r}$. For $0\leq i\leq r$, let $\ell_{i}$ denote the length of the longest path in $\mathcal{P}_{i}$. We shall prove that $\ell_i=k-2r+2i$ and that any path with $\ell_i$ edges in $\mathcal P_i$ contains $\mathcal P_0$ for $0\leq i \leq r$. For $i=0$, the claim is immediate from the definition.

Suppose that the claim holds for some $i$ satisfying $0\leq i <r$. Consider any path of length $\ell_i$ in $\mathcal P_i$. Its endpoints are leaves $v$ and $w$ of $\mathcal P_i$, since otherwise it could be extended to a longer path in $\mathcal P_i$. Extending this path in $\mathcal P_{i+1}$ by adding leaves attached to $v$ and $w$ gives $\ell_{i+1}\geq \ell_i+2$. On the other hand, the endpoints of any path of maximal length in $\mathcal P_{i+1}$ are leaves of $\mathcal P_{i+1}$, and removing the two endpoints of such a path yields a path in $\mathcal P_i$ which has length at most $\ell_i$, so $\ell_{i+1}\leq \ell_i+2$. It follows that $\ell_{i+1}=\ell_i+2=k-2r+2i+2$. Moreover, any path of length $\ell_{i+1}$ in $\mathcal P_{i+1}$ contains a path of length $\ell_i$ in $\mathcal P_i$. By assumption, all paths of length $\ell_i$ in $\mathcal P_i$ contain the skeleton, so it follows that all paths of length $\ell_{i+1}$ in $\mathcal P_{i+1}$ contain the skeleton. This completes the induction step.
\end{proof}

With notation as in Example \ref{ex:blossoming}, we observe that $\mathcal{P}_{0}$ and $\mathcal{P}_{1}$ can also be viewed as special cases of our family of graphs. More precisely, $\mathcal{P}_{0}$ is the graph $\mathcal{H}_{3^{2}}^{0}$, since it is a path on $k-2r=2$ edges. Similarly, $\mathcal{P}_{1}$ is the graph $\mathcal{H}_{3^{4}}^{1}$. Indeed, in this case we have $\left(  k^{\prime},r^{\prime}\right)  =\left(  4,1\right)  $, and so let $\mathcal{P}_{0}^{\prime}$ be the path with $k^{\prime}-2r^{\prime}=2$ edges. In particular, $\mathcal{P}_{0}^{\prime}\cong\mathcal{P}_{0}$. The $3$-blossoming of $\mathcal{P}_{0}$ then produces $\mathcal{P}_{1}$, which is by definition $\mathcal{H}_{3^{4}}^{1}$, upon assing weights of $p$ to its edges. More generally, if $\mathcal{P}_{0},\mathcal{P}_{1},\ldots ,\mathcal{P}_{r}$ is the sequence of graphs used to construct $\mathcal{H}_{p^{k}}^{r}$, then each $\mathcal{P}_{i}$ can be viewed as some $\mathcal{H}_{p^{k^{\prime}}}^{r^{\prime}}$. The following lemma makes this precise.

\begin{lemma}
\label{Lem:specializationsgraphs}Let $\mathcal{P}_{0},\mathcal{P}_{1},\ldots,\mathcal{P}_{r}$ be the sequence of graphs used to construct $\mathcal{H}_{p^{k}}^{r}$. Then, if each edge of $\mathcal{P}_{i}$ is assigned
weight $p$, then
\[
\mathcal{P}_{i}\cong\mathcal{H}_{p^{k-2(r-i)}}^{i}.
\]

\end{lemma}

\begin{proof}
By definition, $\mathcal{P}_{0}$ is a path with $k-2r$ edges, and $\mathcal{P}_{i}$ is obtained from $\mathcal{P}_{0}$ by applying~$i$ successive $p$-blossomings. Now let $\mathcal{P}_{0}^{\prime},\mathcal{P}_{1}^{\prime},\ldots,\mathcal{P}_{i}^{^{\prime}}$ be the sequence of graphs used to construct $\mathcal{H}_{p^{k-2(r-i)}}^{i}$. Then $\mathcal{P}_{0}^{\prime}$ is a path with
\[
k-2\left(  r-i\right)  -2i=k-2r
\]
edges. In particular, $\mathcal{P}_{0}\cong\mathcal{P}_{0}^{\prime}$. Moreover, $\mathcal{P}_{i}^{\prime}$ is obtained from $\mathcal{P}_{0}^{\prime}$ by applying $i$ successive $p$-blossomings, and so $\mathcal{P}_{i}^{\prime}\cong\mathcal{P}_{i}$. The result now follows since $\mathcal{H}_{p^{k-2(r-i)}}^{i}$ is obtained from $\mathcal{P}_{i}^{\prime}$ upon assigning each edge a weight of $p$.
\end{proof}

We next define the notions of bloom depth and bloom, which will aid us in the investigation of the graphs $\mathcal{H}_{p^k}^r$ and $\mathcal{H}_{p^\infty , +}^r$. Throughout this discussion, a tree graph $\mathcal{T}$ may be finite or infinite unless explicitly stated otherwise. We also recall that in a tree, every pair of distinct vertices is connected by exactly one path, and this holds even if the tree is infinite.

\begin{definition}
Let $\mathcal T$ be a tree, and $v$ be any vertex of $\mathcal T$. The \textit{$p$-bloom depth} of $v$ is the largest nonnegative integer $j$ such that every vertex which is at a distance less than $j$ from $v$ is $p$-bloomed. If such an integer does not exist, we say $v$ has $p$-bloom depth $\infty$.
\end{definition} 

When $p$ is clear from context, we will omit it from the notation and only write bloom depth. While bloom depth is defined for any tree, we will only use it for the graphs $\mathcal H_{p^k}^r$ and their infinite analogues, which will be discussed in the next section. In the case of $\mathcal H_{p^k}^r$, the bloom depth has a particularly simple description. 

\begin{lemma}\label{lem:bloomdepthdescription}
    If $r=0$, every vertex of $\mathcal H_{p^k}^r$ has bloom depth $0$. If $r>0$, the bloom depth of a vertex $v$ in $\mathcal H_{p^k}^r$ equals the minimum distance from $v$ to a leaf of $\mathcal H_{p^k}^r$.
\end{lemma}
\begin{proof}
    If $r=0$, none of the vertices are $p$-bloomed, and the claim follows. Suppose $r>0$, and let $v$ be any vertex of $\mathcal H_{p^k}^r$, with $j$ the distance from $v$ to the nearest leaf. Since leaves are not $p$-bloomed, the bloom depth of $v$ is at most $j$. On the other hand, none of the vertices at a distance less than $j$ from $v$ are leaves. Since any vertex which is not a leaf in $\mathcal H_{p^k}^r$ is $p$-bloomed, it follows that the bloom depth of $v$ is at least $j$, and therefore equals $j$.
\end{proof}
\begin{remark}
    Suppose that the $p$-primary graph $\mathcal G_p(E/K)$ is isomorphic to some $\mathcal H_{p^k}^r$. Then the bloom depth of an elliptic curve $E' \in V(\mathcal G_p(E/K))$ is the largest nonnegative integer $j$ such that all $p^j$-isogenies of $E'$ are defined over $K$.
\end{remark}

\begin{definition}
Let $\mathcal{T}$ be a tree containing a fixed distinguished subgraph $\mathcal{S}$, where $\mathcal{S}$ is either a $k$-spine, a line graph, or a ray graph. For a vertex $x\in V(\mathcal{T})$, we define the \textit{distance from }$x$\textit{ to }$\mathcal{S}$, denoted $\operatorname*{dist}_{\mathcal{S}}(x)$, to be the minimum distance from $x$ to a vertex of $\mathcal{S}$. In particular, if $v\in V(\mathcal{S})$, then $\operatorname*{dist}_{\mathcal{S}}(v)=0$.

Since $\mathcal{T}$ is a tree, each vertex $x\in V(\mathcal{T})$ has a unique nearest vertex $v\in V(\mathcal{S})$. We define the map $\pi:\mathcal{T}\rightarrow \mathcal{S}$ by $\pi(x)=v$, and refer to $\pi$ as the \textit{projection of $\mathcal{T}$ onto $\mathcal{S}$}.
\end{definition}

\begin{example}
Let $\mathcal{T}$ be a tree containing a $k$-spine $\mathcal{S}$. If $x\not \in V(\mathcal{S})$ is a vertex adjacent to $v\in V(\mathcal{S})$, then $\operatorname*{dist}_{\mathcal{S}}(x)=1$ and $\pi(x)=v$. More generally, if $x$ lies a distance $d$ from~$\mathcal{S}$, then $\operatorname*{dist}_{\mathcal{S}}(x)=d$ and $\pi(x)$ is the unique vertex on $\mathcal{S}$ such that a path from $x$ to $\mathcal{S}$ first meets $\mathcal{S}$ at $\pi(x)$.
\end{example}

\begin{definition}
Let $\mathcal{T}$ be a tree containing a fixed distinguished subgraph $\mathcal{S}$, where $\mathcal{S}$ is either a $k$-spine, a line graph, or a ray graph, and let $\pi:\mathcal{T}\rightarrow\mathcal{S}$ be the corresponding projection. For a vertex $v\in V(\mathcal{S})$, we define the \textit{bloom at $v$ relative to $\mathcal{S}$}, denoted $\mathcal{B}(v)$, to be the induced subgraph of $\mathcal{T}$ on the vertex set $\pi^{-1}(v)$; that is,
\[
\mathcal{B}(v)=\mathcal{T}[\pi^{-1}(v)].
\]
\end{definition}

\begin{remark}
When the chosen distinguished subgraph $\mathcal{S}$ is clear from context, we omit the phrase ``relative to~$\mathcal{S}$.''
\end{remark}

\begin{example}\label{ex:booms}
Let $\mathcal{T}$ be the tree shown below, and let $\mathcal{S}$ be the $3$-spine with vertex set $\left\{  v_{0},v_{1},v_{2},v_{3}\right\}  $ and edges $\left\{  v_{i},v_{i+1}\right\}  $ for $0\leq i\leq3$. Let $\pi:\mathcal{T}\rightarrow\mathcal{S}$ be the corresponding projection. 
\[
\begin{tikzpicture}[scale=1, transform shape]

    \Vertex[shape=circle, size=0.4, color=white, style={draw=red}, x=0, label=$v_0$]{0}

    \Vertex[shape=circle, size=0.4, color=white, x=2, label = $v_1$]{1}
    \Vertex[shape=circle, size=0.4, color=white, x=4, label = $v_2$]{2}
    \Vertex[shape=circle, size=0.4, color=white, style={draw=blue}, x=6, label = $v_3$]{3}

    \Vertex[shape=circle, size=0.4, color=white, style={draw=red}, x=-0.25, y=1]{011}

    \Vertex[shape=circle, size=0.4, color=white, style={draw=red}, x=0.25, y=1]{012}

    \Vertex[shape=circle, size=0.4, color=white, style={draw=red}, x=-0.75, y=2]{021}

    \Vertex[shape=circle, size=0.4, color=white, style={draw=red}, x=-0.25, y=2]{022}

    \Vertex[shape=circle, size=0.4, color=white, style={draw=red}, x=0.25, y=2]{023}

    \Vertex[shape=circle, size=0.4, color=white, style={draw=red}, x=0.75, y=2]{024}

    \Vertex[shape=circle, size=0.4, color=white, style={draw=red}, x=-1.75, y=3]{031}

    \Vertex[shape=circle, size=0.4, color=white, style={draw=red}, x=-1.25, y=3]{032}    

    \Vertex[shape=circle, size=0.4, color=white, style={draw=red}, x=-.75, y=3]{033}

    \Vertex[shape=circle, size=0.4, color=white, style={draw=red}, x=-.25, y=3]{034}   

    \Vertex[shape=circle, size=0.4, color=white, style={draw=red}, x=.25, y=3]{035}

    \Vertex[shape=circle, size=0.4, color=white, style={draw=red}, x=.75, y=3]{036}   

    \Vertex[shape=circle, size=0.4, color=white, style={draw=red}, x=1.25, y=3]{037}

    \Vertex[shape=circle, size=0.4, color=white, style={draw=red}, x=1.75, y=3]{038} 

    \Vertex[shape=circle, size=0.4, color=white, style={draw=blue}, x=5, y=1]{311}

    \Vertex[shape=circle, size=0.4, color=white, style={draw=blue}, x=6, y=1]{312}

    \Vertex[shape=circle, size=0.4, color=white, style={draw=blue}, x=7, y=1]{313}

    \Vertex[shape=circle, size=0.4, color=white, style={draw=blue}, x=4, y=2]{321}

    \Vertex[shape=circle, size=0.4, color=white, style={draw=blue}, x=4.5, y=2]{322}

    \Vertex[shape=circle, size=0.4, color=white, style={draw=blue}, x=5, y=2]{323}

    \Vertex[shape=circle, size=0.4, color=white, style={draw=blue}, x=5.5, y=2]{324}

    \Vertex[shape=circle, size=0.4, color=white, style={draw=blue}, x=6, y=2]{325}

    \Vertex[shape=circle, size=0.4, color=white, style={draw=blue}, x=6.5, y=2]{326}

    \Vertex[shape=circle, size=0.4, color=white, style={draw=blue}, x=7, y=2]{327}

    \Vertex[shape=circle, size=0.4, color=white, style={draw=blue}, x=7.5, y=2]{328}

    \Vertex[shape=circle, size=0.4, color=white, style={draw=blue}, x=8, y=2]{329}

    \Vertex[shape=circle, size=0.4, color=white, x=1.5, y=1]{111}

    \Vertex[shape=circle, size=0.4, color=white, x=2, y=1]{112}

    \Vertex[shape=circle, size=0.4, color=white, x=2.5, y=1]{113}

    \Vertex[shape=circle, size=0.4, color=white, x=2, y=2]{121}

    \Vertex[shape=circle, size=0.4, color=white, x=2.5, y=2]{122}

    \Vertex[shape=circle, size=0.4, color=white, x=3, y=2]{123}

    \Edge(0)(1)
    \Edge(1)(2)
    \Edge(2)(3)
    \Edge[style={color=red}](0)(011)
    \Edge[style={color=red}](0)(012)
    \Edge[style={color=red}](011)(021)
    \Edge[style={color=red}](011)(022)
    \Edge[style={color=red}](012)(023)
    \Edge[style={color=red}](012)(024)
    \Edge[style={color=red}](021)(031)
    \Edge[style={color=red}](021)(032)
    \Edge[style={color=red}](022)(033)
    \Edge[style={color=red}](022)(034)
    \Edge[style={color=red}](023)(035)
    \Edge[style={color=red}](023)(036)
    \Edge[style={color=red}](024)(037)
    \Edge[style={color=red}](024)(038)

    \Edge[style={color=blue}](3)(311)
    \Edge[style={color=blue}](3)(312)
    \Edge[style={color=blue}](3)(313)

    \Edge[style={color=blue}](311)(321)
    \Edge[style={color=blue}](311)(322)
    \Edge[style={color=blue}](311)(323)
    \Edge[style={color=blue}](312)(324)
    \Edge[style={color=blue}](312)(325)
    \Edge[style={color=blue}](312)(326)
    \Edge[style={color=blue}](313)(327)
    \Edge[style={color=blue}](313)(328)
    \Edge[style={color=blue}](313)(329)

    \Edge(1)(111)
    \Edge(1)(112)
    \Edge(1)(113)
    \Edge(112)(121)
    \Edge(112)(122)
    \Edge(113)(123)

    \node at (3.5, 3) { $\mathcal{T}$};

\end{tikzpicture}
\]
The illustrations below depict the blooms $\mathcal{B}(v_{i})$ for $0\leq i\leq3$. Observe that each bloom is organized into layers according to the distance from the $3$-spine $\mathcal{S}$.

Now observe that for each prime $p$, the $p$-bloom depth of $v_i$ in $\mathcal{T}$ is $0$. However, if we restrict to $\mathcal{B}(v_0)$, then $v_0$ has $2$-bloom depth $3$. Similarly, $v_3$ in $\mathcal{B}(v_3)$ has $3$-bloom depth $2$. However, the vertices $v_1$ and $v_2$ continue to have $p$-bloom depth $0$ for each prime $p$ in their respective blooms.
\[
\begin{tikzpicture}[scale=1, transform shape]

    \Vertex[shape=circle, size=0.4, color=white, style={draw=red}, x=0, label=$v_0$]{0}

    \Vertex[shape=circle, size=0.4, color=white, x=3, label = $v_1$]{1}
    \Vertex[shape=circle, size=0.4, color=white, x=5.5, label = $v_2$]{2}
    \Vertex[shape=circle, size=0.4, color=white, style={draw=blue}, x=8, label = $v_3$]{3}

    \Vertex[shape=circle, size=0.4, color=white, style={draw=red}, x=-0.25, y=1]{011}

    \Vertex[shape=circle, size=0.4, color=white, style={draw=red}, x=0.25, y=1]{012}

    \Vertex[shape=circle, size=0.4, color=white, style={draw=red}, x=-0.75, y=2]{021}

    \Vertex[shape=circle, size=0.4, color=white, style={draw=red}, x=-0.25, y=2]{022}

    \Vertex[shape=circle, size=0.4, color=white, style={draw=red}, x=0.25, y=2]{023}

    \Vertex[shape=circle, size=0.4, color=white, style={draw=red}, x=0.75, y=2]{024}

    \Vertex[shape=circle, size=0.4, color=white, style={draw=red}, x=-1.75, y=3]{031}

    \Vertex[shape=circle, size=0.4, color=white, style={draw=red}, x=-1.25, y=3]{032}    

    \Vertex[shape=circle, size=0.4, color=white, style={draw=red}, x=-.75, y=3]{033}

    \Vertex[shape=circle, size=0.4, color=white, style={draw=red}, x=-.25, y=3]{034}   

    \Vertex[shape=circle, size=0.4, color=white, style={draw=red}, x=.25, y=3]{035}

    \Vertex[shape=circle, size=0.4, color=white, style={draw=red}, x=.75, y=3]{036}   

    \Vertex[shape=circle, size=0.4, color=white, style={draw=red}, x=1.25, y=3]{037}

    \Vertex[shape=circle, size=0.4, color=white, style={draw=red}, x=1.75, y=3]{038} 

    \Vertex[shape=circle, size=0.4, color=white, style={draw=blue}, x=7, y=1]{311}

    \Vertex[shape=circle, size=0.4, color=white, style={draw=blue}, x=8, y=1]{312}

    \Vertex[shape=circle, size=0.4, color=white, style={draw=blue}, x=9, y=1]{313}

    \Vertex[shape=circle, size=0.4, color=white, style={draw=blue}, x=6, y=2]{321}

    \Vertex[shape=circle, size=0.4, color=white, style={draw=blue}, x=6.5, y=2]{322}

    \Vertex[shape=circle, size=0.4, color=white, style={draw=blue}, x=7, y=2]{323}

    \Vertex[shape=circle, size=0.4, color=white, style={draw=blue}, x=7.5, y=2]{324}

    \Vertex[shape=circle, size=0.4, color=white, style={draw=blue}, x=8, y=2]{325}

    \Vertex[shape=circle, size=0.4, color=white, style={draw=blue}, x=8.5, y=2]{326}

    \Vertex[shape=circle, size=0.4, color=white, style={draw=blue}, x=9, y=2]{327}

    \Vertex[shape=circle, size=0.4, color=white, style={draw=blue}, x=9.5, y=2]{328}

    \Vertex[shape=circle, size=0.4, color=white, style={draw=blue}, x=10, y=2]{329}

    \Vertex[shape=circle, size=0.4, color=white, x=2.5, y=1]{111}

    \Vertex[shape=circle, size=0.4, color=white, x=3, y=1]{112}

    \Vertex[shape=circle, size=0.4, color=white, x=3.5, y=1]{113}

    \Vertex[shape=circle, size=0.4, color=white, x=3, y=2]{121}

    \Vertex[shape=circle, size=0.4, color=white, x=3.5, y=2]{122}

    \Vertex[shape=circle, size=0.4, color=white, x=4, y=2]{123}

    \Edge[style={color=red}](0)(011)
    \Edge[style={color=red}](0)(012)
    \Edge[style={color=red}](011)(021)
    \Edge[style={color=red}](011)(022)
    \Edge[style={color=red}](012)(023)
    \Edge[style={color=red}](012)(024)
    \Edge[style={color=red}](021)(031)
    \Edge[style={color=red}](021)(032)
    \Edge[style={color=red}](022)(033)
    \Edge[style={color=red}](022)(034)
    \Edge[style={color=red}](023)(035)
    \Edge[style={color=red}](023)(036)
    \Edge[style={color=red}](024)(037)
    \Edge[style={color=red}](024)(038)

    \Edge[style={color=blue}](3)(311)
    \Edge[style={color=blue}](3)(312)
    \Edge[style={color=blue}](3)(313)

    \Edge[style={color=blue}](311)(321)
    \Edge[style={color=blue}](311)(322)
    \Edge[style={color=blue}](311)(323)
    \Edge[style={color=blue}](312)(324)
    \Edge[style={color=blue}](312)(325)
    \Edge[style={color=blue}](312)(326)
    \Edge[style={color=blue}](313)(327)
    \Edge[style={color=blue}](313)(328)
    \Edge[style={color=blue}](313)(329)

    \Edge(1)(111)
    \Edge(1)(112)
    \Edge(1)(113)
    \Edge(112)(121)
    \Edge(112)(122)
    \Edge(113)(123)

\node at (0, -0.75) { $\mathcal{B}(v_0)$};
\node at (3, -0.75) { $\mathcal{B}(v_1)$};
\node at (5.5, -0.75) { $\mathcal{B}(v_2)$};
\node at (8, -0.75) { $\mathcal{B}(v_3)$};

\end{tikzpicture}
\]
\end{example}

We next record the following statement, which is a direct consequence of the definitions.

\begin{corollary}
\label{cor:isographsblooms}Let $\mathcal{T}$ and $\mathcal{T}^{\prime}$ be trees, and let $\mathcal{S} \subseteq \mathcal{T}$ and $\mathcal{S}^{\prime} \subseteq \mathcal{T}^{\prime}$ be distinguished subgraphs of the same type: either both are $k$-spines, or both are line graphs, or both are ray graphs. Write
\[
V(\mathcal{S})=\left\{  v_{i}\right\}_{i\in I}  \qquad \text{and} \qquad V(\mathcal{S}^{\prime})=\{  v_{i}^{\prime}\}_{i \in I}  ,
\]
where $I\subseteq \mathbb{Z}$ is a common indexing set such that the edges of $\mathcal{S}$ and $\mathcal{S}^{\prime}$ are precisely
\[
\{  v_{i},v_{i+1}\}  \qquad \text{and} \qquad \{  v_{i}^{\prime},v_{i+1}^{\prime}\}  ,
\]
respectively, for all indices $i \in I$ for which $i+1 \in I$. Let $\mathcal{B}(v_{i})$ (resp. $\mathcal{B}(v_{i}^{\prime})$) denote the bloom of $v_{i}$ (resp. $v_{i}^{\prime}$) relative to $\mathcal{S}$ (resp. $\mathcal{S}^{\prime}$). If $\mathcal{B}(v_{i})\cong\mathcal{B}(v_{i}^{\prime})$ for each $i\in I$, then $\mathcal{T}\cong\mathcal{T}^{\prime}$.
\end{corollary}

We next describe the projection $\pi$ onto a $k$-spine of $\mathcal H_{p^k}^r$ in more detail.

\begin{proposition}\label{lem:bloomdepth}\label{prop:uniquetree}
Let $\mathcal{S}$ be a $k$-spine of $\mathcal{H}_{p^{k}}^{r}$ with vertices $\left\{  v_{0},v_{1},\ldots,v_{k}\right\}  $ and edges $\left\{  v_{i},v_{i+1}\right\}$. Let $d_{i}=\min\!\left\{  i,r,k-i\right\}$. Then: 
\begin{enumerate}[label=\textnormal{(\arabic*)}]
\item $\pi^{-1}(v_0)=\left\{  v_0\right\}$ and $\pi^{-1}(v_k)=\{v_k\}$.

\item Let $1\leq i \leq k-1$. Each vertex $x\in\pi^{-1}(v_i)$ with $\operatorname*{dist}_{\mathcal{S}}(x)<d_i$ is $p$-bloomed, while every vertex $x\in\pi^{-1}(v_i)$ with $\operatorname*{dist}_{\mathcal{S}}(x)=d_i$ is a leaf.
\end{enumerate}
Conversely, if $\mathcal T$ is a tree containing a $k$-spine $\mathcal S$ with vertices  $\left\{  v_{0},v_{1},\ldots,v_{k}\right\}  $ and edges $\left\{  v_{i},v_{i+1}\right\}$ such that the projection $\pi$ onto $\mathcal S$ satisfies properties \textnormal{(1)} and \textnormal{(2)}, then $\mathcal T\cong \mathcal H_{p^k}^r$.
\end{proposition}

\begin{proof}
By Lemma \ref{Lem:skel}, $\mathcal{S}$ contains the skeleton of $\mathcal{H}_{p^{k}}^{r}$, which consists of $k+1-2r$ vertices. Since $\mathcal{S}$ is a longest path by Lemma \ref{Lem:longpath}, its endpoints $v_{0}$ and $v_{k}$ are leaves. By construction of $\mathcal{H}_{p^{k}}^{r}$, $v_{0}$ and $v_{k}$ are a distance $r$ from the skeleton, and thus, after a possible relabeling, we may take $v_{r},v_{r+1},\ldots,v_{k-r}$ to be the skeletal vertices. Now let $\mathcal{P}_{0},\mathcal{P}_{1},\ldots,\mathcal{P}_{r}$ be the sequence of graphs used in the construction of $\mathcal{H}_{p^{k}}^{r}$. In particular, $\mathcal{P}_{i}$ is the $p$-blossoming of $\mathcal{P}_{i-1}$ for $0\leq i<r$. Now consider the bloom $\mathcal{B}(v_{i})$ with $r\leq i\leq k-r$. The vertices in $\mathcal{B}(v_{i})$ at distance $j$ from $v_{i}$ are leaves in $\mathcal{P}_{j}$. But for $j<r$, these vertices are $p$-bloomed in $\mathcal{P}_{j+1}$. The vertices at distance $r$ from $v_{i}$ are leaves, which proves the claim. 

Now suppose that $0\leq i<r$ so that $d_i=\min\!\left\{  i,r,k-i\right\}=i$. Then, $\operatorname*{dist}(v_{i},v_{r})=r-i$, and so $v_{i}$ is a leaf in $\mathcal{P}_{r-i}$. Now consider the bloom $\mathcal{B}(v_{i})$. The vertices in $\mathcal{B}(v_{i})$ at a distance $j$ from $v_{i}$ are leaves in $\mathcal{P}_{r-i+j}$. For $j<i$, these vertices are $p$-bloomed in $\mathcal{P}_{r-i+j+1}$, and the vertices at distance $i$ from $v_{i}$ are leaves in $\mathcal{P}_{r}$.

Lastly, suppose that $k-r<i\leq r$ so that $d_i=\min\!\left\{  i,r,k-i\right\}=k-i$. Then $\operatorname*{dist}(v_{i},v_{k-r})=i-k+r$, and thus $v_{i}$ is a leaf in $\mathcal{P}_{i-k+r}$. The vertices in the bloom $\mathcal{B}(v_{i})$ at a distance $j$ from $v_{i}$ are leaves in $\mathcal{P}_{i-k+r+j}$. For $j<k-i$, these vertices are $p$-bloomed in $\mathcal{P}_{i-k+r+j+1}$, and the vertices at a distance $k-i$ from $v_{i}$ are leaves in $\mathcal{P}_{r}$. This proves the first part of the lemma.

The converse follows from Corollary \ref{cor:isographsblooms}. 
\end{proof}

Note that this result gives a different description of the graphs $\mathcal H_{p^k}^r$. Instead of starting with a path on $k-2r$ edges, one can start with any $k$-spine and attach a suitable bloom to each vertex.  This gives a constructive way to describe the graph $\mathcal{H}_{p^{k}}^{r}$ without referring back to the recursive sequence $\mathcal{P}_{0},\mathcal{P}_{1},\ldots,\mathcal{P}_{r}$ used in its construction. Note that the resulting construction does not depend on the initial $k$-spine. We now record the following corollary, which is a direct consequence of Proposition~\ref{lem:bloomdepth}.
\begin{corollary}\label{cor:bloomdepth}
    Let $\mathcal{S}$ be a $k$-spine of $\mathcal{H}_{p^{k}}^{r}$ with vertices $\left\{  v_{0},v_{1},\ldots,v_{k}\right\}  $ and edges $\left\{  v_{i},v_{i+1}\right\}$. Let $d_{i}=\min\!\left\{  i,r,k-i\right\}$ and let $\mathcal{B}(v_i)$ be the bloom of $v_i$ relative to $\mathcal{S}$. Then, $d_i$ is the bloom depth of $v_i$, and if $x\in \mathcal{B}(v_i)$, then $\operatorname{dist}(v_i,x)\leq d_i$. If equality holds, then $x$ is a leaf of $\mathcal{H}_{p^{k}}^{r}$. Otherwise, $x$ is $p$-bloomed.
\end{corollary}

By construction, the graph $\mathcal{H}_{p^{k}}^{r}$ is symmetric with respect to reflection about the midpoint of any $k$-spine $\mathcal{S}$ with vertices $\{v_0, v_1, \ldots, v_k\}$ and edges $\{v_i, v_{i+1}\}$. The vertices $v_0, v_1, v_2 \ldots, v_{r-1}$ have bloom depths $0,1,2,\ldots,r-1$ as one moves toward the skeleton, while the vertices $v_{k-r+1}, \allowbreak v_{k-r+2}, \allowbreak \ldots, v_k$ have bloom depths $r-1,r-2,\ldots,0$ as one moves away from the skeleton. The vertices on the skeleton all have bloom depth $r$. 

Now observe that $\mathcal{H}_{p^{k}}^{0}$ is a $k$-spine by definition, and for a fixed prime $p$ and positive integer $k$, we have the following nested chain of graphs:
\[
\begin{tikzcd}[column sep=large, row sep=large]
\mathcal{H}_{p^{k}}^{0} \arrow[r,hook] \arrow[d,hook] &
\mathcal{H}_{p^{k}}^{1} \arrow[r,hook] \arrow[d,hook] &
\cdots \arrow[r,hook] &
\mathcal{H}_{p^{k}}^{\lfloor k/2\rfloor} \arrow[d,hook] \\
\mathcal{H}_{p^{k+1}}^{0} \arrow[r,hook] &
\mathcal{H}_{p^{k+1}}^{1} \arrow[r,hook] &
\cdots \arrow[r,hook] &
\mathcal{H}_{p^{k+1}}^{\lfloor k/2\rfloor} \arrow[r,hook] &
\mathcal{H}_{p^{k+1}}^{\lfloor (k+1)/2\rfloor}
\end{tikzcd}
\]

\begin{notation}
Since $\mathcal{H}_{p^{k}}^{r}$ is determined by any $k$-spine together with
the associated bloom depths, we introduce a compact notation for depicting
these graphs. Henceforth, we depict $\mathcal{H}_{p^{k}}^{r}$ as a $k$-spine
whose internal vertices are labeled by their bloom depths. In particular, the
bloom depth is given in a star marker, as shown in the figure below,
where the endpoints of the spine are drawn as ordinary vertices, while the
internal vertices are given with their bloom depth marker.
\[
\begin{tikzpicture}[scale=1.2, transform shape]
    \Vertex[shape=circle, size=0.1, color=white, x=0]{1}
    \Vertex[shape=star, size=0.5, color=white, x=1, label=1]{2}
    \Vertex[shape=star, size=0.5, color=white, x=2, label=2]{3}
    \draw[line width=1.5pt, dash pattern=on 5pt off 3pt] (2.2,0) -- (4,0);
    \Vertex[shape=star, size=0.5, color=white, x=3.8, label=r]{4}
    \draw[line width=1.5pt, dash pattern=on 5pt off 3pt] (4,0) -- (5.5,0);
    \Vertex[shape=star, size=0.5, color=white, x=5.8, label=r]{5}
    \draw[line width=1.5pt, dash pattern=on 5pt off 3pt] (6,0) -- (7.5,0);
    \Vertex[shape=star, size=0.5, color=white, x=7.8, label=2]{8}
    \Vertex[shape=star, size=0.5, color=white, x=8.8, label=1]{9}
    \Vertex[shape=circle, size=0.1, color=white, x=9.8]{10}

    \Edge(1)(2)
    \Edge(2)(3)
    \Edge(8)(9)
    \Edge(9)(10)

    \node at (5, -0.75) { $\mathcal{H}_{p^k}^r$};

\end{tikzpicture}
\]
\end{notation}

\begin{example}
In Example~\ref{ex:blossoming}, we constructed $\mathcal{H}_{3^6}^2$. We now depict this graph as
\[
\begin{tikzpicture}[scale=1.2, transform shape]

    \Vertex[ size=0.1, color=white]{1}
    \Vertex[shape=star, size=0.5, color=white, x=1, label=1]{2}
    \Vertex[shape=star, size=0.5, color=white, x=2, label=2]{3}
    \Vertex[shape=star, size=0.5, color=white, x=3, label=2]{4}
    \Vertex[shape=star, size=0.5, color=white, x=4, label=2]{5}
    \Vertex[shape=star, size=0.5, color=white, x=5, label=1]{6}
    \Vertex[size=0.1, color=white, x=6]{7}
   
    \Edge(1)(2)
    \Edge(2)(3)
    \Edge(3)(4)
    \Edge(4)(5)
    \Edge(5)(6)
    \Edge(6)(7)

\node at (3, -0.75) { $\mathcal{H}_{3^6}^2$};
        
\end{tikzpicture}
\]
\end{example}

We now record some combinatorial properties of the graphs $\mathcal{H}_{p^{k}}^{r}$, which will be play a key role in the results that follow.

\begin{lemma}\label{lem:kspineproperty}
Let $\mathcal{P}_{0}$ be the skeleton of $\mathcal{H}_{p^{k}}^{r}$, and for $v\in V(\mathcal{P}_{0})$, let $\mathcal{B}_{0}(v)$ denote the bloom of~$v$ relative to $\mathcal{P}_{0}$. Then,

\begin{enumerate}[label=\textnormal{(\arabic*)}]
\item if $v$ and $v^{\prime}$ are endpoints of the skeleton, with $w$ a leaf of $\mathcal{B}_{0}(v)$ and $w^{\prime}$ a leaf of $\mathcal{B}_{0}(v^{\prime})$, then the unique path from $w$ to $w^{\prime}$ is a $k$-spine, provided that when $v=v^{\prime}$, the unique paths from $v$ to $w$ and from $v$ to $w^{\prime}$ meet only at $v$;

\item A vertex $x\in V(\mathcal{H}_{p^{k}}^{r})$ lies on some $k$-spine if and only if $x$ is a skeletal vertex or $x$ lies in the bloom, relative to the skeleton, of an endpoint of the skeleton.

\item If a vertex $x\in V(\mathcal{H}_{p^{k}}^{r})$ does not lie on any $k$-spine, then $r<\left\lfloor \frac{k}{2} \right \rfloor$ and $x$ lies in the bloom of an internal skeletal vertex $u$. Further, if $v$ and $v^{\prime}$ are the endpoints of the skeleton chosen so that $\operatorname*{dist}(u,v^{\prime})\leq \operatorname*{dist}(u,v)$, then every longest path containing $x$ has one endpoint equal to a leaf of $\mathcal{\mathcal{B}}_{0}(u)$ and the other equal to a leaf of $\mathcal{\mathcal{B}}_{0}(v)$. Moreover, the length of such a path is $2r+\operatorname*{dist}(u,v)<k$.
\end{enumerate}
\end{lemma}

\begin{proof}
First suppose $r=0$. Then, $\mathcal{H}_{p^{k}}^{0}$ is its own skeleton, and thus claims $(1)$, $(2)$, and $(3)$ are immediate. So it remains to consider the case when $r\geq1$. 

Now let $v$ and $v^{\prime}$ be endpoints of the skeleton, and let $w$ and $w^{\prime}$ be leaves of $\mathcal{B}_{0}(v)$ and $\mathcal{B}_{0}(v^{\prime})$, respectively. If $v=v^{\prime}$, i.e., $r=\frac{k}{2}$, suppose that the paths from $v$ to $w$ and from $v$ to $w^{\prime}$ meet only at $v$. Since the bloom depth of each endpoint of the skeleton is $r$, we have $\operatorname*{dist}(w,v)=r$ and $\operatorname*{dist}(v^{\prime},w^{\prime})=r$. By definition, the distance between the endpoints of the skeleton is $k-2r$. Since $\mathcal{H}_{p^{k}}^{r}$ is a tree, the unique path $\mathcal{S}$ from $w$ to $w^{\prime}$ passes through $v$, then along the skeleton from $v$ to $v^{\prime}$, and then through $v^{\prime}$ to $w^{\prime}$. Therefore its length is
\[
\operatorname*{dist}(w,w^{\prime})=\operatorname*{dist}(w,v)+\operatorname*{dist}(v,v^{\prime})+\operatorname*{dist}(v^{\prime},w^{\prime})=r+(k-2r)+r=k.
\]
So the unique path from $w$ to $w^{\prime}$ has $k$ edges, and hence is a $k$-spine, which establishes $(1)$.

Next, suppose that $x\in V(\mathcal{H}_{p^{k}}^{r})$ lies on a $k$-spine $\mathcal{S}$. By Lemma \ref{Lem:skel}, $\mathcal{S}$ contains the skeleton $\mathcal{P}_{0}$ as a subgraph. Since $\mathcal{S}$ is a $k$-spine containing the subpath $\mathcal{P}_{0}$, any vertex of $\mathcal{S}\backslash V(\mathcal{P}_{0})$ must lie in a bloom attached to one of the endpoints of
$\mathcal{P}_{0}$. Hence,
\[
x\in V(\mathcal{B}_{0}(v))\cup V(\mathcal{P}_{0})\cup\mathcal{B}_{0}(v^{\prime}),
\]
where $v$ and $v^{\prime}$ are the endpoints of the skeleton.

Conversely, suppose first that $x\in V(\mathcal{P}_{0})$. Since $\mathcal{H}_{p^{k}}^{r}$ has a $k$-spine and each $k$-spine contains the skeleton by Lemma~\ref{Lem:skel}, we have that $x$ lies on some $k$-spine. 

Now let $v$ and $v^{\prime}$ be the endpoints of the skeleton $\mathcal{P}_{0}$ of $\mathcal{H}_{p^{k}}^{r}$. Note that if $r=\frac{k}{2}$, then $v=v^{\prime}$. Now suppose that $x\in V(\mathcal{B}(v))$, so that $\operatorname*{dist}(x,v)=j\leq r$ by Corollary~\ref{cor:bloomdepth}. Since $\mathcal{B}_{0}(v)$ is a finite tree, we can choose a leaf $w$ of $\mathcal{B}_{0}(v)$ such that $x$ lies on the unique path from $v$ to $w$. Note that since $w$ is a leaf, $\operatorname*{dist}(w,v)=r$. If $v\neq v^{\prime}$, choose any leaf $w^{\prime}\in V(\mathcal{B}_{0}(v^{\prime}))$. If $v=v^{\prime}$, choose a leaf $w^{\prime}\in V(\mathcal{B}_{0}(v))$ whose path from $v$ meets the path $v$ to $w$ only at $v$; such a choice is possible because the bloom at an endpoint has at least two branches. In either case, part (1) shows that the unique path from $w$ to $w^{\prime}$ is a $k$-spine, and by construction, it contains $x$. The case $x\in V(\mathcal{B}(v^{\prime}))$ is identical by symmetry, which proves $(2)$.

Now suppose that $x\in V(\mathcal{H}_{p^{k}}^{r})$ does not lie on some $k$-spine . By (2), $x$ is neither a skeletal vertex nor a vertex lying in the bloom, relative to the skeleton, of an endpoint of the skeleton. Hence $x$ lies in the bloom of an internal vertex $u$. In particular, the skeleton has at least three vertices, and so $k-2r\geq2$, and thus $r<\left\lfloor \frac{k}{2} \right \rfloor$. Now let $v$ and $v^{\prime}$ be the endpoints of the skeleton, and without loss of generality we may assume that $\operatorname*{dist}(u,v^{\prime})\leq\operatorname*{dist}(u,v)$. Since $\mathcal{B}_{0}(u)$ is a finite tree, there exists a leaf $w$ of $\mathcal{B}_{0}(u)$ such that $x$ lies on the unique path from $u$ to $w$. Likewise, let $w^{\prime}$ be any leaf of $\mathcal{B}_{0}(v)$. The unique path from $w$ to $w^{\prime}$ passes through $u$, then along the skeleton from $u$ to $v$, and finally from $v$ to $w^{\prime}$. Since every skeletal vertex has bloom depth $r$ by Corollary~\ref{cor:bloomdepth}, we have that the length of the path is
\[
\operatorname*{dist}(w,w^{\prime})=\operatorname*{dist}(w,u)+\operatorname*{dist}(u,v)+\operatorname*{dist}(v,w^{\prime})=2r+\operatorname*{dist}(u,v).
\]
This path contains $x$ by construction. Now observe that any longest path containing $x$ must pass through $u$. Since each skeletal vertex has bloom depth $r$, the longest path is achieved by traveling from $u$ to the farthest skeletal vertex from $u$. By our choice of $v$, this skeletal vertex is $v$. It follows that the longest path has length $2r+\operatorname*{dist}(u,v)$, which shows (3).
\end{proof}

\begin{lemma}
\label{Lem:bloomvertices}Let $v$ be a vertex on a $k$-spine of the graph $\mathcal{H}_{p^{k}}^{r}$ with bloom depth $d$, and let~$\varphi$ denote Euler's totient function. Then, the number of vertices in $\mathcal{B}(v)$ that are a distance $i\leq d$ from $v$ is $\varphi(p^{i})$. In particular, $\mathcal{B}(v)$ has $p^{d}$ vertices.
\end{lemma}

\begin{proof}
For $0\leq i\leq d$, let
\[
a_{i}=\left\vert w\in\mathcal{B}(v)\mid\operatorname*{dist}(v,w)=i\right\vert.
\]
Then, the result holds for $i=0$ since $a_{0}=1=\varphi(p^{0})$. So suppose the bloom depth is $d\geq1$. Then $v$ is $p$-bloomed, and thus $\deg v=p+1$. However, two of the edges incident to $v$ are on the $k$-spine, so exactly $p-1$ edges from $v$ are in $\mathcal{B}(v)$. Hence, $a_{1}=p-1=\varphi(p)$. Now fix $1\leq i<d$. By definition of bloom depth, every vertex at distance $i$ from $v$ is $p$-bloomed. Such a vertex has one edge to a vertex in the previous layer, and exactly $p$ edges to vertices in the next layer. Thus, each vertex at distance $i$ gives rise to exactly $p$ vertices at distance $i+1$, and so $a_{i+1}=pa_{i}$. By induction on $d$, we then obtain that
\[
a_{i}=p^{i-1}(p-1)=\varphi(p^{i})\qquad\text{for }1\leq i\leq d.
\]
Consequently, the number of vertices in the bloom at $v$ is
\[
\left\vert \mathcal{B}(v)\right\vert =\sum_{i=0}^{d}\varphi(p^{i})=1+\sum_{i=1}^{d}(p-1)p^{i-1}=p^{d}. \qedhere
\]
\end{proof}

\begin{corollary}\label{number_of_vertices}
The number of vertices in the graph $\mathcal{H}_{p^{k}}^{r}$ is
\[
\left(  k+1-2r\right)  p^{r}+2\cdot\frac{p^{r}-1}{p-1}.
\]

\end{corollary}

\begin{proof}
Let $\mathcal{S}$ be a $k$-spine for $\mathcal{H}_{p^{k}}^{r}$ with vertices $\left\{  v_{0},v_{1},\ldots,v_{k}\right\}  $ and edges $\left\{v_{i},v_{i+1}\right\}  $. Then, the bloom depth of $v_{i}$ is $d_{i}=\min \!\left\{  i,r,k-i\right\}  $ by Corollary ~\ref{cor:bloomdepth}. From the definition of bloom, we obtain that
\[
V\!\left(  \mathcal{H}_{p^{k}}^{r}\right)  =\bigcup_{i=0}^{k}V\!\left(\mathcal{B}(v_{i})\right)  .
\]
Since the blooms are pairwise disjoint, we obtain from Lemma~\ref{Lem:bloomvertices} that
\begin{align*}
\left\vert V(\mathcal{H}_{p^{k}}^{r})\right\vert =\sum_{i=0}^{k}\left\vert \mathcal{B}(v_{i})\right\vert  & =\sum_{i=0}^{k}p^{d_{i}}\\
& =\sum_{i=k-r+1}^{k}p^{k-i}+ \sum_{i=0}^{r-1}p^{i}+\sum_{i=r}^{k-r}p^{r}\\
& =\left(  k+1-2r\right)  p^{r}+2\sum_{i=0}^{r-1}p^{i}\\
& =\left(  k+1-2r\right)  p^{r}+2\cdot\frac{p^{r}-1}{p-1}. \qedhere
\end{align*}
\end{proof}

\begin{lemma}\label{pathcounting} 
Let $v_{0}$ be a leaf of $\mathcal{H}_{p^{k}}^{r}$ which lies on a $k$-spine. For $0\leq j\leq k$, the number of vertices which are at a distance $j$ from $v_{0}$ is $p^{\min\left\{  r,\left\lfloor j/2\right\rfloor \right\}  }$.
\end{lemma}

\begin{proof}
Let $\mathcal{S}$ be a $k$-spine with $V(\mathcal{S})=\left\{  v_{0},v_{1},\ldots,v_{k}\right\}  $ and edges $\left\{  v_{i},v_{i+1}\right\}  $ for $0\leq i<k$. By Corollary \ref{cor:bloomdepth}, the bloom depth of $v_{i}$ is $d_{i}=\min\!\left\{  i,r,k-i\right\}  $. Since $\mathcal{H}_{p^{k}}^{r}$ is a tree, each vertex $w\in V(\mathcal{H}_{p^{k}}^{r})$ satisfying $\operatorname*{dist}(v_{0},w)=j$ determines a unique path $\mathcal{P}_{w}$ of length $j$ starting at $v_{0}$ and ending at $w$. It thus suffices to count the number of paths of length $j$.

Let $\mathcal{P}$ be a path of length $j$ starting at $v_{0}$. Let $v_{i}$ be the last vertex of $\mathcal{P}$ on $\mathcal{S}$, and let $t$ be the length of the subpath of $\mathcal{P}$ lying in the bloom $\mathcal{B}(v_{i})$. In particular, $i+t=j$ and $t\leq d_{i}$. We further have that $t\leq i$ and $t\leq r$. The former, together with $i=j-t$ implies that $2t\leq j$. It then follows that $t\leq\left\lfloor j/2\right\rfloor $, and so $t\leq\min\!\left\{r,\left\lfloor j/2\right\rfloor \right\}  $. In particular, $0\leq t\leq \min\!\left\{  r,\left\lfloor j/2\right\rfloor \right\}  $, and the last vertex
on $\mathcal{S}$ is forced to be $v_{j-t}$.

By Lemma \ref{Lem:bloomvertices}, the number of choices for the endpoint of $\mathcal{P}$ in $\mathcal{B}(v_{j-t})$ at distance $t$ from $v_{j-t}$ is $\varphi(p^{t})$. It now follows that 
\begin{align*}
\left\vert \left\{  w\in V(\mathcal{H}_{p^{k}}^{r})\mid\operatorname*{dist}(v_{0},w)=j\right\}  \right\vert  &  =\sum_{t=0}^{\min\!\left\{  r,\left\lfloor \frac{j}{2}\right\rfloor \right\}  }\varphi(p^{t})\\
&  =1+\sum_{t=1}^{\min\!\left\{  r,\left\lfloor \frac{j}{2}\right\rfloor\right\}  }(p-1)p^{t-1}\\
&  =p^{\min\left\{  r,\left\lfloor \frac{j}{2}\right\rfloor \right\}
}.\qedhere
\end{align*}

\end{proof}

\subsection{Classification of finite primary isogeny graphs}\label{subsec:modpGalois}

We now classify, for each prime $p$, the possible finite $p$-primary subgraphs of an elliptic curve over a field $K$ of characteristic zero with $\operatorname{End}_K(E)\cong \mathbb{Z}$. Our aim is to prove that the possible finite $p$-primary subgraphs of isogeny graphs are exactly the graphs $\mathcal H_{p^k}^r$. The proof strategy is as follows. Consider an elliptic curve $E$ in the $p$-primary subgraph $\mathcal G_p(E/K)$ which lies on a path of maximal length. The image of its $p$-adic Galois representation uniquely determines the number of $p^i$-isogenies of $E$ for each~$i$. This corresponds to the number of vertices at a distance $i$ in the $p$-primary subgraph of the isogeny graph. If these numbers coincide with the number of vertices at a distance $i$ from the leaf in $\mathcal H_{p^k}^r$, and if we know that one of the graphs is a subgraph of another, then it follows that the $p$-primary subgraph is isomorphic to $\mathcal H_{p^k}^r$.

We now turn to Galois representations of elliptic curves. The following technical result is a slight generalization of \cite[Proposition 4.4]{bourdon-najman}.

\begin{proposition}\label{prop:imagechangeuponisogeny}
Let $E_{0}$ be an elliptic curve defined over a field $K$ of characteristic $0$, and suppose that $(P,Q)$ is a basis for $E_{0}[p^{k}]$ such that $\left\langle p^{k-l}P\right\rangle $ is $G_{K}$-invariant with $1\le l\leq k$. For $i\in\left\{  0,1,\ldots,l\right\}  $, let $E_{i}\cong E_{0}/\left\langle p^{k-i}P\right\rangle $, and let $\phi_{i}:E_{0}\rightarrow E_{i}$ denote the corresponding $K$-rational $p^{i}$-isogeny. Suppose further that $P_{i}$ is any point in $E_{0}$ such that $p^{i}P_{i}=P$ and set $(R,S):=\left(  \phi_{i}(P_{i}),\phi_{i}(Q)\right)  $. Then, $(R,S)$ is a basis for $E_{i}[p^{k}]$, and if $\sigma\in G_{K}$, then
\[
\rho_{E_{0},p^{k}}(\sigma)=\left(
\begin{array}
[c]{cc}
a & b\\
wp^{l} & c
\end{array}
\right)  \qquad\text{and}\qquad\rho_{E_{i},p^{k}}(\sigma)=\left(
\begin{array}
[c]{cc}
a & bp^{i}\\
wp^{l-i}+tp^{k-i} & c
\end{array}
\right)
\]
in basis $(P,Q)$ and $(R,S)$, respectively, for some $a,b,c,w,t\in\mathbb{Z}/p^{k}\mathbb{Z}$. Further, if $\left\langle p^{k-l}P\right\rangle $ is the largest $G_{K}$-invariant subgroup of $\left\langle P\right\rangle $, then 
\[
\min\!\left\{  v_{p}(w)\left\vert \left(
\begin{array}
[c]{cc}
a & b\\
wp^{l} & c
\end{array}
\right)  \in\rho_{E,p^{k}}(G_{K})\right.  \right\}  =0.
\]

In the special case when $l=k$, then we have that for $\sigma \in G_K$,
\[
\rho_{E_{0},p^{k}}(\sigma)=\left(
\begin{array}
[c]{cc}
a & b\\
0 & c
\end{array}
\right)  \qquad\text{and}\qquad\rho_{E_{i},p^{k}}(\sigma)=\left(
\begin{array}
[c]{cc}
a & bp^{i}\\
wp^{k-i} & c
\end{array}
\right)
\]
for some $a,b,c,w\in\mathbb{Z}/p^{k}\mathbb{Z}$.
\end{proposition}

\begin{proof}
We first show that $(R,S)$ is a basis for $E_{i}[p^{k}]$. To this end, observe that $P_{i}$ has order $p^{k+i}$ since $p^{i}P_{i}=P$ and $P$ has order $p^{k}$. Further, $p^{k-i}P=p^{k-i}(p^{i}P_{i})=p^{k}P_{i}$. Thus, $\ker\phi_{i}=\left\langle p^{k-i}P\right\rangle =\left\langle p^{k}P_{i}\right\rangle $. It then follows from the first isomorphism theorem that
\[
\left\langle \phi_{i}(P_{i})\right\rangle \cong\frac{\left\langle P_{i}\right\rangle }{\left\langle P_{i}\right\rangle \cap\ker\phi_{i}}=\frac{\left\langle P_{i}\right\rangle }{\left\langle p^{k}P_{i}\right\rangle
}.
\]
Consequently, $\phi_{i}(P_{i})$ has order $p^{k}$. Since $\phi_{i}(Q)$ has order $p^{k}$, it suffices to prove that $R=\phi_{i}(P_{i})$ and $S=\phi_{i}(Q)$ are independent. So suppose that for some $\alpha,\beta\in\mathbb{Z}/p^{k}\mathbb{Z}$ one has
\[
\alpha\phi_{i}(P_{i})=\beta\phi_{i}(Q).
\]
Now for $x\in\mathbb{Z}/p^{k}\mathbb{Z}$, let $\widetilde{x}\in\mathbb{Z}/p^{k+i}\mathbb{Z}$ be a lift of $x$ so that $x=\widetilde{x}\ \operatorname{mod}p^{k}$. Then $\widetilde{\alpha}P_{i}-\beta Q\in\ker\phi_{i}$, and so there is a $\gamma\in\mathbb{Z}/p^{k}\mathbb{Z}$ such that
\[
\widetilde{\alpha}P_{i}-\beta Q=\gamma p^{k-i}P=\widetilde{\gamma p^{k}}P_{i}\qquad\Longrightarrow\qquad\left(  \widetilde{\alpha}-\widetilde{\gamma
p^{k}}\right)  P_{i}=\beta Q.
\]
Since $P_{i}$ and $Q$ are independent, it follows that $\beta=0$ and $\widetilde{\alpha}=\widetilde{\gamma p^{k}}=\widetilde{0}$. But then $\alpha=0$, which shows that $\left(  R,S\right)  $ is a basis for $E_{i}[p^{k}]$.

Now let $\sigma\in G_{K}$, and with respect to the basis $(P,Q)$, let
\[
\rho_{E_{0},p^{k}}(\sigma)=\left(
\begin{array}
[c]{cc}
a & b\\
u & c
\end{array}
\right)
\]
for some $a,b,c,u\in\mathbb{Z}/p^{l}\mathbb{Z}$. Since $\left\langle p^{k-l}P\right\rangle $ is $G_{K}$-invariant, we have that
\[
\sigma(p^{k-l}P)=p^{k-l}\sigma(P)=p^{k-l}(aP+uQ)=ap^{k-l}P+p^{k-l} uQ\in\left\langle p^{k-l}P\right\rangle .
\]
Consequently, $u=wp^{l}$ for some $w\in\mathbb{Z}/p^{k}\mathbb{Z}$. Now observe that
\[
\sigma(p^{k-l-v_{p}(w)}P)=ap^{k-l-v_{p}(w)}P+p^{k-l-v_{p}(w)}wp^{l} Q=ap^{k-l-v_{p}(w)}P\in\left\langle p^{k-l-v_{p}(w)}P\right\rangle .
\]
Since $\left\langle p^{k-l-v_{p}(w)}P\right\rangle $ is cyclic of order $l+v_{p}(w)$, we have that if $\left\langle p^{k-l}P\right\rangle $ is the largest $G_{K}$-invariant subgroup of $\left\langle P\right\rangle $, then there is at least one $w$ for which $v_{p}(w)=0$.

With $\rho_{E_{0},p^{k}}(\sigma)=\sm{a &  b  \\ w p^l  & c}$, we now calculate $\sigma(R)$ and $\sigma(S)$. Then,
\[
p^{i}\sigma(R)=\sigma(\phi_{i}(p^{i}P_{i}))=\sigma(\phi_{i}(P))=\phi_{i}(\sigma(P))=\phi_{i}(aP+wp^{l}Q)=a\phi_{i}(p^{i}P_{i})+wp^{l}\phi_{i}(Q).
\]
Thus, $p^{i}\sigma(R)=ap^{i}R+wp^{l}S$, and so $p^{i}(\sigma(R)-aR)=wp^{l}S$. Consequently, $p^{k+i-l}(\sigma(R)-aR)=\mathcal{O}$ and so $\sigma(R)-aR$ has order dividing $p^{k+i-l}$. Hence, $\sigma(R)-aR\in E[p^{k+i-l}]$. Note that $E[p^{k+i-l}]$ can be identified with the subgroup $p^{l-i}E[p^{k}]$ of $E[p^{k}]$, and we deduce that there are $y,z\in \mathbb{Z}/p^{k}\mathbb{Z}$ such that
\[
\sigma(R)=aR+p^{l-i}(yR+zS).
\]
Now observe that
\begin{align*}
wp^{l}S  & =p^{i}(\sigma(R)-aR)=p^{i}(p^{l-i}(yR+zS))=p^{l}(yR+zS)\\
& \Longrightarrow\qquad yp^{l}R+p^{l}(z-w)S=\mathcal{O}.
\end{align*}
Since $R$ and $S$ are independent, we have that $v_{p}(y)\geq k-l$ and $v_{p}(z-w)\geq k-l$. Therefore, $y=sp^{k-l}$ and $z=w+tp^{k-l}$ for some $s,t\in\mathbb{Z}/p^{k}\mathbb{Z}$. This shows that
\[
\sigma(R)=(a+sp^{k-i})R+(wp^{l-i}+tp^{k-i})S.
\]

On the other hand,
\[
\sigma(S)=\sigma(\phi_{i}(Q))=\phi_{i}(\sigma(Q))=\phi_{i}(bP+cQ)=b\phi_{i}(p^{i}P)+c\phi_{i}(Q)=bp^{i}R+cS.
\]
Written as a matrix, the action of $\sigma$ on $(R,S)$ is
\[
\rho_{E_{i},p^{k}}(\sigma)=\left(
\begin{array}
[c]{cc}
a+sp^{k-i} & bp^{i}\\
wp^{l-i}+tp^{k-i} & c
\end{array}
\right)  .
\]
It remains to prove that $sp^{k-i}=0$. In this direction, let $\chi_{p^{k}}:G_{K}\rightarrow(\mathbb{Z}/p^{k}\mathbb{Z})^{\times}$ be the cyclotomic character, and recall that for every elliptic curve $E/K$ and every $\sigma\in G_{K}$,
\[
\det(\rho_{E,p^{k}}(\sigma))=\chi_{p^{k}}(\sigma),
\]
and so the determinant depends only on the cyclotomic action and is independent of $E$. It then follows that
\[
\det(\rho_{E_{0},p^{k}}(\sigma))=\det(\rho_{E_{i},p^{k}}(\sigma))\qquad \Longrightarrow\qquad ac-bwp^{l}=ac+scp^{k-i}-bwp^{l},
\]
and so $scp^{k-i}=0$. Further, since $\det(\rho_{E_{0},p^{k}}(\sigma))\equiv ac\ \operatorname{mod}p$, we have that $v_{p}(ac)=0$. It must then be the case that $v_{p}(s)\geq i$, and so $a+sp^{k-i}=a$. We conclude by observing that the claim for $l=k$ is now immediate.
\end{proof}

We are now ready to classify the possible finite $p$-primary subgraphs. In the results that follow, to ease notation, we omit reference to the $K$-isomorphism class when referring to the vertices of $\mathcal{G}_{p}(E/K)$. Thus, when we write $E^{\prime}\in V(\mathcal{G}_{p}(E/K))$, it is understood that we mean $\left[  E^{\prime}\right]  _{K} \in V(\mathcal{G}_{p}(E/K))$.

\begin{theorem}\label{classificationGpk}
    Let $E$ be an elliptic curve defined over a field $K$ of characteristic $0$ such that $\operatorname*{End}_{K}E\cong\mathbb{Z}$. If for a prime $p$, $\deg\mathcal{G}_{p}(E/K)=p^{k}$ for some nonnegative integer $k$, then
\[
\mathcal{G}_{p}(E/K)\cong\mathcal{H}_{p^{k}}^{r},
\]
for some nonnegative integer $r \leq\frac{k}{2}$. 

Conversely, for any two nonnegative integers $k$ and $r$ such that $ r\leq\frac{k}{2}$, there exists an elliptic curve $E$ defined over a field $K$ of characteristic $0$ such that $\mathcal{G}_{p}(E/K)\cong\mathcal{H}_{p^{k}}^{r}$, unless $(p,r)=(2,0)$ and $k\geq2$.
\end{theorem}
\begin{proof}
    If $\deg\mathcal{G}_{p}(E/K)=1$, then the $p$-primary graph of $E$ is trivial, and thus isomorphic to the trivial graph $\mathcal{H}_{p^{0}}^{0}$. So suppose that $\deg\mathcal{G}_{p}(E/K)=p^{k}$ for some positive integer $k$. By Corollaries~\ref{cor:sameprimary}, there exists some $E_{0}\in V(\mathcal{G}_{p}(E/K))$ such that $E_{0}$ admits a $K$-rational $p^{k}$-isogeny $\phi_{k}$. Let $\ker\phi_{k}=\left\langle P\right\rangle $, and for $0\leq i\leq k$, set $E_{i}=E/\left\langle p^{k-i}P\right\rangle $. Let $\phi_{i}:E_{0}\rightarrow E_{i}$ be the corresponding $K$-rational $p^{i}$-isogeny. In particular, we have the following $k$-spine $\mathcal{S}$ in $\mathcal{G}_{p}(E/K)$:

\[\begin{tikzcd}
	{E_0} & {E_1} & {E_2} & \cdots & {E_k}
	\arrow[no head, from=1-1, to=1-2]
	\arrow[no head, from=1-2, to=1-3]
	\arrow[no head, from=1-3, to=1-4]
	\arrow[no head, from=1-4, to=1-5]
\end{tikzcd}\]
 Further, since $\deg\mathcal{G}_{p}(E/K)=p^{k}$, we have from Corollary \ref{cor:sameprimary} that $E_{0}$ and $E_{k}$ do not admit a $p^{k+1}$-isogeny over $K$. In particular, the endpoints $E_{0}$ and $E_{k}$ of $\mathcal{S}$ are leaves in $\mathcal{G}_{p}(E/K)$. 

Since $E_{0}$ admits a $K$-rational $p^{k}$-isogeny, we may choose a basis $(P,Q)$ for $E_{0}[p^{k}]$ such that $\rho_{E_{0},p^{k}}(G_{K})$ is contained inside the Borel subgroup $B_{0}(p^{k})$ by Lemma \ref{lem:Borel}. 
Next, set
\begin{equation}\label{eq:alphavalue}
\alpha=\min\!\left\{  v_{p}(a-c)\left\vert \left(
\begin{array}
[c]{cc}
a & b\\
0 & c
\end{array}
\right)  \in\rho_{E_{0},p^{k}}(G_{K})\right.  \right\}
\end{equation}
and let $r=\min\!\left\{  \alpha,\left\lfloor k/2\right\rfloor \right\}  $. We now show that $\mathcal{G}_{p}(E/K)\cong\mathcal{H}_{p^{k}}^{r}$. In this direction, let $P_{i}$ be any point in $E_{0}$ such that $p^{i}P_{i}=P$. By Proposition \ref{prop:imagechangeuponisogeny}, $\left(  \phi_{i}(P_{i}),\phi_{i}(Q)\right)  $ is a basis for $E_{i}[p^{k}]$. For any $\sigma\in G_{K}$, let $\rho_{E_{0},p^{k}}(\sigma)=\sm{a & b \\ 0 & c}  $. By loc. cit., for $0\leq i\leq k$, there exists some $w_{i}\in\mathbb{Z}/p^{k}\mathbb{Z}$ such that
\[
\rho_{E_{i},p^{k}}(\sigma)=\left(
\begin{array}
[c]{cc}
a & p^{i}b\\
p^{k-i}w_{i} & c
\end{array}
\right)
\]
in basis $\left(  \phi_{i}(P_{i}),\phi_{i}(Q)\right)  $. Next, for $0\leq i\leq k$, let $d_{i}=\min\!\left\{  i,r,k-i\right\}  $ and observe that
\begin{equation}\label{eq:thmfinitecenter}
\rho_{E_{i},p^{k}}(\sigma)=\left(
\begin{array}
[c]{cc}
a & p^{i}b\\
p^{k-i}w_{i} & c
\end{array}
\right)  \equiv\left(
\begin{array}
[c]{cc}
a & 0\\
0 & c
\end{array}
\right)  \ \operatorname{mod}p^{d_{i}}.
\end{equation}
Consequently, $\rho_{E_{i},p^{d_{i}}}(G_{K})$ is contained in the center of $\operatorname*{GL}\nolimits_{2}(\mathbb{Z}/p^{d_{i}}\mathbb{Z})$. Hence, all $p^{d_{i}}$-isogenies of $E_{i}$ are defined over $K$. Thus, every subgroup of $E_{i}[p^{d_{i}}]$ is $G_{K}$-invariant. In particular, every subgroup of order $p^{j}$ for $0\leq j\leq d_{i}$ is $G_{K}$-invariant, and so $E_{i}$ admits all of its $p^{j}$-isogenies over~$K$. In particular $E_r, E_{r+1}, \ldots, E_{k-r}$ admit all $p^{r}$-isogenies over $K$. Consider the induced subgraph $\mathcal G_p'(E/K)$ whose vertices are $E_r, \ldots, E_{k-r}$ and all  elliptic curves $p^j$-isogenous to those elliptic curves for $1\leq j \leq r$. Then, the argument in Proposition~\ref{prop:uniquetree} can be adapted to show that $\mathcal G_p'(E/K)$ is isomorphic to $\mathcal H_{p^k}^r$, with the isomorphism sending $E_r, E_{r+1}, \ldots, E_{k-r}$ to the skeletal vertices $v_r, v_{r+1},\ldots, v_{k-r}$, respectively, of $\mathcal H_{p^k}^r$. In particular, $\mathcal G_p(E/K)$ contains a subgraph isomorphic to $\mathcal H_{p^k}^r$. We claim that $\mathcal G_p'(E/K)$ is in fact the entire $\mathcal G_p(E/K).$
 
 We prove this by comparing the number of $p^{j}$-isogenies that $E_{0}$ admits for $1\leq j\leq k$ with the number of paths in $\mathcal{G}_{p}^{\prime}(E/K)$ starting at $E_{0}$ of length $j$, which is $p^{\min\left\{  r,[j/2]\right\}  }$ by Lemma \ref{pathcounting}.

From the proof of \cite[Proposition 3.1]{ivan-numberofiso}, we have that the number of $p^{j}$-isogenies admitted by $E_{0}$ over $K$ is either $p^{\min \left\{  \alpha,\left\lfloor j/2\right\rfloor \right\}  }=p^{\min\left\{r,\left\lfloor j/2\right\rfloor \right\}  }$ or $2p^{\alpha}$. By Lemma~\ref{pathcounting}, it suffices to show that $2p^{\alpha}$ is impossible for each $j$ such that $1\leq j \leq k$. By way of contradiction, suppose that for some $j$ with $1\leq j \leq k$, the number of $p^{j}$ isogenies of $E_{0}$ defined over $K$ is $2p^{\alpha}$. At the end of the proof of~\cite[Proposition 3.1]{ivan-numberofiso}, it is shown that in this setting, $\alpha <\frac{j+\beta}{2}$ where
\[
\beta=\min\!\left\{  v_{p}(b)\left\vert \left(
\begin{array}
[c]{cc}
a & b\\
0 & c
\end{array}
\right)  \in\rho_{E_{0},p^{k}}(G_{K})\right.  \right\}  .
\]
Since $E_{0}$ does not admit a $p^{k+1}$-isogeny over $K$, we have that there exists a matrix $\sm{a & b \\ 0 & c}  \in\rho_{E_{0},p^{k}}(G_{K})$ such that $v_{p}(b)=0$, as otherwise $E_{0}$ would have independent cyclic isogenies of degrees $p^{k}$ and $p$. Thus, $\beta=0$ and $\alpha<\frac{j}{2}$. It now follows that $r=\alpha$. We further obtain from the proof of loc. cit. that the system of equations
\[
\left\{  (a-c)y=by^{2}\left\vert \left(
\begin{array}
[c]{cc}
a & b\\
0 & c
\end{array}
\right)  \in\rho_{E_{0},p^{j}}(G_{K}), y\in \mathbb{Z}/p^{j}\mathbb{Z}\right.  \right\}
\]
has a proper solution $y\in \mathbb{Z}/p^{j}\mathbb{Z}$ such that $(a-c)y\neq0$ for some $\sm{a & b \\ 0 & c}  \in\rho_{E_{0},p^{j}}(G_{K})$. It is also shown that such a proper solution satisfies $v_{p}(y)=\alpha-\beta$, and so $v_{p}(y)=r$. Consequently, $v_{p}\!\left(  (a-c)y\right)  =2r<2r+1$. In particular, the system of equations has a proper solution modulo $p^{2r+1}$. It then follows from the proof of loc. cit. that the number of $K$-rational $p^{2r+1}$-isogenies admitted by $E_{0}$ is $2p^{r}$. We claim that any $p^{2r+1}$-isogeny admitted by $E_{0}$ factors though the degree $p^{r+1}$-isogeny $\phi_{r+1}:E_{0}\rightarrow E_{r+1}$. In this direction, let $E^{\prime}$ be an elliptic curve that is $p^{2r+1}$-isogenous to $E_{0}$, with corresponding isogeny $\phi:E_{0}\rightarrow E^{\prime}$. Since $E_{0}$ is a leaf of $\mathcal{G}_{p}(E/K)$, there is a vertex $E_{i}$ on $\mathcal{S}$ such that $\phi$ factors through $E_{i}$, but not $E_{i+1}$. Consequently, $\phi=\phi^{\prime}\circ\phi_{i}$ where $\phi^{\prime}:E_{i}\rightarrow E^{\prime}$ and $E^{\prime}\in V(\mathcal{B}(E_{i}))$. In particular, $\deg\phi^{\prime}=p^{2r+1-i}$. Further, the isogeny $\phi^{\prime}$ admitted by $E_{i}$ is independent from the $p^{k-i}$-isogeny $E_{i}\rightarrow E_{k}$, and so $\mathcal{G}_{p}(E/K)$ contains a path of length $\left(  2r+1-i\right)  +\left(  k-i\right)  =k+\left(  2r+1-2i\right)$ from $E^{\prime}$ to $E_{k}$ as demonstrated in the graph below. 
 \[\begin{tikzcd}[row sep=0.5em]
	{E_0} & {E_1} & {E_2} & \ldots & {E_i} & {E_{i+1}} & \ldots & {E_k} \\
	\\
	\\
	&&&& {E'}
	\arrow[no head, from=1-1, to=1-2]
	\arrow[no head, from=1-2, to=1-3]
	\arrow[no head, from=1-3, to=1-4]
	\arrow[no head, from=1-4, to=1-5]
	\arrow[no head, from=1-5, to=1-6]
	\arrow[no head, "{\text{distance of } 2r+1-i}", from=1-5, to=4-5]
	\arrow[no head, from=1-6, to=1-7]
	\arrow[no head, from=1-7, to=1-8]
\end{tikzcd}\]

Since $\deg\mathcal{G}_{p}(E/K)=p^{k}$, it follows from Corollary~\ref{cor:sameprimary} that $k$ is the maximum length of paths in the graph, and so $2r+1\leq2i$ and hence $r+1\leq i$, which establishes the claim. Now let $\theta:E_{r+1}\rightarrow E_{i}$ be the degree $p^{i-r-1}$-isogeny such that $\phi=\phi^{\prime}\circ\theta\circ\phi_{r+1}$. In particular, $\deg(\phi^{\prime}\circ\theta)=p^{r}$. Thus, to count the number of distinct isogenies $\phi$ of degree $p^{2r+1}$ that factor through $E_{r+1}$, it suffices to count the number of possible compositions $\theta\circ\phi_{r+1}$. Observe that $\theta\circ\phi_{r+1}$ is obtained as the composition of $r$ successive $p$-isogenies. At each step, there are $p$ choices for the next $p$-isogeny, so the total number of such compositions is at most $p^{r}$. Therefore, the maximum number of distinct isogenies $\phi$ of degree $p^{2r+1}$ admitted by $E_{0}$ is $p^{r}$, which is a contradiction of the assumption that $E_{0}$ has $2p^{r}$ isogenies of degree $p^{2r+1}$. We conclude that for $0\leq j\leq k$, the number of $p^{j}$-isogenies of $E_{0}$ defined over $K$ is $p^{\min\left\{  r,\left\lfloor j/2\right\rfloor \right\}}$.

However, by Lemma \ref{pathcounting}, $p^{\min\left\{  r,\left\lfloor j/2\right\rfloor \right\}  }$ is equal to the number of paths of length $j$ starting from $E_{0}$ in $\mathcal{G}_{p}^{\prime}(E/K)$. If $\mathcal{G}_{p}(E/K)$ had more edges, then $E_{0}$ would either have more $p^{j}$-isogenies for some $j\leq k$ or a $p^{k+1}$-isogeny, which is impossible. Thus, $\mathcal{G}_{p}(E/K)=\mathcal{G}_p'(E/K)$, which completes the proof of the forward direction.
    
For the converse direction, note that the value of $\alpha$ as given in \eqref{eq:alphavalue} uniquely determines the value of $r=\min\!\left\{\alpha_,\left\lfloor k/2\right\rfloor \right\}  $, and hence uniquely determines the shape of $\mathcal{G}_{p}(E/K)$. For each nonnegative integer $r \le \frac{k}{2}$, define
\[
H_{p^k}^r := \left\{
\begin{pmatrix} a & b \\ 0 & c \end{pmatrix} \in \GL_2(\mathbb{Z}/p^k\mathbb{Z})
\;\middle|\; v_p(a-c) \ge r
\right\}.
\]
It is easily verified that $H_{p^k}^r \leq \GL_2(\mathbb{Z}/p^k\mathbb{Z})$. It then follows from the Galois correspondence that there is an elliptic curve $E$ defined over a field $K$ of characteristic $0$ such that the following conditions are satisfied:
\begin{enumerate}
\item[(i)] $\rho_{E,p^{k}}(G_{K})$ is conjugate to a subgroup of $H_{p^{k}}^{r}$;

\item[(ii)] if $r+1\leq\frac{k}{2}$, then $\rho_{E,p^{k}}(G_{K})$ is not conjugate to a subgroup of $H_{p^{k}}^{r+1}$;

\item[(iii)] $\rho_{E,p^{k}}(G_{K})$ is not conjugate to a subgroup of $\operatorname*{GL}\nolimits_{2}(\mathbb{Z}/p^{k}\mathbb{Z})$ whose elements reduce to diagonal matrices modulo $p$;

\item[(iv)] $\rho_{E,p^{k+1}}(G_{K})$ is not conjugate to a subgroup of $B_{0}(p^{k+1})$.
\end{enumerate}
Then, $\mathcal{G}_{p}(E/K)\cong\mathcal{H}_{p^{k}}^{r}$. Finally, the case
$\left(  p,r\right)  =\left(  2,0\right)  $ with $k\geq 2$ is impossible since any uppertriangular matrix $\sm{a & b \\ 0 & c} \in \operatorname*{GL}\nolimits_{2}(\mathbb{Z}/2^{k}\mathbb{Z})$ has $v_2(a-c)>0$.
\end{proof}

The proof of Theorem \ref{classificationGpk} introduced the subgroup $H_{p^{k}}^{r}\leq B_{0}(p^{k})$. Our next definition introduces a more general family of groups. 

\begin{definition}\label{def:groupsfinite}
For a prime number $p$ and a positive integer $k$, let $r$ and $i$ be nonnegative integers satisfying $r\leq\frac{k}{2}$ and $i\leq k$. We define the group $H_{p^{k}}^{r}(i)$ to be
\[
H_{p^{k}}^{r}(i):=\left\{  \left.  \left(
\begin{array}
[c]{cc}
a & p^{i}b\\
p^{k-i}w & c
\end{array}
\right)  \in\operatorname*{GL}\nolimits_{2}(\mathbb{Z}/p^{\max\left\{  i,k-i\right\}  }\mathbb{Z})\right\vert v_{p}(a-c)\geq r\right\}  .
\]
We further take $H_{p^{k}}^{r}:=H_{p^{k}}^{r}(0)$ so that
\[
H_{p^{k}}^{r}=\left\{  \left.  \left(
\begin{array}
[c]{cc}
a & b\\
0 & c
\end{array}
\right)  \in B_{0}(p^{k})\right\vert v_{p}(a-c)\geq r\right\}  .
\]

\end{definition}

\begin{remark}
We note that the groups $H_{p^{k}}^{r}(i)$ and $H_{p^{k}}^{r}(k-i)$ are conjugate subgroups inside $\operatorname*{GL}\nolimits_{2}(\mathbb{Z}/p^{\max\left\{  i,k-i\right\}  }\mathbb{Z})$.
\end{remark}

\begin{corollary}\label{cor:subgroup-subgraph-corr}
Let $p$ be a prime number and let $E$ be an elliptic curve over a field $K$ of characteristic $0$. If the $p$-primary graph $\mathcal{G}_{p}(E/K)$ has a subgraph isomorphic to $\mathcal{H}_{p^{k}}^{r}$, then there is a $k$-spine $\mathcal{S}$ in $\mathcal{G}_{p}(E/K)$ with vertices $\left\{  E_{0},E_{1},\ldots,E_{k}\right\}  $ and edges $\left\{  E_{i},E_{i+1}\right\}  $ such that for each nonnegative integer $i\leq k$, $\rho_{E_{i},p^{\max\left\{i,k-i\right\}  }}(G_{K})$ is conjugate to a subgroup of $H_{p^{k}}^{r}(i)$.

Conversely, if $\operatorname*{End}_{K}\! E\cong\mathbb{Z}$ and there exists an elliptic curve $E^{\prime}\in V(\mathcal{G}_{p}(E/K))$ such that $\rho_{E^{\prime},p^{\max\left\{  i,k-i\right\}  }}(G_{K})$ is conjugate to a subgroup of $H_{p^{k}}^{r}(i)$ for some $i\in\left\{0,1,\ldots,k\right\}  $, then $\mathcal{G}_{p}(E/K)$ has a subgraph isomorphic to $\mathcal{H}_{p^{k}}^{r}$.
\end{corollary}

\begin{proof}
Suppose that $\mathcal{G}_{p}(E/K)$ has a subgraph isomorphic to $\mathcal{H}_{p^{k}}^{r}$. Then, there is an embedding $\iota:\mathcal{H}_{p^{k}}^{r}\hookrightarrow\mathcal{G}_{p}(E/K)$. Let $\mathcal{S}$ be a $k$-spine of $\mathcal{H}_{p^{k}}^{r}$ with vertices $\left\{  v_{0},v_{1},\ldots,v_{k}\right\}  $ and edges $\left\{  v_{i},v_{i+1}\right\}  $ for $0\leq i<k$. For $0\leq i\leq k$, let $\iota(v_{i})=E_{i}\in V(\mathcal{G}_{p}(E/K))$. From our assumptions, we have that for $0\leq i\leq k$, there is a $p^{i}$-isogeny $\phi_{i}:E_{0}\rightarrow E_{i}$. Let $\ker\phi_{k}=\left\langle P\right\rangle $, and let $(P,Q)$ be a basis for $E_{0}[p^{k}]$. Then, for any $\sigma\in G_{K}$, let $\rho_{E_{0},p^{k}}(\sigma)=\sm{a & b \\ 0 & c}$. By the proof of Theorem~\ref{classificationGpk}, $v_{p}(a-c)\geq r$. For $0\leq i\leq k$, let $P_{i}$ be a point of $E_{0}$ such that $p^{i}P_{i}=P$. By Proposition~\ref{prop:imagechangeuponisogeny}, $(\phi_{i}(P_{i}),\phi_{i}(Q))$ is a basis for $E_{i}[p^{k}]$, and with respect to this basis, we have that
\[
\rho_{E_{i},p^{k}}(\sigma)=\left(
\begin{array}
[c]{cc}
a & p^{i}b\\
p^{k-i}w_{i} & c
\end{array}
\right)
\]
for some $w_{i}\in\mathbb{Z}/p^{i}\mathbb{Z}$. Reducing modulo $p^{\max\left\{  i,k-i\right\}  }$, we obtain that $\rho_{E_{i},p^{\max\left\{  i,k-i\right\}  }}(G_{K})$ is conjugate to a subgroup of $H_{p^{k}}^{r}(i)$.

Conversely, suppose that $E^{\prime}\in V(\mathcal{G}_{p}(E/K))$ such that $\rho_{E^{\prime},p^{\max\left\{  i,k-i\right\}  }}(G_{K})$ is conjugate to a subgroup of $H_{p^{k}}^{r}(i)$ for some $i\in\left\{  0,1,\ldots,k\right\}  $. Now choose a basis $\left(  P,Q\right)$ of $E^{\prime}[p^{\max\left\{  i,k-i\right\}  }]$ such that $\rho_{E^{\prime},p^{\max\left\{  i,k-i\right\}  }}(G_{K})\leq H_{p^{k}}^{r}(i)$. Without loss of generality, we may assume that $i\leq\frac{k}{2}$ so that
\[
H_{p^{k}}^{r}(i)=\left\{  \left.  \left(
\begin{array}
[c]{cc}
a & p^{i}b\\
0 & c
\end{array}
\right)  \in\operatorname*{GL}\nolimits_{2}(\mathbb{Z}/p^{k-i}\mathbb{Z})\right\vert v_{p}(a-c)\geq r\right\}  .
\]
Now consider $\rho_{E', p^k}(G_K)$. We have $$\rho_{E', p^k}(G_K) \subseteq \left \{ \m{a & p^i b \\ p^{k-i}t & c}\in \GL_2(\Z/p^k\Z) \;\middle|\; v_p(a-c) \ge r
\right\}.$$ Then, by Proposition \ref{prop:imagechangeuponisogeny}, there is an elliptic curve $E_0$ which is $p^j$-isogenous to $E_j$ such that $\rho_{E_0, p^k}(G_K)$ is a subgroup of $H_{p^k}^r(0)$.
    Then $E_0$ has a $p^k$-isogeny, so there is a path  $$(E_0, E_1, \ldots,  E_k)$$ of length $k$ in $\mathcal G_p(E/K)$ such that $E'=E_i$. Furthermore, by Proposition \ref{prop:imagechangeuponisogeny}, if follows that $\rho_{E_j, p^{\max\{j,k-j\}}}(G_K)$ is contained in $H_{p^k}^r(j)$ for $0\leq j \leq k$.  In particular, the images of Galois representations of elliptic curves are all scalar modulo $p^r$. Hence, those elliptic curves have all $K$-rational $p^j$-isogenies for $1\leq j \leq r$. Taking the induced subgraph of $\mathcal G_p(E/K)$ whose vertices are $E_r, \ldots, E_{k-r}$ and the targets of these isogenies, we obtain a copy of $\mathcal H_{p^k}^r$ inside $\mathcal G_p(E/K)$.
\end{proof}

\begin{remark}\label{rem:subgroup-subgraph-generalized}
Let $E$ be an elliptic curve over a field $K$ of characteristic $0$ with $\operatorname*{End}_{K}\! E\cong\mathbb{Z}$. By Corollary~\ref{cor:subgroup-subgraph-corr}, it follows that the image of $\rho_{E,p^{\max(i,k-i)}}(G_K)$ is conjugate to a subgroup of $H_{p^k}^r(i)$ if and only if there exists an embedding $\iota$ of $\mathcal H_{p^k}^r$ into $\mathcal G_p(E/K)$ such that $\iota(v_i)=E$ for some $k$-spine of $\mathcal H_{p^k}^r$ with vertices $\left\{  v_{0},v_{1},\ldots,v_{k}\right\}  $ and edges $\left\{  v_{i},v_{i+1}\right\}  $ for $0\leq i<k$.
\end{remark}


\section{The \texorpdfstring{$p$}{p}-blooming invariant}\label{sec:bloominv}
We now take a closer look at the quantity~$\alpha$ in (\ref{eq:alphavalue}), which appeared in the proofs of \cite[Proposition~3.1]{ivan-numberofiso} and Theorem~\ref{classificationGpk}. In those settings, $\alpha$ was introduced relative to a chosen basis for the $\operatorname{mod}p^{k}$ Galois representation. To proceed towards our classification of infinite primary graphs, we require a $p$-adic version of $\alpha$. With the choice of basis in (\ref{eq:alphavalue}), the matrices are upper-triangular, and thus, their diagonal entries correspond to the eigenvalues of the matrix. Below, we introduce the main object of study in this section, the $p$-blooming invariant of an elliptic curve, and we prove afterward that it is an isogeny class invariant. We then show that for finite primary graphs, the quantity $r$ is determined from $k$ and the blooming invariant (see Corollary~\ref{cor:bloompotinv}). We note that in the next section, the quantity $r$ associated to infinite $p$-primary graphs is the $p$-blooming invariant (see Theorem~\ref{thm:infinitegpks}). We conclude by showing that in fields with a real embedding, the isogeny class degree uniquely determines finite isogeny graphs (see Proposition~\ref{Prop:BloomInvReal} and Corollary~\ref{Cor:BloomInvReal}). To begin, we note that in the definition below, $\Lambda(M)$ denotes the \textit{spectrum} of a square matrix $M$, i.e., the multiset of eigenvalues of $M$.
\begin{definition}\label{def:p-bloominginv}
Let $E$ be an elliptic curve over a field $K$ of characteristic $0$. For a prime $p$, we define the \textit{$p$-blooming invariant} of $E/K$ to be
\[
\mathfrak{I}_{p^{\infty}}(E/K)=\min_{\sigma\in G_{K}}\!\left\{  v_{p}\!\left(\lambda_{1}-\lambda_{2}\right)  \mid\Lambda(\rho_{E,p^{\infty}}(\sigma))=\left\{  \lambda_{1},\lambda_{2}\right\}  \right\}  ,
\]
and we set $\mathfrak{I}(E/K)=\mathfrak{I}_{p^{\infty}}(E/K)$ when $p$ is clear from context. Similarly, we omit reference to $K$ when it is clear from context, and simply write $\mathfrak{I}(E)$.
\end{definition}

We note that the $p$-blooming invariant $\mathfrak{I}(E)\in\mathbb{Z}[\frac{1}{2}]_{\ge0}\cup\{\infty\}$. To see this, fix $\sigma\in G_K$, and let $\Lambda(\rho_{E,p^{\infty}}(\sigma))=\left\{\lambda_{1},\lambda_{2}\right\}  $. Then, $\lambda_{1}$ and $\lambda_{2}$ are roots of the characteristic polynomial
\[
X^{2}-\operatorname*{Tr}(\rho_{E,p^{\infty}}(\sigma))X+\det(\rho_{E,p^{\infty}}(\sigma))\in\mathbb{Z}_{p}[X].
\]
Consequently, each $\lambda_{i}\in\overline{\mathbb{Q}}_{p}$ is integral, and thus $v_{p}(\lambda_{i})\geq0$. Moreover, as $\det(\rho_{E,p^{\infty}}(\sigma))$ does not vanish modulo $p$, we get that $v_{p}(\lambda_{1}\lambda_2)=0$, and so $v_p(\lambda_i)=0$ for $i\in \{1,2\}$. Hence, $v_{p}(\lambda_{1}-\lambda_{2})\geq0$. Further, $\lambda_1$ and $\lambda_2$ are integral elements of a quadratic extension of $\mathbb{Q}_p$, and so $v_p(\lambda_1-\lambda_2)\in \mathbb{Z}[\frac{1}{2}]_{\ge0}\cup\{\infty\}$. 

Next, we show that the $p$-blooming invariant is an isogeny class invariant. This is, in fact, automatic, since isogenous elliptic curves have isomorphic $p$-adic Tate modules. The proof below makes this precise.
\begin{proposition}\label{prop:bloomingpotential}
Let $p$ be a prime and $k$ be a positive integer. If $E_{1}$ and $E_{2}$ are isogenous elliptic curves over some field $K$ of characteristic $0$, then, for $\sigma\in G_{K}$, we have the following equality of spectrums:
\[
\Lambda(\rho_{E_{1},p^{\infty}}(\sigma))=\Lambda(\rho_{E_{2},p^{\infty}}(\sigma)).
\]
In particular, the blooming invariants of $E_{1}$ and $E_{2}$ are equal, and thus, the blooming invariant is an isogeny class invariant.
\end{proposition}

\begin{proof}
From our assumption, we have that there is a $K$-rational isogeny $\phi:E_{1}\rightarrow E_{2}$. Now let $T_{p}(E_{i})$ denote the Tate module of $E_{i}$. By \cite[Lemma~IV.2.2]{milneAV}, we have a short exact sequence
\[
0\longrightarrow T_{p}(E_{1})\longrightarrow T_{p}(E_{2})\longrightarrow
C\longrightarrow0,
\]
with $C$ finite. Since $C$ is finite, tensoring with $\mathbb{Q}_{p}$ yields an isomorphism of $p$-adic Tate modules:
\[
V_{p}(E_{1})=T_{p}(E_{1})\otimes_{\mathbb{Z}_{p}}\mathbb{Q}_{p}\cong T_{p}(E_{2})\otimes_{\mathbb{Z}_{p}}\mathbb{Q}_{p}=V_{p}(E_{2}).
\]
Next, let $\Phi:V_{p}(E_{1})\rightarrow V_{p}(E_{2})$ be the corresponding $G_{K}$-equivariant $\mathbb{Q}_{p}$-linear isomorphism. For $\sigma\in G_{K}$, the induced operators on $V_{p}(E_{i})$ are $\mathbb{Q}_{p}$-linear endomorphisms, and since $\Phi$ is $G_{K}$-equivariant, we obtain
\[
\Phi\circ(\rho_{E_{1},p^{\infty}}(\sigma)\otimes1)=(\rho_{E_{2},p^{\infty}}(\sigma)\otimes1)\circ\Phi.
\]
Thus, $\rho_{E_{1},p^{\infty}}(\sigma)\otimes 1$ and $\rho_{E_{2},p^{\infty}}(\sigma)\otimes 1$ are conjugate over $\mathbb{Q}_{p}$. Consequently, they have the same characteristic polynomial over $\mathbb{Q}_{p}$, and therefore $\Lambda(\rho_{E_{1},p^{\infty}}(\sigma))=\Lambda (\rho_{E_{2},p^{\infty}}(\sigma))$. The result now follows since
\[
\mathfrak{I}(E_{1})=\min_{\sigma\in G_{K}}\!\left\{  v_{p}\!\left(  \lambda _{1}-\lambda_{2}\right)  \mid\Lambda(\rho_{E_{1},p^{\infty}}(\sigma))=\left\{\lambda_{1},\lambda_{2}\right\}  \right\}  =\mathfrak{I}(E_{2}). \qedhere
\]
\end{proof}

To relate the blooming invariant to the quantity $\alpha$ in (\ref{eq:alphavalue}), we note that for a $2\times2$ matrix $M$ with eigenvalues $\lambda_{1}$ and $\lambda_{2}$, the discriminant of the characteristic polynomial satisfies $\operatorname*{Tr}(M)^{2}-4\det(M)=(\lambda_{1}-\lambda_{2})^{2}$. Thus, an alternate definition for the blooming invariant is:
\[
\mathfrak{I}(E)=\min_{\sigma\in G_{K}}\!\left\{  \frac{1}{2}v_{p}\!\left(\operatorname*{Tr}(\rho_{E,p^{\infty}}(\sigma))^{2}-4\det(\rho_{E,p^{\infty}}(\sigma))\right)  \right\}  .
\]
When considering the $\operatorname{mod}p^{k}$ Galois representation, it is this equivalent definition that we want to consider, as $\mathbb{Z}/p^{k}\mathbb{Z}$ is not necessarily a domain. In this direction, we first introduce for $k\in\mathbb{N}\cup\left\{  \infty\right\}  $, the $k$\textit{-truncated }$p$\textit{-adic valuation}, which is defined as
\[
v_{p}^{(k)}(x)=\min\!\left\{  v_{p}(x),k\right\}  .
\]
In particular, $v_{p}^{(\infty)}(x)=v_{p}(x)$. We are now ready to introduce the $\operatorname{mod}p^{k}$ blooming invariant. We first define it, and then state a result that makes the connection between the blooming invariant and the $\operatorname{mod}p^{k}$ blooming invariant explicit.

\begin{definition}\label{def:modpkbloominv}
Let $E$ be an elliptic curve over a field $K$ of characteristic $0$. For a prime $p$ and positive integer $k$, we define the \textit{$\operatorname{mod}p^{k}$ blooming invariant} to be
\[
\mathfrak{I}_{p^{k}}(E/K)=\min_{\sigma\in G_{K}}\!\left\{  \frac{1}{2}v_{p}^{(k)}\!\left(  \operatorname*{Tr}(\rho_{E,p^{k}}(\sigma))^{2}-4\det(\rho_{E,p^{k}}(\sigma))\right)  \right\}  ,
\]
and we set $\mathfrak{I}_{k}(E/K)=\mathfrak{I}_{p^{k}}(E/K)$ when $p$ is clear from context.
\end{definition}

\begin{remark}
As before, when $K$ is clear from context, we simply write $\mathfrak{I}_{p^k}(E)$, or $\mathfrak{I}_{k}(E)$ if $p$ is clear from context.
\end{remark}

\begin{corollary}\label{cor:bloompotprop}
Let $E$ be an elliptic curve over a field $K$ of characteristic $0$, and let $p$ be a prime and $k$ a positive integer. Then,

\begin{enumerate}
\item $\mathfrak{I}_{k}(E/K)$ is independent of basis of $E[p^{k}]$;

\item $\mathfrak{I}_{k}(E/K)=\min\!\left\{  \mathfrak{I}(E/K),\frac{k}{2}\right\}  $ and $\mathfrak{I}(E/K)=\lim_{k\rightarrow\infty}\mathfrak{I}_{k}(E/K)$;

\item $\mathfrak{I}_{k}(E/K)\leq\mathfrak{I}_{k+1}(E/K)$ for each positive
integer $k$;

\item $\mathfrak{I}_{k}(E/K)$ is an isogeny class invariant;

\item if $F/K$ is a field extension, then $\mathfrak{I}_{p^k}(E/F)\geq \mathfrak{I}_{p^k}(E/K)$ for every $k\in\mathbb{N}\cup\left\{  \infty\right\}  $.
\end{enumerate}
\end{corollary}

\begin{proof}
Observe that (1) follows by definition. Now let $\sigma\in G_{K}$, and suppose that $\Lambda(\rho_{E,p^{\infty}}(\sigma))=\left\{  \lambda_{1},\lambda_{2}\right\}  $. Then,
\begin{align*}
(\lambda_{1}-\lambda_{2})^{2}  & =\operatorname*{Tr}(\rho_{E,p^{\infty}}(\sigma))^{2}-4\det(\rho_{E,p^{\infty}}(\sigma))\\
& \equiv\operatorname*{Tr}(\rho_{E,p^{k}}(\sigma))^{2}-4\det(\rho_{E,p^{k}}(\sigma))\ \operatorname{mod}p^{k}\\
\Longrightarrow\qquad\frac{1}{2}v_{p}^{(k)}((\lambda_{1}-\lambda_{2})^{2})  &
=\frac{1}{2}v_{p}^{(k)}\!\left(  \operatorname*{Tr}(\rho_{E,p^{k}}(\sigma))^{2}-4\det(\rho_{E,p^{k}}(\sigma))\right)  .
\end{align*}
Now observe that $\frac{1}{2}v_{p}^{(k)}\!\left(  x^{2}\right)  =\min\{v_{p}(x),\frac{k}{2}\}$ for $x\in \mathbb{Z}_p$. Thus, if $v_{p}(\lambda_{1}-\lambda_{2})>\frac{k}{2}$, then $\frac{1}{2}v_{p}^{(k)}((\lambda_{1}-\lambda_{2})^{2})=\frac{k}{2}$. Similarly, if $v_{p}(\lambda_{1}-\lambda_{2})\leq\frac{k}{2}$, then
\[
\frac{1}{2}v_{p}^{(k)}((\lambda_{1}-\lambda_{2})^{2})=v_{p}(\lambda_{1}-\lambda_{2})\leq\frac{k}{2}.
\]
From this, we now see that
\begin{align*}
\mathfrak{I}_{k}(E/K)  & =\min_{\sigma\in G_{K}}\!\left\{  \frac{1}{2}v_{p}^{(k)}\!\left(  \operatorname*{Tr}(\rho_{E,p^{k}}(\sigma))^{2}-4\det(\rho_{E,p^{k}}(\sigma))\right)  \right\}  \\
& =\min_{\sigma\in G_{K}}\!\left\{  \frac{1}{2}v_{p}^{(k)}((\lambda_{1}-\lambda_{2})^{2})\mid\Lambda(\rho_{E,p^{\infty}}(\sigma))=\left\{\lambda_{1},\lambda_{2}\right\}  \right\}  \\
& =\min\!\left\{  \mathfrak{I}(E/K),\frac{k}{2}\right\}  .
\end{align*}
It now follows that $\mathfrak{I}(E/K)=\lim_{k\rightarrow\infty}\mathfrak{I}_{k}(E/K)$, which concludes (2). Observe that (3) follows from (2).

Next, suppose that $E$ and $E^{\prime}$ are $K$-isogenous. From Proposition
\ref{prop:bloomingpotential}, we have that $\mathfrak{I}(E/K)=\mathfrak{I}(E^{\prime}/K)$. It now follows from (2) that $\mathfrak{I}_{k}(E/K)=\mathfrak{I}_{k}(E^{\prime}/K)$.

Lastly, suppose $F/K$ is a field extension and let $\overline{K}$ be an algebraic closure of $K$ that contains~$F$. We can then identify $G_{F}=\operatorname*{Gal}(\overline{K}/F)$, so that $G_{F}\leq G_{K}$ is a closed subgroup. Since $\rho_{E/F,p^{\infty}} = \rho_{E/K,p^\infty}|_{G_F}$, we deduce that
\begin{align*}
\mathfrak{I}(E/K)  & =\min_{\sigma\in G_{K}}\!\left\{  v_{p}\!\left(\lambda_{1}-\lambda_{2}\right)  \mid\lambda(\rho_{E/K,p^{\infty}}(\sigma))=\left\{  \lambda_{1},\lambda_{2}\right\}  \right\}  \\
& \leq\min_{\sigma\in G_{F}}\!\left\{  v_{p}\!\left(  \lambda_{1}-\lambda_{2}\right)  \mid\lambda(\rho_{E/F,p^{\infty}}(\sigma))=\left\{  \lambda_{1},\lambda_{2}\right\}  \right\}  \\
& =\mathfrak{I}(E/F).
\end{align*}
We now have that (5) follows from (2) since $\mathfrak{I}_{k}(E/K)\leq\mathfrak{I}_{k}(E/F)$ for each positive integer $k$.
\end{proof}

In what follows, unless stated otherwise, we work over a fixed field $K$. We now give an explicit description of $\mathfrak{I}_{p^k}(E)$ for $k\in \mathbb{N} \cup \{\infty\}$ in the case where the Galois image is conjugate to a subgroup of the Borel subgroup.

\begin{lemma}\label{lem:explicitbloompot}
Let $E$ be an elliptic curve over a field $K$ of characteristic $0$, and let $p$ be a prime and $k\in\mathbb{N}\cup\left\{  \infty\right\}  $. Suppose that, after choosing a basis for $E[p^{k}]$ if $k<\infty$ or $T_{p}(E)$ if $k=\infty$, it is the case that $\rho_{E,p^{k}}(G_{K})\leq B_{0}(p^{k})$. Then,
\[
\mathfrak{I}_{p^{k}}(E/K)=\min\!\left\{  \min\!\left\{  v_{p}^{(k)}(a-c),\frac{k}{2}\right\}  \left\vert \left(
\begin{array}
[c]{cc}
a & b\\
0 & c
\end{array}
\right)  \in\rho_{E,p^{k}}(G_{K})\right.  \right\}  .
\]
In particular, $\mathfrak{I}_{p^{\infty}}(E/K)\in\mathbb{Z}_{\geq0}$, and if $\mathfrak{I}_{p^{k}}(E/K)<\frac{k}{2}$ for some positive integer $k$, then $\mathfrak{I}_{p^{k}}(E/K)\in\mathbb{Z}_{\geq0}$.
\end{lemma}

\begin{proof}
As in the proof of Corollary \ref{cor:bloompotprop} (2), we observe that $\frac{1}{2}v_{p}^{(k)}\!\left(  x^{2}\right)  =\min\{v_{p}(x),\frac{k}{2}\}$, with the convention that $\frac{1}{2}v_{p}^{(\infty)}\!\left(  x^{2}\right) =v_{p}\!\left(  x\right)  $. From our assumptions, we then obtain
\begin{align*}
\mathfrak{I}_{p^{k}}(E/K)  & =\min_{\sigma\in G_{K}}\!\left\{  \frac{1}{2}v_{p}^{(k)}\!\left(  \operatorname*{Tr}(\rho_{E,p^{k}}(\sigma))^{2}-4\det(\rho_{E,p^{k}}(\sigma))\right)  \right\}  \\
& =\min\!\left\{  \frac{1}{2}v_{p}^{(k)}\!\left(  (a-c)^{2}\right)  \left\vert
\left(
\begin{array}
[c]{cc}
a & b\\
0 & c
\end{array}
\right)  \in\rho_{E_{0},p^{k}}(G_{K})\right.  \right\}  \\
& =\min\!\left\{  \min\!\left\{  v_{p}^{(k)}(a-c),\frac{k}{2}\right\}  \left\vert
\left(
\begin{array}
[c]{cc}
a & b\\
0 & c
\end{array}
\right)  \in\rho_{E,p^{k}}(G_{K})\right.  \right\}  .
\end{align*}
Now let $M=\sm{a & b  \\ 0  & c}\in\rho_{E,p^{k}}(G_{K})$. If $k=\infty$, then $\min\!\left\{  v_{p}^{(k)}(a-c),\frac{k}{2}\right\}  =v_{p}^{(k)}(a-c)\in\mathbb{Z}_{\geq0}$ since $a,c\in\mathbb{Z}_{p}$. Similarly, if $k<\infty$ and $v_{p}^{(k)}(a-c)<\frac{k}{2}$, then $v_{p}^{(k)}(a-c)$ is an integer since $a,c\in\mathbb{Z}/p^{k}\mathbb{Z}$. The conclusion of the lemma now follows.
\end{proof}

\begin{corollary}\label{cor:bloompotinv}
Let $E$ be an elliptic curve over a field $K$ of characteristic $0$ with $\operatorname*{End}_K(E)\cong\mathbb{Z}$, and let $p$ be a prime. If the $p$-primary graph $\mathcal{G}_{p}(E/K)$ is finite, then $\mathcal{G}_{p}(E/K)\cong\mathcal{H}_{p^{k}}^{r}$ for some nonnegative integer $k$ with
\[
r=\min \!\left\{  \mathfrak{I}(E),\left\lfloor k/2\right\rfloor \right\}=\min \!\left\{  \mathfrak{I}_{k}(E),\left\lfloor k/2\right\rfloor \right\}  .
\]

\end{corollary}

\begin{proof}
By Theorem \ref{classificationGpk}, $\mathcal{G}_{p}(E/K)\cong\mathcal{H}_{p^{k}}^{r}$ for some nonnegative integers $k$ and $r$ with $r\leq\frac{k}{2}$. In particular, the result holds if $k=0$. So suppose that $k$ is a positive integer, and let $E_{0}\in V(\mathcal{G}_{p}(E/K))$ be a leaf that lies on some $k$-spine of $\mathcal{G}_{p}(E/K)$. Then, $E_{0}$ admits a $K$-rational $p^{k}$-isogeny, and without loss of generality, we may assume that $\rho_{E_{0},p^{k}}(G_{K})\leq B_{0}(p^{k})$. By the proof of loc. cit., we have that $r=\min\{\alpha,\left\lfloor k/2\right\rfloor \}$ where
\[
\alpha=\min\!\left\{  v_{p}(a-c)\left\vert \left(
\begin{array}
[c]{cc}
a & b\\
0 & c
\end{array}
\right)  \in\rho_{E_{0},p^{k}}(G_{K})\right.  \right\}  .
\]
By Corollary \ref{cor:bloompotprop} and Lemma \ref{lem:explicitbloompot}, we have that $\mathfrak{I}(E)=\mathfrak{I}(E_{0})$ and thus
\[
\mathfrak{I}_{k}(E)=\min\!\left\{  \mathfrak{I}(E_0),\frac{k}{2}\right\}=\min\!\left\{  \min \!\left\{  v_{p}^{(k)}(a-c),\frac{k}{2}\right\}  \left\vert\left(
\begin{array}
[c]{cc}
a & b\\
0 & c
\end{array}
\right)  \in\rho_{E_{0},p^{k}}(G_{K})\right.  \right\}  .
\]
It follows that if $\mathfrak{I}_{k}(E)<\frac{k}{2}$, then $\alpha=\mathfrak{I}_{k}(E)$. Therefore, $r=\mathfrak{I}_{k}(E)$. So suppose that $\mathfrak{I}_{k}(E)=\frac{k}{2}$. Then, $\alpha\geq\frac{k}{2}$, and so $r=\left\lfloor \frac{k}{2}\right\rfloor =\min\{\mathfrak{I}_{k}(E),\left\lfloor \frac{k}{2}\right\rfloor \}$. Thus, in both cases, we have that
\[
r=\min\!\left\{  \mathfrak{I}_{k}(E),\left\lfloor \frac{k}{2}\right\rfloor\right\}  =\min\!\left\{  \min \!\left\{  \mathfrak{I}(E),\frac{k}{2}\right\},\left\lfloor \frac{k}{2}\right\rfloor \right\}  =\min\!\left\{\mathfrak{I}(E),\left\lfloor \frac{k}{2}\right\rfloor \right\}  . \qedhere
\]
\end{proof}

\begin{corollary}\label{cor:longestpathbp}
Let $E$ be an elliptic curve over a field $K$ of characteristic $0$ with $p$-primary graph $\mathcal{G}_{p}(E/K)\cong\mathcal{H}_{p^{k}}^{r}$ for some positive integer $k$ and nonnegative integer $r\leq\frac{k}{2}$. Suppose further that $j$ is the length of a longest path in $\mathcal{G}_{p}(E/K)$ that contains $E$. If $E$ does not lie on any $k$-spine of $\mathcal{G}_{p}(E/K)$, then $r=\mathfrak{I}(E)=\mathfrak{I}_k(E)\leq\left\lfloor j/2\right\rfloor $.
\end{corollary}

\begin{proof}
By Corollary \ref{cor:bloompotinv}, $r=\min\{\mathfrak{I}(E),\left\lfloor k/2\right\rfloor \}=\min\{\mathfrak{I}_{k}(E),\left\lfloor k/2\right\rfloor \}$. Since our assumptions are the same as those of Lemma~\ref{lem:kspineproperty} (3), we obtain that $r<\frac{k}{2}$ and $2r<j$. Consequently, $r=\mathfrak{I}(E)=\mathfrak{I}_{k}(E)\leq\left\lfloor j/2\right\rfloor $.
\end{proof}

While Corollary \ref{cor:bloompotinv} provides instructions for determining $r$, it does require knowledge of the image of the $\operatorname{mod}p^{k}$ Galois representation. In particular, suppose that $E$ is an elliptic curve with $\mathcal{G}_{p}(E/K)\cong\mathcal{H}_{p^{k}}^{r}$ such that the level $n$ of its Galois representation satisfies $v_{p}(n)<k$. Then, $\mathfrak{I}_{v_{p}(n)}(E)\leq\mathfrak{I}_{k}(E)$. Since $\mathfrak{I}_{v_{p}(n)}(E)\leq\frac {v_{p}(n)}{2}$, it is not necessarily the case that $r=\min\!\left\{\mathfrak{I}_{v_{p}(n)}(E),\left\lfloor k/2\right\rfloor \right\}  $. In this setting, if we want to recover $r$, we require the quantity $\alpha$ in (\ref{eq:alphavalue}). We note that unlike the blooming invariant, $\alpha$ does depend on choice of basis. The following result makes this precise.

\begin{lemma}\label{lem:rvalue}
Let $E$ be an elliptic curve over a field $K$ of characteristic $0$, and suppose that $\mathcal{G}_{p}(E/K)\cong\mathcal{H}_{p^{k}}^{r}$ for some nonnegative integers $k$ and $r$ with $r\leq\frac{k}{2}$. Suppose further that the degrees of the largest independent $K$-rational $p$ power isogenies admitted by $E$ are $p^{i}$ and $p^{j}$, respectively, with $j\leq i$. In particular, there is a basis $(P,Q)$ of $E[p^{i}]$ such that $\left\langle P\right\rangle $ and $\left\langle p^{i-j}Q\right\rangle $ are $G_{K}$-invariant of order $p^{i}$ and $p^{j}$, respectively. With respect to this basis, let
\[
\alpha=\min\!\left\{  v_{p}^{(i)}(a-c)\left\vert \left(
\begin{array}
[c]{cc}
a & b\\
0 & c
\end{array}
\right)  \in\rho_{E,p^{i}}(G_{K})\right.  \right\}  ,
\]
Then, $r=\min\!\left\{  \alpha,\left\lfloor \frac{i+j}{2}\right\rfloor \right\}
=\min\!\left\{  \mathfrak{I}_{k}(E),\left\lfloor \frac{i+j}{2}\right\rfloor
\right\}  $.
\end{lemma}

\begin{proof}
By Corollary~\ref{cor:bloompotinv}, we have that $r=\min\{\mathfrak{I}_{k}(E),\left\lfloor k/2\right\rfloor \}$. Now let $\mathcal{L}$ be a longest path in $\mathcal{G}_{p}(E/K)$ that contains $E$ of length $l$. Our assumptions imply that $i+j=l$. Observe that Lemma \ref{Lem:longpath} gives $l\leq k$. Next, observe that from our assumptions, we have a basis $(P,Q)$ of $E[p^{i}]$ such that $\rho_{E,p^{i}}(G_{K})\leq B_{0}(p^{i})$. We now claim that for each $\sigma\in G_{K}$, if $\rho_{E,p^{i}}(\sigma)=\sm{a & b  \\ 0 & c}$, then $v_{p}(b)\geq j$. Indeed, since $\left\langle p^{i-j}Q\right\rangle $ is $G_{K}$-invariant of order $p^{i}$ we must have that
\[
\sigma(p^{i-j}Q)=p^{i-j}bP+p^{i-j}cQ\in<p^{i-j}Q>\qquad\Longrightarrow\qquad
p^{i-j}bP=\mathcal{O},
\]
and so $v_{p}(b)\geq j$, which establishes the claim. Now let $\pi:\operatorname*{GL}\nolimits_{2}(\mathbb{Z}/p^{k}\mathbb{Z})\rightarrow\operatorname*{GL}\nolimits_{2}(\mathbb{Z}/p^{i}\mathbb{Z})$ be the natural projection and set
\[
H=\pi^{-1}(\rho_{E,p^{i}}(G_{K})).
\]
Now choose a basis $(\widetilde{P},\widetilde{Q})$ of $E[p^{k}]$ satisfying $p^{k-i}\widetilde{P}=P$ and $p^{k-i}\widetilde{Q}=Q$. Then, with respect to the basis $(\widetilde{P},\widetilde{Q})$, every lift of $\sm{a & b  \\ 0 & c}\in\rho_{E,p^{i}}(G_{K})$ has the form
\begin{equation}
\left(
\begin{array}
[c]{cc}
\widetilde{a} & \widetilde{b}\\
wp^{i} & \widetilde{c}
\end{array}
\right)  \in H\qquad\text{where\ }\left(
\begin{array}
[c]{cc}
\widetilde{a} & \widetilde{b}\\
0 & \widetilde{c}
\end{array}
\right)  \equiv\left(
\begin{array}
[c]{cc}
a & b\\
0 & c
\end{array}
\right)  \ \operatorname{mod}p^{i}.\label{eq:matrixlift}
\end{equation}
Now fix $M=\sm{a & b  \\ 0 & c}\in\rho_{E,p^{i}}(G_{K})$ with lift $\widetilde{M}$ as given in (\ref{eq:matrixlift}). Consequently, $v_{p}(\widetilde{b})\geq j$ and $v_{p}^{(k)}(\widetilde{a}-\widetilde{c})\geq v_{p}^{(i)}(a-c)$. In the latter, equality holds if $v_{p}^{(i)}(\widetilde{a}-\widetilde{c})<i$. Otherwise, $v_{p}^{(i)}(a-c)=i\leq v_{p}^{(k)}(\widetilde{a}-\widetilde{c})$. Now set
\[
\alpha^{\prime}=\min\!\left\{  v_{p}^{(k)}(\widetilde{a}-\widetilde{c})\left\vert \left(
\begin{array}
[c]{cc}
\widetilde{a} & \widetilde{b}\\
wp^{i} & \widetilde{c}
\end{array}
\right)  \in\rho_{E,p^{k}}(G_{K})\right.  \right\}  ,
\]
and observe that $\alpha\leq\alpha^{\prime}$, and equality holds if $\alpha^{\prime}<i$. Otherwise we must have that $\alpha=i\leq\alpha^{\prime}$. We now break into cases.

\textbf{Case 1.} Suppose $l=k$. Then, with notation as in (\ref{eq:matrixlift}), we have that $v_{p}^{(k)}(w\widetilde{b}p^{i})=k$. Now observe that
\begin{align*}
\mathfrak{I}_{k}(E)  & =\min_{\sigma\in G_{K}}\!\left\{  \frac{1}{2}v_{p}^{(k)}\!\left(  \operatorname*{Tr}(\rho_{E_{i},p^{k}}(\sigma))^{2}-4\det(\rho_{E_{i},p^{k}}(\sigma))\right)  \right\}  \\
& =\min\!\left\{  \frac{1}{2}v_{p}^{(k)}\!\left(  (\widetilde{a}-\widetilde{c})^{2}\right)  \left\vert \left(
\begin{array}
[c]{cc}
\widetilde{a} & \widetilde{b}\\
wp^{i} & \widetilde{c}
\end{array}
\right)  \in\rho_{E,p^{k}}(G_{K})\right.  \right\}  \\
& =\min\!\left\{  \min\!\left\{  v_{p}^{(k)}\!(\widetilde{a}-\widetilde{c}),\frac{k}{2}\right\}  \left\vert \left(
\begin{array}
[c]{cc}
\widetilde{a} & \widetilde{b}\\
wp^{i} & \widetilde{c}
\end{array}
\right)  \in\rho_{E,p^{k}}(G_{K})\right.  \right\}  .
\end{align*}
It follows that if $\mathfrak{I}_{k}(E)=\frac{k}{2}$, then $r=\left\lfloor k/2\right\rfloor $. But in this case we have that $\alpha^{\prime}\geq\frac{k}{2}$. Since $i\geq j$, we have that $i\geq\left\lfloor k/2\right\rfloor $. So either $\alpha=\alpha^{\prime}\geq\frac{k}{2}$ or $\left\lfloor k/2\right\rfloor \leq i=\alpha\leq\alpha^{\prime}$. In either case we have that $\min\!\left\{  \alpha,\left\lfloor k/2\right\rfloor \right\} =\left\lfloor k/2\right\rfloor =r$. So suppose that $\mathfrak{I}_{k}(E)<\frac{k}{2}$. It follows that $r=\mathfrak{I}_{k}(E)=\alpha^{\prime}$. Further, $\alpha^{\prime}<i$ and so $\alpha=\alpha^{\prime}$. So in both cases, we have that $\min\!\left\{  \alpha,\left\lfloor k/2\right\rfloor \right\}  =\mathfrak{I}_{k}(E)$.

\textbf{Case 2.} Suppose $l<k$. Then Lemma~\ref{lem:kspineproperty} implies that $2r<l$ and $E$ lies in the bloom of an internal skeletal vertex of $\mathcal{G}_{p}(E/K)$ relative to some $k$-spine. We further deduce from the lemma that $i$ is the length of the longest path from $E$ to a leaf on some $k$-spine and $j$ is the distance to a leaf in the aforementioned bloom. Since the bloom depth of a skeletal vertex is $r$, we obtain that $j<r$. Since $\mathcal{L}$ is a longest path containing $E$, and $i$ is the distance on this path from $E$ to a leaf on some $k$-spine, we have that $i>\left\lfloor k/2\right\rfloor $. We note that from Corollary~\ref{cor:longestpathbp} we now obtain that if $E$ does not lie on any $k$-spine of $\mathcal{G}_{p}(E/K)$, then
\begin{equation}
j<r=\mathfrak{I}_{k}(E)\leq\left\lfloor \frac{l}{2}\right\rfloor
.\label{valueokbloominvalpa}
\end{equation}
Now observe that with notation as in (\ref{eq:matrixlift}), we have that
$v_{p}^{(k)}(w\widetilde{b}p^{i})\geq l$. Now observe that
\begin{align*}
\mathfrak{I}_{k}(E)  & =\min_{\sigma\in G_{K}}\!\left\{  \frac{1}{2}v_{p}
^{(k)}\!\left(  \operatorname*{Tr}(\rho_{E_{i},p^{k}}(\sigma))^{2}-4\det
(\rho_{E_{i},p^{k}}(\sigma))\right)  \right\}  \\
& =\min\!\left\{  \frac{1}{2}v_{p}^{(k)}\!\left(  (\widetilde{a}-\widetilde{c}
)^{2}+4p^{i}\widetilde{b}w\right)  \left\vert \left(
\begin{array}
[c]{cc}
\widetilde{a} & \widetilde{b}\\
wp^{i} & \widetilde{c}
\end{array}
\right)  \in\rho_{E,p^{k}}(G_{K})\right.  \right\}  .
\end{align*}
In particular, if $r=\mathfrak{I}_{k}(E)=\left\lfloor l/2\right\rfloor $, we
must have that $\alpha^{\prime}\geq\left\lfloor l/2\right\rfloor $. Now
observe that  $\left\lfloor l/2\right\rfloor <i$, and so if $\alpha^{\prime
}<i$, then $\alpha\geq\left\lfloor l/2\right\rfloor $. If $\alpha^{\prime}\geq
i$, then $\alpha=i>\left\lfloor l/2\right\rfloor $ and so $\min\!\left\{
\alpha,\left\lfloor l/2\right\rfloor \right\}  =\left\lfloor l/2\right\rfloor
=r$.

So suppose that $r=\mathfrak{I}_{k}(E)<l/2$. Then,
\begin{align*}
\mathfrak{I}_{k}(E)  & =\min\!\left\{  \min\!\left\{  v_{p}^{(k)}\!(\widetilde{a}-\widetilde{c}),\frac{k}{2}\right\}  \left\vert \left(
\begin{array}
[c]{cc}
\widetilde{a} & \widetilde{b}\\
wp^{i} & \widetilde{c}
\end{array}
\right)  \in\rho_{E,p^{k}}(G_{K})\right.  \right\}  \\
& =\min\!\left\{  v_{p}^{(k)}\!(\widetilde{a}-\widetilde{c})\left\vert \left(
\begin{array}
[c]{cc}
\widetilde{a} & \widetilde{b}\\
wp^{i} & \widetilde{c}
\end{array}
\right)  \in\rho_{E,p^{k}}(G_{K})\right.  \right\}  \\
& =\alpha^{\prime}
\end{align*}
Since $\frac{l}{2}<i$, we have that $r=\mathfrak{I}_{k}(E)=\alpha=\min\!\left\{  \alpha,\left\lfloor l/2\right\rfloor \right\}  $, which concludes the proof.
\end{proof}

As an application of our results, we now show that elliptic curves defined over a field admitting a real embedding have smallest $p$-blooming invariant. Before stating the result, we note that in the next section, we complete the proof of Theorem \ref{mathmclass} for infinite graphs (see Theorem~\ref{thm:infinitegpks}). In the process, it is determined that $r=\mathfrak{I}_{p}(E/F)$. The statement below assumes this fact.

\begin{proposition}\label{Prop:BloomInvReal}
Let $E$ be an elliptic curve over a field $K$. If $K$ admits a real embedding, then $\operatorname*{End}_{K}E\cong\mathbb{Z}$ and the $p$-blooming invariant of $E$ is
\[
\mathfrak{I}_{p}(E/K)=\left\{
\begin{array}
[c]{cl}
1 & \text{if }p=2,\\
0 & \text{if }p\geq3.
\end{array}
\right.
\]
In particular, if $\deg\mathcal{G}_{p}(E/K)=p^{k}$ with $k\in\mathbb{Z}_{\geq0}\cup\{\infty\}$, then $\mathcal{G}_{p}(E/K)$ is isomorphic to either $\mathcal{H}_{p^{k}}^{r}$ or $\mathcal{H}_{p^{\infty,+}}^{r}$ with $r=\min\{\mathfrak{I}_{p}(E/K),\left\lfloor k/2\right\rfloor \}$.
\end{proposition}

\begin{proof}
We first show that $\operatorname*{End}_{K}E\cong\mathbb{Z}$. By way of contradiction, suppose this is not the case. Since $K$ has characteristic $0$, it follows that $\operatorname*{End}_{K}E\otimes_{\mathbb{Z}}\mathbb{Q}\cong L$ for some imaginary quadratic field $L$. By Theorem \cite[Theorem 2.3]{GrepsCM}, the field of definition of $\operatorname*{End}E=\operatorname*{End}_{K}E$ is $L(j(E))$. In particular, $L\subseteq L(j(E))\subseteq K$. But then, composing a real embedding of $K$ with the restriction of the inclusion $L\hookrightarrow K$ to $L$ yields a real embedding of $L$, which is
impossible. Thus, $\operatorname*{End}_{K}E\cong\mathbb{Z}$.

By assumption, we have a real embedding $\iota:K\hookrightarrow\mathbb{R}$. Next, we choose an embedding $\overline{\iota}:\overline{K}\hookrightarrow\mathbb{C}$ such that the following commutative diagram holds:
\[
\begin{tikzcd}
\overline{K} \arrow[r, "\overline{\iota}", hook] & \mathbb{C}                 \\
K \arrow[r, "\iota", hook] \arrow[u, hook]       & \mathbb{R} \arrow[u, hook]
\end{tikzcd}
\]
We consequently obtain an injection $\operatorname*{Gal}(\mathbb{C}/\mathbb{R})\hookrightarrow G_{K}$, which depends on the choice of $\overline{\iota}$. Thus, we get a complex conjugation $c\in G_{K}$ which is well defined up to conjugation. Now let $p$ be a prime and $k$ a positive integer. If $p=2$, suppose $k\geq2$. Then, $c(\zeta_{p^{k}})=\zeta_{p^{k}}^{-1}$, where $\zeta_{p^{k}}$ denotes a primitive $p^{k}$-th root of unity. Consequently, the $p^{k}$-th cyclotomic character $\chi_{p^{k}}:G_{K}\rightarrow(\mathbb{Z}/p^{k}\mathbb{Z})^{\times}$ satisfies $\chi_{p^{k}}(c)\equiv-1\ \operatorname{mod}p^{k}$.

Next, fix a basis of $E[p^{k}]$ and consider the $\operatorname{mod}p^{k}$ Galois representation $\rho_{E,p^{k}}:G_{K}\rightarrow\operatorname*{GL}\nolimits_{2}(\mathbb{Z}/p^{k}\mathbb{Z})$ attached to $E$. Since the cyclotomic action is independent of elliptic curves, we have that
\[
\det(\rho_{E,p^{k}}(c))=\chi_{p^{k}}(c)=-1.
\]
Moreover, since $c^{2}$ is the identity, we obtain that $\rho_{E,p^{k}}(c)^2\equiv\sm{1 & 0  \\ 0  & 1}\ \operatorname{mod}p^{k}$. It then follows that $\operatorname*{Tr}\rho_{E,p^{k}}(c)=0$ \cite[Lemma 6.5]{GrepsCM}. Therefore,
\[
\operatorname*{Tr}(\rho_{E,p^{k}}(c))^{2}-4\det\rho_{E,p^{k}}(c)=4.
\]
From Definition \ref{def:modpkbloominv}, we deduce that
\[
\mathfrak{I}_{p^{k}}(E/K)=\left\{
\begin{array}
[c]{cl}
1 & \text{if }p=2\text{ and }k\geq2,\\
0 & \text{if }p\geq3\text{ and }k\geq1.
\end{array}
\right.
\]
We then Corollary \ref{cor:bloompotprop} implies that
\[
\mathfrak{I}_{p}(E/K)=\lim_{k\rightarrow\infty}\mathfrak{I}_{p^{k}}(E/K)=\left\{
\begin{array}
[c]{cl}
1 & \text{if }p=2,\\
0 & \text{if }p\geq3.
\end{array}
\right.
\]
The proof now follows since $r=\min\{\mathfrak{I}_{p}(E/K),\left\lfloor k/2\right\rfloor \}$ by Corollary~\ref{cor:bloompotinv} and Theorem~\ref{thm:infinitegpks}.
\end{proof}

As consequence of Proposition \ref{Prop:BloomInvReal}, we obtain the following result which gives that the isogeny class degree uniquely determines the isogeny graphs for fields admitting a real embedding.

\begin{corollary}\label{Cor:BloomInvReal}
Let $K$ be a field admitting a real embedding. Then, the finite isogeny graphs $\mathcal{G}(E/K)$ that occur over $K$ are uniquely determined by the isogeny class degree. In particular, if $\deg\mathcal{G}(E/K)=\prod_{p}p^{k_{p}}$, then
\[
\mathcal{G}(E/K)\cong\left\{
\begin{array}
[c]{cl}
\underset{p}{\square}\mathcal{H}_{p^{k_{p}}}^{0} & \text{if }k_{2}\leq1,\\
\underset{p\neq2}{\square}\mathcal{H}_{p^{k_{p}}}^{0}\square\mathcal{H}
_{2^{k_{2}}}^{1} & \text{if }k_{2}\geq2.
\end{array}
\right.
\]
\end{corollary}

\begin{proof}
Let $E$ be an elliptic curve defined over $K$ such that $\mathcal{G}(E/K)$ is finite. Since $K$ admits a real embedding, Proposition~\ref{Prop:BloomInvReal} implies that $\operatorname*{End}_{K}E\cong\mathbb{Z}$. The result now follow loc. cit., together with Theorems~\ref{mainthmCarPro} and~\ref{mathmclass}~(1).
\end{proof}

We conclude this section with the following result, which showcases how a large blooming invariant measures the potential for an isogeny graph to further blossom over a field extension.

\begin{lemma}
Let $E$ be an elliptic curve over a field $K$ of characteristic $0$ with $\operatorname*{End}_{K}\! E\cong\mathbb{Z}$, and let $p$ be a prime such that $\mathcal{G}_p(E/K)$ is finite. Suppose further that $j$ is a nonnegative integer such that $\left\lfloor j/2\right\rfloor \leq\mathfrak{I}(E/K)$. If $F/K$ is a field extension such that the base change $E/F$ admits a $F$-rational $p^{j}$-isogeny, then
\[
\mathcal{H}_{p^{j}}^{\left\lfloor j/2\right\rfloor }\hookrightarrow \mathcal{G}_{p}(E/F).
\]

\end{lemma}

\begin{proof}
By Theorem \ref{classificationGpk} and Corollary \ref{cor:bloompotinv}, $\mathcal{G}_{p}(E/F)\cong\mathcal{H}_{p^{k}}^{r}$ for some positive integer $k$ and $r=\min\!\left\{  \mathfrak{I}(E/F),\left\lfloor k/2\right\rfloor \right\}$. Since $E$ admits a $F$-rational $p^{j}$-isogeny, we have that $j\leq k$. By Corollary \ref{cor:bloompotprop},
\[
\mathfrak{I}(E/K)\leq\mathfrak{I}(E/F).
\]
It then follows that $\left\lfloor j/2\right\rfloor \leq r$. From the discussion following Corollary~\ref{cor:bloomdepth}, we have that
\[
\mathcal{H}_{p^{j}}^{\left\lfloor j/2\right\rfloor }\hookrightarrow \mathcal{H}_{p^{k}}^{\left\lfloor j/2\right\rfloor }\hookrightarrow \mathcal{H}_{p^{k}}^{r}\cong\mathcal{G}_{p}(E/F). \qedhere
\]
\end{proof}


\section{Infinite \texorpdfstring{$p$}{p}-primary isogeny graphs}\label{sec:infgraphs}
In this section, we consider the case when the $p$-primary isogeny graphs are infinite. The layout of this section is similar to Section \ref{sec:gps}. Namely, we first introduce the infinite graphs $\mathcal{H}_{p^{\infty}}^{r}$ and $\mathcal{H}_{p^{\infty,+}}^{r}$, and then proceed with a careful study of the $p$-adic Galois representations attached to the elliptic curves to establish the main result of this section, Theorem \ref{thm:infinitegpks}. As a consequence, we obtain Theorem \ref{mathmclass} (2), and the converse in the cases corresponding to the graphs $\mathcal{H}_{p^{\infty}}^{r}$ and $\mathcal{H}_{p^{\infty,+}}^{r}$. We now begin by introducing the infinite graphs $\mathcal{H}_{p^{\infty}}^{r}$ and $\mathcal{H}_{p^{\infty,+}}^{r}$.

\subsection{The families of infinite graphs \texorpdfstring{$\mathcal H_{p^\infty}^r$}{Hrpinfinity} and \texorpdfstring{$\mathcal H_{p^{\infty,+}}^r$}{Hrpinfinity,+}}\label{subsec:infinitegraphs}
The construction of the finite family of graphs $\mathcal{H}_{p^{k}}^{r}$ proceeded by applying a finite sequence of $p$-blossomings to a path graph. The infinite analogues are built in the same spirit, but now the starting graph is either a line or a ray. By Kőnig' infinity lemma (Lemma~\ref{koenig}), every infinite tree contains a ray as a subgraph. In our setting, this distinction separate the two infinite families: the graphs $\mathcal{H}_{p^{\infty}}^{r}$ contain a line as a subgraph, while the graphs $\mathcal{H}_{p^{\infty,+}}^{r}$ do not contain a line, but do contain a ray as a subgraph. We now begin the construction of these infinite families.

\begin{definition}
Let $p$ be a prime number and let $r\in\mathbb{Z}_{\geq0}\cup\left\{  \infty\right\}  $. We define the infinite graphs $\mathcal{H}_{p^{\infty}}^{r}$ as follows.

\begin{enumerate}
\item For $r=0$, let $\mathcal{H}_{p^{\infty}}^{0}$ be the line graph with vertex set $\left\{  v_{i}\right\}_{i \in\mathbb{Z}}  $ and edges $\left\{  v_{i},v_{i+1}\right\}  $ for each $i\in \mathbb{Z}$. In addition, each edge is assigned weight~$p$.

\item For $0<r<\infty$, define $\mathcal{H}_{p^{\infty}}^{r}$ to be the $p$-blossom of $\mathcal{H}_{p^{\infty}}^{r-1}$, with each edge assigned weight~$p$.

\item For $r=\infty$, define
\[\mathcal{H}_{p^{\infty}}^{\infty}:=\bigcup_{r\geq0}\mathcal{H}_{p^{\infty}}^{r},
\]
where the union is taken over the increasing sequence of graphs
\[
\mathcal{H}_{p^{\infty}}^{0}\subseteq\mathcal{H}_{p^{\infty}}^{1}\subseteq\mathcal{H}_{p^{\infty}}^{2}\subseteq\cdots.
\]

\end{enumerate}
\noindent The vertices $\left\{  v_{i}\right\}_{i\in\mathbb{Z}}  $ of the initial line graph are called the \textit{skeletal vertices} of $\mathcal{H}_{p^{\infty}}^{r}$, and the set of all skeletal vertices is called the \textit{skeleton} of $\mathcal{H}_{p^{\infty}}^{r}$.
\end{definition}

\begin{example}\label{Ex:Hpinfinitygraphs}
We now demonstrate the construction of $\mathcal{H}_{3^{\infty}}^{2}$. We begin with the line graph with vertex set $\left\{  v_{i}\right\}  _{i\in\mathbb{Z}}$ and edges $\left\{  v_{i},v_{i+1}\right\}  $ for each $i\in\mathbb{Z}$. This is the skeleton of $\mathcal{H}_{3^{\infty}}^{2}$, and it is also $\mathcal{H}_{3^{\infty}}^{0}$ after assigning weight $3$ to each edge. To obtain $\mathcal{H}_{3^{\infty}}^{2}$, we apply two successive $3$-blossomings to $\mathcal{H}_{3^{\infty}}^{0}$. The first $3$-blossoming produces the bi-infinite caterpillar graph $\mathcal{H}_{3^{\infty}}^{1}$. The graphs $\mathcal{H}_{3^{\infty}}^{0}$ and $\mathcal{H}_{3^{\infty}}^{1}$ are shown below, and the vertices and edges of $\mathcal{H}_{3^{\infty}}^{0}$ are drawn in blue, while the new vertices and edges obtained in the $3$-blossoming are drawn in red.
\[
\begin{tikzpicture}[scale=.9, transform shape]

    \draw[line width=1.5pt, style={draw=blue}, dash pattern=on 5pt off 3pt] (-1.45,0) -- (-.25,0);
    \Vertex[shape=circle, size=0.6, color=white, style={draw=blue},label =$v_{-2}$, x=0]{m2}
    \Vertex[shape=circle, size=0.6, color=white, style={draw=blue},label =$v_{-1}$, x=1]{m1}
    \Vertex[shape=circle, size=0.6, color=white, style={draw=blue},label =$v_0$, x=2]{0}
    \Vertex[shape=circle, size=0.6, color=white, style={draw=blue},label =$v_1$, x=3]{1}
    \Vertex[shape=circle, size=0.6, color=white, style={draw=blue},label =$v_2$, x=4]{2}
    \draw[line width=1.5pt, style={draw=blue}, dash pattern=on 5pt off 3pt] (4.3,0) -- (5.5,0);

    \Edge[style={color=blue}](m2)(m1)
    \Edge[style={color=blue}](m1)(0)
    \Edge[style={color=blue}](0)(1)
    \Edge[style={color=blue}](1)(2)

        \node at (2, -1.75) { $\mathcal{H}_{3^\infty}^0$};

\end{tikzpicture}
\qquad \qquad
\begin{tikzpicture}[scale=.9, transform shape]

    \draw[line width=1.5pt, style={draw=blue}, dash pattern=on 5pt off 3pt] (-1.45,0) -- (-.25,0);
    \Vertex[shape=circle, size=0.6, color=white, style={draw=blue},label =$v_{-2}$, x=0]{m2}
    \Vertex[shape=circle, size=0.6, color=white, style={draw=blue},label =$v_{-1}$, x=1]{m1}
    \Vertex[shape=circle, size=0.6, color=white, style={draw=blue},label =$v_0$, x=2]{0}
    \Vertex[shape=circle, size=0.6, color=white, style={draw=blue},label =$v_1$, x=3]{1}
    \Vertex[shape=circle, size=0.6, color=white, style={draw=blue},label =$v_2$, x=4]{2}
    \draw[line width=1.5pt, style={draw=blue}, dash pattern=on 5pt off 3pt] (4.3,0) -- (5.5,0);

    \Vertex[shape=circle, size=0.5,  color=white, style={draw=red}, x=0, y=1]{m21}
    \Vertex[shape=circle, size=0.5,  color=white, style={draw=red}, x=0, y=-1]{m22}
    \Vertex[shape=circle, size=0.5,  color=white, style={draw=red}, x=1, y=1]{m11}
    \Vertex[shape=circle, size=0.5,  color=white, style={draw=red}, x=1, y=-1]{m12}
    \Vertex[shape=circle, size=0.5,  color=white, style={draw=red}, x=2, y=1]{01}
    \Vertex[shape=circle, size=0.5,  color=white, style={draw=red}, x=2, y=-1]{02}
    \Vertex[shape=circle, size=0.5,  color=white, style={draw=red}, x=3, y=1]{11}
    \Vertex[shape=circle, size=0.5,  color=white, style={draw=red}, x=3, y=-1]{12}
    \Vertex[shape=circle, size=0.5,  color=white, style={draw=red}, x=4, y=1]{21}
    \Vertex[shape=circle, size=0.5,  color=white, style={draw=red}, x=4, y=-1]{22}

    \Edge[style={color=blue}](m2)(m1)
    \Edge[style={color=blue}](m1)(0)
    \Edge[style={color=blue}](0)(1)
    \Edge[style={color=blue}](1)(2)
    \Edge[style={color=red}](m2)(m21)
    \Edge[style={color=red}](m2)(m22)
    \Edge[style={color=red}](m1)(m11)
    \Edge[style={color=red}](m1)(m12)
    \Edge[style={color=red}](0)(01)
    \Edge[style={color=red}](0)(02)
    \Edge[style={color=red}](1)(11)
    \Edge[style={color=red}](1)(12)
    \Edge[style={color=red}](2)(21)
    \Edge[style={color=red}](2)(22)

        \node at (2, -1.75) { $\mathcal{H}_{3^\infty}^1$};

\end{tikzpicture}
\]
Observe that the skeletal vertices are $3$-bloomed in $\mathcal{H}_{3^{\infty}}^{1}$. In particular, relative to $\mathcal{H}_{3^{\infty}}^{0}$, each skeletal vertex has bloom depth $1$ in $\mathcal{H}_{3^{\infty}}^{1}$. Applying the second $3$-blossoming, or equivalently, $3$-blossoming $\mathcal{H}_{3^{\infty}}^{1}$, we obtain $\mathcal{H}_{3^{\infty}}^{2}$, which is shown below. We keep the same color scheme to distinguish $\mathcal{H}_{3^{\infty}}^{0}$ and $\mathcal{H}_{3^{\infty}}^{1}$ as subgraphs of $\mathcal{H}_{3^{\infty}}^{2}$, and draw the new edges and vertices in black.
\[
\begin{tikzpicture}[scale=.9, transform shape]

    \draw[line width=1.5pt, style={draw=blue}, dash pattern=on 5pt off 3pt] (-1.45,0) -- (-.25,0);
    \Vertex[shape=circle, size=0.6, color=white, style={draw=blue},label =$v_{-2}$, x=0]{m2}
    \Vertex[shape=circle, size=0.6, color=white, style={draw=blue},label =$v_{-1}$, x=2]{m1}
    \Vertex[shape=circle, size=0.6, color=white, style={draw=blue},label =$v_0$, x=4]{0}
    \Vertex[shape=circle, size=0.6, color=white, style={draw=blue},label =$v_1$, x=6]{1}
    \Vertex[shape=circle, size=0.6, color=white, style={draw=blue},label =$v_2$, x=8]{2}
    \draw[line width=1.5pt, style={draw=blue}, dash pattern=on 5pt off 3pt] (8.3,0) -- (9.5,0);

    \Vertex[shape=circle, size=0.5,  color=white, style={draw=red}, x=0, y=1]{m21}
    \Vertex[shape=circle, size=0.5,  color=white, style={draw=red}, x=0, y=-1]{m22}
    \Vertex[shape=circle, size=0.5,  color=white, style={draw=red}, x=2, y=1]{m11}
    \Vertex[shape=circle, size=0.5,  color=white, style={draw=red}, x=2, y=-1]{m12}
    \Vertex[shape=circle, size=0.5,  color=white, style={draw=red}, x=4, y=1]{01}
    \Vertex[shape=circle, size=0.5,  color=white, style={draw=red}, x=4, y=-1]{02}
    \Vertex[shape=circle, size=0.5,  color=white, style={draw=red}, x=6, y=1]{11}
    \Vertex[shape=circle, size=0.5,  color=white, style={draw=red}, x=6, y=-1]{12}
    \Vertex[shape=circle, size=0.5,  color=white, style={draw=red}, x=8, y=1]{21}
    \Vertex[shape=circle, size=0.5,  color=white, style={draw=red}, x=8, y=-1]{22}

    \Vertex[shape=circle, size=0.5,  color=white, x=-.65, y=2]{m211}
    \Vertex[shape=circle, size=0.5,  color=white, x=0, y=2]{m212}
    \Vertex[shape=circle, size=0.5,  color=white, x=.65, y=2]{m213}
    \Vertex[shape=circle, size=0.5,  color=white, x=-.65, y=-2]{m221}
    \Vertex[shape=circle, size=0.5,  color=white, x=0, y=-2]{m222}
    \Vertex[shape=circle, size=0.5,  color=white, x=.65, y=-2]{m223}

    \Vertex[shape=circle, size=0.5,  color=white, x=1.35, y=2]{m111}
    \Vertex[shape=circle, size=0.5,  color=white, x=2, y=2]{m112}
    \Vertex[shape=circle, size=0.5,  color=white, x=2.65, y=2]{m113}
    \Vertex[shape=circle, size=0.5,  color=white, x=1.35, y=-2]{m121}
    \Vertex[shape=circle, size=0.5,  color=white, x=2, y=-2]{m122}
    \Vertex[shape=circle, size=0.5,  color=white, x=2.65, y=-2]{m123}

    \Vertex[shape=circle, size=0.5,  color=white, x=3.35, y=2]{011}
    \Vertex[shape=circle, size=0.5,  color=white, x=4, y=2]{012}
    \Vertex[shape=circle, size=0.5,  color=white, x=4.65, y=2]{013}
    \Vertex[shape=circle, size=0.5,  color=white, x=3.35, y=-2]{021}
    \Vertex[shape=circle, size=0.5,  color=white, x=4, y=-2]{022}
    \Vertex[shape=circle, size=0.5,  color=white, x=4.65, y=-2]{023}

    \Vertex[shape=circle, size=0.5,  color=white, x=5.35, y=2]{111}
    \Vertex[shape=circle, size=0.5,  color=white, x=6, y=2]{112}
    \Vertex[shape=circle, size=0.5,  color=white, x=6.65, y=2]{113}
    \Vertex[shape=circle, size=0.5,  color=white, x=5.35, y=-2]{121}
    \Vertex[shape=circle, size=0.5,  color=white, x=6, y=-2]{122}
    \Vertex[shape=circle, size=0.5,  color=white, x=6.65, y=-2]{123}

    \Vertex[shape=circle, size=0.5,  color=white, x=7.35, y=2]{211}
    \Vertex[shape=circle, size=0.5,  color=white, x=8, y=2]{212}
    \Vertex[shape=circle, size=0.5,  color=white, x=8.65, y=2]{213}
    \Vertex[shape=circle, size=0.5,  color=white, x=7.35, y=-2]{221}
    \Vertex[shape=circle, size=0.5,  color=white, x=8, y=-2]{222}
    \Vertex[shape=circle, size=0.5,  color=white, x=8.65, y=-2]{223}

    \Edge[style={color=blue}](m2)(m1)
    \Edge[style={color=blue}](m1)(0)
    \Edge[style={color=blue}](0)(1)
    \Edge[style={color=blue}](1)(2)
    \Edge[style={color=red}](m2)(m21)
    \Edge[style={color=red}](m2)(m22)
    \Edge[style={color=red}](m1)(m11)
    \Edge[style={color=red}](m1)(m12)
    \Edge[style={color=red}](0)(01)
    \Edge[style={color=red}](0)(02)
    \Edge[style={color=red}](1)(11)
    \Edge[style={color=red}](1)(12)
    \Edge[style={color=red}](2)(21)
    \Edge[style={color=red}](2)(22)

    \Edge(m21)(m211)
    \Edge(m21)(m212)
    \Edge(m21)(m213)
    \Edge(m22)(m221)
    \Edge(m22)(m222)
    \Edge(m22)(m223)

    \Edge(m11)(m111)
    \Edge(m11)(m112)
    \Edge(m11)(m113)
    \Edge(m12)(m121)
    \Edge(m12)(m122)
    \Edge(m12)(m123)

    \Edge(01)(011)
    \Edge(01)(012)
    \Edge(01)(013)
    \Edge(02)(021)
    \Edge(02)(022)
    \Edge(02)(023)

    \Edge(11)(111)
    \Edge(11)(112)
    \Edge(11)(113)
    \Edge(12)(121)
    \Edge(12)(122)
    \Edge(12)(123)

    \Edge(21)(211)
    \Edge(21)(212)
    \Edge(21)(213)
    \Edge(22)(221)
    \Edge(22)(222)
    \Edge(22)(223)

        \node at (4, -2.75) { $\mathcal{H}_{3^\infty}^2$};

\end{tikzpicture}
\]
This example illustrates the recursive nature of the construction: the first step $3$-blooms each skeletal vertex, and each subsequent step adds one new layer to every existing leaf.
\end{example}

Example \ref{Ex:Hpinfinitygraphs} highlights the key feature of the graphs $\mathcal{H}_{p^{\infty}}^{r}$: the construction proceeds by $r$ successive $p$-blossomings, with the first step making each skeletal vertex being $p$-bloomed. Each subsequent step adds one new layer to every existing leaf by $p$-blossoming the current leaves. In particular, for $r<\infty$, each skeletal vertex will have bloom depth $r$. This is made precise in Lemma~\ref{Lem:linebloomdepth}.

\begin{figure}
\includegraphics[scale=0.60]{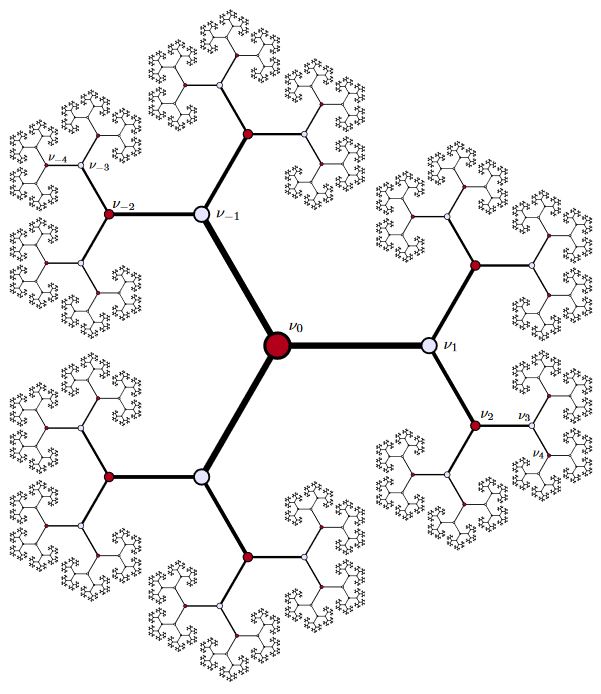}
\caption{The graph of $\mathcal{H}_{2^{\infty}}^\infty$. Figure reproduced from \cite{Casselman2019}.}
\label{fig:BruhatTits}
\end{figure}

Next, consider the increasing sequence $\mathcal{H}_{p^{\infty}}^{0}\subseteq\mathcal{H}_{p^{\infty}}^{1}\subseteq\mathcal{H}_{p^{\infty}}^{2}\subseteq\cdots$. This leads to
\[
\mathcal{H}_{p^{\infty}}^{\infty}=\bigcup_{r\geq0}\mathcal{H}_{p^{\infty}}^{r}.
\]
Necessarily, each skeletal vertex has bloom depth $\infty$. This graph is precisely the $p$-Bruhat-Tits tree, in which every vertex is $p$-bloomed. These graphs occur as the primary isogeny graph of non-CM elliptic curves over an algebraically closed field. See Figure~\ref{fig:BruhatTits} for an illustration of the graph $\mathcal{H}_{2^{\infty}}^{\infty}$.

We now define the infinite graphs $\mathcal{H}_{p^{\infty,+}}^{r}$.

\begin{definition}
Let $p$ be a prime number and let $r\in\mathbb{Z}_{\geq0}\cup\left\{  \infty\right\}  $. We define the infinite graphs $\mathcal{H}_{p^{\infty,+}}^{r}$ as follows:

\begin{enumerate}
\item For $r\in\mathbb{Z}_{\geq0}$, let $\mathcal{I}_{0}$ be the ray graph with vertex set $\left\{v_{i}\right\}_{ i\geq r}  $ and edges $\left\{  v_{i},v_{i+1}\right\}  $ for $i\geq r$. For $1\leq j\leq r$, let $\mathcal{I}_{j}$ be the $p$-blossom of $\mathcal{I}_{j-1}$. We define $\mathcal{H}_{p^{\infty,+}}^{r}$ to be the graph $\mathcal{I}_{r}$, with the additional property that each edge is assigned weight $p$.

\item For $r=\infty$, let $\mathcal{I}_{0}$ be the ray graph with vertex set $\left\{  v_{i}\right\}_{ i\geq0}  $ and edges $\left\{  v_{i},v_{i+1}\right\}  $ for $i\geq0$. For each $j\geq1$, define $\mathcal{I}_{j}$ to be the graph obtained from $\mathcal{I}_{j-1}$ by $p$-blossoming every vertex whose projection onto $\mathcal{I}_{0}$ lies in the tail $\left\{  v_{i}\right\}_{ i\geq j}  $. In this way, we obtain an increasing sequence of graphs $\mathcal{I}_{0}\subseteq\mathcal{I}_{1}\subseteq\mathcal{I}_{2}\subseteq\cdots$, and we define
\[
\mathcal{H}_{p^{\infty,+}}^{\infty}:=\bigcup_{j\geq0}\mathcal{I}_{j}.
\] 
\end{enumerate}
\noindent 
\end{definition}

\begin{remark}
The construction of $\mathcal H_{p^{\infty, +}}^\infty$ differs from our previous constructions in that the graph is not obtained from just considering a sequence of successive $p$-blossomings applied to its skeleton. For this reason, we take the convention that $\mathcal H_{p^{\infty, +}}^\infty$ does not have a skeleton.
\end{remark}

\begin{example}\label{Ex:Hprayinfinitygraphs}
We now demonstrate the construction of $\mathcal{H}_{3^{\infty,+}}^{2}$. We begin with a ray $\mathcal{I}_{0}$ with vertices $\left\{  v_{i}\right\}  _{i\geq2}$ and edges $\left\{  v_{i},v_{i+1}\right\}  $ for $i\geq2$. We then obtain $\mathcal{I}_{1}$ by $3$-blossoming $\mathcal{I}_{0}$. Consequently, each vertex $v_{i}$ of $\mathcal{I}_{0}$ becomes $3$-bloomed, and thus $v_{2}$ acquires three new adjacent vertices. In the graph below, we label one of these new adjacent vertices by $v_{1}$. The vertices and edges of $\mathcal{I}_{0}$ are drawn in blue, while the new vertices and edges obtained in the $3$-blossoming are drawn in red.
\[
\begin{tikzpicture}[scale=.9, transform shape]

    \Vertex[shape=circle, size=0.6, color=white, style={draw=blue},label =$v_{2}$, x=0]{m2}
    \Vertex[shape=circle, size=0.6, color=white, style={draw=blue},label =$v_{3}$, x=1]{m1}
    \Vertex[shape=circle, size=0.6, color=white, style={draw=blue},label =$v_4$, x=2]{0}
    \Vertex[shape=circle, size=0.6, color=white, style={draw=blue},label =$v_5$, x=3]{1}
    \Vertex[shape=circle, size=0.6, color=white, style={draw=blue},label =$v_6$, x=4]{2}
    \draw[line width=1.5pt, style={draw=blue}, dash pattern=on 5pt off 3pt] (4.3,0) -- (5.5,0);

    \Edge[style={color=blue}](m2)(m1)
    \Edge[style={color=blue}](m1)(0)
    \Edge[style={color=blue}](0)(1)
    \Edge[style={color=blue}](1)(2)

        \node at (2, -1.75) { $\mathcal{I}_{0}$};

\end{tikzpicture}
\qquad \qquad
\begin{tikzpicture}[scale=.9, transform shape]

    \Vertex[shape=circle, size=0.6, color=white, style={draw=red},label =$v_{1}$, x=0]{m2}
    \Vertex[shape=circle, size=0.6, color=white, style={draw=blue},label =$v_{2}$, x=1]{m1}
    \Vertex[shape=circle, size=0.6, color=white, style={draw=blue},label =$v_3$, x=2]{0}
    \Vertex[shape=circle, size=0.6, color=white, style={draw=blue},label =$v_4$, x=3]{1}
    \Vertex[shape=circle, size=0.6, color=white, style={draw=blue},label =$v_5$, x=4]{2}
    \Vertex[shape=circle, size=0.6, color=white, style={draw=blue},label =$v_6$, x=5]{3}
    \draw[line width=1.5pt, style={draw=blue}, dash pattern=on 5pt off 3pt] (5.3,0) -- (6.5,0);

    \Vertex[shape=circle, size=0.5,  color=white, style={draw=red}, x=1, y=1]{m11}
    \Vertex[shape=circle, size=0.5,  color=white, style={draw=red}, x=1, y=-1]{m12}
    \Vertex[shape=circle, size=0.5,  color=white, style={draw=red}, x=2, y=1]{01}
    \Vertex[shape=circle, size=0.5,  color=white, style={draw=red}, x=2, y=-1]{02}
    \Vertex[shape=circle, size=0.5,  color=white, style={draw=red}, x=3, y=1]{11}
    \Vertex[shape=circle, size=0.5,  color=white, style={draw=red}, x=3, y=-1]{12}
    \Vertex[shape=circle, size=0.5,  color=white, style={draw=red}, x=4, y=1]{21}
    \Vertex[shape=circle, size=0.5,  color=white, style={draw=red}, x=4, y=-1]{22}
    \Vertex[shape=circle, size=0.5,  color=white, style={draw=red}, x=5, y=1]{31}
    \Vertex[shape=circle, size=0.5,  color=white, style={draw=red}, x=5, y=-1]{32}

    \Edge[style={color=red}](m2)(m1)
    \Edge[style={color=blue}](m1)(0)
    \Edge[style={color=blue}](0)(1)
    \Edge[style={color=blue}](1)(2)
    \Edge[style={color=blue}](2)(3)
    \Edge[style={color=red}](m1)(m11)
    \Edge[style={color=red}](m1)(m12)
    \Edge[style={color=red}](0)(01)
    \Edge[style={color=red}](0)(02)
    \Edge[style={color=red}](1)(11)
    \Edge[style={color=red}](1)(12)
    \Edge[style={color=red}](2)(21)
    \Edge[style={color=red}](2)(22)
    \Edge[style={color=red}](3)(31)
    \Edge[style={color=red}](3)(32)

        \node at (3, -1.75) { $\mathcal{I}_{1}$};

\end{tikzpicture}
\]
Now consider the ray $\left\{  v_{i}\right\}  _{i\geq1}$ in $\mathcal{I}_{1}$. The vertices $v_{i}$ with $i\geq2$ have bloom depth $1$, while $v_{1}$ has bloom depth $0$ since it is a leaf. We then obtain $\mathcal{H}_{3^{\infty,+}}^{2}$ by performing another $3$-blossoming and assigning weight $3$ to each edge. In particular, $\mathcal{I}_{1}$ is a subgraph of $\mathcal{H}_{3^{\infty,+}}^{2}$, and the new vertices are leaves. Below, we keep the same color scheme to distinguish $\mathcal{I}_{0}\subseteq\mathcal{I}_{1}\subseteq\mathcal{H}_{3^{\infty,+}}^{2}$, with the new edges and vertices in $\mathcal{H}_{3^{\infty,+}}^{2}$ drawn in black.
\[
\begin{tikzpicture}[scale=.9, transform shape]

    \Vertex[shape=circle, size=0.6, color=white,label =$v_{0}$, x=0]{m2}
    \Vertex[shape=circle, size=0.6, color=white, style={draw=red},label =$v_{1}$, x=2]{m1}
    \Vertex[shape=circle, size=0.6, color=white, style={draw=blue},label =$v_2$, x=4]{0}
    \Vertex[shape=circle, size=0.6, color=white, style={draw=blue},label =$v_3$, x=6]{1}
    \Vertex[shape=circle, size=0.6, color=white, style={draw=blue},label =$v_4$, x=8]{2}
    \Vertex[shape=circle, size=0.6, color=white, style={draw=blue},label =$v_5$, x=10]{3}
    \Vertex[shape=circle, size=0.6, color=white, style={draw=blue},label =$v_6$, x=12]{4}
    \draw[line width=1.5pt, style={draw=blue}, dash pattern=on 5pt off 3pt] (12.3,0) -- (13.5,0);

    \Vertex[shape=circle, size=0.5,  color=white, x=2, y=1]{m11}
    \Vertex[shape=circle, size=0.5,  color=white, x=2, y=-1]{m12}
    \Vertex[shape=circle, size=0.5,  color=white, style={draw=red}, x=4, y=1]{01}
    \Vertex[shape=circle, size=0.5,  color=white, style={draw=red}, x=4, y=-1]{02}
    \Vertex[shape=circle, size=0.5,  color=white, style={draw=red}, x=6, y=1]{11}
    \Vertex[shape=circle, size=0.5,  color=white, style={draw=red}, x=6, y=-1]{12}
    \Vertex[shape=circle, size=0.5,  color=white, style={draw=red}, x=8, y=1]{21}
    \Vertex[shape=circle, size=0.5,  color=white, style={draw=red}, x=8, y=-1]{22}
    \Vertex[shape=circle, size=0.5,  color=white, style={draw=red}, x=10, y=1]{31}
    \Vertex[shape=circle, size=0.5,  color=white, style={draw=red}, x=10, y=-1]{32}
    \Vertex[shape=circle, size=0.5,  color=white, style={draw=red}, x=12, y=1]{41}
    \Vertex[shape=circle, size=0.5,  color=white, style={draw=red}, x=12, y=-1]{42}

    \Vertex[shape=circle, size=0.5,  color=white, x=3.35, y=2]{011}
    \Vertex[shape=circle, size=0.5,  color=white, x=4, y=2]{012}
    \Vertex[shape=circle, size=0.5,  color=white, x=4.65, y=2]{013}
    \Vertex[shape=circle, size=0.5,  color=white, x=3.35, y=-2]{021}
    \Vertex[shape=circle, size=0.5,  color=white, x=4, y=-2]{022}
    \Vertex[shape=circle, size=0.5,  color=white, x=4.65, y=-2]{023}

    \Vertex[shape=circle, size=0.5,  color=white, x=5.35, y=2]{111}
    \Vertex[shape=circle, size=0.5,  color=white, x=6, y=2]{112}
    \Vertex[shape=circle, size=0.5,  color=white, x=6.65, y=2]{113}
    \Vertex[shape=circle, size=0.5,  color=white, x=5.35, y=-2]{121}
    \Vertex[shape=circle, size=0.5,  color=white, x=6, y=-2]{122}
    \Vertex[shape=circle, size=0.5,  color=white, x=6.65, y=-2]{123}

    \Vertex[shape=circle, size=0.5,  color=white, x=7.35, y=2]{211}
    \Vertex[shape=circle, size=0.5,  color=white, x=8, y=2]{212}
    \Vertex[shape=circle, size=0.5,  color=white, x=8.65, y=2]{213}
    \Vertex[shape=circle, size=0.5,  color=white, x=7.35, y=-2]{221}
    \Vertex[shape=circle, size=0.5,  color=white, x=8, y=-2]{222}
    \Vertex[shape=circle, size=0.5,  color=white, x=8.65, y=-2]{223}

    \Vertex[shape=circle, size=0.5,  color=white, x=9.35, y=2]{311}
    \Vertex[shape=circle, size=0.5,  color=white, x=10, y=2]{312}
    \Vertex[shape=circle, size=0.5,  color=white, x=10.65, y=2]{313}
    \Vertex[shape=circle, size=0.5,  color=white, x=9.35, y=-2]{321}
    \Vertex[shape=circle, size=0.5,  color=white, x=10, y=-2]{322}
    \Vertex[shape=circle, size=0.5,  color=white, x=10.65, y=-2]{323}

    \Vertex[shape=circle, size=0.5,  color=white, x=11.35, y=2]{411}
    \Vertex[shape=circle, size=0.5,  color=white, x=12, y=2]{412}
    \Vertex[shape=circle, size=0.5,  color=white, x=12.65, y=2]{413}
    \Vertex[shape=circle, size=0.5,  color=white, x=11.35, y=-2]{421}
    \Vertex[shape=circle, size=0.5,  color=white, x=12, y=-2]{422}
    \Vertex[shape=circle, size=0.5,  color=white, x=12.65, y=-2]{423}

    \Edge(m2)(m1)
    \Edge[style={color=red}](m1)(0)
    \Edge[style={color=blue}](0)(1)
    \Edge[style={color=blue}](1)(2)
    \Edge[style={color=blue}](2)(3)
    \Edge[style={color=blue}](3)(4)
    \Edge(m1)(m11)
    \Edge(m1)(m12)
    \Edge[style={color=red}](0)(01)
    \Edge[style={color=red}](0)(02)
    \Edge[style={color=red}](1)(11)
    \Edge[style={color=red}](1)(12)
    \Edge[style={color=red}](2)(21)
    \Edge[style={color=red}](2)(22)
    \Edge[style={color=red}](3)(31)
    \Edge[style={color=red}](3)(32)
    \Edge[style={color=red}](4)(41)
    \Edge[style={color=red}](4)(42)

    \Edge(01)(011)
    \Edge(01)(012)
    \Edge(01)(013)
    \Edge(02)(021)
    \Edge(02)(022)
    \Edge(02)(023)

    \Edge(11)(111)
    \Edge(11)(112)
    \Edge(11)(113)
    \Edge(12)(121)
    \Edge(12)(122)
    \Edge(12)(123)

    \Edge(21)(211)
    \Edge(21)(212)
    \Edge(21)(213)
    \Edge(22)(221)
    \Edge(22)(222)
    \Edge(22)(223)

    \Edge(31)(311)
    \Edge(31)(312)
    \Edge(31)(313)
    \Edge(32)(321)
    \Edge(32)(322)
    \Edge(32)(323)

    \Edge(41)(411)
    \Edge(41)(412)
    \Edge(41)(413)
    \Edge(42)(421)
    \Edge(42)(422)
    \Edge(42)(423)

        \node at (8.25, -2.75) { $\mathcal{H}_{3^{\infty,+}}^2$};

\end{tikzpicture}
\]
Now observe that in $\mathcal{H}_{3^{\infty,+}}^{2}$, the vertex $v_{1}$ acquires three new edges, and we label one of these new adjacent vertices by $v_{0}$. The vertices $v_{i}$ with $i\geq2$ have bloom depth $2$, while $v_{1}$ and $v_{0}$ have bloom depths $1$ and $0$, respectively.
\end{example}

\begin{example}\label{Ex:Hprayinfinity2graphs}
We now consider the construction of $\mathcal{H}_{2^{\infty,+}}^{\infty}$. We begin with a ray $\mathcal{I}_{0}$ having vertices $\left\{  v_{i}\right\}  _{i\geq0}$ and edges $\left\{v_{i},v_{i+1}\right\}  $ for $i\geq0$. We construct $\mathcal{I}_{1}$ by $2$-blossoming every vertex $v_{i}$ with $i\geq1$. The graphs $\mathcal{I}_{0}$ and $\mathcal{I}_{1}$ are shown below, with the vertices and edges of $\mathcal{I}_{0}$ drawn in blue and the new vertices and edges in $\mathcal{I}_{1}$ drawn in red.
\[
\begin{tikzpicture}[scale=.9, transform shape]

    \Vertex[shape=circle, size=0.5, color=white, style={draw=blue},label =$v_{0}$, x=0]{0}
    \Vertex[shape=circle, size=0.5, color=white, style={draw=blue},label =$v_{1}$, x=1]{1}
    \Vertex[shape=circle, size=0.5, color=white, style={draw=blue},label =$v_2$, x=2]{2}
    \Vertex[shape=circle, size=0.5, color=white, style={draw=blue},label =$v_3$, x=3]{3}
    \Vertex[shape=circle, size=0.5, color=white, style={draw=blue},label =$v_4$, x=4]{4}
    \draw[line width=1.5pt, style={draw=blue}, dash pattern=on 5pt off 3pt] (4.25,0) -- (5.5,0);

    \Edge[style={color=blue}](0)(1)
    \Edge[style={color=blue}](1)(2)
    \Edge[style={color=blue}](2)(3)
    \Edge[style={color=blue}](3)(4)

        \node at (2, -.75) { $\mathcal{I}_{0}$};

\end{tikzpicture}
\qquad \qquad
\begin{tikzpicture}[scale=.9, transform shape]

    \Vertex[shape=circle, size=0.5, color=white, style={draw=blue},label =$v_{0}$, x=0]{0}
    \Vertex[shape=circle, size=0.5, color=white, style={draw=blue},label =$v_{1}$, x=1]{1}
    \Vertex[shape=circle, size=0.5, color=white, style={draw=blue},label =$v_2$, x=2]{2}
    \Vertex[shape=circle, size=0.5, color=white, style={draw=blue},label =$v_3$, x=3]{3}
    \Vertex[shape=circle, size=0.5, color=white, style={draw=blue},label =$v_4$, x=4]{4}
    \Vertex[shape=circle, size=0.5, color=white, style={draw=red}, x=1,y=1]{11}
    \Vertex[shape=circle, size=0.5, color=white, style={draw=red}, x=2,y=1]{21}
    \Vertex[shape=circle, size=0.5, color=white, style={draw=red}, x=3,y=1]{31}
    \Vertex[shape=circle, size=0.5, color=white, style={draw=red}, x=4,y=1]{41}
    \draw[line width=1.5pt, style={draw=blue}, dash pattern=on 5pt off 3pt] (4.25,0) -- (5.5,0);

    \Edge[style={color=blue}](0)(1)
    \Edge[style={color=blue}](1)(2)
    \Edge[style={color=blue}](2)(3)
    \Edge[style={color=blue}](3)(4)
    \Edge[style={color=red}](1)(11)
    \Edge[style={color=red}](2)(21)
    \Edge[style={color=red}](3)(31)
    \Edge[style={color=red}](4)(41)

        \node at (2, -.75) { $\mathcal{I}_{1}$};

\end{tikzpicture}
\]    
To obtain $\mathcal{I}_{2}$, we $2$-blossom every vertex of $\mathcal{I}_{1}$ whose projection onto $\mathcal{I}_{0}$ is a vertex $v_{i}$ with $i\geq2$. The graph $\mathcal{I}_{2}$ is given below, and we maintain the color scheme $\mathcal{I}_{0},\mathcal{I}_{1}\subseteq\mathcal{I}_{2}$, drawing the new vertices and edges in teal. We then construct $\mathcal{I}_{3}$ by $2$-blossoming every vertex of $\mathcal{I}_{2}$ whose projection onto $\mathcal{I}_{0}$ is a vertex $v_{i}$ with $i\geq3$. The graph $\mathcal{I}_{3}$ is shown below, maintaining the same color scheme to distinguish $\mathcal{I}_{0},\mathcal{I}_{1}$, and $\mathcal{I}_{2}$ as subgraphs, and the draw the new edges and vertices in black. Of note is that in the $\mathcal{I}_{j}$ considered, the vertex $v_{i}$ has bloom depth $\min\!\left\{  i,j\right\}  $.
\[
\begin{tikzpicture}[scale=.8, transform shape]

    \Vertex[shape=circle, size=0.5, color=white, style={draw=blue},label =$v_{0}$, x=0]{0}
    \Vertex[shape=circle, size=0.5, color=white, style={draw=blue},label =$v_{1}$, x=1.25]{1}
    \Vertex[shape=circle, size=0.5, color=white, style={draw=blue},label =$v_2$, x=2.5]{2}
    \Vertex[shape=circle, size=0.5, color=white, style={draw=blue},label =$v_3$, x=3.75]{3}
    \Vertex[shape=circle, size=0.5, color=white, style={draw=blue},label =$v_4$, x=5]{4}
    \Vertex[shape=circle, size=0.5, color=white, style={draw=red}, x=1.25,y=1]{11}
    \Vertex[shape=circle, size=0.5, color=white, style={draw=red}, x=2.5,y=1]{21}
    \Vertex[shape=circle, size=0.5, color=white, style={draw=red}, x=3.75,y=1]{31}
    \Vertex[shape=circle, size=0.5, color=white, style={draw=red}, x=5,y=1]{41}
    \draw[line width=1.5pt, style={draw=blue}, dash pattern=on 5pt off 3pt] (5.25,0) -- (6.5,0);
    \Vertex[shape=circle, size=0.5, color=white, style={draw=teal}, x=2.2,y=2]{211}
    \Vertex[shape=circle, size=0.5, color=white, style={draw=teal}, x=2.8,y=2]{212}
    \Vertex[shape=circle, size=0.5, color=white, style={draw=teal}, x=3.45,y=2]{311}
    \Vertex[shape=circle, size=0.5, color=white, style={draw=teal}, x=4.05,y=2]{312}
    \Vertex[shape=circle, size=0.5, color=white, style={draw=teal}, x=4.7,y=2]{411}
    \Vertex[shape=circle, size=0.5, color=white, style={draw=teal}, x=5.3,y=2]{412}

    \Edge[style={color=blue}](0)(1)
    \Edge[style={color=blue}](1)(2)
    \Edge[style={color=blue}](2)(3)
    \Edge[style={color=blue}](3)(4)
    \Edge[style={color=red}](1)(11)
    \Edge[style={color=red}](2)(21)
    \Edge[style={color=red}](3)(31)
    \Edge[style={color=red}](4)(41)
    \Edge[style={color=teal}](21)(211)
    \Edge[style={color=teal}](21)(212)
    \Edge[style={color=teal}](31)(311)
    \Edge[style={color=teal}](31)(312)
    \Edge[style={color=teal}](41)(411)
    \Edge[style={color=teal}](41)(412)

        \node at (2.5, -.75) { $\mathcal{I}_{2}$};

\end{tikzpicture}
\quad \quad
\begin{tikzpicture}[scale=.8, transform shape]

    \Vertex[shape=circle, size=0.5, color=white, style={draw=blue},label =$v_{0}$, x=0]{0}
    \Vertex[shape=circle, size=0.5, color=white, style={draw=blue},label =$v_{1}$, x=2]{1}
    \Vertex[shape=circle, size=0.5, color=white, style={draw=blue},label =$v_2$, x=4]{2}
    \Vertex[shape=circle, size=0.5, color=white, style={draw=blue},label =$v_3$, x=6]{3}
    \Vertex[shape=circle, size=0.5, color=white, style={draw=blue},label =$v_4$, x=8]{4}
    \Vertex[shape=circle, size=0.5, color=white, style={draw=red}, x=2,y=1]{11}
    \Vertex[shape=circle, size=0.5, color=white, style={draw=red}, x=4,y=1]{21}
    \Vertex[shape=circle, size=0.5, color=white, style={draw=red}, x=6,y=1]{31}
    \Vertex[shape=circle, size=0.5, color=white, style={draw=red}, x=8,y=1]{41}
    \draw[line width=1.5pt, style={draw=blue}, dash pattern=on 5pt off 3pt] (8.25,0) -- (9.5,0);
    \Vertex[shape=circle, size=0.5, color=white, style={draw=teal}, x=3.6,y=2]{211}
    \Vertex[shape=circle, size=0.5, color=white, style={draw=teal}, x=4.4,y=2]{212}
    
    \Vertex[shape=circle, size=0.5, color=white, style={draw=teal}, x=5.6,y=2]{311}
    \Vertex[shape=circle, size=0.5, color=white, style={draw=teal}, x=6.4,y=2]{312}
    \Vertex[shape=circle, size=0.5, color=white, style={draw=teal}, x=7.6,y=2]{411}
    \Vertex[shape=circle, size=0.5, color=white, style={draw=teal}, x=8.4,y=2]{412}

    \Vertex[shape=circle, size=0.4, color=white, x=5.1,y=3]{3111}
    \Vertex[shape=circle, size=0.4, color=white, x=5.6,y=3]{3112}
    \Vertex[shape=circle, size=0.4, color=white, x=6.1,y=3]{3121}
    \Vertex[shape=circle, size=0.4, color=white, x=6.6,y=3]{3122}

    \Vertex[shape=circle, size=0.4, color=white, x=7.1,y=3]{4111}
    \Vertex[shape=circle, size=0.4, color=white, x=7.6,y=3]{4112}
    \Vertex[shape=circle, size=0.4, color=white, x=8.1,y=3]{4121}
    \Vertex[shape=circle, size=0.4, color=white, x=8.6,y=3]{4122}

    \Edge[style={color=blue}](0)(1)
    \Edge[style={color=blue}](1)(2)
    \Edge[style={color=blue}](2)(3)
    \Edge[style={color=blue}](3)(4)
    \Edge[style={color=red}](1)(11)
    \Edge[style={color=red}](2)(21)
    \Edge[style={color=red}](3)(31)
    \Edge[style={color=red}](4)(41)
    \Edge[style={color=teal}](21)(211)
    \Edge[style={color=teal}](21)(212)
    \Edge[style={color=teal}](31)(311)
    \Edge[style={color=teal}](31)(312)
    \Edge[style={color=teal}](41)(411)
    \Edge[style={color=teal}](41)(412)
    \Edge(311)(3111)
    \Edge(311)(3112)
    \Edge(312)(3121)
    \Edge(312)(3122)
    \Edge(411)(4111)
    \Edge(411)(4112)
    \Edge(412)(4121)
    \Edge(412)(4122)

        \node at (4, -.75) { $\mathcal{I}_{3}$};

\end{tikzpicture}
\]
Continuing inductively, we obtain an increasing sequence of graphs $\mathcal{I}_{0}\subseteq\mathcal{I}_{1}\subseteq\mathcal{I}_{2}\subseteq\cdots$, and define
\[
\mathcal{H}_{2^{\infty,+}}^{\infty}=\bigcup_{j\geq0}\mathcal{I}_{j}.
\]
The construction shows that the vertex $v_{i}$ stabilizes after the $i$-th stage: once the $i$-step is reached, no further branching occurs above $v_{i}$. In other words, the bloom $\mathcal{B}(v_{i})$ of $v_{i}$ in $\mathcal{H}_{2^{\infty,+}}^{\infty}$ is contained in $\mathcal{I}_{t}$ for each $t\geq i$. Consequently, the bloom depth of $v_{i}$ is $i$.
\end{example}

Examples \ref{Ex:Hprayinfinitygraphs} and \ref{Ex:Hprayinfinity2graphs} show the pattern that will later be formalized: the graphs $\mathcal{H}_{p^{\infty,+}}^{r}$ contain a ray $\left\{  v_{i}\right\}  _{i\geq0}$ such that relative to this ray, the vertex $v_{i}$ has bloom depth $i$ in $\mathcal{H}_{p^{\infty,+}}^{r}$. This is made precise in Lemma \ref{Lem:raybloomdepth}.

The next results establish the bloom depth profiles and corresponding uniqueness statements for the graphs $\mathcal{H}_{p^{\infty}}^{r}$ and $\mathcal{H}_{p^{\infty,+}}^{r}$. In particular, these results provide constructive descriptions of the infinite graphs in terms of a distinguished line or ray, together with the associated bloom depths, without referring back to the recursive sequences used in their construction.

\begin{lemma}
    If $r=0$, then every vertex of $\mathcal H_{p^\infty}^r$  has bloom depth $0$. If $r>0$, the bloom depth of a vertex $v$ in $\mathcal H_{p^\infty}^r$ equals the minimum distance from $v$ to a leaf. The same holds for the family of graphs $\mathcal H_{p^{\infty,+}}^r$.
\end{lemma}
\begin{proof}
    The proof is exactly the same as the proof of Lemma~\ref{lem:bloomdepthdescription}. If $r=0$, then there are no $p$-bloomed vertices, so every vertex has bloom depth $0$. If $r>0$, every vertex that is not a leaf is $p$-bloomed, and the claim follows.
\end{proof}

\begin{lemma}\label{Lem:linebloomdepth}
Let $\mathcal{S}$ be the skeleton of $\mathcal{H}_{p^{\infty}}^{r}$, where $r\in\mathbb{Z}_{\geq0}\cup\left\{  \infty\right\}  $, with vertices $\left\{  v_{i}\right\}_{i\in\mathbb{Z}}$ and edges $\left\{  v_{i},v_{i+1}\right\}  $ for each $i\in\mathbb{Z}$. Then, each skeletal vertex has $p$-bloom depth $r$ and all vertices in the bloom $\mathcal B(v_i)$ at a distance $d_i$ from $v_i$ are leaves.  

Conversely, if $\mathcal{I}$ is a line graph with vertices $\left\{  v_{i}\right\}_{i\in\mathbb{Z}}$ and edges $\left\{  v_{i},v_{i+1}\right\}  $ for each $i\in\mathbb{Z}$, then $\mathcal{H}_{p^{\infty}}^{r}$ is, up to isomorphism, the unique tree containing $\mathcal{I}$ such that each vertex $v_{i}\in V(\mathcal{I})$ has bloom depth $r$ and each vertex in the bloom $\mathcal B(v_i)$ at a distance $r$ from $v_i$ is a leaf, with each edge assigned weight $p$.
\end{lemma}

\begin{proof}
If $r=0$, then $\mathcal{S}=\mathcal{H}_{p^{\infty}}^{0}$ is a line graph, and so each skeletal vertex has bloom depth $0$. Now suppose that $0<r<\infty$. By definition, $\mathcal{H}_{p^{\infty}}^{r}$ is obtained from the line graph $\mathcal{H}_{p^{\infty}}^{0}=\mathcal{S}$ by applying $r$ successive $p$-blossomings. Each blossoming step adds one new layer of leaves by $p$-blooming the leaves of the previous layer. Consequently, after the $r$-th step, every skeletal vertex $v\in V(\mathcal{S})$ has exactly $r$ layers in its bloom $\mathcal{B}(v)$. In particular, if $x\in V(\mathcal{B}(v))$, then $\operatorname*{dist}(v,x)\leq r$, and if equality holds, then $x$ is a leaf. Otherwise, $x$ is $p$-bloomed. Thus each skeletal vertex has $p$-bloom depth $r$.

Now suppose $r=\infty$. By definition,
\[
\mathcal{H}_{p^{\infty}}^{\infty}=\bigcup_{r\geq0}\mathcal{H}_{p^{\infty}}^{r},
\]
where the union is taken over the increasing sequence of graphs
\[
\mathcal{S}=\mathcal{H}_{p^{\infty}}^{0}\subseteq\mathcal{H}_{p^{\infty}}^{1}\subseteq\mathcal{H}_{p^{\infty}}^{2}\subseteq\cdots
\]
Now fix a skeletal vertex $v\in V(\mathcal{S})$, and consider its bloom $\mathcal{B}(v)$. If $x\in V(\mathcal{B}(v))$ and $\operatorname*{dist}(v,x)=j$, then $x$ is a leaf in $\mathcal{H}_{p^{\infty}}^{j}$, and $x$ is $p$-bloomed in $\mathcal{H}_{p^{\infty}}^{t}$ for each $t>j$. Since this holds for each $j$, the bloom $\mathcal{B}(v)$ contains no leaves. Consequently, every vertex of $\mathcal{B}(v)$ is $p$-bloomed, and so the bloom depth of $v$ is $\infty$. The converse follows from the above and Corollary~\ref{cor:isographsblooms}.
\end{proof}

\begin{lemma}
\label{Lem:raybloomdepth}Let $\mathcal{S}$ be a ray in $\mathcal{H}_{p^{\infty,+}}^{r}$, where $r\in\mathbb{Z}_{\geq0}\cup\left\{  \infty\right\}  $, with vertices $\left\{  v_{i}\right\}_{i\geq0}$ and edges $\left\{  v_{i},v_{i+1}\right\}  $ for each $i\geq0$. Suppose further that $\mathcal{S}$ contains the skeleton of $\mathcal{H}_{p^{\infty,+}}^{r}$, and that $v_{0}$ is a leaf of $\mathcal{H}_{p^{\infty,+}}^{r}$. Then, the $p$-bloom depth of $v_{i}$ is $d_{i}=\min\!\left\{i,r\right\}  $ and all vertices in the bloom $\mathcal B(v_i)$ at a distance $d_i$ from $v_i$ are leaves.  

Conversely, let $\mathcal{I}$ be a ray graph with vertices $\left\{  v_{i}\right\}_{i\geq0}$ and edges $\left\{  v_{i},v_{i+1}\right\}  $ for each $i\geq0$, and set $d_{i}=\min\{r,i\}$. Then $\mathcal{H}_{p^{\infty,+}}^{r}$ is, up to isomorphism, the unique tree containing $\mathcal{I}$ such that each vertex $v_{i}\in V(\mathcal{I})$ has bloom depth $d_{i}$ and each vertex in the bloom $\mathcal B(v_i)$ at a distance $d_i$ from $v_i$ is a leaf, with each edge assigned weight $p$.
\end{lemma}

\begin{proof}
First suppose that $r<\infty$. By construction, $\mathcal{H}_{p^{\infty,+}}^{r}$ is obtained from a ray by applying $r$ successive $p$-blossomings. Since $\mathcal{S}$ contains the skeleton of $\mathcal{H}_{p^{\infty,+}}^{r}$ and $v_{0}$ is a leaf of $\mathcal{H}_{p^{\infty,+}}^{r}$, it follows that $\operatorname*{dist}_{\mathcal{S}}(v_{0})=r$. In particular, the vertices $v_{0},v_{1},\ldots,v_{r-1}$ are the non-skeletal vertices of $\mathcal{S}$, while the vertices $v_{i}$ with $i\geq r$ are skeletal vertices. 

Now suppose that $i\geq r$. Since $\mathcal{H}_{p^{\infty,+}}^{r}$ is obtained by $r$ successive $p$-blossomings, every skeletal vertex acquires exactly $r$ layers in its bloom. Hence, the bloom depth of $v_{i}$ is $r$.

Next, suppose $i<r$. Then $v_{i}$ first appears as a leaf after the first $r-i$ blossoming steps, and the remaining $i$ blossoming steps add one layer to its bloom at each stage. Consequently, the bloom $\mathcal{B}(v_{i})$ has the property that every vertex $x\in V(\mathcal{B}(v_{i}))$ with $\operatorname*{dist}(v_{i},x)<i$ is $p$-bloomed, while every vertex $x\in V(\mathcal{B}(v_{i}))$ with $\operatorname*{dist}(v_{i},x)=i$ is a leaf. Thus, the bloom depth of $v_{i}$ is $i$. In sum, for every $v_{i}\in V(\mathcal{S})$, the bloom depth of $v_{i}$ is $d_{i}=\min\!\left\{  i,r\right\}  $.

Now suppose that $r=\infty$, and let $\mathcal{I}_{0}\subseteq\mathcal{I}_{1}\subseteq\mathcal{I}_{2}\subseteq\cdots$ be the increasing sequence of graphs in the construction of $\mathcal{H}_{p^{\infty,+}}^{\infty}$. Since $\mathcal{S}$ is a ray containing the ray $\mathcal{I}_{0}$, and $\mathcal{I}_{0}$ has exactly one leaf, it follows that $\mathcal{S} =\mathcal{I}_{0}$. By construction, $\mathcal{I}_{j}$ is obtained from $\mathcal{I}_{j-1}$ by $p$-blossoming every vertex whose projection onto $\mathcal{S}$ lies in the tail $\left\{  v_{i}\right\}  _{i\geq j}$. In particular, the vertex $v_{i}$ acquires one additional bloom layer with each blossoming, and there are exactly $i$ successive blossoming steps before $v_{i}$ becomes part of the stable tail. Hence the bloom $\mathcal{B}(v_{i})$ is contained in $\mathcal{I}_{i}$, and every vertex $x\in V(\mathcal{B}(v_{i}))$ satisfies $\operatorname*{dist}(v_{i},x)\leq i$. Moreover, if equality holds, then $x$ is a leaf, while otherwise, $x$ is $p$-bloomed. It follows that the bloom depth of $v_{i}$ is $i=\min\!\left\{  i,\infty\right\}  $. Now observe that by the above, and Corollary~\ref{cor:isographsblooms}, we obtain the converse statement.
\end{proof}

\begin{notation}
Since the infinite graphs $\mathcal{H}_{p^{\infty}}^{r}$ and $\mathcal{H}_{p^{\infty,+}}^{r}$ are determined by a distinguished line or ray together with the associated bloom depths, we introduce a compact notation for depicting these graphs. The essence of this notation is the same as in the finite case case: the distinguished line or ray is drawn as the underlying spine, and each vertex is labeled by its bloom depth. In particular, the bloom depth is indicated by a star-shaped marker, as shown in the figures below. For a line, all vertices are displayed in order along the bi-infinite spine; for a ray, the initial vertex is drawn as a leaf and the remaining vertices are labeled according to their bloom depths.
\[
\begin{tikzpicture}[scale=1.2, transform shape]
\draw[line width=1.5pt, dash pattern=on 5pt off 3pt] (-1.2,0) -- (0,0);
    \Vertex[shape=star, size=0.5, color=white, x=0, label=r]{1}
    \Vertex[shape=star, size=0.5, color=white, x=1, label=r]{2}
    \Vertex[shape=star, size=0.5, color=white, x=2, label=r]{3}
    \Vertex[shape=star, size=0.5, color=white, x=3, label=r]{4}
    \Vertex[shape=star, size=0.5, color=white, x=4, label=r]{5}
    \Vertex[shape=star, size=0.5, color=white, x=5, label=r]{6}
    \Vertex[shape=star, size=0.5, color=white, x=6, label=r]{7}

\draw[line width=1.5pt, dash pattern=on 5pt off 3pt] (6.2,0) -- (7.2,0);
    \Edge(1)(2)
    \Edge(2)(3)
    \Edge(3)(4)
    \Edge(4)(5)
    \Edge(5)(6)
    \Edge(6)(7)

      \node at (-2, 0) { $\mathcal{H}_{p^\infty}^r$};
\end{tikzpicture}
\]
\vspace{-.5em}
\[
\begin{tikzpicture}[scale=1.2, transform shape]
\draw[line width=1.5pt, dash pattern=on 5pt off 3pt] (-1.2,0) -- (0,0);
    \Vertex[shape=star, size=0.5, color=white, x=0, label=$\infty$]{1}
    \Vertex[shape=star, size=0.5, color=white, x=1, label=$\infty$]{2}
    \Vertex[shape=star, size=0.5, color=white, x=2, label=$\infty$]{3}
    \Vertex[shape=star, size=0.5, color=white, x=3, label=$\infty$]{4}
    \Vertex[shape=star, size=0.5, color=white, x=4, label=$\infty$]{5}
    \Vertex[shape=star, size=0.5, color=white, x=5, label=$\infty$]{6}
    \Vertex[shape=star, size=0.5, color=white, x=6, label=$\infty$]{7}

\draw[line width=1.5pt, dash pattern=on 5pt off 3pt] (6.2,0) -- (7.2,0);
    \Edge(1)(2)
    \Edge(2)(3)
    \Edge(3)(4)
    \Edge(4)(5)
    \Edge(5)(6)
    \Edge(6)(7)

      \node at (-2, 0) { $\mathcal{H}_{p^\infty}^\infty$};
\end{tikzpicture}
\]
\vspace{-.5em}
\[
\begin{tikzpicture}[scale=1.2, transform shape]
    \Vertex[shape=circle, size=0.1, color=white, x=0]{1}
    \Vertex[shape=star, size=0.5, color=white, x=1, label=1]{2}
    \Vertex[shape=star, size=0.5, color=white, x=2, label=2]{3}
    \Vertex[shape=star, size=0.5, color=white, x=3, label=3]{4}
    \draw[line width=1.5pt, dash pattern=on 5pt off 3pt] (3.2,0) -- (5.2,0);
    \Vertex[shape=star, size=0.5, color=white, x=5, label=r]{5}
    \Vertex[shape=star, size=0.5, color=white, x=6, label=r]{6}
    \Vertex[shape=star, size=0.5, color=white, x=7, label=r]{7}
    \draw[line width=1.5pt, dash pattern=on 5pt off 3pt] (7.2,0) -- (8.2,0);
     
    \Edge(1)(2)
    \Edge(2)(3)
    \Edge(3)(4)
    \Edge(5)(6)
    \Edge(6)(7)

    \node at (-1, 0) { $\mathcal{H}_{p^{\infty,+}}^r$};
\end{tikzpicture}
\]
\vspace{-.5em}
\[
\begin{tikzpicture}[scale=1.2, transform shape]
    \Vertex[shape=circle, size=0.1, color=white, x=0]{1}
    \Vertex[shape=star, size=0.5, color=white, x=1, label=1]{2}
    \Vertex[shape=star, size=0.5, color=white, x=2, label=2]{3}
    \Vertex[shape=star, size=0.5, color=white, x=3, label=3]{4}
    \Vertex[shape=star, size=0.5, color=white, x=4, label=4]{5}
    \Vertex[shape=star, size=0.5, color=white, x=5, label=5]{6}
    \Vertex[shape=star, size=0.5, color=white, x=6, label=6]{7}
    \Vertex[shape=star, size=0.5, color=white, x=7, label=7]{8}
    \draw[line width=1.5pt, dash pattern=on 5pt off 3pt] (7.2,0) -- (8.2,0);
     
    \Edge(1)(2)
    \Edge(2)(3)
    \Edge(3)(4)
    \Edge(4)(5)
    \Edge(5)(6)
    \Edge(6)(7)
    \Edge(7)(8)

    \node at (-1, 0) { $\mathcal{H}_{p^{\infty,+}}^\infty$};
\end{tikzpicture}
\]
\end{notation}

\begin{lemma}\label{lem:infinitecountline}
Consider the graph $\mathcal{H}_{p^{\infty}}^{r}$ with $r<\infty$. For a skeletal vertex $v$, the number of vertices $w\in V(\mathcal{H}_{p^{\infty}}^{r})$ a distance  $j\geq1$ from $v$ is
\[
\left\vert \left\{  w\in V(\mathcal{H}_{p^{\infty}}^{r})\mid \operatorname*{dist}(v,w)=j\right\}  \right\vert =\left\{
\begin{array}
[c]{cl}
p^j+p^{j-1} & \text{if }1\leq j\leq r,\\
2p^{r} & \text{if }j>r,
\end{array}
\right.
\]
where $\varphi$ denotes Euler's totient function.
\end{lemma}

\begin{proof}
By Lemma \ref{Lem:linebloomdepth}, each skeletal vertex of $\mathcal{H}_{p^{\infty}}^{r}$ has bloom depth $r$.  Now let $v$ be a skeletal vertex. Suppose that $1\leq j\leq r$. Then a path of length $j$ starting at $v$ may begin in any of the $p+1$ edges incident to $v$: two edges along the skeleton, and the $p-1$ edges entering the bloom of $v$ relative to the skeleton. Since each of the skeletal vertices has bloom depth $r$, we have that after the first step, each of the remaining $j-1$ steps has exactly $p$ choices. Consequently, the number of vertices a distance $j$ from $v$ is $(p+1)p^{j-1}=p^j+p^{j-1}$.

It remains to consider the case $j>r$. In this direction, let the set of skeletal vertices be $\left\{  v_{i}\right\}  _{i\in\mathbb{Z}}$ with edges $\left\{  v_{i},v_{i+1}\right\}  _{i\in\mathbb{Z}}$ and $v=v_{0}$. Next, let $\mathcal{I}$ be a path of length $j$, and let $v_{i}$ be the last vertex of $\mathcal{I}$ on the skeleton, and let $t$ be the length of the subpath of $\mathcal{I}$ lying in the bloom $\mathcal{B}(v_{i})$ relative to the skeleton. In particular, $\left\vert i\right\vert +t=j$ and $t\leq r$. Note that $\left\vert i\right\vert \geq1$ and thus for each fixed $r$, there are two choices for $v_{i}$. Moreover, by Lemma~\ref{Lem:bloomvertices}, the number of vertices in $\mathcal{B}(v_{i})$ at distance $t$ from $v_{i}$ is $\varphi(p^{t})$. Thus, the number of vertices at distance $j$ from $v_{0}$ is
\[
2\sum_{t=0}^{r}\varphi(p^{t})=2p^{r},
\]
which concludes the proof.
\end{proof}

\begin{corollary}\label{cor:infinitecount}
Let $r\in\mathbb{Z}_{\geq0}\cup\left\{  \infty\right\}  $, and let $\mathcal{L}$ be a ray with leaf $v_{0}$ that contains the skeleton of $\mathcal{H}_{p^{\infty,+}}^{r}$. Then, for each nonnegative integer $j$, the number of vertices which are at a distance $j$ from $v_{0}$ is $p^{\min\left\{  r,\left\lfloor j/2\right\rfloor\right\}  }$.
\end{corollary}

\begin{proof}
Suppose that $\mathcal{L}$ has vertex set $\left\{  v_{i}\right\}  _{i\geq0}$ and edges $\left\{  v_{i},v_{i+1}\right\}  $ for $i\geq0$. By Lemma~\ref{Lem:raybloomdepth}, the bloom depth of $v_{i}$ is $d_{i}=\min\!\left\{i,r\right\}  $. Now observe that upon fixing $j$, and setting $k=2j$, the
proof of Lemma \ref{pathcounting} shows that
\[
\left\vert \left\{  w\in V(\mathcal{H}_{p^{\infty,+}}^{r})\mid\operatorname*{dist}(v_{0},w)=j\right\}  \right\vert =p^{\min\left\{r,\left\lfloor j/2\right\rfloor \right\}  }. \qedhere
\]
\end{proof}

\subsection{Classification of infinite primary isogeny graphs}\label{subsec:padicGalois}
We conclude this section with the classification of infinite $p$-primary graphs. To achieve this, we require a $p$-adic version of Proposition~\ref{prop:imagechangeuponisogeny}. Recall that if $E$ is an elliptic curve over a field $K$ of characteristic $0$, the elements of its Tate module $T_p(E)$ are sequences $\mathcal P=(P_i)_{i \ge 1}$ such that $P_i\in E[p^i](\overline K)$ and $pP_{i+1}=P_i$. The Tate module is a free $\Z_p$-module of dimension $2$. A pair $(\mathcal P, \mathcal Q)$ is a basis for $T_p(E)$  if $(P_i, Q_i)$ is a basis for $E[p^i]$ for each $i$. 
    
\begin{proposition}\label{prop:infiniteimagechange}
Let $p$ be a prime, and suppose that $E_{0}$ is an elliptic curve defined over a field $K$ of characteristic $\neq p$, and suppose that $(\mathcal{P},\mathcal{Q})$ is a basis for the Tate module $T_{p}(E_{0})$ with $\mathcal{P}=\left(  P_{j}\right)  _{j\geq1}$ and $\mathcal{Q}=\left(  Q_{j}\right)  _{j\geq1}$. Suppose further that for some positive integers $k$ and $i$ with $i\leq k$, $\left\langle P_{k}\right\rangle$ is $G_{K}$-invariant, and let $E_{i}=E_{0}/\left\langle p^{k-i}P_{k}\right\rangle $ be the corresponding $K$-rational $p^{i}$-isogeny $\phi_{i}:E_{0}\rightarrow E_{i}$. Now set $\mathcal{R}=\left(  \phi_{i}(P_{i+j})\right)  _{j\geq1}$ and $\mathcal{S}=(\phi_{i}(Q_{j}))_{j\geq1}$. Then, $(\mathcal{R},\mathcal{S})$ is a basis for $T_{p}(E_{i})$, and if $\sigma\in G_{K}$, then
\[
\rho_{E_{0},p^{\infty}}(\sigma)=\left(
\begin{array}
[c]{cc}
a & b\\
wp^{k} & c
\end{array}
\right)  \qquad\text{and}\qquad\rho_{E_{i},p^{\infty}}(\sigma)=\left(
\begin{array}
[c]{cc}
a & bp^{i}\\
wp^{k-i} & c
\end{array}
\right)
\]
in basis $(\mathcal{P},\mathcal{Q})$ and $(\mathcal{R},\mathcal{S})$, respectively, for some $a,b,c,w\in\mathbb{Z}_{p}$.
    \end{proposition}
    
\begin{proof}
Let $\pi_{k}:T_{p}(E)\rightarrow E[p^{k}]$ be the reduction modulo $p^{k}$. Now let $\sigma\in G_{K}$ and let $\rho_{E_{0},p^{\infty}}(\sigma)=\sm{a &  b  \\ u  & c}$ for some $a,b,c,u\in\mathbb{Z}_{p}$ in basis $(\mathcal{P},\mathcal{Q})$. Since $\left\langle P_{k}\right\rangle $ is $G_{K}$-invariant, we obtain that
\[
\pi_{k}(\sigma(\mathcal{P}))=\overline{a}P_{k}+\overline{u}Q_{k}\in\left\langle P_{k}\right\rangle ,
\]
where $\overline{a}$ and $\overline{u}$ are the corresponding reductions modulo $p^{k}$. Consequently, $v_{p}(u)\geq k$, and thus $u=wp^{k}$ for some $w\in\mathbb{Z}_{p}$.

We next show that $\mathcal{(\mathcal{R}},\mathcal{S)}$ is a basis for $T_{p}(E_{i})$. In this direction, observe that for each $j\geq1$, $P_{i+j}=p^{k-i}P_{k+j}$ and $P_{k}=p^{j}P_{k+j}$. Thus,
\[
\ker\phi_{i}\cong\left\langle p^{k-i}P_{k}\right\rangle \cong\left\langle p^{k-i}(p^{j}P_{k+j})\right\rangle \cong\left\langle p^{k+j-i}P_{k+j}
\right\rangle .
\]
We then deduce from the first isomorphism theorem that
\[
\left\langle \phi_{i}(P_{i+j})\right\rangle \cong\frac{\left\langle P_{i+j}\right\rangle }{\left\langle P_{i+j}\right\rangle \cap\ker\phi_{i}}\cong\frac{\left\langle p^{k-i}P_{k+j}\right\rangle }{\left\langle p^{k+j-i}P_{k+j}\right\rangle }.
\]
Since $p^{k-i}P_{k+j}$ and $p^{k+j-i}P_{k+j}$ have order $p^{i+j}$ and $p^{i}$, respectively, we conclude that $\phi_{i}(P_{i+j})$ has order $p^{j}$ for each $j\geq1$. Since $\phi_{i}(Q_{j})$ also has order $p^{j}$, we have that $\mathcal{R},\mathcal{S}\in T_{p}(E_{i})$. It suffices to prove that they are independent. Suppose that for some $\alpha,\beta\in\mathbb{Z}_{p}$ one has
\[
\alpha\mathcal{R}=\beta\mathcal{S}\qquad\Longrightarrow\qquad\alpha\phi_{i}(P_{i+j})=\beta\phi_{i}(Q_{j})
\]
for each $j\geq1$. In particular, for each $j\ge 1$, we have that $\alpha P_{i+j}-\beta Q_{j}\in\ker\phi_{i}$, so
\[
\alpha P_{i+j}-\beta Q_{j}=\gamma_j p^{k-i}P_{k}
\]
for some $\gamma_j \in\mathbb{Z}_{p}$. A multiplication-by-$p^{i}$ then shows that $\alpha p^{i}P_{i+j}=\beta p^{i}Q_{j}$. Consequently, $p^i (\alpha P_{i+j}-\beta Q_{j})=\mathcal{O}$ for each $j\ge 1$. But $T_p(E)$ is torsion free, and so it must be the case that $\alpha P_{i+j}=\beta Q_{j}$. Since $P_{i+j}$ and $Q_{j}$ are independent for each $j\geq1$, we conclude that $\alpha=\beta=0$. Hence, $\mathcal{(\mathcal{R}},\mathcal{S)}$ is a basis for $T_{p}(E_{i})$.

Now set $R_{j}=\phi_{i}(P_{i+j})$ and $S_{j}=\phi_{i}(Q_{j})$ for each $j\geq1$, and let $\sigma\in G_{K}$ with $\rho_{E_{0},p^{\infty}}(\sigma)=\sm{a &  b  \\ w p^k  & c}$ for some $a,b,c,w\in\mathbb{Z}_{p}$. We now calculate $\sigma(R_{j})$ and $\sigma(S_{j})$ for each $j\geq1$. In the calculations that follow, for $x\in\mathbb{Z}_{p}$, we will write $\overline{x}=x\ \operatorname{mod}p^{i+j}$ and $\widetilde{x}=x\ \operatorname{mod}p^{j}$. One has
\[
\sigma(S_{j})=\sigma(\phi_{i}(Q_{j}))=\phi_{i}(\sigma(Q_{j}))=\phi_{i}(\overline{b}P_{j}+\overline{c}Q_{j})=\phi_{i}(\overline{bp^{i}} P_{i+j}+\overline{c}Q_{j})=\widetilde{bp^{i}}R_{j}+\widetilde{c}S_{j}.
\]
On the other hand,
\[
\sigma(R_{j})=\phi_{i}(\sigma(P_{i+j}))=\phi_{i}(\overline{a}P_{i+j}+\overline{wp^{k}}Q_{i+j})=\phi_{i}(\overline{a}P_{i+j}+\overline{wp^{k-i}}Q_{j})=\widetilde{a}R_{j}+\widetilde{wp^{k-i}}S_{j}.
\]
Written as a matrix, the action of $\sigma$ with respect to the basis
$(\mathcal{R},\mathcal{S})$ is
\[
\rho_{E_{i},p^{\infty}}(\sigma)=\left(
\begin{array}
[c]{cc}
a & bp^{i}\\
wp^{k-i} & c
\end{array}
\right)  . \qedhere
\]
\end{proof}

We now use Proposition~\ref{prop:infiniteimagechange} to establish the lemma below, which will be used in the classification of infinite primary graphs.

\begin{lemma}\label{lem:raycounting}
Let $E$ be an elliptic curve over a field $K$ of characteristic $0$ with $\operatorname{End}_K(E)\cong \mathbb{Z}$, and let $p$ be a prime. Suppose further that $\mathcal{G}_{p}(E/K)$ is infinite, and does not contain a line graph. Then, there exist a leaf $E_{0}\in V(\mathcal{G}_{p}(E/K))$ for which the number of $K$-rational $p^{j}$-isogenies admitted by $E_{0}$ is $p^{\min\{\mathfrak{I}(E),\left\lfloor j/2\right\rfloor \}}$.
\end{lemma}

\begin{proof}
    Since $\mathcal{G}_{p}(E/K)$ is an infinite tree and does not contain a line graph, we have by Lemma~\ref{koenig} that $\mathcal{G}_{p}(E/K)$ contains a ray. In particular, $\mathcal{G}_{p}(E/K)$ must have at least one leaf, as otherwise $\mathcal{G}_{p}(E/K)$ would have a line as a subgraph. So suppose that $\left(  E_{i}\right)  _{i\geq0}$ is a ray $\mathcal{L}$ in $\mathcal{G}_{p}(E/K)$ with $E_{0}$ a leaf.  Thus, for each $i\geq1$, there is a $p^{i}$-isogeny $\phi_{i}:E_{0}\rightarrow E_{i}$. Let $P_{i}=\ker\phi_{i}$, and let $\mathcal{P}=(P_{i})_{i\geq1}\in T_{p}(E)$. Choose $\mathcal{Q}=(Q_{i})_{i\geq1}\in T_{p}(E)$ such that $(\mathcal{P},\mathcal{Q})$ is a basis for $T_{p}(E)$. For $i\geq1$, set $\mathcal{R}_{i}=(\phi_{i}(P_{i+j}))_{j\geq1}$ and $\mathcal{S}_{i}=(\phi_{i}(Q_{j}))_{j\geq1}$. Then, by Proposition \ref{prop:infiniteimagechange}, $(\mathcal{R}_{i},\mathcal{S}_{i})$ is a basis for $T_{p}(E_{i})$, and since each $\left\langle P_{i}\right\rangle $ is $G_{K}$-invariant, the proof of loc. cit. shows that for $\sigma\in G_{K}$,
\begin{equation}
\rho_{E_{0},p^{\infty}}(\sigma)=\left(
\begin{array}
[c]{cc}
a & b\\
0 & c
\end{array}
\right)  \qquad\text{and}\qquad\rho_{E_{i},p^{\infty}}(\sigma)=\left(
\begin{array}
[c]{cc}
a & bp^{i}\\
0 & c
\end{array}
\right)  \label{changeofbasisTp}
\end{equation}
in basis $(\mathcal{P},\mathcal{Q})$ and $(\mathcal{R}_{i},\mathcal{S}_{i})$, respectively, for some $a,b,c\in\mathbb{Z}_{p}$. Further, since $E_{0}$ is a leaf, there is a $\sigma\in G_{K}$ for which $v_{p}(b)=0$. By Proposition~\ref{prop:bloomingpotential} and Lemma~\ref{lem:explicitbloompot}, we have that
\begin{equation}
\mathfrak{I}(E)=\mathfrak{I}(E_0)=\min\!\left\{  v_{p}(a-c)\left\vert \left(
\begin{array}
[c]{cc}
a & b\\
0 & c
\end{array}
\right)  \in\rho_{E_{0},p^{\infty}}(G_{K})\right.  \right\}  .\label{alphaZp}
\end{equation}
In what follows, we set $r=\mathfrak I(E_0)$.  

If $r=0$, the ray $(E_i)_{i\geq 0}$ is the entire isogeny graph, because it follows from the proof of \cite[Proposition 3.1]{ivan-numberofiso} that the number of $p^j$-isogenies of $E_0$ for any $j$ equals $1$. 

If $r=\infty$, then it again follows from the proof of \cite[Proposition 3.1]{ivan-numberofiso} that the number of $K$-rational $p^j$-isogenies of $E_0$ is $p^{\lfloor j/2 \rfloor}$ for any $j$ and we are done.

Suppose $0<r<\infty$, and let $\pi$ be the projection from $\mathcal{G}_p(E/K)$ onto the ray $\left(  E_{i}\right)  _{i\geq r}$. For $i \ge r$, let $\mathcal B(E_i)=\mathcal{G}_p(E/K)[\pi^{-1}(E_i)]$ denote the bloom of $E_i$ relative to the ray $\left(  E_{i}\right)  _{i\geq r}$. These are finite trees because $\mathcal{G}_p(E/K)$ does not contain a line. Now consider the bloom $\mathcal B(E_r)$. Note that $\mathcal B(E_r)$ contains $E_0, E_1, \ldots, E_r$ by construction. Let $F_0$ be a leaf in this bloom whose distance to $E_r$ is maximal. Since $E_r$ attains all $p^r$-isogenies, it follows that the distance from $F_0$ to $E_r$ is at least $r$. Denote this distance by $r+t$ for $t\geq 0$,  and let $(F_i)_{i \ge0}$ be the unique ray starting at $F_0$ such that $F_{r+t+i}=E_{r+i}$ for $i\geq 0$. In particular $F_{r+t}=E_r$.  

    We will now prove that the number of $K$-rational $p^j$-isogenies of $F_0$ is $p^{\min\{r, \lfloor j/2 \rfloor \}}$. By Proposition~\ref{prop:bloomingpotential} and Lemma~\ref{lem:explicitbloompot}, we have that
    $$r=\mathfrak{I}(E)=\mathfrak{I}(F_0)=\min\!\left\{  v_{p}(a-c)\left\vert \left(
\begin{array}
[c]{cc}
a & b\\
0 & c
\end{array}
\right)  \in\rho_{F_{0},p^{\infty}}(G_{K})\right.  \right\}  .$$ 
Consequently, the elliptic curves $F_i$ attains all of its $K$-rational $p^{d_i}$-isogenies where $d_i=\min\{i,r\}$. Further, we claim that for $i\leq r$, the elliptic curve $F_i$ does not admit a $p^{i+1}$-isogeny which is independent from the isogeny $F_i\to F_{i+1}$. Indeed, suppose that there exists such an isogeny $F_i\to F'$. Then $F' \in V(\mathcal B(E_r))$ and 
    $$\dist(F', E_r)=\dist(F', F_{r+t})=\dist(F', F_i)+\dist(F_i, F_{r+t})=i+1+(r+t-i)=r+t+1,$$
    which is a contradiction with $r+t$ being the maximal distance from a vertex in $\mathcal B(E_r)$ to $E_r$.

        We are now ready to count the number of $p^j$-isogenies admitted by $F_0$. First, take $j\leq 2r+1$. We claim that any $p^{j}$-isogeny $F_0\to F$ factors through the $p^{\lceil j/2 \rceil}$-isogeny $F_0\to F_{\lceil j/2 \rceil}$. Namely, let $i$ be the largest integer such that the isogeny $F_0\to F$ factors through the isogeny $F_0\to F_i$. In particular, we have the following diagram:
\[
\begin{tikzcd}[row sep=0.5em]
                       &                           & F \arrow[d, no head, dashed]      &                            &                           & E_0 \arrow[d, no head]         &                            &        \\
                       &                           & \vdots \arrow[d, no head, dashed] &                            &                           & \vdots \arrow[d, no head]      &                            &        \\
F_0 \arrow[r, no head] & \cdots \arrow[r, no head] & F_i \arrow[r, no head]            & F_{i+1} \arrow[r, no head] & \cdots \arrow[r, no head] & F_{r+t}=E_r \arrow[r, no head] & E_{r+1} \arrow[r, no head] & \cdots
\end{tikzcd}
\]
    Then $F_i\to F$ is an isogeny of degree $p^{j-i}$ independent from $F_i \to F_{i+1}$. If $i<j/2$, then $i\leq r$ and $j-i\geq i+1$, which is a contradiction with the fact that for $i\leq r$, the elliptic curve $F_i$ has no $p^{i+1}$-isogenies independent from $F_i\to F_{i+1}$. Thus, $i\geq \lceil j/2\rceil$. This means that the first $\lceil j/2 \rceil$ steps of any $p^j$-isogeny are fixed, and go from $F_0$ to $F_{\lceil j/2 \rceil}$. The elliptic curve $F_{\lceil j/2 \rceil}$ attains all $p^{\lfloor j/2 \rfloor}$-isogenies, so there are $p$ choices for each of the remaining $\lfloor j/2 \rfloor$ steps.  In conclusion, $F_0$ has $p^{\lfloor j/2 \rfloor}$ isogenies for $j\leq 2r+1$. 

    With $r$ as before, the proof of \cite[Proposition 3.1]{ivan-numberofiso} shows that for each $j\geq1$, the number of $p^{j}$-isogenies of $F_{0}$ is either $p^{\min\left\{  r,\left\lfloor j/2\right\rfloor \right\}  }$ or $2p^{r}$. Furthermore, if this number is $2p^{r}$ for some $j$, then $j>2r$, and in fact, the number of $p^{2r+1}$-isogenies admited by $F_{0}$ also equals $2p^{r}$. Since we have just proven that the number of $p^{2r+1}$-isogenies of $F_0$ equals $p^{\left \lfloor \frac{2r+1}{2} \right \rfloor}=p^r$, it follows that the number of $p^j$-isogenies for any $j>2r$ equals $p^{r}$, which proves the claim. 
\end{proof}

\begin{theorem}\label{thm:infinitegpks}
Let $E$ be an elliptic curve defined over a field $K$ of characteristic $0$ such that $\operatorname*{End}_{K}E\cong\mathbb{Z}$. If for a prime $p$, the $p$-primary graph $\mathcal{G}_{p}(E/K)$ is infinite, then
\[
\mathcal{G}_{p}(E/K)\cong\mathcal{H}_{p^{\infty}}^{r}\qquad\text{or}\qquad\mathcal{G}_{p}(E/K)\cong\mathcal{H}_{p^{\infty,+}}^{r}
\]
with $r=\mathfrak{I}_p(E/K)  $.

Conversely, for any $r\in\mathbb{Z}_{\geq0}\cup\left\{  \infty\right\}  $, each of the graphs $\mathcal{H}_{p^{\infty}}^{r}$ and $\mathcal{H}_{p^{\infty,+}}^{r}$ is realized as $\mathcal{G}_{p}(E/K)$ for some elliptic curve $E$ defined over a field $K$ of characteristic $0$, unless $(p,r)=(2,0)$.
\end{theorem}

\begin{proof}
 By Lemma \ref{koenig}, $\mathcal{G}_{p}(E/K)$ contains a ray $\left(  E_{i}\right)  _{i\geq0}$. We then have for each $i\geq1$, a $p^{i}$-isogeny $\phi_{i}:E_{0}\rightarrow E_{i}$ with $\ker\phi_{i}=\left\langle P_{i}\right\rangle $ a $G_{K}$-invariant subgroup. Thus $\mathcal{P}=(P_{i})_{i\geq1}\in T_{p}(E)$. Now choose $\mathcal{Q}=(Q_{i})_{i\geq1}\in T_{p}(E)$ such that $(\mathcal{P},\mathcal{Q})$ is a basis for $T_{p}(E)$. Proceeding as in the proof of Lemma \ref{lem:raycounting}, we obtain via by Proposition \ref{prop:infiniteimagechange}, a basis $(\mathcal{R}_{i},\mathcal{S}_{i})$ of $T_{p}(E_{i})$ such that for $\sigma\in G_{K}$,
\begin{equation}
\rho_{E_{0},p^{\infty}}(\sigma)=\left(
\begin{array}
[c]{cc}
a & b\\
0 & c
\end{array}
\right)  \qquad\text{and}\qquad\rho_{E_{i},p^{\infty}}(\sigma)=\left(
\begin{array}
[c]{cc}
a & bp^{i}\\
0 & c
\end{array}
\right)  \label{eq:thmrepn}
\end{equation}
in basis $(\mathcal{P},\mathcal{Q})$ and $(\mathcal{R}_{i},\mathcal{S}_{i})$, respectively, for some $a,b,c\in\mathbb{Z}_{p}$. Next, set $r=\mathfrak{I}(E)$. By Proposition~\ref{prop:bloomingpotential} and Lemma~\ref{lem:explicitbloompot},
$$r=\mathfrak I(E_0)=\min\!\left\{  v_{p}(a-c)\left\vert \left(
\begin{array}
[c]{cc}
a & b\\
0 & c
\end{array}
\right)  \in\rho_{E_{0},p^{\infty}}(G_{K})\right.  \right\}.$$ We proceed by cases, depending on whether $\mathcal{G}_{p}(E/K)$ contains a line graph or not.

\textbf{Case 1.} Suppose $\mathcal{G}_{p}(E/K)$ does not contain a line graph. By Lemma \ref{lem:raycounting}, we may further assume that $E_{0}$ is a leaf with the property that for each nonnegative integer $j$, the number of vertices at a distance $j$ from $v_{0}$ is $p^{\min\left\{  r,\left\lfloor j/2\right\rfloor \right\}  }$. Now set $d_{i}=\min\!\left\{  i,r\right\}  $, and observe that from (\ref{eq:thmrepn}), we obtain for $i\geq1$ that
\[
\rho_{E_{i},p^{\infty}}(\sigma)=\left(
\begin{array}
[c]{cc}
a & bp^{i}\\
0 & c
\end{array}
\right)  \equiv\left(
\begin{array}
[c]{cc}
a & 0\\
0 & a
\end{array}
\right)  \ \operatorname{mod}p^{d_{i}}.
\]
Now consider the proof of Theorem \ref{classificationGpk}. Observe that the argument following (\ref{eq:thmfinitecenter}) implies in our setting that for $i\geq1$, each elliptic curve at a distance $m<d_{i}$ from $E_{i}$ is $p$-bloomed. We then further obtain from the argument, together with Lemma~\ref{Lem:linebloomdepth}, that $\mathcal{G}_{p}(E/K)$ contains a subgraph $\mathcal{G}_{p}^{\prime}(E/K)$ such that $\mathcal{G}_{p}^{\prime}(E/K)\cong\mathcal{H}_{p^{\infty,+}}^{r}$. Note that by construction of the isomorphism, $\mathcal{G}_{p}^{\prime}(E/K)$ contains the ray $\left(E_{i}\right)  _{i\geq0}$, together with those elliptic curves a distance at most $d_{i}$ from $E_{i}$ for $i\geq0$. By Corollary~\ref{cor:infinitecount}, the number of vertices in $\mathcal{G}_{p}^{\prime}(E/K)$ which are at a distance $j$ from $E_{0}$ is $p^{\min\left\{  r,\left\lfloor j/2\right\rfloor \right\}  }$. This number agrees with our choice of $E_{0}$ via Lemma~\ref{lem:raycounting}, and so we conclude that if $\mathcal{G}_{p}(E/K)$ had more edges than $\mathcal{G}_{p}^{\prime}(E/K)$, then $E_{0}$ would admit more $p^{j}$-isogenies for some $j\geq1$, which is impossible. Thus, $\mathcal{G}_{p}(E/K)=\mathcal{G}_{p}^{\prime}(E/K)$.

\textbf{Case 2.} Suppose $\mathcal{G}_{p}(E/K)$ contains a line graph. Then, upon possibly changing our choice of ray $\left(  E_{i}\right)  _{i\geq0}$, we may assume that the ray $\left(  E_{i}\right)  _{i\geq0}$ lies on a line $\left(  E_{i}\right)  _{i\in\mathbb{Z}}$. We may then choose $\mathcal{Q}=(Q_{i})_{i\geq1}\in T_{p}(E)$ such that $E_{-i}=E_{0}/\left\langle Q_{i}\right\rangle $ for $i\geq1$. Thus, each $\left\langle Q_{i}\right\rangle $ is $G_{K}$-invariant and thus for $\sigma\in G_{K}$, we have that
\[
\rho_{E_{i},p^{\infty}}(\sigma)=\left(
\begin{array}
[c]{cc}
a & 0\\
0 & c
\end{array}
\right)
\]
for some $a,c\in\mathbb{Z}_{p}$.

If $r=\infty$, then $\rho_{E_{i},p^{\infty}}(G_{K})$ is in the center of $\operatorname*{GL}\nolimits_{2}(\mathbb{Z}_{p})$, which is equivalent to each $E_{i}$ admitting all of its $p^{k}$-isogenies for each positive integer $k$. Thus, $\mathcal{G}_{p}(E/K)$ is the entire $p$-Bruhat-Tits tree $\mathcal{H}_{p^{\infty}}^{\infty}$, and so $\mathcal{G}_{p}(E/K)\cong\mathcal{H}_{p^{\infty}}^{\infty}$.

Now suppose that $r<\infty$. Then, for each $i\in \mathbb{Z}$,
\[
\rho_{E_{i},p^{\infty}}(\sigma)=\left(
\begin{array}
[c]{cc}
a & 0\\
0 & c
\end{array}
\right)  \equiv\left(
\begin{array}
[c]{cc}
a & 0\\
0 & a
\end{array}
\right)  \ \operatorname{mod}p^{r}.
\]
Consequently, each $E_{i}$ admits all of its $p^{j}$-isogenies over $K$ for $j\leq r$. Further, the argument in the proof of Theorem~\ref{classificationGpk} quoted in Case 1 shows that for a fixed $i$, each elliptic curve a distance $j<r$ from $E_{i}$ is $p$-bloomed. Let $\mathcal{G}_{p}^{\prime}(E/K)$ be the subgraph of $\mathcal{G}_{p}(E/K)$ consisting of those elliptic curves a distance at most $r$ from the line $\left(  E_{i}\right)  _{i\in\mathbb{Z}}$. From Lemma~\ref{Lem:raybloomdepth}, we obtain that $\mathcal{G}_{p}^{\prime}(E/K)\cong\mathcal{H}_{p^{\infty}}^{r}$. By the proof of \cite[Proposition 3.1]{ivan-numberofiso}, the number of $K$-rational $p^{j}$-isogenies for $j\geq1$ admitted by $E_{i}$ is
\[
\left\vert \left\{  E^{\prime}\in V(\mathcal{H}_{p^{\infty}}^{r})\mid\operatorname*{dist}(E_{i},E^{\prime})\right\}  =j\right\vert =\left\{
\begin{array}
[c]{cl}
p^j+p^{j-1} & \text{if }1\leq j\leq r,\\
2p^{r} & \text{if }j>r.
\end{array}
\right.
\]
Note that this number agrees with the number of vertices a distance $j$ from a skeletal vertex, which is the context of Lemma \ref{lem:infinitecountline}. Therefore if $\mathcal{G}_{p}(E/K)$ had more edges than $\mathcal{G}_{p}^{\prime}(E/K)$, then some vertex $E_{i}$ would admit more $p^{j}$-isogenies for some $j\geq1$, which is impossible. Thus, $\mathcal{G}_{p}(E/K)=\mathcal{G}_{p}^{\prime}(E/K)$, and so $\mathcal{G}_{p}(E/K)\cong\mathcal{H}_{p^{\infty}}^{r}$.

It remains to prove that each graph can be achieved. By Galois theory, any closed subgroup of the image of a $p$-adic Galois representation can be realized after restricting to the fixed field of its preimage. It suffices to find a suitable closed subgroup for each of the mentioned graphs.

For $r\in\mathbb{Z}_{\geq0}\cup\left\{  \infty\right\}  $, let
{\small \begin{align*}
H_{p^{\infty,+}}^{r} &  :=\left\{  \left.  \left(
\begin{array}
[c]{cc}
a & b\\
0 & c
\end{array}
\right)  \in\operatorname*{GL}\nolimits_{2}(\mathbb{Z}_{p})\right\vert
v_{p}(a-c)\geq r\right\}  =\underset{k}{\underleftarrow{\lim}}\left\{  \left.
\left(
\begin{array}
[c]{cc}
a & b\\
0 & c
\end{array}
\right)  \in\operatorname*{GL}\nolimits_{2}(\mathbb{Z}/p^{k}\mathbb{Z}
)\right\vert v_{p}(a-c)\geq r\right\}  ,\\
H_{p^{\infty}}^{r} &  :=\left\{  \left.  \left(
\begin{array}
[c]{cc}
a & 0\\
0 & c
\end{array}
\right)  \in\operatorname*{GL}\nolimits_{2}(\mathbb{Z}_{p})\right\vert
v_{p}(a-c)\geq r\right\}  =\underset{k}{\underleftarrow{\lim}}\left\{  \left.
\left(
\begin{array}
[c]{cc}
a & 0\\
0 & c
\end{array}
\right)  \in\operatorname*{GL}\nolimits_{2}(\mathbb{Z}/p^{k}\mathbb{Z}
)\right\vert v_{p}(a-c)\geq r\right\}  .
\end{align*}}
In particular, $H_{p^{\infty,+}}^{r}$ and $H_{p^{\infty}}^{r}$ are closed subgroups of $\operatorname*{GL}\nolimits_{2}(\mathbb{Z}_{p})$. From the forward direction, we have that if $E$ is an elliptic curve over a field $K$ of characteristic $0$ with $\operatorname*{End}_{K}\! E\cong\mathbb{Z}$ and $\rho_{E,p^{\infty}}(G_{K})$ is isomorphic to $H_{p^{\infty,+}}^{r}$ (resp. $H_{p^{\infty}}^{r}$), then $\mathcal{G}_{p}(E/K)$ is isomorphic to $\mathcal{H}_{p^{\infty,+}}^{r}$ (resp. $\mathcal{H}_{p^{\infty}}^{r}$). Finally, observe that if $\sm{a &  b  \\ 0  & c}\in\operatorname*{GL}\nolimits_{2}(\mathbb{Z}_{2})$, then $v_{2}(a-c)\geq 1$. Consequently, the case $\left(  p,r\right)  =\left(  2,0\right)  $ is impossible.
\end{proof}

The proof of Theorem \ref{thm:infinitegpks} introduced the subgroups~$H_{p^{\infty}}^{r}$ and~$H_{p^{\infty,+}}^{r}$ of $B_{0}(p^{\infty})$, which are defined below.

\begin{definition}\label{def:groupsinfinite}
For a prime number $p$ and $r\in\mathbb{Z}_{\geq0}\cup\left\{  \infty\right\}  $, define
\begin{align*}
H_{p^{\infty}}^{r}  & :=\left\{  \left.  \left(
\begin{array}
[c]{cc}
a & 0\\
0 & c
\end{array}
\right)  \in\operatorname*{GL}\nolimits_{2}(\mathbb{Z}_{p})\right\vert v_{p}(a-c)\geq r\right\}  ,\\
H_{p^{\infty,+}}^{r}  & :=\left\{  \left.  \left(
\begin{array}
[c]{cc}
a & b\\
0 & c
\end{array}
\right)  \in\operatorname*{GL}\nolimits_{2}(\mathbb{Z}_{p})\right\vert v_{p}(a-c)\geq r\right\}  .
\end{align*}

\end{definition}

Next, we observe that proceeding via a similar argument as that considered in the proof of Corollary \ref{cor:subgroup-subgraph-corr}, yields the following infinite analogue:

\begin{corollary}\label{cor:infiniteimageconjugate}
Let $p$ be a prime number and $r\in\mathbb{Z}_{\geq0}\cup\left\{  \infty\right\}  $, and let $E$ be an elliptic curve over a field $K$ of characteristic $0$. Then, the $p$-primary graph $\mathcal{G}_{p}(E/K)$ has a subgraph isomorphic to $\mathcal{H}_{p^{\infty}}^{r}$ (resp. $\mathcal{H}_{p^{\infty,+}}^{r}$) if and only if there exists an elliptic curve $E_{0}\in V(\mathcal{G}_{p}(E/K))$ such that $\rho_{E_{0},p^{\infty}}(G_{K})$ is conjugate to a subgroup of $H_{p^{\infty}}^{r}$ (resp. $H_{p^{\infty,+}}^{r}$).
\end{corollary}

\begin{remark}
    As in the finite case, one can also define subgroups $$H_{p^{\infty, +}}^r(j)=\left\{
\begin{pmatrix} a & p^jb \\ 0 & c \end{pmatrix} \in \GL_2(\mathbb{Z}_p)
\;\middle|\; v_p(c-a) \ge r
\right\}.$$ Elliptic curves whose $p$-adic Galois image is contained in $H_{p^{\infty, +}}^r(j)$ for some $j$ also contain $\mathcal H_{p^{\infty, +}}^r$ as a subgraph of $\mathcal G_p(E/K)$.
\end{remark}



\section{Potential complex multiplication}\label{sec:pCM}
In this section, we study the isogeny graphs of complex elliptic curves with potential complex multiplication, and conclude with a classification of their isogeny graphs (see Theorem~\ref{thm:genkwon}). As a consequence, we obtain Corollary~\ref{genKwonBC}, which generalizes theorems of Kwon~\cite{Kwon} and Bourdon--Clark~\cite{BourdonClark1}. In this direction, let $E$ be a complex elliptic curve with $\operatorname*{End}E\cong\mathcal{O}$, where $\mathcal{O}$ is an order in an imaginary quadratic field $K$, and let $F=\mathbb{Q}(j(E))$. A key observation is that $F$ admits a real embedding \cite[Remark 5.2]{SilverbergFoD}. Consequently, Proposition~\ref{Prop:BloomInvReal} applies to elliptic curves defined over $F$, providing strong restrictions on their $p$-primary graphs, as seen in Corollary~\ref{Cor:BloomInvReal}.

We now recall some facts, and in the process describe why $F$ admits a real embedding. Let~$\Delta$ denote the discriminant of $\mathcal{O}$. Then $j(E)$ is a root of the Hilbert class polynomial $H_{\Delta}(X)$ \cite[\S 13]{CoxPrimes},\cite[Chapter~2]{MR1312368}. Since $H_{\Delta}(X)\in \mathbb{Z}[X]$ is irreducible, we have that $F\cong\mathbb{Q}[X]/(H_{\Delta}(X))$. Moreover, $H_{\Delta}(X)$ possesses a real root \cite[(5.4.3)]{Shimura}, and so $F$ admits a real embedding. As $K$ is an imaginary quadratic field, it follows that $K \not \subseteq F $, and hence $\operatorname*{End}_{F}E\cong\mathbb{Z}$. Equivalently, the field of definition of $\operatorname*{End}E$ is $KF$~\cite[Theorem 2.3]{GrepsCM},\cite[Remark 5.2]{SilverbergFoD}. 

Now observe that by Proposition~\ref{Prop:BloomInvReal}, it is the case that for each odd prime $p$, there is no elliptic curve in the $p$-primary graph $\mathcal{G}_{p}(E/F)$ that is $p$-bloomed. We similarly have that no vertex in $\mathcal{G}_{2}(E/F)$ has $2$-bloom depth $\geq2$. Further, these conclusions continue to hold if we base change to any field with a real embedding. This leads to the following observation: $K$ must be contained in the splitting field of the level $p^{k}$ modular polynomial with $p$ a prime and~$k$ a positive integer such that $k\geq2$ if $p=2$. We now demonstrate this observation by using Lozano-Robledo's characterization \cite{GrepsCM} of the Galois representations of elliptic curves attached to CM elliptic curves. In particular, the following result generalizes~\cite[Lemma~3.15]{bourdonclarkstan}, which shows that $K\subseteq F(E[n])$ for each integer $n\geq3$.

\begin{proposition}\label{Prop:JsinQ}
Suppose that $E$ is a complex elliptic curve with $\operatorname*{End}E\otimes_{\mathbb{Z}}\mathbb{Q}\cong K$ an imaginary quadratic field. For $n\geq3$, let $\Phi_{n}(X,Y)$ denote the level $n$ modular polynomial, and set
\begin{equation}
J_{E}[n]=\left\{  j\in\overline{\mathbb{Q}}\mid\Phi_{n}X,j(E))=0\right\}
.\label{jinvfld}
\end{equation}
If $E$ is defined over $F=\mathbb{Q}(j(E))$, then $K\subseteq F(J_{E}[n])$.
\end{proposition}

\begin{proof}
Let $\Delta_{K}$ and $\Delta$ denote the discriminant of $K$ and $\operatorname*{End}E$, respectively. In particular, $\Delta=\mathfrak{f}^{2}\Delta_{K}$, where $\mathfrak{f}$ is the conductor of $\operatorname*{End}E$. Define $\delta,\phi$ by
\[
\left(  \delta,\phi\right)  =\left\{
\begin{array}
[c]{cl}
\left(  \frac{\Delta}{4},0\right)   & \text{if }\Delta\equiv 0\ \operatorname{mod}4\text{ or }n\text{ odd,}\\
\left(  \frac{\Delta_{K}-1}{4}\mathfrak{f}^{2},\mathfrak{f}\right)   & \text{if }\Delta\equiv1\ \operatorname{mod}4\text{ and }n\text{ even.}
\end{array}
\right.
\]
Now consider the Cartan subgroup $\mathcal{C}_{\delta,\phi}(n)\leq \operatorname*{GL}\nolimits_{2}(\mathbb{Z}/n\mathbb{Z})$ defined by
\[
\mathcal{C}_{\delta,\phi}(n)=\left\langle \left.  \left(
\begin{array}
[c]{cc}
a+b\phi & b\\
\delta b & a
\end{array}
\right)  \right\vert a,b\in\mathbb{Z}/n\mathbb{Z},a^{2}+ab\phi-\delta b^{2}\in(\mathbb{Z}/n\mathbb{Z})^{\times}\right\rangle .
\]
In particular, we have that
\begin{equation}
\left\langle \left.  \left(
\begin{array}
[c]{cc}
a & 0\\
0 & a
\end{array}
\right)  \right\vert a\in(\mathbb{Z}/n\mathbb{Z})^{\times}\right\rangle \leq
\mathcal{C}_{\delta,\phi}(n).\label{CenterGL2CM}
\end{equation}
By \cite[Theorem 1.1]{GrepsCM}, there is a basis of $E[n]$ such that the image of $\rho_{E,n}:G_{F}\rightarrow\operatorname*{GL}\nolimits_{2}(\mathbb{Z}/n\mathbb{Z})$ satisfies
\[
\rho_{E,n}(G_{F})\leq\mathcal{N}_{\delta,\phi}(n)=\left\langle \mathcal{C}
_{\delta,\phi}(n),\left(
\begin{array}
[c]{cc}
-1 & 0\\
\phi & 1
\end{array}
\right)  \right\rangle .
\]
In particular, $\mathcal{C}_{\delta,\phi}(n)$ is an index $2$ subgroup of $\mathcal{N}_{\delta,\phi}(n)$.

Next, following Section 6.2 of \cite{GrepsCM}, let $c\in G_{F}$ be a complex conjugation such that $\left\langle c|_{KF}\right\rangle \cong\operatorname*{Gal}(KF/F)$ acts on $\operatorname*{Gal}(K(E[n])/K)$ as complex conjugation. Further,
\[
\operatorname*{Gal}\!\left(  F(E[n])/F\right)  \cong\operatorname*{Gal}\!\left(  F(E[n])/KF\right)  \rtimes\left\langle c\right\rangle .
\]
From Theorem 6.7 of loc. cit., we obtain that the image of a complex conjugation via $\rho_{E,n}$ is contained in $\mathcal{N}_{\delta,\phi}(n)\backslash\mathcal{C}_{\delta,\phi}(n)$. This, paired with the discussion in Section 6.2 of loc. cit., shows that the fixed field of $C=\rho_{E,n}(G_{F})\cap\mathcal{C}_{\delta,\phi}(n)$ is $KF$. Now let
\[
Z=\rho_{E,n}(G_{F})\cap\left\{  \left(
\begin{array}
[c]{cc}
a & 0\\
0 & a
\end{array}
\right)  \in\operatorname*{GL}\nolimits_{2}(\mathbb{Z}/n\mathbb{Z})\right\}  .
\]
Then, $F(J_{E}[n])$ is the fixed field of $Z$. From (\ref{CenterGL2CM}), we deduce that $Z\leq C$. It follows from the Galois correspondence that
\[
K\subseteq KF\subseteq F(J_{E}[n]). \qedhere
\]
\end{proof}

By Corollary~\ref{Cor:BloomInvReal}, the $p$-primary graphs of elliptic curves $E$ with potential multiplication over $\mathbb{Q}(j(E))$ are uniquely determined by their isogeny class degree. Our next result extends Kwon's characterization~\cite[Corollary 4.2]{Kwon} to determine the $p$-primary graphs $\mathcal{G}_p(E/\mathbb{Q}(j(E)))$ from $\operatorname{End}(E)$.

\begin{lemma}\label{lem:potCM}
Let $\mathcal{O}$ be an order with discriminant $\Delta<-4$ in a quadratic field $K$. Suppose that $E$ is an elliptic curve such that $\operatorname*{End}E\cong\mathcal{O}$. If $E$ is defined over $F=\mathbb{Q}(j(E))$, then $\mathcal{G}(E/F)\cong\square_{p}\mathcal{G}_{p}(E/F)$, where
\[
\mathcal{G}_{p}(E/F)\cong\left\{
\begin{array}
[c]{cl}
\mathcal{H}_{p^{v_{p}(\Delta)}}^{0} & \text{if }(i)\ p\text{ is odd or }(ii)\ p=2\text{ and }v_{2}(\Delta)=0\text{ with }\Delta\equiv 5\ \operatorname{mod}8,\\
\mathcal{H}_{2}^{0} & \text{if }p=2\text{ and either }\Delta\equiv1\ \operatorname{mod}8\text{ or }v_{2}(\Delta)\in\{2,3\},\\
\mathcal{H}_{2^{v_{2}(\Delta)-2}}^{1} & \text{if }p=2\text{ and }v_{2}(\Delta)\geq4.
\end{array}
\right.
\]
Furthermore, the isogeny class degree and isogeny class size of $E/F$ are:
\begin{align*}
\deg\mathcal{G}(E/F)  & =2^{d_{2}}\prod_{\underset{p\text{ odd}}{p|\Delta}}p^{v_{p}(\Delta)}\qquad\text{and}\qquad\left\vert V(\mathcal{G}(E/F))\right\vert =s_{2}\prod_{\underset{p\text{ odd}}{p|\Delta}}(v_{p}(\Delta)+1),\\
\text{where }\left(  d_{2},s_{2}\right)    & =\left\{
{\renewcommand{\arraystretch}{1.3}\begin{array}
[c]{cl}
\left(  0,1\right)   & \text{if }\mathcal{G}(E/F)\cong\mathcal{H}_{2^{0}}^{0},\\
\left(  1,2\right)   & \text{if }\mathcal{G}(E/F)\cong\mathcal{H}_{2}^{0},\\
\left(  v_{2}(\Delta)-2,2(v_{2}(\Delta)-2)\right)   & \text{if }\mathcal{G}(E/F)\cong\mathcal{H}_{2^{v_{2}(\Delta)-2}}^{1}.
\end{array}}
\right.
\end{align*}
\end{lemma}

\begin{proof}
By the discussion at the start of the section, $F$ admits at least one real embedding. Consequently, Corollary \ref{Cor:BloomInvReal} implies that
\[
\mathcal{G}(E/F)\cong\left\{
\begin{array}
[c]{cl}
\underset{p}{\square}\mathcal{H}_{p^{k_{p}}}^{0} & \text{if }k_{2}\leq1,\\
\underset{p\text{ odd}}{\square}\mathcal{H}_{p^{k_{p}}}^{0}\square\mathcal{H}
_{2^{k_{2}}}^{1} & \text{if }k_{2}\geq2.
\end{array}
\right.
\]
In particular, $\deg\mathcal{G}_{p}(E/F)$ uniquely determines the $p$-primary graph. Thus, by Corollary~\ref{cor:sameprimary}, it suffices to determine the largest $p$-power isogeny occuring in $\mathcal{G}_{p}(E/F)$. In this direction, we recall Kwon's characterization \cite[Corollary~4.2]{Kwon} for the existence of an $n$-isogeny over $F=\mathbb{Q}(j(E))$. Now write $\Delta=\mathfrak{f}^{2}\Delta_{K}$, where $\mathfrak{f}$ denotes the conductor of $\mathcal{O}$ and $\Delta_{K}$ the discriminant of $K$. The table below summarizes Kwon's result: the third column gives the necessary and sufficient conditions for the existence of such an isogeny, while the fourth column gives the number of distinct $F$-rational $n$-isogenies admitted by $E$.
{\begingroup 
\renewcommand{\arraystretch}{1.18}
 \begin{longtable}{cccl}
 	\caption{Necessary and sufficient conditions for $E/\mathbb{Q}(j(E))$ to admit an $n$-isogeny, and the number $m $ of such isogenies.}\label{ta:kwon}\\
	\hline
	\(v_2(\mathfrak{f})\) & Congruence on \(\Delta_K\) & Divisibility condition & \# $m$ of $n$-isogenous curves  \\
	\hline

	\endfirsthead
	\hline
	\(\mathfrak{f}\) & Congruence on  \(\Delta_K\) & Divisibility condition  & \# $m$ of $n$-isogenous curves  \\
	\hline
	\endhead
	\hline

	\multicolumn{4}{r}{\emph{continued on next page}}
	\endfoot
	\hline
	\endlastfoot

$\ge 1$  & $\Delta_K \equiv 0 \mod 4$ & \multirow{2}{*}{$n \mid \dfrac{\mathfrak{f}^2\Delta_K}{4}$} & $2$ if  $n\equiv 0 \mod 4$\\\cmidrule{4-4}
$\ge 2$ & $\Delta_K \equiv 1 \mod 4$ & & $1$ otherwise\\
\midrule
$ 1$ & $\Delta_K \equiv 1 \mod 4$ & \multirow{2}{*}{$\begin{array}{c}\mbox{$n=n'$ or $2n'$, $n'$ odd} \\
 \mbox{$n' \mid \mathfrak{f}^2\Delta_K$}\end{array}$} &  \multirow{2}{*}{$1$} \\
$0$ & $\Delta_K \not\equiv 5 \mod 8$ & & \\
\midrule
$0$ & $\Delta_K \equiv 5 \mod 8$ & $n \mid \mathfrak{f}^2\Delta_K$ & $1$

\end{longtable}
\endgroup}

We first consider the case of an odd prime. If $v_{p}(\Delta)=0$, then loc. cit. implies that $E$ does not admit a $F$-rational $p$-isogeny. Consequently, $\mathcal{G}_{p}(E/F)$ is trivial. Now suppose that $p$ is an odd prime and let $s$ be a positive integer. We now consider the necessary and sufficient conditions in the above table to deduce when $E$ admits a $p^{s}$-isogeny. Since $p\neq2$, the conditions in the third column become equivalent to $p^{s}|\Delta$. Consequently, by loc. cit., $E$ admits a $F$-rational $p^{s}$-isogeny if and only if $s\leq v_{p}(\Delta)$. Moreover, in each of these cases, the table shows that there is a unique elliptic curve that is $p^{s}$-isogenous to $E$ over $F$. Hence $\mathcal{G}_{p}(E/F)=\mathcal{H}_{p^{v_{p}(\Delta)}}^{0}$.

It remains to analyze the prime $2$. In this direction, let $s$ be a positive integer. We now proceed by considering the cases appearing in the above table separately, to deduce when $E$ admits a $2^{s}$-isogeny via Kwon's criterion.

\textbf{Case 1.} Suppose that $(i)$ $v_{2}(\mathfrak{f})\geq1$ and $\Delta_{K}\equiv0\ \operatorname{mod}4$ or $(ii)$ $v_{2}(\mathfrak{f})\geq2$ and $\Delta_{K}\equiv1\ \operatorname{mod}4$. Then, $v_{2}(\Delta)\geq4$, and the condition $2^{s}|\frac{\mathfrak{f}\Delta_{K}}{4}$ is equivalent to $s\leq v_{2}(\Delta)-2$. It follows from Kwon's criterion, that when $s=1$, there is a unique elliptic curve that is $2$-isogenous to $E$ over~$F$. Whereas for every $s\geq2$ satisfying $s\leq v_{2}(\Delta)$, it is the case that there are exactly two such elliptic curves. These are precisely the incidence relations of $\mathcal{H}_{2^{v_{2}(\Delta)-2}}^{1}$. Therefore, in this case, $\mathcal{G}_{2}(E/F)\cong\mathcal{H}_{2^{v_{2}(\Delta)-2}}^{1}$.

\textbf{Case 2.} Suppose that $(i)$ $v_{2}(\mathfrak{f})=1$ and $\Delta_{K}\equiv1\ \operatorname{mod}4$ or $(ii)$ $v_{2}(\mathfrak{f})=0$ and $\Delta_{K}\not \equiv 5\ \operatorname{mod}8$. Then, our assumptions are equivalent to either $\Delta\equiv1\ \operatorname{mod}8$ or $v_{2}(\Delta)\in\{2,3\}$. From Kwon's criterion, we obtain that $E$ admits a $F$-rational $2^{s}$-isogeny if and only if $s=1$. Moreover, the corresponding elliptic curve is unique. Thus, the $2$-primary graph is $\mathcal{G}_{2}(E/F)\cong\mathcal{H}_{2}^{0}$.

\textbf{Case 3.} Suppose $v_{2}(\mathfrak{f})=0$ and $\Delta_{K}\equiv5\ \operatorname{mod}8$. Then, $\Delta\equiv5\ \operatorname{mod}8$, and thus $v_{2}(\Delta)=0$. Applying Kwon's criterion yields that $E$ does not admit a $2$-isogeny over $F$. Hence $\mathcal{G}_{2}(E/F)$ is trivial.
Finally, we observe that the claim regarding the isogeny class degree and isogeny class size follow from Proposition~\ref{prop:maxmat} and Corollary~\ref{number_of_vertices}, respectively.
\end{proof}

\begin{example}\label{ex:CMvolc}
Let $\theta=\sqrt{6+2\sqrt{21}}$, and consider the number field $F=\mathbb{Q}\left(  \theta\right)  $ with LMFDB label \href{https://www.lmfdb.org/NumberField/4.2.1323.1}{4.2.1323.1}. Let $E$ be the elliptic curve defined over $F$ with Weierstrass model
\[
E:x^{3}-(53865\theta^{3}-26460\theta^{2}-159705\theta-120015)x-4614624\theta
^{3}+2055942\theta^{2}+13780368\theta+10139472.
\]
We then check that $\operatorname{End}E\otimes_\mathbb{Z} \mathbb{Q} \cong K=\mathbb{Q}(\sqrt{-7})$ and $\operatorname*{End}E\cong\mathcal{O}=\mathbb{Z}+9\mathcal{O}_{K}$ with~$\mathcal{O}_{K}$ denoting the ring of integers of $K$. Further, $F\cong\mathbb{Q}(j(E))$ and $KF=K(j(E))$ is the ring class field of $\mathcal{O}$. In particular, the class number of $\mathcal{O}$ is $4$ since $\left[KF:K\right]  =\left[  F:\mathbb{Q}\right]  =4$. Since the hypothesis of Lemma~\ref{lem:potCM} apply, and the discriminant of $\mathcal{O}$ is $\Delta=\mathfrak{f}^{2}\Delta_{K}=-63 \equiv 1\ \operatorname{mod}8$, we obtain that
\[
\mathcal{G}(E/F)\cong\mathcal{H}_{2}^{0}\square\mathcal{H}_{9}^{0}\square\mathcal{H}_{7}^{0}.
\]
In particular, the isogeny class $V(\mathcal{G}(E/F))$ consists of $12$ elliptic curves. Computations in \textsc{SageMath} show that four distinct endomorphism rings occur among these curves. We label the vertices of $\mathcal{G}(E/F)$ by the $F$-isomorphism classes $[E_{D,m}]_{F}$, where $E_{D,m}$ denotes an elliptic curve whose endomorphism ring has discriminant $D$. With this notation, $E\cong E_{-63,1}$.
\[
 \adjustbox{scale=0.75}{\begin{tikzcd}
{[E_{-63,1}]_F} \arrow[rd, "7", no head, red] \arrow[dd, "2", no head] \arrow[rr, "3", no head, blue] &                                                         & {[E_{-7,1}]_F} \arrow[rd, "7", no head, red] \arrow[rr, "3", no head, blue] \arrow[dd, no head] &                                                        & {[E_{-63,3}]_F} \arrow[rd, "7", no head, red] \arrow[dd, no head] &                                          \\
                                                                                           & {[E_{-63,2}]_F} \arrow[rr, no head, blue] \arrow[dd, no head] &                                                                                      & {[E_{-7,2}]_F} \arrow[rr, no head, blue] \arrow[dd, no head] &                                                              & {[E_{-63,4}]_F} \arrow[dd, "2", no head] \\
{[E_{-252,1}]_F} \arrow[rd, "7", no head, red] \arrow[rr, no head, blue]                              &                                                         & {[E_{-28,1}]_F} \arrow[rr, no head, blue] \arrow[rd, "7", no head, red]                         &                                                        & {[E_{-252,3}]_F} \arrow[rd, "7", no head, red]                    &                                          \\
                                                                                           & {[E_{-252,2}]_F} \arrow[rr, "3", no head, blue]               &                                                                                      & {[E_{-28,2}]_F} \arrow[rr, "3", no head, blue]               &                                                              & {[E_{-252,4}]_F}                        
\end{tikzcd}}
\]
Among these elliptic curves, $E_{-7,1},E_{-7,2},E_{-28,1},$ and $E_{-28,2}$ are the base changes of \href{https://www.lmfdb.org/EllipticCurve/Q/49/a/4}{49.a4}, \href{https://www.lmfdb.org/EllipticCurve/Q/49/a/2}{49.a2}, \href{https://www.lmfdb.org/EllipticCurve/Q/49/a/3}{49.a3}, and \href{https://www.lmfdb.org/EllipticCurve/Q/49/a/1}{49.a1}, respectively. We now base change the isogeny class to the Hilbert class field~$KF$. Over $KF$, the following pairs of elliptic curves becomes isomorphic:
\[
\adjustbox{scale=0.95}{$\begin{array}
[c]{rclcrclcrcl}
\left[  E_{-7,1}\right]  _{KF} & = & \left[  E_{-7,2}\right]  _{KF} & \quad & \left[  E_{-63,1}\right]  _{KF} & = & \left[  E_{-63,4}\right]  _{KF} & \quad & \left[  E_{-252,1}\right]  _{KF} & = & \left[  E_{-252,4}\right]  _{KF}\\
\left[  E_{-28,1}\right]  _{KF} & = & \left[  E_{-28,2}\right]  _{KF} &  & \left[  E_{-63,2}\right]  _{KF} & = & \left[  E_{-63,3}\right]  _{KF} &  & \left[  E_{-252,2}\right]  _{KF} & = & \left[  E_{-252,3}\right]  _{KF}
\end{array}$}
\]
Further, $V(\mathcal{G}(E/KF))$ consists of ten elliptic curves, and we denote the four elliptic curves that do not descent to $F$ by $E_{-63,5},E_{-63,6},E_{-252,5},$ and $E_{-252,6}$.  For $D\in\{-63,-252\}$, the four $j$-invariants corresponding to the elliptic curves $E_{D,m}$, where $m\in\{1,2,5,6\}$, are precisely the roots of the Hilbert class polynomial $H_{D}(X)$.

Since $KF$ is the Hilbert class field of the order of discriminant $-256$, the determination of the full isogeny graph $\mathcal{G}(E/F)$ now falls under the theory of isogeny volcanoes~\cite{Clarkvolcanoes,DrewVolcanoes}. As in Example~\ref{Ex:CMnotCarPro}, elliptic curves with isomorphic endomorphism rings are joined by infinitely many horizontal isogenies. The graph below represents each infinite family of horizontal isogenies by a single dashed edge labeled with the least prime degree occurring between the corresponding vertices.
\[
 \adjustbox{scale=0.75}{\begin{tikzcd}
                                                                                                                                                     &  &                                                                                                       & {[E_{-7,1}]_{KF}} \arrow[llld, "3"', no head, bend right,blue] \arrow[rrrd, "3", no head, bend left,blue] \arrow[ld, "3"', no head, bend right,blue] \arrow[rd, "3", no head, bend left,blue] \arrow[ddddd, no head] &                                                                    &  &                                                \\
{[E_{-63,1}]_{KF}} \arrow[rr, "2", no head,dashed] \arrow[rrrrrr, no head, bend right,dashed, "\hspace{1em}2"] \arrow[rrrr, "7", no head, bend right,red] \arrow[ddddd, "2"', no head] &  & {[E_{-63,5}]_{KF}} \arrow[rr, "\hspace{1em} 2"', no head,dashed] \arrow[rrrr, "7", no head, bend right,red] \arrow[ddddd, no head] &                                                                                                                                                                                                  & {[E_{-63,2}]_{KF}} \arrow[rr, "2", no head,dashed] \arrow[ddddd, no head] &  & {[E_{-63,6}]_{KF}} \arrow[ddddd, "2", no head] \\
                                                                                                                                                     &  &                                                                                                       &                                                                                                                                                                                                  &                                                                    &  &                                                \\
                                                                                                                                                     &  &                                                                                                       &                                                                                                                                                                                                  &                                                                    &  &                                                \\
                                                                                                                                                     &  &                                                                                                       &                                                                                                                                                                                                  &                                                                    &  &                                                \\
                                                                                                                                                     &  &                                                                                                       & {[E_{-28,1}]_{KF}} \arrow[llld, "3"', no head, bend right,blue] \arrow[rrrd, "3", no head, bend left,blue] \arrow[rd, "3", no head, bend left,blue] \arrow[ld, "3"', no head, bend right,blue]                       &                                                                    &  &                                                \\
{[E_{-252,1}]_{KF}} \arrow[rr, "11", no head,teal,dashed] \arrow[rrrrrr, "11"', no head, bend right,teal,dashed] \arrow[rrrr, "7", no head, bend right,red]                     &  & {[E_{-252,5}]_{KF}} \arrow[rr, "11", no head,teal,dashed] \arrow[rrrr, "7", no head, bend right,red]                 &                                                                                                                                                                                                  & {[E_{-252,2}]_{KF}} \arrow[rr, "11", no head,teal,dashed]                      &  & {[E_{-252,6}]_{KF}}                           
\end{tikzcd}}
\]
In particular, $\mathcal{G}_{2}(E_{-63,1}/KF)$ is the following $2$-volcano of depth $1$, whose crater consists of the elliptic curves with endomorphism ring of discriminant $-63$.
\[
 \adjustbox{scale=0.8}{\begin{tikzcd}
                                                                                                                      & {[E_{-63,5}]_{KF}} \arrow[rd, "2"', no head, bend left] \arrow[d, "2", no head] &                                            \\
{[E_{-63,1}]_{KF}} \arrow[ru, "2"', no head, bend left] \arrow[rd, "2", no head, bend right] \arrow[d, "2"', no head] & {[E_{-252,5}]_{KF}}                                                             & {[E_{-63,2}]_{KF}} \arrow[d, "2", no head] \\
{[E_{-252,1}]_{KF}}                                                                                                   & {[E_{-63,6}]_{KF}} \arrow[ru, "2", no head, bend right] \arrow[d, "2", no head] & {[E_{-252,2}]_{KF}}                        \\
                                                                                                                      & {[E_{-252,6}]_{KF}}                                                             &                                           
\end{tikzcd}}
\]
\end{example}

The only cases not covered by lemma~\ref{lem:potCM} are those with $\Delta \in \{-3,-4\}$. These cases are exceptional, since the corresponding elliptic curves admit twists that are not quadratic. They have already been treated in \cite[Table~5]{MR4203041}. More precisely, the authors obtained the following classification for their isogeny graphs.
\renewcommand{\arraystretch}{1.3}
\begin{center}
\begin{tabular}{ccc}

\(\Delta\) & \(E\) & \(\mathcal{G}(E/\mathbb{Q})\) \\
\hline
\multirow{3}{*}{\(-3\)}
& \(y^2=x^3+t^3\) & \(\mathcal{H}^0_2 \square \mathcal{H}^0_3\) \\\cmidrule{2-3}
& \(y^2=x^3+16t^3\) & \(\mathcal{H}^0_{27}\) \\\cmidrule{2-3}
& \(y^2=x^3+s,\ s\neq t^3,16t^3\) & \(\mathcal{H}^0_3\) \\
\hline
\multirow{2}{*}{\(-4\)}
& \(y^2=x^3\pm t^2x\) &  \(\mathcal{H}^1_4\) \\\cmidrule{2-3}
& \(y^2=x^3+sx,\ s\neq \pm t^2\) & \(\mathcal{H}^0_2\)  \\
\hline
\end{tabular}
\end{center}
In addition, the distribution of the $p$-adic valuations of the conductors corresponding to the endomorphism rings occurring in the possible $p$-primary graphs are well-known. The lemma below summarizes this information.

\begin{lemma}\label{lem:specialdiscpotcm}
Let $\mathcal{O}$ be an order with discriminant $\Delta\in\{-3,-4\}$. Let $E$ be an elliptic curve such that $\operatorname*{End}E\cong\mathcal{O}$. In particular, $j(E)\in\{0,1728\}$. Suppose further that $E$ is defined over~$\mathbb{Q}$. For a prime $p$, Table~\ref{ta:potcmendvp34} gives a graph isomorphic to the $p$-primary graph $\mathcal{G}_{p}(E/\mathbb{Q})$. The table further gives, for each vertex $E_{i}$ of $\mathcal{G}_{p}(E/\mathbb{Q})$, the $p$-adic valuation of $\mathfrak{f}_{i}$, where $\mathfrak{f}_{i}$ denotes the conductor of $\operatorname*{End}E_{i}$.
\end{lemma}
{\begingroup \small
\renewcommand{\arraystretch}{1.3}
 \begin{longtable}{ccccc}
 	\caption{For $E/\mathbb{Q}$ with $\operatorname*{End}E$ having discriminant $\Delta\in\{-3,-4\}$, the table below gives the possible $p$-primary graphs $\mathcal{G}_{p}(E/\mathbb{Q})$ for each prime $p$. In addition, for each vertex $E_{i}$ in $\mathcal{G}_{p}(E/\mathbb{Q})$, let $\mathfrak{f}_{i}$ denote the conductor of $\operatorname*{End}E_{i}$. Then, $(v_{p}(\mathfrak{f}_{i}))_{i}$ is as given with respect to the chosen orientation in $\mathcal{G}_{p}(E/\mathbb{Q})$.}\label{ta:potcmendvp34}\\
	\hline
	$\Delta$ & $p$ & \multicolumn{2}{c}{$\mathcal{G}_{p}(E/\mathbb{Q})$} & $\left(  v_{p}(\mathfrak{f}_{i})\right)  _{i}$ \\
	\hline

	\endfirsthead
	\hline
	$\Delta$ & $p$ & \multicolumn{2}{c}{$\mathcal{G}_{p}(E/\mathbb{Q})$} & $\left(  v_{p}(\mathfrak{f}_{i})\right)  _{i}$  \\
	\hline
	\endhead
	\hline

	\multicolumn{5}{r}{\emph{continued on next page}}
	\endfoot
	\hline
	\endlastfoot

$-3$ & $\geq5$ & $\mathcal{H}_{1}^{0}$ & \adjustbox{scale=1}{\begin{tikzcd}
E_0
\end{tikzcd}} & $\left(  0\right)  $\\\cmidrule{2-5}
& $3$ & $\mathcal{H}_{3}^{0}$ & \adjustbox{scale=1}{\begin{tikzcd}
E_0 \arrow[r, "3", no head] & E_1
\end{tikzcd}} & $\left(  0,0\right)  $\\\cmidrule{3-5}
&  & $\mathcal{H}_{27}^{0}$ &  
\adjustbox{scale=1}{\begin{tikzcd}
E_1 \arrow[r, "3", no head] \arrow[d, "3"', no head] & E_2 \arrow[d, "3", no head] \\
E_0                                                  & E_3                        
\end{tikzcd}}
& $\left(  1,0,0,1\right)  $\\\cmidrule{2-5}
& $2$ & $\mathcal{H}_{2}^{0}$ &
\adjustbox{scale=1}{\begin{tikzcd}
E_1                         \\
E_0 \arrow[u, "2", no head]
\end{tikzcd}}
& $\left(  1,0\right)  $\\\hline
$-4$ & $\geq3$ & $\mathcal{H}_{1}^{0}$ & \adjustbox{scale=1}{\begin{tikzcd}
E_0
\end{tikzcd}} & $\left(  0\right)  $\\\cmidrule{2-5}
& $2$ & $\mathcal{H}_{2}^{0}$ & 
\adjustbox{scale=1}{\begin{tikzcd}
E_0 \arrow[r, "2", no head] & E_1
\end{tikzcd}}
& $\left(  0,0\right)  $\\\cmidrule{3-5}
& 
& $\mathcal{H}_{4}^{1}$ & 
\adjustbox{scale=1}{\begin{tikzcd}
E_1 \arrow[r, "2", no head] \arrow[d, "2"', no head] \arrow[rd, "2", no head] & E_2 \\
E_0                                                                           & E_3
\end{tikzcd}}
& $\left(  1,0,0,1 \right)  $

\end{longtable}
\endgroup}

Our next result extends Lemma~\ref{lem:potCM} by determining the distribution of the $p$-adic valuations of the conductors corresponding to the endomorphism rings occurring in $\mathcal{G}_p(E/F)$ along a $k$-spine $(E_i)_{i=0}^k$ in $\mathcal{G}_p(E/F)$ with $E\cong E_0$.

\begin{theorem}\label{thm:CMdiscclass}
Let $\mathcal{O}$ be an order with discriminant $\Delta<-4$ in a quadratic field $K$, and write $\Delta=\mathfrak{f}^{2}\Delta_{K}$ where $\mathfrak{f}$ and $\Delta_{K}$ denote the conductor of $\mathcal{O}$ and discriminant of $K$, respectively. Let $E$ be an elliptic curve such that $\operatorname*{End} E\cong\mathcal{O}$. Suppose further that $E$ is defined over $F=\mathbb{Q}(j(E))$. For a prime $p$, Table~\ref{ta:potcmendvp} gives a graph isomorphic to the $p$-primary graph $\mathcal{G}_{p}(E/F)$ of $E$. Further, with $k=v_{p}(\deg\mathcal{G}(E/F))$ and $(E_{i})_{i=0}^{k}$  a $k$-spine in $\mathcal{G}_{p}(E/F)$ with $E\cong E_{0}$, loc. cit. gives the value of $v_{p}(\mathfrak{f}_{i})$ for $0\leq i\leq k$, where $\mathfrak{f}_{i}$ denotes the conductor of the order $\operatorname*{End}E_{i}$.
\end{theorem}

{\begingroup \small
\renewcommand{\arraystretch}{1.9}
 \begin{longtable}{cccc}
 	\caption{For a $k$-spine $\left(  E_{i}\right)  _{i=0}^{k}$ of $\mathcal{G}_{p}(E/F)\cong\mathcal{H}_{p^{k}}^{r}$ with $E\cong E_{0}$, the conditions below determine $\mathcal{G}_p(E/F)$ and $v_{p}(\mathfrak{f}_{i})$, where $\mathfrak{f}_{i}$ denotes the conductor of $\operatorname*{End}E_{i}$.}\label{ta:potcmendvp}\\
	\hline
	$p$ & Conditions  & $\mathcal{G}_{p}(E/F)$ & $v_{p}(\mathfrak{f}_{i})$  \\
	\hline

	\endfirsthead
	\hline
	$p$ & Conditions & $\mathcal{G}_{p}(E/F)$ & $v_{p}(\mathfrak{f}_{i})$  \\
	\hline
	\endhead
	\hline

	\multicolumn{4}{r}{\emph{continued on next page}}
	\endfoot
	\hline
	\endlastfoot

$p\ge 3$ & $v_p(\Delta)\ge 0$ & $\mathcal{H}_{p^{v_{p}(\Delta)}}^{0}$ & $\max\!\left\{
\left\lfloor \frac{v_{p}(\Delta)}{2}\right\rfloor -i,i-\left\lceil \frac
{v_{p}(\Delta)}{2}\right\rceil \right\}  $\\\hline
$2$ & $v_{2}(\mathfrak{f})=0$ and $\Delta_{K}\equiv5\ \operatorname{mod}8$ &
 $\mathcal{H}_{1}^{0}$ & $0$\\\cmidrule{2-4}
& $v_{2}(\mathfrak{f})=1$ and $\Delta_{K}\equiv1\ \operatorname{mod}4$  &
$\mathcal{H}_{2}^{0}$ & $1-i$\\\cmidrule{2-4}
& $v_{2}(\mathfrak{f})=0$ and $\Delta_{K}\equiv1\ \operatorname{mod}8$ & 
$\mathcal{H}_{2}^{0}$ & $i$\\\cmidrule{2-4}
& $v_{2}(\mathfrak{f})=0$ and $\Delta_{K}\equiv0\ \operatorname{mod}4$ &
$\mathcal{H}_{2}^{0}$ & $0$\\\cmidrule{2-4}
& $v_{2}(\Delta)\geq4$ and $\Delta_{K}\equiv1\ \operatorname{mod}4$ & $\mathcal{H}_{2^{v_{2}(\Delta)-2}}^{1}$ & $\max\!\left\{  \frac
{v_{2}(\Delta)}{2}-i,i+\frac{4-v_{2}(\Delta)}{2}\right\}  $\\\cmidrule{2-4}
& $v_{2}(\Delta)\geq4$ and $\Delta_{K}\equiv0\ \operatorname{mod}4$ &
 $\mathcal{H}_{2^{v_{2}(\Delta)-2}}^{1}$ & $\max\!\left\{  \left\lfloor
\frac{v_{p}(\Delta)-2}{2}\right\rfloor -i,i-\left\lceil \frac{v_{p}(\Delta)-2}{2}\right\rceil \right\}  $

\end{longtable}
\endgroup}

\begin{proof}
By Lemma~\ref{lem:potCM}, the third column of Table~\ref{ta:potcmendvp} is as claimed. It thus suffices to show that the fourth column is as claimed. We note that the cases $\left(  i\right)  $ $p$ odd and $(ii)$ $p=2$ with $v_{2}(\mathfrak{f})=0$ and $\Delta_{K}\equiv5\ \operatorname{mod}8$ are automatic since $\mathcal{G}_{p}(E/F)$ is trivial in these cases. Further, the proof of loc. cit. shows that $E$ is a leaf in $\mathcal{G}_{p}(E/F)$ since it does not admit two independent $p$-isogenies. To establish the remaining cases, we begin with a review of $p$-volcanoes, and assume the terminology in~\cite[pp. 1-7]{DrewVolcanoes}. Since $\operatorname*{End}E$ is an order in an imaginary quadratic field, we have that $\mathcal{G}_{p}(E/\overline{F})$ is a $p$-volcano and $E$ is $p$-bloomed in $\mathcal{G}_{p}(E/\overline{F})$. Each $p$-isogeny $\phi:E\rightarrow E^{\prime}$ is either horizontal or vertical. In what follows, let $\mathfrak{f}_{E}$ (resp. $\mathfrak{f}_{E^{\prime}}$) denote the conductor of $\operatorname*{End}E$ (resp. $\operatorname*{End}E^{\prime}$).

An isogeny $\phi$ is horizontal if $\mathfrak{f}_{E}=\mathfrak{f}_{E^{\prime}}$, and $E/\overline{F}$ admits $h$ horizontal $p$-isogenies, where
\begin{equation}
h=\left\{
\begin{array}
[c]{cl}
0 & \text{if }p\text{ is inert in }K\text{ or }v_{p}(\mathfrak{f}_{E})\geq1,\\
1 & \text{if }p\text{ is ramified in }K\text{ and }v_{p}(\mathfrak{f}_{E})=0,\\
2 & \text{if }p\text{ splits in }K\text{ and }v_{p}(\mathfrak{f}_{E})=0.
\end{array}
\right.  \label{eq:hvalcoro}
\end{equation}
In particular, if $E$ admits a horizontal $p$-isogeny, then it is necessary that $\operatorname*{End}E$ be maximal at~$p$. That is, $v_{p}(\left[\mathcal{O}_{K}:\operatorname*{End}E\right]  )=0$. We further note that if $p$ is ramified in $K$ and $v_{p}(\mathfrak{f})=0$, then the unique horizontal $p$-isogeny admitted by $E$ is $F$-rational. 

Now suppose that $\phi$ is a vertical isogeny. Then, it is either an ascending or descending. We say that $\phi$ is ascending if $\left[  \operatorname*{End}E^{\prime}:\operatorname*{End}E\right]  =\frac{\mathfrak{f}_{E}}{\mathfrak{f}_{E^{\prime}}}=p$. Similarly, we say that $\phi$ is descending if $\left[  \operatorname*{End}E:\operatorname*{End}E^{\prime}\right]=\frac{\mathfrak{f}_{E^{\prime}}}{\mathfrak{f}_{E}}=p$. If $v_{p}(\mathfrak{f}_{E})=0$, then the number vertical isogenies admitted by $E/\overline{F}$ is $p+1-h$. So suppose instead that $v_{p}(\mathfrak{f}_{E})\geq1$. By~(\ref{eq:hvalcoro}), each $p$-isogeny admitted by $E$ must be vertical. In fact, $E$ admits exactly one ascending $p$-isogeny \cite[Lemma~6]{DrewVolcanoes}, and it is $F$-rational. In sum, if $d$ denotes the number of descending $p$-isogenies admitted by $E/\overline{F}$, then
\[
p+1=\left\{
\begin{array}
[c]{cl}
h+d & \text{if }v_{p}(\mathfrak{f}_{E})=0,\\
1+d & \text{if }v_{p}(\mathfrak{f}_{E})\geq1.
\end{array}
\right.
\]

We now return to considering the $p$-primary graph $\mathcal{G}_{p}(E/F)$ for the remaining cases. Since $\operatorname*{End}_{F}E\cong\mathbb{Z}$, we have that $\operatorname*{End}_{F}E^{\prime}\cong\mathbb{Z}$ for each $E^{\prime}\in V(\mathcal{G}_{p}(E/F))$. Therefore, by the above discussion, it follows that for $E^{\prime}\in V(\mathcal{G}_{p}(E/F))$, the number $h_{F}(E^{\prime})$ of $F$-rational horizontal $p$-isogenies admitted by $E^{\prime}$ is
\begin{equation}
h_{F}(E^{\prime})=\left\{
\begin{array}
[c]{cl}
0 & \text{if }v_{p}(\Delta_{K})=0,\\
1 & \text{if }v_{p}(\Delta_{K})\geq1\text{ and }v_{p}(\mathfrak{f}_{E^{\prime
}})=0.
\end{array}
\right.  \label{eq:valofh0}
\end{equation}
Here we note that if $p$ splits in $K$, the two corresponding horizontal~$p$-isogenies admitted by $E/\overline{F}$ have field of definition $KF$, and are not $F$-rational.

For each of the remaining cases, we have that $k\geq1$ and $\mathcal{G}_{p}(E/F)\cong\mathcal{H}_{p^{k}}^{r}$ with $r\in\{0,1\}$ are as given in the third and fourth column, respectively, of Table~\ref{ta:potcmendvp}. By the proof of Lemma~\ref{lem:potCM}, $E$ admits a $F$-rational $p^{k}$-isogeny. In particular, there exists a $k$-spine $\left(  E_{i}\right)  _{i=0}^{k}$ in $\mathcal{G}_{p}(E/F)$ such that $E\cong E_{0}$. Next, for $0\leq i<k$, let $\phi_{i}:E_{i}\rightarrow E_{i+1}$ be the corresponding $F$-rational $p$-isogeny. Now set $\mathfrak{f}_{i}=\mathfrak{f}_{E_i}$, so that $\Delta=\mathfrak{f}_{0}^{2}\Delta_{K}$. We now proceed by cases.

\textbf{Case 1.} Suppose that $p$ is odd. If $v_{p}(\Delta)=1$, then $\mathcal{G}_{p}(E/F)\cong\mathcal{H}_{p}^{0}$, $v_{p}(\mathfrak{f}_{0})=0$, and $p$ is ramified in $K$. In particular, by (\ref{eq:valofh0}), $\phi_{0}$ must be the unique horizontal $p$-isogeny. Hence, $v_{p}(\mathfrak{f}_{0})=v_{p}(\mathfrak{f}_{1})=0$. So suppose that $k=v_{p}(\Delta)\geq2$ so that $\mathcal{G}_{p}(E/F)\cong\mathcal{H}_{p^{k}}^{0}$. Since $p$ is odd, $v_{p}(\Delta_{K})\leq1$, and so $v_{p}(\mathfrak{f}_{0})\geq1$. This implies that $\phi_{0}$ is the unique ascending $p$-isogeny. Now set
\begin{equation}
t=\min\left\{  i\mid v_{p}(\mathfrak{f}_{i})=0\right\}
.\label{eg:valtinpotcm}
\end{equation}
Then, for $0\leq i<t$, we have that $\phi_{i}$ is the unique ascending $p$-isogeny admitted by $E_{i}$ since $v_{p}(\mathfrak{f}_{i})\geq1$. It follows that $\mathfrak{f}_{i}=p\mathfrak{f}_{i+1}$, and thus $v_{p}(\mathfrak{f}_{i})=t-i$ for $0\leq i\leq t$.\ In particular, $v_{p}(\mathfrak{f}_{0})=t$, and so
\[
k=v_{p}(\Delta)=\left\{
\begin{array}
[c]{cl}
2t & \text{if }v_{p}(\Delta_{K})=0,\\
2t+1 & \text{if }v_{p}(\Delta_{K})=1.
\end{array}
\right.
\]
In particular, $t=\left\lfloor k/2\right\rfloor <k$. Now observe that since $v_{p}(\mathfrak{f}_{t})=0$, $E_{t}$ admits $h_{F}(E_{t})$ horizontal $p$-isogenies. Consequently, $\phi_{t}$ is descending (resp. horizontal) if $v_{p}(\Delta_{K})=0$ (resp. $1$). Now set
\begin{equation}
t^{\ast}=\left\{
\begin{array}
[c]{cl}
t & \text{if }v_{p}(\Delta_{K})=0,\\
t+1 & \text{if }v_{p}(\Delta_{K})=1,
\end{array}
\right.  \label{eq:valtstar}
\end{equation}
so that $v_{p}(\mathfrak{f}_{t})=v_{p}(\mathfrak{f}_{t^{\ast}})$. By construction, we have that for $t^{\ast}\leq i<k$, the $p$-isogeny $\phi_{i}$ is descending. Indeed, the unique ascending $p$-isogeny admitted by $E_{i}$ is the dual $\widehat{\phi}_{i-1}$ of $\phi_{i-1}$. It follows that $p\mathfrak{f}_{i}=\mathfrak{f}_{i+1}$ for $t^{\ast}\leq i\leq k$. This, in turn, implies that $v_{p}(\mathfrak{f}_{i})=i-t^{\ast}$ for $t^{\ast}\leq i\leq k$. Since $k-t^{\ast}=t$, we observe that $\left\lfloor k/2\right\rfloor -i=v_{p}(\mathfrak{f}_{i})=v_{p}(\mathfrak{f}_{k-i})$ for $0\leq i\leq\frac {k}{2}$. Below, we illustrate the two possible $k$-spines for $v_{p}(\Delta_{K})\in\{0,1\}$.
\[
 \adjustbox{scale=0.72}{\begin{tikzcd}
                                                                 & v_p(\Delta_K)=0                     &                                              & v_p(\mathfrak{f}_i)                     &                                                                                  & v_p(\Delta_K)=1 &                                                                         \\
                                                                 & E_{\frac{k}{2}} \arrow[rd, no head] &                                              & 0                                       & E_{\left\lfloor \frac{k}{2}\right\rfloor} \arrow[d, no head] \arrow[rr, no head] &                 & E_{\left\lfloor \frac{k}{2}\right\rfloor+1} \arrow[d, no head]          \\
E_{\frac{k}{2}-1} \arrow[ru, no head] \arrow[d, no head, dotted] &                                     & E_{\frac{k}{2}+1} \arrow[d, no head, dotted] & 1 \arrow[d, no head, dotted]            & E_{\left\lfloor \frac{k}{2}\right\rfloor-1} \arrow[d, no head, dotted]           &                 & E_{\left\lfloor \frac{k}{2} \right\rfloor+2} \arrow[d, no head, dotted] \\
E_1 \arrow[d, no head]                                           &                                     & E_{k-1} \arrow[d, no head]                   & \left\lfloor \frac{k}{2}\right\rfloor-1 & E_1 \arrow[d, no head]                                                           &                 & E_{k-1} \arrow[d, no head]                                              \\
E_0                                                              &                                     & E_k                                          & \left\lfloor \frac{k}{2}\right\rfloor   & E_0                                                                              &                 & E_k                                                                    
\end{tikzcd}}
\]
The case now follows from the following observation:
\[
v_{p}(\mathfrak{f}_{i})=\left\{
\begin{array}
[c]{cl}
\left\lfloor \frac{k}{2}\right\rfloor -i & \text{if }0\leq i<t,\\
0 & \text{if }t\in\{t,t^{\ast}\}\\
i-\left\lceil k/2\right\rceil  & \text{if }t^{\ast}<i\leq k
\end{array}
\right.  =\max\left\{  \left\lfloor \frac{v_{p}(\Delta)}{2}\right\rfloor
-i,i-\left\lceil \frac{v_{p}(\Delta)}{2}\right\rceil \right\}  .
\]

\textbf{Case 2.} Suppose $p=2$ with $v_{2}(\mathfrak{f}_{0})=1$ and $\Delta_{K}\equiv1\ \operatorname{mod}4$. Then, $\phi_{0}$ must be an ascending $2$-isogeny since $E_{0}$ is a leaf. The case now follows since $v_{2}(\mathfrak{f}_{1})=0$ and $\mathcal{G}_{p}(E/F)\cong\mathcal{H}_{2}^{0}$.

\textbf{Case 3.} Suppose $p=2$ with $v_{2}(\mathfrak{f})=0$ and $\Delta_{K}\equiv1\ \operatorname{mod}8$. Since $2$ does not ramify in $K$ and $\operatorname*{End}E_{0}$ is maximal at $p$, it follows that $\phi_{0}$ must be a descending $2$-isogeny. The case is thus established, since $v_{2}(\mathfrak{f}_{1})=1$ and $\mathcal{G}_{p}(E/F)\cong\mathcal{H}_{2}^{0}$.

\textbf{Case 4.} Suppose $p=2$ with $v_{2}(\mathfrak{f})=0$ and $\Delta_{K}\equiv0\ \operatorname{mod}4$. Since $2$ is ramified in $K$, $\phi_{0}$ is a horizontal $2$-isogeny, and thus $v_{p}(\mathfrak{f}_{0})=v_{p}(\mathfrak{f}_{1})=0$. Since, $\mathcal{G}_{p}(E/F)\cong\mathcal{H}_{2}^{0}$ this concludes this case.

\textbf{Case 5.} Suppose $p=2$ with $v_{2}(\Delta)\geq4$ and $\Delta_{K}\equiv1\ \operatorname{mod}4$. Since $v_{2}(\Delta)=2v_{2}(\mathfrak{f}_{0})$, we have that $v_{2}(\mathfrak{f}_{0})\geq2$. In this setting we have that $\mathcal{G}_{2}(E/F)\cong\mathcal{H}_{2^{k}}^{1}$ with $k=v_{2}(\Delta)-2$. Now let
\[
t=\min\left\{  i\mid v_{p}(\mathfrak{f}_{i})=1\right\}  .
\]
By an argument similar to that given following (\ref{eg:valtinpotcm}), we deduce that $v_{2}(\mathfrak{f}_{i})=t+1-i$ for $0\leq i\leq t$. In particular, $v_{p}(\mathfrak{f}_{0})=t+1$ and thus $k=v_{2}(\Delta)-2=2t$. In particular, $t\geq2$ and the isogeny from $E_{t}$ to $E_{k}$ has degree $2^{t}$. Consequently, $E_{t}$ admits three $2$-isogenies, one of which is the descending isogeny $\widehat{\phi}_{i-1}$. Thus, $\phi_{t}$ is either ascending or descending. We claim that it is descending. By way of contradiction, suppose $\phi_{t}$ is ascending so that $v_{p}(\mathfrak{f}_{t+1})=0$. Since $2$ is unramified in $K$, it follows that $\phi_{t+1}$ is a descending $2$-isogeny. It follows that for $t+1\leq i\leq2t=k$, the isogeny $\phi_{i}$ is descending. Consequently, $v_{2}(\mathfrak{f}_{k})=t-1$, and thus $\mathfrak{f}_{0}=4\mathfrak{f}_{k}$. Now observe that the order $\operatorname*{End}E_{k}$ has discriminant $\mathfrak{f}_{k}^{2}\Delta_{K}$ and satisfies $v_{2}(\mathfrak{f}_{k}^{2}\Delta_{K})=2(t-1)\geq2$.

Now consider the proof of the forward direction of \cite[Theorem~4.1~(1)]{Kwon} with $a=\frac{\mathfrak{f}_{0}}{\mathfrak{f}_{k}}=4$ and $b=\frac{2^{2t}}{a}=2^{2t-2}$. With notation as in loc. cit., we have that the existence of the $2^{2t}$-isogeny $\phi:E_{0}\rightarrow E_{k}$ implies that $b\in S_{\mathfrak{f}_{k}}^{\prime}$, where $S_{\mathfrak{f}_{k}}^{\prime}$ is the set of all integers which are indexes of primitive proper ideals of $\operatorname*{End}E_{k}$ invariant under complex conjugation. For a prime $p$, let
\begin{equation}
A_{p}=\left\{
\begin{array}
[c]{cl}
p^{v_{p}(\mathfrak{f}_{k}^{2}\Delta_{K})} & \text{if }p\text{ is odd,}\\
\left\{  4,2^{v_{2}(\mathfrak{f}_{k}^{2}\Delta_{K})-2}\right\}   & \text{if }p=2\text{ with }16|\mathfrak{f}_{k}^{2}\Delta_{K}\text{,}\\
\left\{  2\right\}   & \text{if }p=2\text{ with }16\nmid\mathfrak{f}_{k}^{2}\Delta_{K}\text{, }\mathfrak{f}_{k}\text{ odd, and }2|\Delta_{K},\\
\varnothing & \text{otherwise.}
\end{array}
\right.  \label{eq:KwonAp1}
\end{equation}
By \cite[Proposition 3.1]{Kwon},
\begin{equation}
S_{\mathfrak{f}_{k}}^{\prime}=\left\{  \prod_{p|\mathfrak{f}_{k}^{2}\Delta_{K}}n_{p}\mid n_{p}\in\left\{  1\right\}  \cup A_{p}\right\}
.\label{eq:KwonAp}
\end{equation}
Now recall that $b=2^{2t-2}\in S_{\mathfrak{f}_{k}}^{\prime}$ with $t\geq2$. If $t=2$, then $v_{2}(\mathfrak{f}_{k}^{2}\Delta_{K})=2$ and so $b=1$. But this contradicts that $ab=4b=2^{2t}$. So suppose that $t\geq3$ so that $v_{2}(\mathfrak{f}_{k}^{2}\Delta_{K})=2(t-1)\geq4$. Then, $b\in \{1,4,2^{2(t-2)}\}$. But this is impossible since $b=2^{2t-2}\geq16$. Consequently, by \cite[Theorem 4.1 (1)]{Kwon}, there can be no $2^{2t}$-isogeny $\phi:E_{0}\rightarrow E_{k}$. This establishes the claim, and so $\phi_{t}$ must be descending. Now observe that the argument following (\ref{eq:valtstar}), implies that $\phi_{i}$ is descending for $t\leq i<k$. In particular, $v_{2}(\mathfrak{f}_{i})=1+i-t$ for $t\leq i\leq k=2t$. In sum,
\begin{align*}
v_{2}(\mathfrak{f}_{i})  & =\left\{
\begin{array}
[c]{cl}
t+1-i & \text{if }0\leq i\leq t,\\
1+i-t & \text{if }t\leq i\leq k,
\end{array}
\right.  \\
& =\left\{
\begin{array}
[c]{cl}
\frac{v_{2}(\Delta)}{2}-i & \text{if }0\leq i\leq t,\\
i+\frac{4-v_{2}(\Delta)}{2} & \text{if }t\leq i\leq k,
\end{array}
\right.  \\
& =\max\!\left\{  \frac{v_{2}(\Delta)}{2}-i,i+\frac{4-v_{2}(\Delta)}
{2}\right\}  .
\end{align*}
Below, we illustrate the graph of $\mathcal{G}_{2}(E/F)$, and the $2$-adic valuations of the corresponding conductors at each stage. We note that $E'_i$ is the $p$-isogenous elliptic curve to $E_i$ for $0<i<k$ such that the corresponding $p$-isogeny is neither $\phi_i$ or $\widehat{\phi}_{i-1}$.
\[
\adjustbox{scale=0.82}{\begin{tikzcd}
v_2(\mathfrak{f}_i)          &                    &                                                           &                                                        &                                                                   &                    \\
0                            &                    &                                                           & E'_{\frac{k}{2}}                                       &                                                                   &                    \\
1                            &                    &                                                           & E_{\frac{k}{2}} \arrow[u, no head] \arrow[rd, no head] &                                                                   &                    \\
2                            &                    & E_{\frac{k}{2}-1} \arrow[ld, no head] \arrow[ru, no head] &                                                        & E_{\frac{k}{2}+1} \arrow[dd, no head, dotted] \arrow[rd, no head] &                    \\
3 \arrow[d, no head, dotted] & E'_{\frac{k}{2}-1} &                                                           &                                                        &                                                                   & E'_{\frac{k}{2}+1} \\
\frac{k}{2}-1                &                    & E_1 \arrow[ld, no head] \arrow[uu, no head, dotted]       &                                                        & E_{k-1} \arrow[d, no head] \arrow[rd, no head]                    &                    \\
\frac{k}{2}                  & E'_1               & E_0 \arrow[u, no head]                                    &                                                        & E_k                                                               & E'_{k-1}          
\end{tikzcd}}
\]

\textbf{Case 6.} Suppose $p=2$ with $v_{2}(\Delta)\geq4$ and $\Delta_{K}\equiv0\ \operatorname{mod}4$. Note that $v_{2}(\Delta)=2v_{2}(\mathfrak{f}_{0})+v_{2}(\Delta_{K})$ where $v_{2}(\Delta_{K})\in\{2,3\}$. In particular, $v_{2}(\Delta)$ is even if and only if $v_{2}(\Delta_{K})=2$. Also observe that $\mathcal{G}_{2}(E/F)\cong\mathcal{H}_{2^{k}}^{1}$ where $k=2v_{2}(\Delta)-2\geq2$. Now let $t$ be as given in~(\ref{eg:valtinpotcm}), and observe that the argument following the equation holds verbatim in our setting, to imply that~$\phi_{i}$ for $0\leq i\leq t-1$ is ascending and thus $v_{p}(\mathfrak{f}_{i})=t-i$ for $0\leq i\leq t-1$. Since $v_{p}(\mathfrak{f}_{0})=t$, we have that
\begin{equation}
k=v_{2}(\Delta)-2=2(t-1)+v_{2}(\Delta_{K}).\label{eq:valkpotcmc6}
\end{equation}
Note that since $v_{2}(\Delta)\geq4$ by assumption, we have that $t\geq2$. Next, we observe that since $v_{p}(\mathfrak{f}_{t-1})=1$, it is the case that $\phi_{t-1}$ is either ascending or descending. Since $2$ is ramified in $K$, the elliptic curve in the codomain of the ascending isogeny must admit a horizontal isogeny. In particular, it is not a leaf. It follows that the elliptic curve in the codomain of the descending isogeny must be a leaf, and so $\phi_{t-1}$ must be an ascending isogeny.

Next, we consider the isogeny $\phi_{t}$, which must either be horizontal or descending. We first suppose that $\phi_{t}$ is descending. Then, $\phi_{i}$ for $t\leq i<k$ is descending. It is this checked that $v_{p}(\mathfrak{f}_{i})=i-t$ for $t\leq i\leq k$, and thus $v_{p}(\mathfrak{f}_{k})=k-t=t+v_{2}(\Delta)-2$ by~(\ref{eq:valkpotcmc6}). Next, let
\[
a=\frac{\mathfrak{f}_{k}}{\mathfrak{f}_{0}}=p^{v_{2}(\Delta)-2}\qquad \text{and}\qquad b=\frac{2^{k}}{a}=2^{2t}.
\]
Since $E_{k}$ admits a $p^{k}$-isogeny onto $E_{0}$, we have by \cite[Theorem~4.1 (1)]{Kwon}, that $b\in S_{\mathfrak{f}_{0}}^{\prime}$. Since
\[
v_{2}(\mathfrak{f}_{0}^{2}\Delta_{K})-2=2t-2+v_{2}(\Delta_{K})\geq2,
\]
we have from (\ref{eq:KwonAp1}) and (\ref{eq:KwonAp}) that
\[
b\in\left\{  1,4,2^{2t-2+v_{2}(\Delta_{K})}\right\}  .
\]
Since $t\geq2$, we have that $v_{2}(b)\geq4$ and so $b=2^{2t-2+v_{2}(\Delta_{K})}$. Consequently, $\phi_{t}$ is descending if and only if $v_{2}(\Delta_{K})=2$. We have thus shown that if $v_{2}(\Delta_{K})=2$, then $k=2t$ and
\[
v_{p}(\mathfrak{f}_{i})=\left\{
\begin{array}
[c]{cl}
t-i & \text{if }0\leq i\leq t\\
i-t & \text{if }t<i\leq k
\end{array}
\right.  =\max\!\left\{  \frac{k}{2}-i,i-\frac{k}{2}\right\}  .
\]
It remains to consider the case when $v_{2}(\Delta_{K})=3$, in which case
$\phi_{t}$ is horizontal. Then, $v_{2}(\mathfrak{f}_{t+1})=~0$ and thus,
$E_{t+1}$ admits a unique horizontal $p$-isogeny. However, this unique isogeny
is simply $\widehat{\phi}_{t}$, and so $\phi_{t+1}$ is descending.
Consequently, each isogeny $\phi_{i}$ for $t<i<k$ is descending. The graphs
below illustrate the two possible scenarios corresponding to $v_{2}(\Delta
_{K})\in\{2,3\}$, where $E_{i}^{\prime}$ denotes the $p$-isogenous elliptic
curve to $E_{i}$ for $0<i<k$ such that the corresponding $p$-isogeny is
neither $\phi_{i}$ or $\widehat{\phi}_{i-1}$.
\[
\adjustbox{scale=0.595}{\begin{tikzcd}
                   &                                                                                       & v_2(\Delta_K)=2                                        &                                                                   &                    & v_2(\mathfrak{f}_i)                      &                                                                                      &                                                                                  & v_2(\Delta_K)=3 &                                                                                    &                                                                                              \\
                   &                                                                                       & E_{\frac{k}{2}} \arrow[rd, no head] \arrow[r, no head] & E'_{\frac{k}{2}}                                                  &                    & 0                                        &                                                                                      & E_{\left\lfloor \frac{k}{2}\right\rfloor} \arrow[d, no head] \arrow[rr, no head] &                 & E_{\left\lfloor \frac{k}{2}\right\rfloor+1} \arrow[rd, no head] \arrow[d, no head] &                                                                                              \\
                   & E_{\frac{k}{2}-1} \arrow[ru, no head] \arrow[dd, no head, dotted] \arrow[ld, no head] &                                                        & E_{\frac{k}{2}+1} \arrow[rd, no head] \arrow[dd, no head, dotted] &                    & 1                                        & E_{\left\lfloor \frac{k}{2} \right\rfloor-1} \arrow[rd, no head] \arrow[ru, no head] & E'_{\left\lfloor \frac{k}{2} \right\rfloor}                                      &                 & E'_{\left\lfloor \frac{k}{2}\right\rfloor+1}                                       & E_{\left\lfloor \frac{k}{2} \right\rfloor+2} \arrow[ld, no head] \arrow[dd, no head, dotted] \\
E'_{\frac{k}{2}-1} &                                                                                       &                                                        &                                                                   & E'_{\frac{k}{2}+1} & 2 \arrow[d, no head, dotted]             &                                                                                      & E'_{\left\lfloor \frac{k}{2} \right\rfloor-1}                                    &                 & E'_{\left\lfloor \frac{k}{2} \right\rfloor+2}                                      &                                                                                              \\
                   & E_1 \arrow[ld, no head] \arrow[d, no head]                                            &                                                        & E_{k-1} \arrow[rd, no head] \arrow[d, no head]                    &                    & \left\lfloor \frac{k}{2} \right\rfloor-1 & E_1 \arrow[rd, no head] \arrow[uu, no head, dotted]                                  &                                                                                  &                 &                                                                                    & E_{k-1} \arrow[ld, no head] \arrow[d, no head]                                               \\
E'_1               & E_0                                                                                   &                                                        & E_{k}                                                             & E'_{k-1}           & \left\lfloor \frac{k}{2}\right\rfloor    & E_0 \arrow[u, no head]                                                               & E'_1                                                                             &                 & E'_{k-1}                                                                           & E_{k}                                                                                       
\end{tikzcd}}
\]
Continuing, since each $\phi_{i}$ for $t<i<k$ is descending, we obtain that $v_{p}(\mathfrak{f}_{i})=i-t-1$ for $t+1\leq i\leq k$. Since $v_{2}(\Delta_{K})=3$, we have that $k=2t+1$. The proof now follows since
\begin{align*}
v_{p}(\mathfrak{f}_{i})  & =\left\{
\begin{array}
[c]{cl}
t-i & \text{if }0\leq i\leq t,\\
i-t-1 & \text{if }t<i\leq k,
\end{array}
\right.  \\
& =\left\{
\begin{array}
[c]{cl}
\frac{k-1}{2}-i & \text{if }0\leq i\leq\frac{k}{2},\\
i-\frac{k+1}{2} & \text{if }\frac{k}{2}<i\leq k,
\end{array}
\right.  \\
& =\max\!\left\{  \left\lfloor \frac{k}{2}\right\rfloor -i,i-\left\lceil \frac{k}{2}\right\rceil \right\}  . \qedhere
\end{align*}
\end{proof}

While Theorem~\ref{thm:CMdiscclass} is restricted to $F=\mathbb{Q}(j(E))$, our next theorem shows that knowledge over $F$ suffices to determine the isogeny graph over any field $F$ with a real embedding.

\begin{theorem}\label{thm:genkwon}
Let $E$ be a complex elliptic curve with $\operatorname*{End}E\otimes_{\mathbb{Z}}\mathbb{Q}\cong K$ an imaginary quadratic field. Suppose further that $E$ is defined over a field $F$ such that $\operatorname*{End}_{F}E\cong\mathbb{Z}$. For each $E^{\prime}\in V(\mathcal{G}(E/F))$, let $\mathfrak{f}_{E^{\prime}}$ denote the conductor of the order $\operatorname*{End}E^{\prime}$. For each prime $p$, let
\[
\mathfrak{F}_{p}=\sup\left\{  v_{p}(\mathfrak{f}_{E^{\prime}})\mid E^{\prime
}\in V(\mathcal{G}_{p}(E/F))\right\}  .
\]
Then, the isogeny graph $\mathcal{G}(E/F)\cong\square_{p}(\mathcal{G}_{p}(E/F),[E]_{F})$ where $\mathcal{G}_{p}(E/F)$ is determined as follows:

\begin{enumerate}
\item if $\mathfrak{F}_{p}=\infty$, then
\[
\mathcal{G}_{p}(E/F)\cong\left\{
\begin{array}
[c]{cl}
\mathcal{H}_{2^{\infty}}^{1} & \text{if }p=2,\\
\mathcal{H}_{p^{\infty}}^{0} & \text{if }p\geq3;
\end{array}
\right.
\]

\item if $\mathfrak{F}_{p}<\infty$, then there exists an elliptic curve $E_{0}\in V(\mathcal{G}_{p}(E/F))$ such that $\operatorname*{End}E_{0}$ has discriminant $\Delta=\mathfrak{f}_{E_{0}}^{2}\Delta_{K}$, where $\Delta_{K}$ is the discriminant of $K$ and $v_{p}(\mathfrak{f}_{E_{0}})=\mathfrak{F}_{p}$. Further, let $\mathcal{G}_{p}(E_{0}/\mathbb{Q}(j(E_{0}))$  be as given by Lemma~\ref{lem:specialdiscpotcm} or Theorem~\ref{thm:CMdiscclass}. Then,
\[
\mathcal{G}_{p}(E/F)\cong\mathcal{G}_{p}(E_{0}/\mathbb{Q}(j(E_{0}))).
\]
\end{enumerate}
\end{theorem}

\begin{proof}
By Theorem \ref{mainthmCarPro}, it suffices to show that $\mathcal{G}_{p}(E/F)$ is as claimed. By \cite[Theorem 2.3]{GrepsCM}, the field of definition of $\operatorname*{End}_{F}E$ is $KF$. In particular, $K\not \subseteq F$ since $\operatorname*{End}_{F}E\cong \mathbb{Z}$. From \cite[Lemma~6.5]{GrepsCM}, we have that there is a complex conjugation $c\in G_{F}$ such that $\left\langle c|_{KF}\right\rangle \cong\operatorname*{Gal}(KF/F)$ acts on $\operatorname*{Gal}(K(E[n])/K)$ with the property that for each positive integer~$k$ with $k\geq2$ if $p=2$, it is the case that $\det \rho_{E,p^{k}}(c)=-1$ and $\operatorname*{Tr}\rho_{E,p^{k}}(c)=0$. We then obtain by an identical argument to the one given in the proof of Proposition \ref{Prop:BloomInvReal} that
\begin{equation}
\mathfrak{I}_{p}(E/F)=\left\{
\begin{array}
[c]{cl}
1 & \text{if }p=2,\\
0 & \text{if }p\geq3.
\end{array}
\right.  \label{eq:bloominvpotcm}
\end{equation}

Now suppose that $\mathfrak{F}_{p}=\infty$. Then, there must be infinitely many $p$-isogenies, and so $\mathcal{G}_{p}(E/F)$ is infinite. We claim that $\mathcal{G}_{p}(E/F)$ has a line as a subgraph. In this direction, suppose on the contrary that $\mathcal{G}_{p}(E/F)$ does not contain a line. Then, by Theorem~\ref{mathmclass} and (\ref{eq:bloominvpotcm}), we deduce that $\mathcal{G}_{p}(E/F)\cong\mathcal{H}_{p^{\infty,+}}^{\mathfrak{I}_{p}(E/F)}$. By Lemma \ref{lem:raycounting}, there exists a leaf $E_{0}\in V(\mathcal{G}_{p}(E^{\prime}/F))$ such that the number of $F$-rational $p^{j}$-isogenies admitted by $E_{0}$ is $p^{\min\{\mathfrak{I}_{p}(E/F),\left\lfloor j/2\right\rfloor \}}$. Since $E_{0}$ is defined over $F$, we have that $\mathbb{Q}(j(E_{0}))\subseteq F$. Next, identify $\mathcal{G}_{p}(E/\mathbb{Q}(j(E_{0})))$ with its corresponding isomorphic graph $\mathcal{H}_{p^{k}}^{r}$ as given in Theorem~\ref{thm:CMdiscclass}. Note that by loc. cit., we may assume that $\left(  E_{i}\right)  _{i=0}^{k}$ is a $k$-spine in $\mathcal{H}_{p^{k}}^{r}$ such that $E\cong E_{0}$. Since $\mathcal{H}_{p^{k}}^{r}\hookrightarrow\mathcal{G}_{p}(E/F)$, we may view, after base changing, $\left(  E_{i}\right)  _{i=0}^{k}$ as a $k$-spine of $\mathcal{G}_{p}(E/F)$. If $E_{k}$ is a leaf in $\mathcal{G}_{p}(E/F)$, then $\left(  p,r\right) =\left(  2,1\right)  $. In particular, the $k-1$-spine $\left(  E_{i}\right)_{i=0}^{k-1}$ can be extended, after possibly replacing $E_{k}$, to a ray $\left(  E_{i}\right)  _{i=0}^{\infty}$. Note that we may still suppose that $\left(  E_{i}\right)  _{i=0}^{k}$ is a $k$-spine in $\mathcal{H}_{p^{k}}^{r}$. Next, let $\mathfrak{f}_{i}$ denote the conductor of $\operatorname*{End}E_{i}$. By loc. cit., we have that either
\[
v_{p}\!\left(  \mathfrak{f}_{\left\lfloor \frac{k}{2}\right\rfloor }\right)=0\qquad\text{or}\qquad v_{p}\!\left(  \mathfrak{f}_{\left\lceil \frac{k}{2}\right\rceil }\right)  =0.
\]
Now choose $t\in\{\left\lfloor k/2\right\rfloor ,\left\lceil k/2\right\rceil \}$ such that $v_{p}(\mathfrak{f}_{t})=0$. Further, $E_{t}$ admits either a horizontal or a descending $p$-isogeny. In particular, if
\[
t^{\ast}=\left\{
\begin{array}
[c]{cl}
t & \text{if }v_{p}(\mathfrak{f}_{t})=1,\\
t+1 & \text{if }v_{p}(\mathfrak{f}_{t})=0,
\end{array}
\right.
\]
then $E_{t^{\ast}}$ admits a descending isogeny. In fact, each isogeny $E_{i}\rightarrow E_{i+1}$ for $i\geq t^{\ast}$ is descending. Consequently, $v_{p}(\mathfrak{f}_{i})=i-t^{\ast}$ for $i\geq t^{\ast}$. 

Now observe that $\mathbb{Q}(j(E_{i}))\subseteq F$ for each $i$. Now let $m\geq4$, and consider the elliptic curve $E_{k+m}$. Now identify $\mathcal{G}_{p}(E/\mathbb{Q}(j(E_{k+m})))$ with its corresponding isomorphic graph $\mathcal{H}_{p^{k^{\prime}}}^{r^{\prime}}$ given in loc. cit. Moreover, since $t^{\ast}<k+m$, we have that
\begin{equation}
v_{p}(\mathfrak{f}_{k+m})=k+m-t^{\ast}\geq k-\left\lceil \frac{k}{2}\right\rceil +m\geq\left\lfloor \frac{k}{2}\right\rfloor +4\geq
4.\label{eq:ascendingthm}
\end{equation}
By loc. cit., $E_{k+m}$ is a leaf on a $k^{\prime}$-spine of $\mathcal{H}_{p^{k^{\prime}}}^{r^{\prime}}$. In particular, $E_{k+m}$ admits a $p^{k^{\prime}}$-isogeny $\phi$. We further obtain from loc. cit. and~(\ref{eq:ascendingthm}), that $k^{\prime}\geq2(k+m-t^{\ast}-1)$. Upon base changing back to $F$, we have that $\phi$ remains $F$-rational. In particular, since $\mathcal{G}_{p}(E/F)\cong\mathcal{H}_{p^{\infty,+}}^{\mathfrak{I}_{p}(E/F)}$ with $\mathfrak{I}_{p}(E/F)\in\left\{  0,1\right\}  $, we must have that $\phi$ factors through the isogeny $E_{k+m}\rightarrow E_{t^{\ast}+1}$. Now observe that the isogeny $E_{k+m}\rightarrow E_{0}$ has degree $p^{k+m}$ by construction. Since $k+m<k^{\prime}$ and $E_{0}$ is a leaf in $\mathcal{H}_{p^{\infty,+}}^{\mathfrak{I}_{p}(E/F)}$, we have that $\phi$ factors through $E_{1}$. In particular, the existence of $\phi$ requires that $E_{1}$ be $p$-bloomed. But this is impossible if $p$ is odd since $\mathfrak{I}_{p}(E/F)=0$. Consequently, $p=2$, and each elliptic curve $E'$ through which $\phi$ factors such that $E'$ is $p^l$-isogenous to $E_{k+m}$ with $0<l<k'$ is $2$-bloomed. In particular, we have the following subgraph of $\mathcal{H}_{2^{\infty,+}}^{1}$:
\[
\begin{tikzpicture}[scale=.9, transform shape]

    \draw[line width=1.5pt, style={draw=blue}, dash pattern=on 5pt off 3pt] (-1.45,0) -- (-.25,0);
    \Vertex[shape=circle, size=0.6, color=white, style={draw=blue}, x=0]{m2}
    \Vertex[shape=circle, size=0.6, color=white, style={draw=blue}, x=2]{m1}
    \Vertex[shape=circle, size=0.6, color=white, style={draw=blue},label =$E_1$, x=4]{0}
    \Vertex[shape=circle, size=0.6, color=white, style={draw=blue},label =$E_2$, x=6]{1}
    \Vertex[shape=circle, size=0.6, color=white, style={draw=blue},label =$E_3$, x=8]{2}
    \draw[line width=1.5pt, style={draw=blue}, dash pattern=on 5pt off 3pt] (8.3,0) -- (9.5,0);

    \Vertex[shape=circle, size=0.6,  color=white, style={draw=red}, x=0, y=2]{m21}
    \Vertex[shape=circle, size=0.6,  color=white, style={draw=red}, x=2, y=2]{m11}
    \Vertex[shape=circle, size=0.6,  color=white, style={draw=red}, label =$E_0$, x=4, y=2]{01}
    \Vertex[shape=circle, size=0.6,  color=white, style={draw=red}, x=6, y=2]{11}
    \Vertex[shape=circle, size=0.6,  color=white, style={draw=red}, x=8, y=2]{21}

    \Edge[style={color=blue}](m2)(m1)
    \Edge[style={color=blue}](m1)(0)
    \Edge[style={color=blue}](0)(1)
    \Edge[style={color=blue}](1)(2)
    \Edge[style={color=red}](m2)(m21)
    \Edge[style={color=red}](m1)(m11)
    \Edge[style={color=red}](0)(01)
    \Edge[style={color=red}](1)(11)
    \Edge[style={color=red}](2)(21)

\end{tikzpicture}
\]
Now observe that $E_{0}$ is $8$-isogenous to four elliptic curves. But, this yields our desired contradiction, since $E_{0}$ was chosen such that the number of $F$-rational isogenies admitted by $E_{0}$ is $2=2^{\min\{1,\left\lfloor 3/2\right\rfloor \}}$. With the claim established, it follows that $\mathcal{G}_{p}(E/F)$ has a line as a subgraph. We then obtain from Theorem~\ref{mathmclass} and~(\ref{eq:bloominvpotcm}) that $\mathcal{G}_{p}(E/F)\cong\mathcal{H}_{p^{\infty}}^{\mathfrak{I}_{p}(E/F)}$.

We now consider the case $\mathfrak{F}_{p}<\infty$. In this direction, we begin by observing that there is an elliptic curve $E_{0}\in V(\mathcal{G}_{p}(E/F))$ such that $v_{p}(\mathfrak{f}_{E_{0}})=\mathfrak{F}_{p}$. Next, identify $\mathcal{G}_{p}(E/\mathbb{Q}(j(E_{0})))$ with its corresponding isomorphic graph $\mathcal{H}_{p^{k}}^{r}$ as given in Lemma~\ref{lem:specialdiscpotcm} and Theorem~\ref{thm:CMdiscclass}. Now let $E^{\prime}\in V(\mathcal{G}_{p}(E/\mathbb{Q}(j(E_{0}))))$. From the discussion preceding and after (\ref{eq:hvalcoro}), we have that if $E^{\prime}/F$ admits a horizontal isogeny, then its field of definition $\mathbb{Q}(j(E_{0}))$. We similarly have that if $v_{p}(\mathfrak{f}_{E^{\prime}})\geq 1$, then $E^{\prime}$ admits a unique ascending isogeny with field of definition $\mathbb{Q}(j(E_{0}))$. Consequently, upon base changing to $F$, we have that if $E^{\prime}$ admits a $F$-rational isogeny with field of definition $L\not \supseteq \mathbb{Q}(j(E_{0}))$, then the isogeny must be descending. Now consider $\mathcal{H}_{p^{k}}^{r}\hookrightarrow\mathcal{G}_{p}(E/F)$. It suffices to show that this inclusion is in fact an isomorphism. Towards a contradiction, suppose that $\mathcal{G}_{p}(E/F)\not\cong\mathcal{H}_{p^{k}}^{r}$. Then, there exists an elliptic curve $E^{\prime}$ in $\mathcal{H}_{p^{k}}^{r}$ that admits a $F$-rational descending $p$-isogeny $\phi:E^{\prime}\rightarrow \widetilde{E}$ such that its field of definition $\mathbb{Q}(\phi)\not \subseteq\mathbb{Q}(j(E_{0}))$. In particular, $\widetilde{E}\not \in V(\mathcal{H}_{p^{k}}^{r})$ and $\mathfrak{f}_{\widetilde{E}}=p\mathfrak{f}_{E^{\prime}}$. In addition, since $E^{\prime}$ and $\widetilde{E}$ are defined over $F$, we have that  $\mathbb{Q}(j(E^{\prime}),j(\widetilde{E}),\phi)\subseteq F$. Next, following Clark~\cite{Clarkvolcanoes}, we say that $j(E^{\prime})$ and $j(\widetilde{E})$ are coreal if the number field $\mathbb{Q}(j(E^{\prime}),j(\widetilde{E}))$ admits a real embedding. 

We now suppose that $\Delta_{K}<-4$, and claim that $j(E^{\prime})$ and $j(\widetilde{E})$ are coreal. Suppose on the contrary that they are not coreal. Since $v_{p}(\mathfrak{f}_{\widetilde{E}})=1+v_{p}(\mathfrak{f}_{E^{\prime}})$, we obtain from \cite[Theorem~4.3]{Clarkvolcanoes}. that $\mathbb{Q}(j(E^{\prime}),j(\widetilde{E}))=K(j(\widetilde{E}))$. But then, $K\subseteq F$, which is our desired contradiction. Consequently, $j(E^{\prime})$ and $j(\widetilde{E})$ are coreal, and another application of loc. cit. yields that $\mathbb{Q}(j(E^{\prime}),j(\widetilde{E}))\cong\mathbb{Q}(j(\widetilde{E}))$. Since
\begin{equation}
\left[\mathbb{Q}(j(E^{\prime}),j(\widetilde{E})):\mathbb{Q}\right]  =\left[\mathbb{Q}(j(\widetilde{E})):\mathbb{Q}\right]  \label{eq:fldindimp}
\end{equation}
and $\mathbb{Q}(j(\widetilde{E}))\subseteq\mathbb{Q}(j(E^{\prime}),j(\widetilde{E}))$, we deduce that $\mathbb{Q}(j(E^{\prime}),j(\widetilde{E}))=\mathbb{Q}(j(\widetilde{E}))$. Moreover, since $\phi$ is a descending $p$-isogeny, it follows from \cite[Theorem 6.2 (c)]{Clarkvolcanoes} that
\[\mathbb{Q}(\phi)=\mathbb{Q}(j(E^{\prime}),j(\widetilde{E}))=\mathbb{Q}
(j(\widetilde{E}))\subseteq F.
\]
Now observe that by construction, there is a $F$-rational $p$-power-isogeny
$\psi:E_{0}\rightarrow\widetilde{E}$. Since $\mathbb{Q}(\psi)\subseteq F$ and $\operatorname*{End}_{F}E\cong\mathbb{Z}$, we have from loc. cit. that $j(E_{0})$ and $j(\widetilde{E})$ are coreal and $\mathbb{Q}(\psi)=\mathbb{Q}(j(E_{0}),j(\widetilde{E}))$. Now recall that by assumption, $\mathfrak{F}_{p}=v_{p}(\mathfrak{f}_{E_{0}})\geq v_{p}(\mathfrak{f}_{\widetilde{E}})$. We then obtain from \cite[Theorem 4.3]{Clarkvolcanoes} that $\mathbb{Q}(j(E_{0}),j(\widetilde{E}))\cong\mathbb{Q}(j(E_{0}))$. Arguing similar to the argument following (\ref{eq:fldindimp}), we deduce that $\mathbb{Q}(j(E_{0}),j(\widetilde{E}))=\mathbb{Q}(j(E_{0}))$. But this leads to our desired contradiction since
\begin{equation}
\mathbb{Q}(j(E_{0}))\not \supseteq\mathbb{Q}(\phi)=\mathbb{Q}(j(\widetilde{E}))\subseteq\mathbb{Q}(j(E_{0})).\label{eq:contradstatement}
\end{equation}
This establishes the case when $\Delta_{K}<-4$.

Now suppose that $\Delta_{K}\in\{-3,-4\}$. Since $\mathfrak{f}_{\widetilde{E}
}=p\mathfrak{f}_{E^{\prime}}$, we have that $v_{p}(\mathfrak{f}_{\widetilde{E}
})\geq1$. In particular, $\mathfrak{f}_{\widetilde{E}}^{2}\Delta_{K}<-4$. We
begin by claiming that $j(E^{\prime})$ and $j(\widetilde{E})$ are coreal. In
this direction, first suppose that
\[
\mathfrak{f}_{E^{\prime}}\in\mathcal{D}=\left\{  -3,-4,-12,-16,-27\right\}  .
\]
Then, \cite[Proposition 2.2]{Clarkvolcanoes2}, implies that $\mathbb{Q}(j(E^{\prime}),j(\widetilde{E}))\cong\mathbb{Q}(j(\widetilde{E}))$. A similar argument as above then gives that  $\mathbb{Q}(j(E^{\prime}),j(\widetilde{E}))=\mathbb{Q}(j(\widetilde{E}))$. Now suppose that $\mathfrak{f}_{E^{\prime}}\not \in \mathcal{D}$. Then, $\mathfrak{f}_{\widetilde{E}}\not \in \mathcal{D}$ and so $p$ divides $\gcd(\mathfrak{f}_{E^{\prime}},\mathfrak{f}_{\widetilde{E}})$. Therefore, $\mathbb{Q}(j(E^{\prime}),j(\widetilde{E}))=\mathbb{Q}(j(\widetilde{E}))$ by another application of loc. cit. Both cases thus yield that $j(E^{\prime})$ and $j(\widetilde{E})$ are coreal. By \cite[Theorem~5.1]{Clarkvolcanoes}, we then obtain that $\mathbb{Q}(\phi)=\mathbb{Q}(j(E^{\prime}),j(\widetilde{E}))=\mathbb{Q}(j(\widetilde{E}))$. As in the previous case, we now consider the $F$-rational $p$-power isogeny $\psi:E_{0}\rightarrow\widetilde{E}$. We now claim that $j(E_{0})$ and $j(\widetilde{E})$ are coreal. Since  $\mathfrak{F}_{p}=v_{p}(\mathfrak{f}_{E_{0}})\geq v_{p}(\mathfrak{f}_{\widetilde{E}})$, we have that $\mathfrak{f}_{E_{0}}=p^{l}\mathfrak{f}_{\widetilde{E}}$ for some nonnegative integer $l$. If $\mathfrak{f}_{\widetilde{E}}\in\mathcal{D}$, we obtain from \cite[Proposition 2.2]{Clarkvolcanoes2}, that $\mathbb{Q}(j(E_{0}),j(\widetilde{E}))\cong\mathbb{Q}(j(E_{0}))$, and so $\mathbb{Q}(j(E_{0}),j(\widetilde{E}))=\mathbb{Q}(j(E_{0}))$. If $\mathfrak{f}_{\widetilde{E}}\not \in \mathcal{D}$, then $\mathfrak{f}_{E_{0}}\not \in \mathcal{D}$ and $p$ divides $\gcd(\mathfrak{f}_{E^{\prime}},\mathfrak{f}_{\widetilde{E}})$. Thus, loc. cit. gives that $\mathbb{Q}(j(E_{0}),j(\widetilde{E}))=\mathbb{Q}(j(E_{0}))$, which establishes the claim. By \cite[Theorem 5.1]{Clarkvolcanoes}, $\phi(\psi)=\mathbb{Q}(j(E_{0}),j(\widetilde{E}))=\mathbb{Q}(j(E_{0}))$. In particular, (\ref{eq:contradstatement}) holds, yielding our desired contradiction.
\end{proof}

Theorem~\ref{thm:genkwon} gives a complete description of the isogeny graph $\mathcal{G}(E/F)$ of an elliptic curve $E$ with potential CM. In the finite case, it shows that every $p$-primary graph is determined by the isogeny graph $\mathcal{G}_p(E_0/\mathbb{Q}(j(E_0)))$, where $E_0$ is an elliptic curve in $\mathcal{G}_p(E/F)$ whose endomorphism order has maximal $p$-power conductor. Moreover, if the $p$-primary graph enlarges after base change, then the new graph is obtained by extending this graph downward along the corresponding chain of descending $p$-isogenies. Consequently, Theorem~\ref{thm:CMdiscclass} determines the $p$-primary graph over every field of definition.

As an application, we obtain a common generalization of results of Kwon \cite[Theorem 4.1 and Corollary 4.2]{Kwon} and Bourdon--Clark \cite[Theorem 5.3]{BourdonClark1}. More precisely, let $E$ be an elliptic curve with potential CM and let $\Delta$ denote the discriminant of the order $\operatorname{End}E$. When $\Delta <-4$, Kwon gave necessary and sufficient conditions for when $E$ admits an $n$-isogeny over $\mathbb{Q}(j(E))$, as seen in the proof of Lemma~\ref{lem:potCM}. Bourdon and Clark subsequently extended Kwon's Theorem to $\Delta \in \{-3,-4\}$, and proved that Kwon's necessary conditions remain valid over any number field~$F$ that neither contains the CM field $K$ nor a subfield isomorphic to $\mathbb{Q}\lbrack X]/(H_{p^{2}\Delta}(X))$. In our setting, the assumption $\operatorname{End}_F \! E \cong \mathbb{Z}$ is equivalent to $K\not \subseteq F$, so the latter hypothesis is automatic. Our classification of isogeny graphs shows that the remaining conditions are also sufficient. Furthermore, the proof does not require $F$ to be a number field.

\begin{corollary}\label{genKwonBC}
Let $E$ be a complex elliptic curve defined over a field $F$ with $\operatorname*{End}_{F}\!E\cong\mathbb{Z}$, and suppose that $\operatorname*{End}E$ is an order of discriminant $\Delta$ in an imaginary quadratic field $K$. Let $n$ be a positive integer, and suppose that for every prime $p|n$, the field $F$ does not contain a subfield isomorphic to $\mathbb{Q}\lbrack X]/(H_{p^{2}\Delta}(X))$, where $H_{p^{2}\Delta}(X)$ denotes the Hilbert class polynomial. Then, $E$ admits an $n$-isogeny if and only if one of the conditions in Table \ref{ta:kwon} holds. Moreover, Table \ref{ta:kwon} gives the number $m$ of distinct $n$-isogenies admitted by $E$.
\end{corollary}

\begin{proof}
For $E^{\prime}\in V(\mathcal{G}(E/F\mathcal{)})$, let $\mathfrak{f}_{E^{\prime}}$ denote the conductor of the order $\operatorname*{End}E^{\prime}$. We begin by observing that the assumption that $F$ does not contain a subfield isomorphic to $\mathbb{Q}\lbrack X]/(H_{p^{2}\Delta}(X))$ for each prime $p|n$, implies that
\[
\mathfrak{F}_{p}=\sup\!\left\{  v_{p}(\mathfrak{f}_{E^{\prime}})\mid
E^{\prime}\in V(\mathcal{G}_{p}(E/F))\}=v_{p}(\mathfrak{f}_{E})\right\}  .
\]
We then obtain from Theorem \ref{thm:genkwon} that for each prime $p|n$,
\[
\mathcal{G}_{p}(E/F)\cong\mathcal{G}_{p}(E/\mathbb{Q}(j(E))).
\]
From Theorem \ref{mainthmCarPro} and Lemma \ref{Lem:IsoCartGraphs}, we then
obtain that
\[
\mathcal{G}(E/F)\cong\left(  \underset{p|n}{\square}\left(  \mathcal{G}_{p}(E/\mathbb{Q}(j(E))),[E]_{\mathbb{Q}(j(E))}\right)  \right)  \square\left(  \underset{p\nmid n}{\square}\left(\mathcal{G}_{p}(E/F),[E]_{F}\right)  \right)  .
\]
In particular, the argument in \ref{prop:maxmat}, specifically, the one preceding (\ref{kronmatrix}) extends to show that any $n$-isogeny in $\mathcal{G}(E/F)$ lifts from an $n$-isogeny in $\underset{p|n}{\square}\left(  \mathcal{G}_{p}(E/\mathbb{Q}(j(E))),[E]_{\mathbb{Q}(j(E))}\right)  $, even if $\mathcal{G}(E/F)$ is infinite. We further have that
\[
\mathcal{G}_{p}(E/\mathbb{Q}(j(E)))\cong\left(  \underset{p|n}{\square}\mathcal{G}_{p}(E/\mathbb{Q}(j(E)))\right)  \square\left(  \underset{p\nmid n}{\square}\mathcal{G}_{p}(E/\mathbb{Q}(j(E)))\right)  ,
\]
with an $n$-isogeny lifting from $\underset{p|n}{\square}\mathcal{G}_{p}(E/\mathbb{Q}(j(E)))$. In particular, Kwon's necessary and sufficient conditions~\cite[Corollary 4.2]{Kwon} for an $n$-isogeny apply, and Table \ref{ta:kwon} provides these conditions, as well as the number $m$ of distinct $n$-isogenies admitted by $E$.
\end{proof}


\section{Determining the isogeny graph from the adelic Galois representation}\label{sec:algorithm image->graph}

In this section, we present an algorithm for determining the isogeny graph $\mathcal{G}(E/K)$ of an elliptic curve $E$ defined over a field $K$ of characteristic $0$ from the image of its associated adelic Galois representation
\[
\rho_{E}:G_{K}\longrightarrow\operatorname*{GL}\nolimits_{2}(\widehat{\mathbb{Z}}).
\]
Throughout, we assume that $\operatorname*{End}_{K}\! E\cong\mathbb{Z}$. For simplicity, we also assume that $\rho_{E}$ has finite level $N$. The input to our algorithm consists of a positive integer $N$, the level of $\rho_{E}$, together with a subset $S\subseteq\operatorname*{GL}\nolimits_{2}(\mathbb{Z}/N \mathbb{Z})$ generating the image $\rho_{E,N}(G_{K})$ with respect to a chosen basis of $E[N]$. From this data, the algorithm outputs $\mathcal{G}(E/K)$. The following result makes this precise.

\begin{theorem}\label{thm:algorithm}
Let $E$ be an elliptic curve over a field $K$ of characteristic $0$ with $\operatorname*{End}_{K}\! E\cong~\mathbb{Z}$. Then, there exists an algorithm which outputs the pointed graph $\left(  \mathcal{G}(E/K),E\right)  $ from the inputs:

\begin{itemize}
\item a positive integer $N$, the level of $\rho_{E}$;

\item a subset $S\subseteq\operatorname*{GL}\nolimits_{2}(\mathbb{Z}/n\mathbb{Z})$ such that $\rho_{E,n}(G_{K})\cong \langle S\rangle$ with respect to some basis of~$E[n]$.
\end{itemize}
\end{theorem}

\begin{proof}
Let $P$ denote the set of prime factors of $n$. By Theorem~\ref{mainthmCarPro}, 
\[
\left(  \mathcal{G}(E/K),E\right)  \cong\square_{p\in P}(\mathcal{G}_{p}(E/K),E).
\]
We now focus on determining the pointed graph $(\mathcal{G}_{p}(E/K),E)$ for each prime $p$. In this direction, fix $p\in P$ and let $v=v_{p}(n)$. Let $S_{p}=\left\{M\ \operatorname{mod}p^{v}\mid M\in S\right\}  $. Denote by $R$ the subgroup of $\operatorname*{GL}\nolimits_{2}(\mathbb{Z}/p^{v}\mathbb{Z})$ generated by $S_{p}$. Then, $\rho_{E,p^{v}}(G_{K})=R$ in some basis. 

\textbf{Step 1.} \textit{Finding the largest isogeny degree}. The first step of our algorithm is to replace $R$ by a suitable conjugate, so that it is easier to read the number of $p^{i}$-isogenies of $E$ for each $i$. In other words, we want to determine the value of the largest $i$ for which $R$ is conjugate to a subgroup of $B_{0}(p^{i})$.

Let $M=\sm{a &  b  \\ c  & d}\in R$, and let $\sm{\alpha &  \beta  \\ \gamma  & \delta}\in\operatorname*{GL}\nolimits_{2}(\mathbb{Z}/p^{v}\mathbb{Z})$. Then a direct calculation shows that the conjugate of $M$ by $\sm{\alpha &  \beta  \\ \gamma  & \delta}$ is upper-triangular modulo $p^{i}$ if and only if
\[
b\gamma^{2}+\left(  d-a\right)  \gamma\delta-c\delta^{2}\equiv 0\ \operatorname{mod}p^{i}.
\]
In other words, $\left(  \gamma,\delta\right)  $ is a zero modulo $p^{i}$ of the quadratic form
\[
f_{M}(x,y)=bx^{2}+(d-a)xy-cy^{2}.
\]
Thus, $R$ is conjugate to a subgroup of upper-triangular matrices modulo
$p^{i}$ if and only if there exists a point $(\gamma:\delta)\in\mathbb{P}^{1}(\mathbb{Z}/p^{i}\mathbb{Z})$ which is a common zero modulo $p^{i}$ of all quadratic forms $f_{M}(x,y)$ for $M\in R$. Since $S_{p}$ generates $R$, it is equivalent to finding a common zero modulo $p^{i}$ for each $M\in S_{p}$. In particular, we need to find the largest exponent $i$ for which the system has a common zero modulo $p^{i}$.

The next step is to simplify the system of equations. Let $S_{p}=\left\{M_{1},M_{2},\ldots,M_{m}\right\}  $, where $M_{j}=\sm{a_j &  b_j  \\ c_j  & d_j}$ for $1\leq j\leq m$. Consider the matrix
\[
\left(
\begin{array}
[c]{ccc}
b_{1} & d_{1}-a_{1} & c_{1}\\
b_{2} & d_{2}-a_{2} & c_{2}\\
\vdots & \vdots & \vdots\\
b_{m} & d_{m}-a_{m} & c_{m}
\end{array}
\right)  .
\]
In particular, the rows of the matrix encode the quadratic forms. We may perform row operations on this matrix to obtain an equivalent system of equations. In other words, we can replace this matrix by a version of its Smith normal form.

The Smith normal form algorithm begins by choosing a pivot element whose $p$-adic valuation is minimal among all entries of the matrix. One then subtracts suitable multiples of the row containing this pivot from the remaining rows in order to clear the rest of the pivot column. 

After performing row operations, we are left with at most $3$ quadratic equations. We need to find the largest exponent $i\leq v$ such that the equations are simultaneously solvable modulo $p^i$. Since we are looking for solutions $(x:y)\in \mathbb P^1(\Z/p^i\Z)$, we separately look for solutions of the form $(x:1)$ and $(1:y)$. The set of $\Z_p$-solutions of a quadratic equation of the form $$Ax^2+Bx+C\equiv 0 \mod{p^i}$$
is a union of at most two sets of the form $\{x\in \Z_p \mid  x\equiv x_0 \mod p^w\}$ where $x_0$ is some solution. These sets are $p$-adic balls inside $\Z_p$. For each of the (at most) three equations and for each~$i$, we calculate the sets of solutions modulo $p^i$. Then, we determine whether the sets have an intersection. This turns out to be simple, because an intersection of $p$-adic balls is either empty or a $p$-adic ball itself, and any point of a $p$-adic ball is its center. 

In the end, we are left with the largest $i$ for which there is a solution, and we also obtain the solution $(\gamma: \delta)$ itself. Finally, we conjugate all elements of $S_p$ by an invertible matrix of the form $\sm{\alpha & \beta \\ \gamma & \delta}.$  In fact, it turns out that the reductions of elements of $S_p$ modulo $p^i$ are enough to determine the entire graph, so we also reduce all elements of $S_p$ modulo $p^i$.

By abuse of notation, in what follows $S_p$ will denote the resulting set. In particular, $S_p$ is a set of matrices in $\GL_2(\Z/p^i\Z)$ which are uppertriangular.

\textbf{Step 2.} \textit{Finding the largest degree of an isogeny independent to the one of maximal degree.} 
We have now determined the largest $i$ for which $E$ has a $p^i$-isogeny. The next step of the algorithm is to determine the largest integer $l$ such that $E$ has a $K$-rational $p^{l}$-isogeny independent from the known $K$-rational $p^{i}$-isogeny. Recalling that $R\leq\operatorname*{GL}\nolimits_{2}(\mathbb{Z}/p^{v}\mathbb{Z})$, we let $\widetilde{R}\equiv R\ \operatorname{mod}p^{i}$. Now fix a basis $(P,Q)$ for $E[p^{i}]$ so that $\rho_{E,p^{i}}(G_{K})\cong\widetilde{R}$. In particular, $\left\langle P\right\rangle $ is $G_{K}$-invariant. If $E$ admits a $K$-rational $p^{l}$-isogeny for some positive integer $l$, then we must have a $G_{K}$-invariant cyclic group of the form $\left\langle Q+tP\right\rangle $ for some $t\in\mathbb{Z}/p^{i}\mathbb{Z}$. Thus, the finding of the largest integer $l$ for which $E$ admits a $K$-rational $p^{l}$-isogeny is equivalent to finding the largest integer $l$ for which there is a basis $(P,Q+tP)$ for some $t\in\mathbb{Z}/p^{i}\mathbb{Z}$ such that with respect to this basis,
\begin{equation}\label{eq:rhoimagealgo}
\rho_{E,p^{i}}(G_{K})\subseteq\left\{  \left(
\begin{array}
[c]{cc}
\ast & p^{l}\cdot\ast\\
0 & \ast
\end{array}
\right)  \in\operatorname*{GL}\nolimits_{2}(\mathbb{Z}/p^{i}\mathbb{Z})\right\}  .
\end{equation}
 Now write $S_{p}=\left\{  \widetilde{M}_{1},\widetilde{M}_{2},\ldots,\widetilde{M}_{m}\right\}  $, where $\widetilde{M}_{j}=\sm{a_j &  b_j  \\ 0  & c_j}$ for $1\leq j\leq m$. Then, upon changing basis to $(P,Q+tP)$, the corresponding new matrix for $\widetilde{M}_{j}$ is
\[
\left(
\begin{array}
[c]{cc}
1 & t\\
0 & 1
\end{array}
\right)  \left(
\begin{array}
[c]{cc}
a_{j} & b_{j}\\
0 & c_{j}
\end{array}
\right)  \left(
\begin{array}
[c]{cc}
1 & t\\
0 & 1
\end{array}
\right)^{-1} =\left(
\begin{array}
[c]{cc}
a_{j} & b_{j}-(c_{j}-a_{j})t\\
0 & c_{j}
\end{array}
\right)  .
\]
It thus suffices to find the largest integer $l$ for which there exists a $t\in\mathbb{Z}/p^{i}\mathbb{Z}$ such that
\[
(c_{j}-a_{j})t-b_{j}\equiv0\ \operatorname{mod}p^{l}
\]
for each $1\leq j\leq m$. This is a system of linear equations. As before, we may consider the matrix
\[
\left(
\begin{array}
[c]{cc}
c_{1}-a_{1} & -b_{1}\\
c_{2}-a_{2} & -b_{2}\\
\vdots & \vdots\\
c_{m}-a_{m} & -b_{m}
\end{array}
\right)  ,
\]
and then use row operations to reduce the problem to at most two linear equations. Similarly to the previous step, we can solve each of the linear equations separately, and then intersect the solution sets. We obtain the maximal exponent $l$ and the corresponding solution $t$.

As before, we now conjugate each element of $S_p$ by $\sm{1 & t \\ 0 & 1}$. 

\textbf{Step 3.} \textit{Determining the primary graph.} From Step 1 we have that the largest $K$-rational isogeny that $E$ admits is of degree $p^{i}$. Step $2$ gives us that the largest $K$-rational isogeny that $E$ admits that is independent of the known $p^{i}$-isogeny is of degree $p^{l}$. Step 2 further provides instructions to determine a basis $(P,Q)$ of $E[p^{i}]$ such that (\ref{eq:rhoimagealgo}) holds. By abuse of notation, with respect to the basis discussed, we set $R=\rho_{E,p^{i}}(G_{K})$ and let $S_{p}=\left\{M\ \operatorname{mod}p^{i}\mid M\in S\right\}  $ be a set that generates $R$. Let $S_{p}=\left\{  M_{1},M_{2},\ldots,M_{m}\right\}  $, where $M_{j}=\sm{a_j & b_j  \\ 0  & c_j}$ for $1\leq j\leq m$ and $v_{p}(b_{j})\geq l$ for $1\leq j\leq m$. Further, at least one $b_{j}$ satisfies $v_{p}(b_{j})=l$ since by construction, $p^{l}$ is the largest degree for an isogeny independent of the known $p^{i}$-isogeny. Next, let
\[
\alpha=\min\left\{  v_{p}^{(i)}(a-c)\left\vert \left(
\begin{array}
[c]{cc}
a & b\\
0 & c
\end{array}
\right)  \in S_{p}\right.  \right\}  .
\]
Since our assumptions are exactly those of Lemma~\ref{lem:rvalue}, we obtain
\[
r=\min\!\left\{  \alpha,\left\lfloor \frac{i+l}{2}\right\rfloor \right\}=\min\!\left\{  \mathfrak{I}_{k}(E),\left\lfloor \frac{i+l}{2}\right\rfloor
\right\}  .
\]
In particular, if $\alpha<\left\lfloor \frac{i+l}{2}\right\rfloor $, then $r=\alpha=\mathfrak{I}_{k}(E)$. Otherwise, $\alpha$ and $\mathfrak{I}_{k}(E)$ are both greater or equal to $\left\lfloor \frac{i+l}{2}\right\rfloor $, in which case $r=\left\lfloor \frac{i+l}{2}\right\rfloor $. It remains to determine the value of $k$ and the location of $E$ in $\mathcal{H}_{p^{k}}^{r}$. We now split into cases depending on the size of $\mathfrak{I}_{k}(E)$.

\qquad\textbf{Case A.} Suppose $\mathfrak{I}_{k}(E)\leq l$ or $\frac{i+l}{2}<\mathfrak{I}_{k}(E)$. By the contrapositive of the statement preceding and including~(\ref{valueokbloominvalpa}), we obtain that $E$ lies on some $k$-spine. Thus, $k=i+l$. We now see that $R$ is conjugate to a subgroup of $H_{p^{k}}^{r}(k-i)$. In particular, there is some $k$-spine $\mathcal{S}$ of $\mathcal{G}_{p}(E/K)$ with vertices $\left\{  E_{0},E_{1},\ldots,E_{k}\right\}  $ and edges $\left\{  E_{n},E_{n+1}\right\}  $ for $0\leq n<k$ such that $E$ is $K$-isomorphic to $E_{l}=E_{k-i}$. In particular, if we construct $\mathcal{H}_{p^{i+l}}^{r}$ with respect to $\mathcal{S}$ as in Proposition~\ref{prop:uniquetree}, we obtain an isomorphism of pointed graphs $(\mathcal{G}_{p}(E/K),E)\cong(\mathcal{H}_{p^{i+l}}^{r},E_{l})$.

\qquad\textbf{Case B.} Suppose $l<\mathfrak{I}_{k}(E)\leq\frac{i+l}{2}$. Since $\frac{i+l}{2}\leq\frac{k}{2}$, we obtain that $r=\mathfrak{I}_{k}(E)$. By assumption, $E$ is on a longest path of length $i+l$. Let $\mathcal{L}$ be the $i+l$-spine with vertices $\left\{  E_{0},E_{1},\ldots,E_{i+l}\right\}  $ and edges $\left\{  E_{n},E_{n+1}\right\}  $ for $0\leq n<i+l$ that contains $E$. In particular, $E_{0}$ and $E_{i+l}$ are leaves. Since a leaf of $\mathcal{L}$ is at a distance $i$ or $l$ from $E$, we may assume without loss of generality that $E_{l}\cong E$. From Lemma \ref{lem:kspineproperty}, we have that each skeletal vertex has bloom depth $r$. Since $l<r$, we have that $E_{0}$ and $E_{l}$ are contained in the bloom of some skeletal vertex. In particular, $E_{r}$ is the skeletal vertex such that its bloom $\mathcal{B}(E_{r})$ contains $E_{0}$ and $E_{l}$. We claim that $E_{i+l}$ lies on some $k$-spine. For a contradiction, suppose $E_{i+l}$ is not on any $k$-spine. By Lemma \ref{lem:kspineproperty}, $E_{i+l}$ must be in the bloom of some internal skeletal vertex and $E_{r}$ must be an endpoint of the skeleton. But then, another application of the lemma yields that $E_{l}$ lies on a $k$-spine. By construction, $\mathcal{L}$ is a longest path containing $E_{l}$ and $E_{i+l}$, and therefore $E_{i+l}$ lies on a $k$-spine, which is our desired contradiction.

Since $E_{i+l}$ lies on a $k$-spine, there is a $k$-spine $\mathcal{S}$ with vertices $\left\{  F_{0},F_{1},\ldots,F_{k}\right\}  $ and edges $\left\{E_{n},E_{n+1}\right\}  $ for $0\leq n<k$ such that $E_{i+l}\cong F_{k}$. Further, by Lemma \ref{Lem:skel}, $\mathcal{S}$ contains the skeleton, and so $E_{r}\cong F_{d}$ for some integer $d$. In particular, we have the following diagram.

\[
\begin{tikzcd}
	&& {E_0} && \\
	&& E=E_l \\
	{F_0} && {F_d=E_r} && {E_{i+l}}
	\arrow["{p^l}", no head, from=1-3, to=2-3]
	\arrow["{p^{r-l}}", no head, from=2-3, to=3-3]
	\arrow["{p^d}", no head, from=3-1, to=3-3]
	\arrow["{p^{i+l-r}}", no head, from=3-3, to=3-5]
\end{tikzcd}
\]

Observe that we now have that $k=d+i+l-r$. Now observe that since $F_{0}$ is a leaf, it is at a distance at least $r$ from the skeleton, and so $r\leq d$. Since $\mathcal{L}$ is a longest path containing $E_{l}$, we have that $d\leq i+l-r$. Since $l-r<0$, we have that $d<i$. We now determine the value of $d$ given that $r\leq d<i$. To determine $d$, we first count the number of $p^{r-l+t}$-isogenies admitted by $E$ for all $t$ satisfying $r\leq t\leq i+l-r$. Note that since $E_{r}$ is a skeletal vertex, we have that for each $t\geq r$, a $p^{r-l+t}$-isogeny from $E$ must factor through $E_{r}$. 

When $t=r$ we have that the number of $p^{2r-l}$-isogenies of $E$ is $p^{r}$ regardless of $d$. Since the isogeny factors through $E_{r}$, there are $p$ ways to choose for each of the remaining $r$ steps. 

Now consider $t$ with $r<t\leq d$. The number of $p^{r-l+t}$-isogenies from $E$ is $2p^{r}$ by a similar argument as that considered in Lemma~\ref{lem:infinitecountline}. Indeed, the isogeny factors through $E\rightarrow E_{r}$, and then, since $t>\alpha$, we must move along the skeleton in one of the two possible directions for the next $t-r$ steps. Finally, we have $p$ choices for each of the last $r$ steps.

Now consider $t$ with $d<t\leq i+l-r$. The number of $p^{i+l-r}$-isogenies of $E$ is $p^{r}$; namely, the isogeny first factors through $E\rightarrow E_{r}$, and then, since $t>d$, we must move along the skeleton towards the right side of the above diagram for the next $t-r$ steps. Since we have $p$ choices for each of the last $r$ steps, we conclude that there are a total of $p^{r}$ isogenies of degree $p^{i+l-r}$ admitted by $E$.

Thus, $d$ can be determined from the number of $p^{r-l+t}$-isogenies admitted by $E$ for all $t$ satisfying $r\leq t\leq i+l-r<i$: the number $d$ is the largest integer greater than $r$ for which the number of $p^{r-l+d}$-isogenies is equal to $2p^{r}$. We note that if no such integer $d$ exists, then $d=r$ and $i+l=k$.

We now argue similarly to the proof of Theorem \ref{classificationGpk} and reference the proof of \cite[Proposition~3.1]{ivan-numberofiso}. We note that in our setting, their value of $\alpha$ is the same as the quantity $\alpha$ considered at the start of this step. In particular, $r=\min\{\alpha,\left\lfloor \frac{i+l}{2}\right\rfloor \}$. Now observe that the proof of loc. cit. establishes that the number of $p^{r-l+d}$-isogenies equals to $2p^{r}$ if and only if there exists a nonzero element $y\in \mathbb{Z}/p^{r-l+d}\mathbb{Z}$ such that
\[
(c_{j}-a_{j})y\equiv b_{j}y^{2}\ \operatorname{mod}p^{r-l+d}\qquad\text{for
each }\sm{a_j &  b_j  \\ 0  & c_j}\in S_{p}.
\]
Further, for such a proper solution $y$, the proof shows, together with~(\ref{eq:rhoimagealgo}), that
\[
v_{p}(y)=r-\min\!\left\{  v_{p}(b_{j})\mid1\leq j\leq m\right\}  =r-l.
\]
So we may write $y=p^{r-l}z$ for some unit $z\in\mathbb{Z}/p^{d}\mathbb{Z}$. We thus obtain the system of equations
\[
(c_{j}-a_{j})p^{r-l}z-b_{j}p^{2r-2l}z^{2}\equiv0\ \operatorname{mod}
p^{r-l+d}\qquad\text{for }1\leq j\leq m.
\]
We may divide by $zp^{r-l}$ to obtain
\[
(c_{j}-a_{j})-b_{j}p^{r-l}z\equiv0\ \operatorname{mod}p^{d}\qquad\text{for
}1\leq j\leq m.
\]
It remains to find the largest $d>r$ for which such a unit $z$ exists; if no such $d$ exists, we take $d=r$. 

Now observe that $v_{p}(c_{j}-a_{j})\geq r$ by Lemma~\ref{lem:explicitbloompot}. Since $v_{p}(b_{j})\geq l$, we also have that $v_{p}(b_{j}p^{r-l})\geq r$. Now set $C_{j}=\frac{c_{j}-a_{j}}{p^{r}}$ and $B_{j}=\frac{b_{j}p^{r-l}}{p^{r}}=\frac{b_{j}}{p^{l}}$. Thus, the system of equations becomes
\[
C_{j}-B_{j}z\equiv0\ \operatorname{mod}p^{d-r}\qquad\text{for }1\leq j\leq m.
\]
Proceeding as before, we consider the following matrix corresponding to the system of equations:
\[
\left(
\begin{array}
[c]{cc}
C_{1} & -B_{1}\\
C_{2} & -B_{2}\\
\vdots & \vdots\\
C_{m} & -B_{m}
\end{array}
\right)  .
\]
We once again apply the Smith normal form algorithm, and search for an element with minimal $p$-adic valuation. By our scaling of the matrix, the minimal $p$-adic valuation occurring in the matrix is $0$. After permuting rows if necessary, we may assume that a pivot with $p$-adic valuation $0$ lies in the first row. Now observe that if $C_{1}B_{1}$ is not a unit, then there is no unit $z$ for which $C_{1}\equiv B_{1}z\ \operatorname{mod}p^{d-r}$. In this case, we take $d=r$. So suppose that $C_{1}B_{1}$ is a unit. If $m=1$, then $z=C_{1}B_{1}^{-1}$ and so the congruence is satisfied modulo $p^{h}$ for every $h$. Thus, $d=i+l-r$ and $k=2(i+l-r)$. So suppose that $m\geq2$. Then substituting $z=C_{1}B_{1}^{-1}$ yields
\[
B_{1}C_{j}-B_{j}C_{1}\equiv0\ \operatorname{mod}p^{d-r}\qquad\text{for }2\leq
j\leq m.
\]
Then our desired response is given by
\[
w=\min\!\left\{  v_{p}(B_{1}C_{j}-B_{j}C_{1})\mid2\leq j\leq m\right\}  .
\]
Consequently, $d=r+w$ and $k=i+l+w$.

By the above, we have determined the value of $k$ and $d$ in all cases. Now construct $\mathcal{H}_{p^{k}}^{r}$ with respect to $\mathcal{S}$ as in Proposition \ref{prop:uniquetree}. The above shows that $E\cong E_{l}$ is in the bloom $\mathcal{B}(F_{d})$ of $F_{d}$. In particular, it is in the $r-l$-th layer of the bloom, and upon identifying the corresponding vertex isomorphic to $E_{l}$ in $\mathcal{B}(F_{d})$ we obtain the isomorphism of pointed graphs $(\mathcal{G}_{p}(E/K),E)\cong(\mathcal{H}_{p^{k}}^{r},E_{l})$.
\end{proof}



\section{Modular curves associated to isogeny graphs}\label{sec:ModularCurves}
In the classification of finite $p$-primary graphs, we introduced the groups $H_{p^{k}}^{r}\leq B_{0}(p^{k})$ for $k\in\mathbb{N}$ and $r$ a nonnegative integer with $r\leq\frac{k}{2}$ (see Definition~\ref{def:groupsfinite}). In this section, we introduce the modular curves associated to isogeny graphs. In the case of prime power level, we consider the modular curves $X_{H_{p^{k}}^{r}}$ corresponding to $H_{p^{k}}^{r}$. To ease notation, we set $X_0^{r}(p^k)=X_{H_{p^{k}}^{r}}$. This notation is chosen so that it lines up with the classical modular curve $X_0(n)$. Indeed, $X_0(p^k)$ corresponds to the Borel subgroup $B_0(p^k)=H_{p^k}^0$, and thus $X_0^{0}(p^k)= X_0(p^k)$. We begin by investigating the field of definition of $X_0^{0}(p^k)$, thereby completing the last remaining proof of Theorem~\ref{mathmclass}. We then extend the definition to consider all modular curves associated to finite isogeny graphs, and classify those modular curves of genus $0$ and $1$. 

To begin, we observe that thus far, we have classified all the possible $p$-primary graphs. We further showed that each $p$-primary graph is realizable. It remains to prove that all possible combinations of $p$-primary graphs can simultaneously occur for some field $K$ of characteristic 0 and some elliptic curve $E/K$. Our first result establishes this fact.

\begin{proposition}\label{prop:allgraphsconverse}
For each prime $p$, let $k_{p}\in\mathbb{Z}_{\geq0}\cup\{\infty,(\infty,+)\}$ and $r_{p}\in\mathbb{Z}_{\geq0}\cup\{\infty\}$ be such that $r_{p}\leq\frac{k_{p}}{2}$. Suppose further that if $p=2$ and $k_{2}\geq2$, then $r_{2}\geq1$. Then, there exists a field~$K$ of characteristic $0$ and an elliptic curve $E/K$ such that $\mathcal{G}_{p}(E/K)\cong\mathcal{H}_{p^{k_{p}}}^{r_{p}}$ for each prime $p$.
\end{proposition}
\begin{proof}
    In Definitions \ref{def:groupsfinite} and \ref{def:groupsinfinite}, we introduced the groups $H_{p^{k_{p}}}^{r_{p}}$. In the former case, we have that $k_{p}<\infty$ and $H_{p^{k_{p}}}^{r_{p}}\leq\operatorname*{GL}\nolimits_{2}(\mathbb{Z}/p^{k}\mathbb{Z})$ is of level $p^{k_{p}}$. In what follows, we consider $H_{p^{k_{p}}}^{r_{p}}\leq\operatorname*{GL}\nolimits_{2}(\mathbb{Z}_{p})$ by considering its lift. In particular, we have that
\[
H_{p^{k_{p}}}^{r_{p}}=\left\{
\begin{array}
[c]{rl}
\left\{  \left.  \left(
\begin{array}
[c]{cc}
a & b\\
p^{k_{p}}w & c
\end{array}
\right)  \in\operatorname*{GL}\nolimits_{2}(\mathbb{Z}_{p})\right\vert v_{p}(a-c)\geq r_{p}\right\}   & \text{if }k_{p}\in\mathbb{Z}_{\geq0}\cup\{(\infty,+)\},\\
\left\{  \left.  \left(
\begin{array}
[c]{cc}
a & 0\\
0 & c
\end{array}
\right)  \in\operatorname*{GL}\nolimits_{2}(\mathbb{Z}_{p})\right\vert v_{p}(a-c)\geq r_{p}\right\}   & \text{if }k_{p}=\infty.
\end{array}
\right.
\]
Next, let $H=\prod_{p}H_{p^{k_{p}}}^{r_{p}}\leq\operatorname*{GL}\nolimits_{2}(\widehat{\mathbb{Z}})$. Since each $H_{p^{k_{p}}}^{r_{p}}$ is closed in $\operatorname*{GL}\nolimits_{2}(\mathbb{Z}_{p})$, it follows that $H$ is closed in $\operatorname*{GL}\nolimits_{2}(\widehat{\mathbb{Z}})$.
    
By the work of Greicius~\cite[Theorem 1.5]{greicius}, there exist a number field $K$ and an elliptic curve $E/K$ with surjective adelic Galois representation. By infinite Galois theory, there exists an extension $L/K$ such that the adelic image of $E/L$ equals $H$. This proves the claim. 
\end{proof}

\begin{remark}
    In fact, Zywina proved in \cite{zywinamaximal} that if $K\neq \Q$ is linearly disjoint from cyclotomic extensions, the density of elliptic curves over $K$ with surjective adelic image equals $1$. 

    The same does not hold over $\Q$. In fact, let $E/\Q$ be a non-CM elliptic curve and denote by $D$ the squarefree part of the discriminant of $E$. Serre proved in \cite[Proposition 22]{serre} that any element $g$ of the adelic image $\rho_{E}(G_\Q)$ satisfies $$\epsilon(g_2)=\chi_{D}
    (\det(g_{4D})),$$ where $g_n$ denotes the reduction $g \mod n$. This implies that the index of the adelic image is always divisible by $2$, as the above equality defines an index $2$ subgroup of $\GL_2(\widehat{\Z})$.
\end{remark}

Motivated by the above, we now extend our definition of $H_{p^{k}}^{r}$. To ease notation, we restrict ourselves to the finite-level setting in what follows.

\begin{definition}
Let $n>1$ be an integer with prime factorization $n=\prod_{p}p^{k_{p}}$, and let $\mathbf{r}=(r_{p})_{p|n}$ be a tuple such that each $r_{p}$ is a nonnegative integer satisfying $r_{p}\leq\frac{k_{p}}{2}$. Then the \textit{groups associated to isogeny graph of level $n$ and tuple $\mathbf{r}$}, denoted $H_{n}^{\mathbf{r}}$, is defined to be the group
\[
H_{n}^{\mathbf{r}}:=\left\{  \left.  \left(
\begin{array}
[c]{cc}
a & b\\
0 & c
\end{array}
\right)  \in\operatorname*{GL}\nolimits_{2}(\mathbb{Z}/n\mathbb{Z})\right\vert v_{p}(a-c)\geq r_p\text{ for each prime }p|n\right\}  .
\]
If the tuple $\mathbf{r}=(r)_{p|n}$ for some nonnegative integer $r\leq \frac{k_{p}}{2}$ for each prime $p$, then we define $H_{n}^{r}=H_{n}^{(r)_{p|n}}$.
\end{definition}

\begin{remark}
Observe that if $n=\prod_{p}p^{k_{p}}$ and $\mathbf{r}=(r_{p})_{p|n}$, then
\[
H_{n}^{\mathbf{r}}\cong\prod_{p}H_{p^{k_{p}}}^{r_{p}}.
\]
In particular, $H_{n}^{0}=B_{0}(n)$, where we recall that $H_{2^{k}}^{0}=H_{2^{k}}^{1}$ if $k\geq2$.
\end{remark}

As a consequence of Theorem \ref{mainthmCarPro}, Corollary~\ref{cor:subgroup-subgraph-corr}, and Proposition~\ref{prop:allgraphsconverse} we obtain:

\begin{corollary}
Let $n>1$ be an integer and for each prime $p|n$, let $r_{p}$ be a nonnegative integer such that $r_{p}\leq\frac{v_{p}(n)}{2}$ with the additional requirement that $r_{2}\geq1$ if $v_{2}(n)\geq2$. Suppose further that $E$ is an elliptic curve over a field $K$ of characteristic~$0$. If the isogeny graph $\mathcal{G}(E/K)$ has a subgraph isomorphic to $\square_{p|n}\mathcal{H}_{p^{v_{p}(n)}}^{r_{p}}$, then there exists $E^{\prime}\in V(\mathcal{G}(E/K))$ such that $\rho_{E^{\prime},n}(G_{K})$ is conjugate to a subgroup of $H_{n}^{\mathbf{r}}$ where $\mathbf{r}=(r_{p})_{p|n}$.

Conversely, if $\operatorname*{End}_{K}E\cong\mathbb{Z}$ and there exists an elliptic curve $E^{\prime}\in V(\mathcal{G}(E/K))$ such that $\rho_{E^{\prime
},n}(G_{K})$ is conjugate to a subgroup of $H_{n}^{\mathbf{r}}$ with
$\mathbf{r}=(r_{p})_{p|n}$ such that $r_{2}\geq1$ if $v_{2}(n)\geq2$, then
$\mathcal{G}(E/K)$ has a subgraph that is isomorphic $\square_{p|n}
\mathcal{H}_{p^{v_{p}(n)}}^{r_{p}}$.
\end{corollary}

Now observe that from Theorems~\ref{mainthmCarPro} and~\ref{mathmclass}, Lemma~\ref{lem:modularcurveimage}, and Corollary~\ref{cor:subgroup-subgraph-corr}, every finite isogeny graph $\mathcal{G}$ has an associated modular curve $X_{H}$ for some open subgroup $H\leq\operatorname*{GL}\nolimits_{2}(\widehat{\mathbb{Z}})$ of level $n=\deg\mathcal{G}(E/K)$. Now suppose that $\mathcal{G}(E/K)\cong\square_{p}\mathcal{H}_{p^{v_{p}(n)}}^{r_{p}}$. If we identify $H$ with its corresponding projection in $\operatorname*{GL}\nolimits_{2}(\mathbb{Z}/n\mathbb{Z})$, we obtain that $X_{H}$ is isomorphic to the following fiber product over $X(1)$:
\[
X_{H}\cong\bigtimes_{p}X_{0}^{r_{p}}(p^{v_{p}(n)}).
\]
This leads us to our next definition.

\begin{definition}\label{def:Modcurveisogra}
Let $X_{H}$ be a modular curve associated to an open subgroup $H\leq \operatorname*{GL}\nolimits_{2}(\widehat{\mathbb{Z}})$ of level $n$. We say that the \textit{modular curve }$X_{H}$\textit{ is associated to an isogeny graph} if $X_{H}$ is isomorphic to a fiber product over $X(1)$ of the form
\[
X_{H}\cong\bigtimes_{p}X_{0}^{r_{p}}(p^{v_{p}(n)}),
\]
where $r_{p}$ is a nonnegative integer satisfying $r_{p}\leq\frac{v_{p}(n)}{2}$ for each prime $p$. With $\mathbf{r}=(r_{p})_{p|n}$, we introduce the shorthand notation
\[
X_{0}^{\mathbf{r}}(n)=\bigtimes_{p|n}X_{0}^{r_{p}}(p^{v_{p}(n)}).
\]
If $\mathbf{r}$ is the constant tuple $(r)_{p}$, then we set $X_0^{r}(n)=X_{0}^{(r)_{p}}(n)$.
\end{definition}

Note that by definition, $X_{0}^{r}(n)\cong\bigtimes_{p}X_{0}^{r}(p^{v_{p}(n)})$, where the fiber product is taken over $X(1)$. Also note that $X_{0}^{0}(2^{k})=X_{0}^{1}(2^{k})$ if $k\geq2$ since with this assumption it is the case that $H_{2^{k}}^{0}=H_{2^{k}}^{1}$. Since $H_{n}^{0}=B_{0}(n)$, and the classical modular curve $X_{0}(n)$ is the modular curve associated to $B_{0}(n)$, we have the following isomorphism of modular curves $X_{0}^{0}(n)\cong X_{0}(n)$. Further, the modular curve associated to the group $H_{n}^{\mathbf{r}}$ is $X_{0}^{\mathbf{r}}(n)$. In particular, if $n=p^{k}$ for some prime $p$ and positive integer $k$ and $r$ is a nonnegative integer with $r\leq\frac{k}{2}$, then the modular curve associated to $H_{p^{k}}^{r}$ is $X_{0}^{r}(p^{k})$.

\begin{lemma}\label{lem:fieldofdefmodcurve}
Let $p$ be a prime number, $k$ a positive integer, and $r$ a nonnegative integer such that $r\leq\frac{k}{2}$. Then, the modular curve $X_{0}^{r}(p^{k})$ has field of definition
\[
K(X_{0}^{r}(p^{k}))=\left\{
\begin{array}
[c]{cl}
\mathbb{Q} & \text{if }r=0\text{ or }p=2\text{ with }r\leq1,\\
\mathbb{Q}(i) & \text{if }(p,r)=(2,2),\\
\mathbb{Q}(\zeta_{8}) & \text{if }p=2\text{ and }r\geq3,\\
\mathbb{Q}(\sqrt{p^{\ast}}) & \text{if }p\neq2\text{ and }r\geq1,\text{ where
}p^{\ast}=(-1)^{(p-1)/2}p.
\end{array}
\right.
\]

\end{lemma}

\begin{proof}
The modular curve $X_{0}^{r}(p^{k})$ is associated to the group $H_{p^{k}}^{r}\leq\operatorname*{GL}\nolimits_{2}(\mathbb{Z}/p^{k}\mathbb{Z})$. In particular,
\[
K(X_{0}^{r}(p^{k}))=\mathbb{Q}(\zeta_{p^{k}})^{\det H_{p^{k}}^{r}}.
\]
By definition of $H_{p^{k}}^{r}$ we have that if $\sm{a &  b  \\ 0  & c}\in H_{p^{k}}^{r}$, then $a\equiv c\ \operatorname{mod}p^{r}$. Since
\[
\det H_{p^{k}}^{r}=\left\{  ac\ \operatorname{mod}p^{k}\left\vert \left(
\begin{array}
[c]{cc}
a & b\\
0 & c
\end{array}
\right)  \in H_{p^{k}}^{r}\right.  \right\}  ,
\]
we deduce that $x\in\det H_{p^{k}}^{r}$ if and only if there is an $s\in(\mathbb{Z}/p^{r}\mathbb{Z})^{\times}$ such that $x\equiv s^{2}\ \operatorname{mod}p^{r}$. It now follows that $K(X_{0}^{r}(p^{k}))=\mathbb{Q}$ if and only if either $r=0$ or $p=2$ with $r\leq1$. Indeed, in these cases, we have that $H_{p^{k}}^{r}$ has full determinant.

Now suppose that $p=2$. If $r=2$, then $x\in\det H_{2^{k}}^{2}$ if and only if $x\equiv1\ \operatorname{mod}4$. Thus, $\det H_{2^{k}}^{2}$ corresponds to an index $2$ subgroup of $(\mathbb{Z}/2^{k}\mathbb{Z})^{\times}$, thus yielding
\[
K(X_{0}^{2}(2^{k}))=\mathbb{Q}(\zeta_{2^{k}})^{\det H_{2^{k}}^{2}}=\mathbb{Q}(i).
\]
Similarly, if $r\geq3$, then $x\in\det H_{2^{k}}^{r}$ if and only if $x\equiv1\ \operatorname{mod}8$. Hence, $\det H_{2^{k}}^{r}$ corresponds to an index $4$ subgroup of $(\mathbb{Z}/2^{k}\mathbb{Z})^{\times}$, and thus
\[
K(X_{0}^{r}(2^{k}))=\mathbb{Q}(\zeta_{2^{k}})^{\det H_{2^{k}}^{r}}=\mathbb{Q}(\zeta_{8}).
\]

Finally, if $p\geq3$ and $r\geq1$, then $x\in\det H_{p^{k}}^{r}$ if and only if $x$ is a square modulo $p^{r}$. Since $p$ is odd, we deduce that $\det H_{p^{k}}^{r}$ is an index $2$ subgroup of $(\mathbb{Z}/p^{k}\mathbb{Z})^{\times}$. It follows that
\[
K(X_{0}^{r}(p^{k}))=\mathbb{Q}(\zeta_{p^{k}})^{\det H_{p^{k}}^{r}}=\mathbb{Q}(\sqrt{p^{\ast}}).\qedhere
\]
\end{proof}

Our next result takes a finer look at the field of definition of $X_{0}^{r}(p^{k})$ for small $r$. Note that the result implies the remaining portion of Theorem~\ref{mathmclass}, which concludes its proof.

\begin{corollary}\label{cor:mainthmfieldofdef} 
Let $p$ be a prime and $k$ a positive integer. If $E$ is an elliptic curve over a field $K$ of characteristic not equal to $p$ such that $\rho_{E,p^{k}}(G_{K})$ is conjugate to a subgroup of $H_{p^{k}}^{r}$, but not conjugate to a subgroup of $H_{p^{k}}^{r+1}$, for some nonnegative integer $r\leq\frac{k}{2}$, then

\begin{enumerate}
\item if $p=2$, then

\begin{enumerate}
\item $r\leq1$ if and only if either $i\not \in K$ or $k\leq3$;

\item $r=2$ if and only if $i\in K$ with $k\geq4$ and either $\zeta_{8}\not \in K$ or $k\in\{4,5\}$;

\item $r\geq3$ if and only if $\zeta_{8}\in K$ and $k\geq6$; 
\end{enumerate}

\item if $p=3$, then

\begin{enumerate}
\item $r=0$ if and only if either $\zeta_{3}\not \in K$ or $k=1$;

\item $r\geq1$ if and only if $\zeta_{3}\in K$ and $k\geq2$; 
\end{enumerate}

\item for $p\geq5$, if $r\geq1$, then $\sqrt{p^{\ast}}\in K$. 
\end{enumerate}
\end{corollary}

\begin{proof}
The forward direction of each statement is a consequence of Lemma~\ref{lem:fieldofdefmodcurve}, together with Lemma~\ref{lem:modularcurveimage} and Corollary~\ref{cor:subgroup-subgraph-corr}. We now consider the converse statements for $p\in\{2,3\}$. In this direction, we note that by assumption we may choose a basis of $E[p^{k}]$ such that $\rho_{E,p^{k}}(G_{K})\leq H_{p^{k}}^{r}$. Now suppose $p=2$. If $i\not \in K$, then $r\leq1$ by Lemma~\ref{lem:fieldofdefmodcurve}. Similarly, if $k\leq3$, then $r\leq 1=\left\lfloor k/2\right\rfloor $. So suppose that $i\in K$ with $k\geq4$. By \cite[Chapter 5.2]{Adelmann},
\[
\rho_{E,2^{k}}(G_{K})\leq H_{2^{k}}^{r}\cap\left\{  M\in\operatorname*{GL}
\nolimits_{2}(\mathbb{Z}/2^{k}\mathbb{Z})\mid\det M\equiv1\ \operatorname{mod}4\right\}  .
\]
Consequently, if $\sm{a &  b  \\ 0  & c}\in\rho_{E,2^{k}}(G_{K})$, then $a\equiv c\ \operatorname{mod}4$. Equivalently, $v_{2}(a-c)\geq2$ and so $r\geq2$. In particular, if $k\in\{4,5\}$, then $2\leq r\leq\frac{k}{2}$ implies that $r=2$. Now suppose $k\geq6$. Now let $\chi_{2^{k}}:G_{K}
\rightarrow(\mathbb{Z}/2^{k}\mathbb{Z})^{\times}$ be the cyclotomic character, and recall that for each $\sigma\in G_{K}$ it is the case that
\begin{equation}
\det(\rho_{E,2^{k}}(\sigma))=\chi_{2^{k}}(\sigma),\label{eq:Galoisactionat2}
\end{equation}
since the determinant depends only on the cyclotomic action and is independent of $E$. Now observe that if $\zeta_{8}\in K$, then $G_{K}$ acts trivially on $\mu_{8}=\left\langle \zeta_{8}\right\rangle $, and so $\chi_{2^{k}}(\sigma)\equiv1\ \operatorname{mod}8$ for each $\sigma\in G_{K}$. In particular, if $\zeta_{8}\not \in K$, then there exists a $\sigma\in G_{K}$ such that $\chi_{2^{k}}(\sigma)\not \equiv 1\ \operatorname{mod}8$. By~\eqref{eq:Galoisactionat2}, $\det(\rho_{E,2^{k}}(\sigma))\not \equiv 1\ \operatorname{mod}8$. Consequently, if $\rho_{E,2^{k}}(\sigma)=\sm{a &  b  \\ 0  & c}$, then $a\not \equiv c\ \operatorname{mod}8$. Thus, $v_{2}(a-c)=2$ and so $r=2$. Observe that the above argument thus shows that if $\zeta_{8}\in K$ and $\sm{a &  b  \\ 0  & c}\in\rho_{E,2^{k}}(G_{K})$, then $a\equiv c\ \operatorname{mod}8$. Therefore, $v_{2}(a-c)\geq3$, which shows that $r\geq3$.

Now suppose that $p=3$ and $k\geq2$ with $\zeta_{3}\in K$. By \cite[Chapter~5.2]{Adelmann},
\[
\rho_{E,3^{k}}(G_{K})\leq H_{3^{k}}^{r}\cap\left\{  M\in\operatorname*{GL}\nolimits_{2}(\mathbb{Z}/3^{k}\mathbb{Z})\mid\det M\equiv1\ \operatorname{mod}3\right\}  .
\]
Thus, if $\sm{a &  b  \\ 0  & c}\in\rho_{E,3^{k}}(G_{K})$, then $a\equiv c\ \operatorname{mod}3$. It now follows that $r\geq1$ since $v_{3}(a-c)\geq1$.
\end{proof}

The converse to Corollary~\ref{cor:mainthmfieldofdef} (3), i.e., the case when $p\geq5$, is not necessarily true. In fact, we have the following as a direct consequence of Proposition~\ref{Prop:BloomInvReal}:
\begin{corollary}\label{cor:flddefsqrt}
Let $p \equiv 1 \mod 4$ be a prime. Then, for each integer $k\ge 2$ and each positive integer $r\le k/2$, the collection of $\mathbb{Q}(\sqrt{p})$-rational non-cuspidal points on the modular curve $X_{0}^{r}(p^k)$ is empty.
\end{corollary}
In Section~\ref{sec:expcalssgenus0}, we revisit Corollary~\ref{cor:flddefsqrt} in the special case of $X_{0}^{1}(25)$. In particular, we recover via explicit methods that  $X_{0}^{1}(25)$ has no non-cuspidal $\mathbb{Q}(\sqrt{5})$-rational points (see Remark~\ref{rmk:Qsqrt5flddef}).

We now pivot towards classifying the genus $0$ and $1$ modular curves associated with isogeny graphs. Towards this, we first recall that the genus of $X_{0}^{0}(n)=X_{0}(n)$ satisfies
\begin{equation}
\operatorname*{genus}(X_{0}^{0}(n))=\left\{
\begin{array}
[c]{cl}
0 & \text{if }n=1,2,\ldots,10,12,13,16,18,25,\\
1 & \text{if }n=11,14,15,17,19,20,21,24,27,32,36,49.
\end{array}
\right.  \label{lowgenusX0}
\end{equation}
By Faltings' Theorem~\cite{FaltingsThm}, if $n$ is a positive integer not appearing in (\ref{lowgenusX0}), then there are finitely many $j$-invariants over a fixed number field $K$ for which $E$ admits a $K$-rational $n$-isogeny. We now build on the above and complete the classification of genus $0$ and $1$ modular curves $X_{H}$ associated with isogeny graphs.

\begin{lemma}\label{lem:gen01classmodcurves}
Let $X_{H}$ be a modular curve associated to an isogeny graph. Then,

\begin{enumerate}
\item $\operatorname*{genus}(X_{H})=0$ if and only if

\begin{enumerate}
\item\label{1a} $X_{H}\cong X_{0}^{0}(n)$ with $n=1,2,\ldots,10,12,13,16,18,25$;

\item $X_{H}$ is isomorphic to a modular curve in the set
\[
\left\{  X_{0}^{1}(9),X_{0}^{2}(16),X_{0}^{0}(2)\times_{X(1)}X_{0}
^{1}(9),X_{0}^{1}(25)\right\}  .
\]

\end{enumerate}

\item $\operatorname*{genus}(X_{H})=1$ if and only if

\begin{enumerate}
\item\label{2a} $X_{H}\cong X_{0}^{0}(n)$ with $n=11,14,15,17,19,20,21,24,27,32,36,49$;

\item $X_{H}$ is isomorphic to a modular curve in the set
\[
\left\{  X_{0}^{1}(27),X_{0}^{2}(32),X_{0}^{1}(36)\right\}  .
\]

\end{enumerate}
\end{enumerate}
\end{lemma}

\begin{proof}
Since $H_{n}^{r}\leq B_{0}(n)$, we obtain a morphism of modular curves $X_{0}^{r}(n)\rightarrow X_{0}(n)$. From Riemann-Hurwitz, we deduce that $\operatorname*{genus}(X_{0}^{r}(n))\geq\operatorname*{genus}(X_{0}(n))$. From~(\ref{lowgenusX0}), it follows that if $\operatorname*{genus}(X_{0}^{r}(n))\leq1$, then
\begin{equation}
n\in\left\{  1,2,\ldots,21,24,25,27,32,36,49\right\}
\label{eq:modcurngenus01}
\end{equation}
Next, we observe that statements \ref{1a} and \ref{2a} are automatic from~(\ref{lowgenusX0}). It remains to consider those modular curves $X_{H}$ such that $H\cong H_{n}^{\mathbf{r}}$ with $n$ satisfying (\ref{eq:modcurngenus01}) and $\mathbf{r}=(r_{p})_{p|n}$ such that $r_{2}\geq2$ or $r_{p}\neq0$ for some odd prime $p|n$. Since $r_{p}\leq\frac{v_{p}(n)}{2}$, it follows that $H$ must isomorphic to one of the following:
\[
H_{9}^{1},H_{16}^{2},H_{2}^{0}\times H_{9}^{1},H_{25}^{1},H_{27}^{1}
,H_{32}^{2},H_{36}^{1},H_{49}^{1}.
\]
In \cite{GitHubIsogenies}, we provide \texttt{Magma} code\footnote{The implementation makes use of Zywina's repository
\url{https://github.com/davidzywina/Modular}, which is based on the methods and computations described in \cite{ZywinaLowGonality}.}
that verifies the claimed genus of all modular curves considered in this work, with the sole exception of $X_0^1(49)$. In this exceptional case, our computations show that $\operatorname{genus}(X_0^1(49))=3$. For the remaining modular curves, Table~\ref{tab:modularcurvemodels} gives explicit models for $X_H$. Observe that all genus $0$ modular curves are isomorphic to $\mathbb{P}^1_{K_H}$, except for $X_0^1(25)$, which has no $K_H=\Q(\sqrt{5})$-rational points and hence is not isomorphic to $\mathbb{P}^1_{K_H}$. Furthermore, each genus $1$ modular curve $X_H$ has no non-cuspidal non-CM points defined over the corresponding field $K_H$.

\begin{table}[ht]
\centering
\renewcommand{\arraystretch}{1.2}
\caption{Models of $X_H$ defined over $K_H$.}
\begin{tabular}{cc}
\hline
$X_H$ & Model of $X_H/K_H$\\
\hline
$X_0^1(9)$ & $\mathbb{P}^1_{\Q(\sqrt{-3})}$ \\
$X_0^2(16)$ & $\mathbb{P}^1_{\Q(i)}$ \\
$X_0^{0}(2)\times X_0^1(9)$ & $\mathbb{P}^1_{\Q(\sqrt{-3})}$ \\
$X_0^1(25)$ &
$\sqrt{5}x^2-3x^2+\sqrt{5}xy-3xy+\sqrt{5}x-3x-2y^2+\sqrt{5}y-y-2=0$ \\
\hline
$X_0^1(27)$ & $y^2+y=x^3-7$ \\
$X_0^2(32)$ & $y^2=x^3+4x$ \\
$X_0^1(36)$ & $y^2=x^3+1$ \\
\hline
\end{tabular}
\label{tab:modularcurvemodels}
\end{table}
\end{proof}

\begin{remark}
In Corollary \ref{cor:subgroup-subgraph-corr}, we defined subgroups $$H_{p^k}^r(j)=\left\{ \m{a & p^jb \\  p^{k-j}w & c}\in \GL_2(\Z_p) \mid v_p(a-c)\geq r\right\}$$ for $j\in \{0, 1\ldots, k\}$, and proved that elliptic curves whose images of $p$-adic Galois representations are contained in $H_{p^k}^r(j)$ also have $p$-primary subgraphs which contain $\mathcal H_{p^k}^r$. In fact, we defined them as subgroups of $\GL_2(\Z/p^{\max\{j,k-j\}}\Z)$, and in what follows we use $H_{p^k}^r(j)$ to refer to the lift to $\GL_2(\Z_p)$. 

Instead of considering modular curves associated to $H_{p^k}^r,$ one could instead consider modular curves associated to $H_{p^k}^r(j)$ for any $j\in \{0,1, \ldots, k\}$. However, the corresponding modular curves are isomorphic. Indeed, note that $$H_{p^k}^r(j)=gH_{p^k}^r(j+1)g^{-1} \text{ for } g=\m{1 & 0 \\ 0 & p}.$$ By work of Deligne and Rapoport \cite[IV-3.11.-3.19.]{Deligne-Rapoport}, modular curves corresponding to subgroups of $\GL_2(\widehat{\Z})$ which are conjugate by elements of $\GL_2(\Q)$ are isomorphic, and the isomorphism is modular, in the sense that they can be described in terms of level structures. 

More concretely, let $(E, \alpha)$ be a point on the modular curve associated with $H_{p^k}^r(j)$. The level structure $\alpha$ determines a $k$-spine in the $p$-primary isogeny graph of $E$, and $E$ occurs as its $j$-th vertex. Let $E'$ be the $(j+1)$-th vertex on this $k$-spine. The corresponding $p$-isogeny $E\to E'$ transfers the level structure on $E$ to a level structure on $E'$, which we denote by $g\cdot \alpha$.  Then $(E', g\cdot \alpha)$ is a point on the modular curve corresponding to $H_{p^k}^r(j+1)$. The map $(E, \alpha) \mapsto (E', g\cdot \alpha)$ is an isomorphism, because one can construct the inverse in a similar way.
\end{remark}

\section{Explicit classification of genus \texorpdfstring{$0$}{0} isogeny graphs}\label{sec:expcalssgenus0}
In this section, we establish that every genus $0$ modular curve associated to an isogeny graph admits an explicit parameterization by families of elliptic curves. These parameterizations simultaneously describe the $K$-rational non-cuspidal points on the modular curves and recover a distinguished portion of the corresponding isogeny graph. Since the cases corresponding to $X_{0}(n)$ were treated in \cite{Bariso}, it remains to consider modular curves corresponding to the four groups $H_{9}^{1},H_{16}^{2},H_{2}^{0}\times H_{9}^{1}$, and $H_{25}^{1}$ by Lemma~\ref{lem:gen01classmodcurves}. To ease notation in this section, we define $H_{n}$ by
\begin{equation}
{\renewcommand{\arraystretch}{1.2}
\begin{array}
[c]{c|cccc}
n & 9 & 16 & 18 & 25\\\hline
H_{n} & H_{9}^{1} & H_{16}^{2} & H_{2}^{0}\times H_{9}^{1} & H_{25}^{1}\\\hline
\end{array}
}\label{eq:Hndef}
\end{equation}
We proceed by first extending the families $\mathcal{C}_{n,j}(t,d)$ in \cite[Table 5]{Bariso} to larger families $\mathcal{C}_{n,m}^{1}(t,d)$ for $n=9,16,18,25$. For each such $n$, the extended family contains the previously constructed family $\mathcal{C}_{n,j}(t,d)$ while having the additional elliptic curves needed to parameterize the modular curves associated with the four groups $H_{n}$ in \eqref{eq:Hndef}. We begin by restating \cite[Theorem~1]{Bariso} in the language of Galois representations used throughout this section, and then introduce the extended families and their properties. We conclude this section by applying these parameterizations to explicitly determine the $n$-division field for $n\in\{3,4,5\}$ for elliptic curves $E/K$ such that $\rho_{E,n}(G_{K})$ is conjugate to a subgroup of the split Cartan subgroup of $\operatorname*{GL}\nolimits_{2}(\mathbb{Z}/n\mathbb{Z})$.

\begin{theorem}[{\cite[Theorem 1]{Bariso}}] \label{thm:barisothm1}
Let $n>1$ be an integer such that $X_{0}(n)$ has genus $0$. Let $E$ be an elliptic curve over a field $K$ of characteristic $0$ or relatively prime to $6n$, and suppose that $\rho_{E,n}(G_{K})$ is conjugate to the subgroup $H_{n}^{0}\leq\operatorname*{GL}\nolimits_{2}(\mathbb{Z}/n\mathbb{Z})$. Suppose further that the $j$-invariant of $E$ is not equal to $0$ or $1728$ if $n$ is prime. Then there exists a $t\in K$ and $d\in K^{\times}/(K^{\times})^{2}$ such that the following hold:

\begin{enumerate}
\item $\left\{  \left[  \mathcal{C}_{n,j}(t,d)\right]  _{K}\right\}  _{j} \subseteq V(\mathcal{G}(E/K))$, where $\mathcal{C}_{n,j}(t,d)$ is as defined in~\cite[Table~5]{Bariso};

\item The elliptic curve $E$ is $K$-isomorphic to $\mathcal{C}_{n,k}(t,d)$, where
\[
{\renewcommand{\arraystretch}{1.1}\begin{array}
[c]{c|ccccc}
n & 2,3,5,6,7,9,10,13,18,25 & 4 & 8 & 12 & 16\\\hline
k & 1 & 4 & 3 & 5 & 2\\\hline
\end{array}}
\]

\item if $\operatorname*{End}_{K}E\cong\mathbb{Z}$, then the partial isogeny graph of $\left\{  \left[  \mathcal{C}_{n,j}(t,d)\right]  _{K}\right\}  _{j}$ given in Table \ref{isographsbar} is a subgraph of $\mathcal{G}(E/K)$.
\end{enumerate}
\end{theorem}

{\begingroup
\renewcommand{\arraystretch}{1.3}
 \begin{longtable}{cccc}
 	\caption{Isogeny graph of $\{[\mathcal{C}_{n,j}(t,d)]_{K}\}_{j}$}\label{isographsbar}\\
	\hline
	$n$  & Isogeny Graph & $n$  & Isogeny graph  \\
	\hline

	\endfirsthead
	\hline
	$n$  & Isogeny Graph & $n$  & Isogeny graph \\
	\hline
	\endhead
	\hline

	\multicolumn{4}{r}{\emph{continued on next page}}
	\endfoot
	\hline
	\endlastfoot

$\begin{array}
[c]{c}
2,3,5,\\
7,13
\end{array}$ & \begin{tikzcd}
{\mathcal{C}_{n,1}} \arrow[r, "n", no head] & {\mathcal{C}_{n,2}}
\end{tikzcd}
& $9,25$  & \begin{tikzcd}
{\mathcal{C}_{n,1}} \arrow[r, "\sqrt{n}", no head] & {\mathcal{C}_{n,2}} \arrow[r, "\sqrt{n}", no head] & {\mathcal{C}_{n,3}}
\end{tikzcd} \\\hline

$6,10$ & \begin{tikzcd}
{\mathcal{C}_{n,1}} \arrow[r, "2", no head] \arrow[d, "\frac{n}{2}", no head] & {\mathcal{C}_{n,2}} \arrow[d, "\frac{n}{2}", no head] \\
{\mathcal{C}_{n,3}} \arrow[r, "2", no head]                                    & {\mathcal{C}_{n,4}}                                  
\end{tikzcd}
& $4$ & \begin{tikzcd}
                                            & {\mathcal{C}_{4,2}}                                                  &                     \\
{\mathcal{C}_{4,3}} \arrow[r, "2", no head] & {\mathcal{C}_{4,1}} \arrow[u, "2"', no head] \arrow[r, "2", no head] & {\mathcal{C}_{4,4}}
\end{tikzcd}\\\hline

$8$ & \adjustbox{scale=0.93}{\begin{tikzcd}
                                            & {\mathcal{C}_{8,2}} \arrow[d, "2", no head] & {\mathcal{C}_{8,5}}                                                 &                     \\
{\mathcal{C}_{8,3}} \arrow[r, "2", no head] & {\mathcal{C}_{8,1}} \arrow[r, "2", no head]          & {\mathcal{C}_{8,4}} \arrow[r, "2", no head] \arrow[u, "2", no head] & {\mathcal{C}_{8,6}}
\end{tikzcd}}
& $12$ & \adjustbox{scale=0.75}{\begin{tikzcd}
                                                                      & {\mathcal{C}_{12,3}} \arrow[r, "3", no head]                                                          & {\mathcal{C}_{12,4}}                          &                                                \\
                                                                      & {\mathcal{C}_{12,1}} \arrow[r, "3", no head] \arrow[rd, "\text{ }\text{ }\text{ }2", no head] \arrow[u, "2", no head] & {\mathcal{C}_{12,2}} \arrow[u, "2"', no head] &                                                \\
{\mathcal{C}_{12,5}} \arrow[ru, "2", no head] \arrow[r, "3", no head] & {\mathcal{C}_{12,6}} \arrow[ru, "2\text{ }\text{ }\text{ }\text{ }", no head]                                                 & {\mathcal{C}_{12,7}} \arrow[r, "3", no head]  & {\mathcal{C}_{12,8}} \arrow[lu, "2"', no head]
\end{tikzcd}} \\\hline

$16$ & \adjustbox{scale=0.75}{\begin{tikzcd}
                                             & {\mathcal{C}_{16,2}} \arrow[d, "2", no head] & {\mathcal{C}_{16,5}} \arrow[d, "2", no head] & {\mathcal{C}_{16,7}}                                                  &                      \\
{\mathcal{C}_{16,3}} \arrow[r, "2", no head] & {\mathcal{C}_{16,1}} \arrow[r, "2", no head] & {\mathcal{C}_{16,4}} \arrow[r, "2", no head] & {\mathcal{C}_{16,6}} \arrow[u, "2"', no head] \arrow[r, "2", no head] & {\mathcal{C}_{16,8}}
\end{tikzcd}}
& $18$ & \begin{tikzcd}
{\mathcal{C}_{18,1}} \arrow[r, "3", no head] \arrow[d, "2", no head] & {\mathcal{C}_{18,3}} \arrow[d, "2", no head] \arrow[r, "3", no head] & {\mathcal{C}_{18,5}} \arrow[d, "2", no head] \\
{\mathcal{C}_{18,2}} \arrow[r, "3", no head]                         & {\mathcal{C}_{18,4}} \arrow[r, "3", no head]                         & {\mathcal{C}_{18,6}}                        
\end{tikzcd}

\end{longtable}
\endgroup}

We now turn to the remaining four genus $0$ modular curves. For each of the corresponding groups $H_n$ in \eqref{eq:Hndef}, we introduce families of elliptic curves $\{\mathcal{C}_{n,m}^{1}(t,d)\}_m$, whose Weierstrass models are given in Table~\ref{ta:curves}. These families are constructed so as to extend the families $\{\mathcal{C}_{n,j}(t,d)\}_j$ from \cite{Bariso}, thereby completing the explicit parameterization of the genus $0$ modular curves that parameterize isogeny graphs.

As in loc. cit., the first step is to establish that, for each $n$, the curves in the family are linked by prescribed isogenies and realize the expected partial isogeny graph. In this direction, our first result is an analogue of \cite[Proposition 3.1]{Bariso}.

\begin{proposition}
\label{prop:expclassg0}Suppose $n\in\left\{  9,16,18,25\right\}  $, and let $j_{n,m}(t)$ and $\mathcal{C}_{n,m}^{1}(t,d)$ be as defined in Tables~\ref{ta:jinv} and~\ref{ta:curves}, respectively. For each $n$, define $k_{1},k_{2},$ and $S_{n}$ as:
\begin{equation}
{\renewcommand{\arraystretch}{1.9}
\begin{array}
[c]{c|ccc}
n & 9,25 & 16 & 18\\\hline
(k_{1},k_{2}) & (1,3) & (1,8) & (1,6)
\end{array}
}\quad\text{and}\quad{\renewcommand{\arraystretch}{1.9}
\begin{array}
[c]{c|ccc}
n & 9,18 & 16 & 25\\\hline
S_{n} & \{\zeta_{3}\} & \{i\} &  \left\{ \sqrt{5},\sqrt{-2(5+\sqrt{5})(t^{2}+4)}\right\} 
\end{array}
}\label{eq:explicitKWnHn}
\end{equation}
Suppose further that $K$ is a field of characteristic $0$ or relatively prime to $6n$ such that $S_{n}\subset K$. Then, the elliptic curves $\mathcal{C}_{n,m}^{1}(t,d)$ are defined over $K(t,d)$ and satisfy the following properties:

\begin{enumerate}
\item The $j$-invariants of $\mathcal{C}_{n,k_{l}}^{1}(t,d)$ for $l=1,2$ are
given by the Fricke parameterizations~$j_{n,l}(t)$ (see Table~\ref{ta:jinv});

\item The elliptic curves $\mathcal{C}_{n,m}^{1}(t,d)$ are pairwise
non-isomorphic over $K(t,d)$;

\item The isogeny class of $\mathcal{C}_{n,1}^{1}(t,d)$ over $K(t,d)$ contains
the set $\{[\mathcal{C}_{n,m}^{1}(t,d)]_{K(t,d)}\}_{m}$;

\item The partial isogeny graph associated to $\{[\mathcal{C}_{n,m}
^{1}(t,d)]_{K(t,d)}\}_{m}$ is given in Table~\ref{isographs}. 
\end{enumerate}
\end{proposition}

{\begingroup
\renewcommand{\arraystretch}{1.3}
 \begin{longtable}{cccc}
 	\caption{Isogeny graph of $\{[\mathcal{C}_{n,m}^{1}(t,d)]_{K(t,d)}\}_{m}$}\label{isographs}\\
	\hline
	$n$  & Isogeny graph & $n$  & Isogeny graph  \\
	\hline

	\endfirsthead
	\hline
	$n$  & Isogeny graph & $n$  & Isogeny graph  \\
	\hline
	\endhead
	\hline

	\multicolumn{4}{r}{\emph{continued on next page}}
	\endfoot
	\hline
	\endlastfoot

$9$& \adjustbox{scale=0.8}{\begin{tikzcd}
                                              & {\mathcal{C}^1_{9,4}} \arrow[d, "3", no head]                         &                       \\
{\mathcal{C}^1_{9,1}} \arrow[r, "3", no head] & {\mathcal{C}^1_{9,2}} \arrow[r, "3", no head] \arrow[d, "3", no head] & {\mathcal{C}^1_{9,3}} \\
                                              & {\mathcal{C}^1_{9,5}}                                                 &                      
\end{tikzcd}} 	

& $16$ & \adjustbox{scale=0.8}{\begin{tikzcd}
                                               & {\mathcal{C}^1_{16,9}}                                                 &                                                                           & {\mathcal{C}^1_{16,10}}                                                &                        \\
                                               & {\mathcal{C}^1_{16,3}}                                                 & {\mathcal{C}^1_{16,5}} \arrow[lu, "2"', no head] \arrow[ru, "2", no head] & {\mathcal{C}^1_{16,7}}                                                 &                        \\
{\mathcal{C}^1_{16,1}} \arrow[r, "2", no head] & {\mathcal{C}^1_{16,2}} \arrow[r, "2", no head] \arrow[u, "2", no head] & {\mathcal{C}^1_{16,4}} \arrow[r, "2", no head] \arrow[u, "2", no head]    & {\mathcal{C}^1_{16,6}} \arrow[r, "2", no head] \arrow[u, "2", no head] & {\mathcal{C}^1_{16,8}}
\end{tikzcd}}   \\\hline

 $25$ & \adjustbox{scale=0.8}{\begin{tikzcd}
{\mathcal{C}^1_{25,4}} \arrow[rd, "5", no head] &                                                                                                                           & {\mathcal{C}^1_{25,5}} \\
{\mathcal{C}^1_{25,1}} \arrow[r, "5", no head]  & {\mathcal{C}^1_{25,2}} \arrow[r, "5", no head] \arrow[ru, "5", no head] \arrow[rd, "5", no head] \arrow[ld, "5", no head] & {\mathcal{C}^1_{25,3}} \\
{\mathcal{C}^1_{25,6}}                          &                                                                                                                           & {\mathcal{C}^1_{25,7}}
\end{tikzcd}}
&  $18$ & \adjustbox{scale=0.8}{\begin{tikzcd}
                                                                            &                                                                                                                         &                                                    &  &                                                & {\mathcal{C}^1_{18,2}} \arrow[d, "3", no head]                         &                         \\
                                                                            & {\mathcal{C}^1_{18,1}} \arrow[d, "3"', no head] \arrow[rrrru, "2", no head]                                             &                                                    &  & {\mathcal{C}^1_{18,8}} \arrow[r, "3", no head] & {\mathcal{C}^1_{18,4}} \arrow[r, "3", no head] \arrow[d, "3", no head] & {\mathcal{C}^1_{18,10}} \\
{\mathcal{C}^1_{18,7}} \arrow[r, "3"', no head] \arrow[rrrru, "2", no head] & {\mathcal{C}^1_{18,3}} \arrow[d, "3"', no head] \arrow[r, no head] \arrow[r, "3"', no head] \arrow[rrrru, "2", no head] & {\mathcal{C}^1_{18,9}} \arrow[rrrru, "2", no head] &  &                                                & {\mathcal{C}^1_{18,6}}                                                 &                         \\
                                                                            & {\mathcal{C}^1_{18,5}} \arrow[rrrru, "2", no head]                                                                      &                                                    &  &                                                &                                                                        &                        
\end{tikzcd}}

\end{longtable}
\endgroup}

\begin{proof}
The proof relies on computer verification, which was done in \textsc{SageMath}. We include a GitHub \cite{GitHubIsogenies} which contains code that verifies the various claims made in the argument below. Specifically, each $n$ is treated separately in the four Jupyter Notebooks \texttt{n=9.ipynb}, \texttt{n=16.ipynb}, \texttt{n=18.ipynb}, and \texttt{n=25.ipynb}. It is then checked that $\left\{  [\mathcal{C}_{n,j}(t,d)]_{K(t,d)}\right\}  _{j}\subseteq\left\{  \lbrack\mathcal{C}_{n,m}^{1}(t,d)]_{K(t,d)}\right\}  _{m}$ since the following isomorphism hold over $K(t,d)$:
\begin{equation}\label{eq:isofamilies}
\mathcal{C}_{n,j}(t,d)\cong\left\{
\begin{array}
[c]{cl}
\mathcal{C}_{16,2}^{1}(t,d) & \text{if }(n,j)=(16,1),\\
\mathcal{C}_{16,1}^{1}(t,d) & \text{if }(n,j)=(16,2),\\
\mathcal{C}_{n,j}^{1}(t,d) & \text{otherwise.}
\end{array}
\right.
\end{equation}
It then follows from \cite[Proposition 3.1]{Bariso} that $j(\mathcal{C}_{n,k_{l}}^{1}(t,d))=j_{n,l}(t)$, which establishes (1). Further, (2) follows since the $j$-invariants of the family $\mathcal{C}_{n,m}^{1}(t,d)$ are distinct. We now proceed by cases and show that the isogeny graph is as claimed, thereby establishing (3) and (4).

\textbf{Case 1.} Suppose $n=9$, and let
\[
T =\frac{-18\zeta_{3}-3t}{3-\zeta_{3}^{2}t}  .
\]
Then, the following isomorphisms hold over $K(t,d)$:
\[
\mathcal{C}_{9,1}^{1}(T,-3d)\cong\mathcal{C}_{9,4}^{1}(t,d),\qquad
\mathcal{C}_{9,2}^{1}(T,-3d)\cong\mathcal{C}_{9,2}^{1}(t,d),\qquad
\mathcal{C}_{9,3}^{1}(T,-3d)\cong\mathcal{C}_{9,5}^{1}(t,d).
\]
Consequently, loc. cit. implies that the following isogenies hold over $K$:
\[
\begin{tikzcd}
{\mathcal{C}_{9,1}^1(t,d)} \arrow[r, "3", no head] & {\mathcal{C}_{9,2}^1(t,d)} \arrow[r, "3", no head] & {\mathcal{C}_{9,3}^1(t,d)} &  & {\mathcal{C}_{9,4}^1(t,d)} \arrow[r, "3", no head] & {\mathcal{C}_{9,2}^1(t,d)} \arrow[r, "3", no head] & {\mathcal{C}_{9,5}^1(t,d)}
\end{tikzcd}
\]
It now follows that the isogeny graph is as claimed.

\textbf{Case 2.} Suppose $n=16$, and let
\[
 T =\frac{4+2it}{2i+t}  .
\]
Then, it is verified that the $K(t,d)$-isomorphisms below hold:
\[
\mathcal{C}_{16,m}^{1}(T,-d)\cong\left\{
\begin{array}
[c]{cl}
\mathcal{C}_{16,m}^{1}(t,d) & \text{if }1\leq m\leq4,\\
\mathcal{C}_{16,6}^{1}(t,d) & \text{if }m=5,\\
\mathcal{C}_{16,5}^{1}(t,d) & \text{if }m=6,\\
\mathcal{C}_{16,9}^{1}(t,d) & \text{if }m=7,\\
\mathcal{C}_{16,10}^{1}(t,d) & \text{if }m=8.
\end{array}
\right.
\]
It follows from loc. cit. that the isogeny graph is as claimed.

\textbf{Case 3.} Suppose $n=18$, and let
\[
T_{w}=\left\{
\begin{array}
[c]{cl}
\zeta_{3}t & \text{if }w=1,\\
\zeta_{3}^{2}t & \text{if }w=2.
\end{array}
\right.
\]
Then, the following $K(t,d)$-isomorphisms hold
\[
\mathcal{C}_{18,m}^{1}(T_{1},d)\cong\left\{
\begin{array}
[c]{cl}
\mathcal{C}_{18,7}^{1}(t,d) & \text{if }m=1,\\
\mathcal{C}_{18,8}^{1}(t,d) & \text{if }m=2,\\
\mathcal{C}_{18,m}^{1}(t,d) & \text{if }3\leq m\leq6,
\end{array}
\right.  \hspace{-.3em}\text{ and }\hspace{.1em}\mathcal{C}_{18,m}^{1}(T_{2},d)\cong\left\{
\begin{array}
[c]{cl}
\mathcal{C}_{18,9}^{1}(t,d) & \text{if }m=1,\\
\mathcal{C}_{18,10}^{1}(t,d) & \text{if }m=2,\\
\mathcal{C}_{18,m}^{1}(t,d) & \text{if }3\leq m\leq6,
\end{array}
\right.
\]
From loc. cit., we then obtain the following partial isogeny graphs, from which the desired isogeny graph is automatic.

\[
\adjustbox{scale=0.93}{\begin{tikzcd}
{\mathcal{C}_{18,7}^1(t,d)} \arrow[r, "3", no head] \arrow[d, "2", no head] & {\mathcal{C}_{18,3}^1(t,d)} \arrow[r, "3", no head] \arrow[d, "2", no head] & {\mathcal{C}_{18,5}^1(t,d)} \arrow[d, "2", no head] &  & {\mathcal{C}_{18,9}^1(t,d)} \arrow[r, "3", no head] \arrow[d, "2", no head] & {\mathcal{C}_{18,3}^1(t,d)} \arrow[r, "3", no head] \arrow[d, "2", no head] & {\mathcal{C}_{18,5}^1(t,d)} \arrow[d, "2", no head] \\
{\mathcal{C}_{18,8}^1(t,d)} \arrow[r, "3", no head]                         & {\mathcal{C}_{18,4}^1(t,d)} \arrow[r, "3", no head]                         & {\mathcal{C}_{18,6}^1(t,d)}                         &  & {\mathcal{C}_{18,10}^1(t,d)} \arrow[r, "3", no head]                        & {\mathcal{C}_{18,4}^1(t,d)} \arrow[r, "3", no head]                         & {\mathcal{C}_{18,6}^1(t,d)}                        
\end{tikzcd}}
\]

\textbf{Case 4.} Suppose $n=25$, and let $s=\sqrt{-2(5+\sqrt{5})(t^{2}+4)}$. Next, for $w\in\left\{  1,2,3,4\right\}  $, define $T_{w}$ and $D_{w}$ as given below:
\[
\left(  T_{w},D_{w}\right)  =\left\{
\begin{array}
[c]{cl}
\left(  \frac{(\sqrt{5}-1)t-s}{4},d\left(  \sqrt{5}-1+\frac{(10+2\sqrt{5})t}{s}\right)  ^{-1}\right)   & \text{if }w=1,\\
\left(  \frac{(\sqrt{5}-1)t+s}{4},d\left(  \sqrt{5}-1-\frac{(10+2\sqrt{5})t}{s}\right)  ^{-1}\right)   & \text{if }w=2,\\
\left(  \frac{-2t(\sqrt{5}+1)+s\sqrt{5}-s}{8},-d\left(  \sqrt{5}+1+\frac{4t\sqrt{5}}{s}\right)  ^{-1}\right)   & \text{if }w=3,\\
\left(  \frac{-2t(\sqrt{5}+1)-s\sqrt{5}+s}{8},-d\left(  \sqrt{5}+1-\frac{4t\sqrt{5}}{s}\right)  ^{-1}\right)   & \text{if }w=4.
\end{array}
\right.
\]
We then obtain, for each $w\in\{1,2,3,4\}$, that the following $K(t,d)$-isomorphisms hold:
\[
\mathcal{C}_{25,m}^{1}(T_{w},D_{w})\cong\left\{
\begin{array}
[c]{cl}
\mathcal{C}_{25,w+3}^{1}(t,d) & \text{if }m=1,\\
\mathcal{C}_{25,m}^{1}(t,d) & \text{if }m=2,3.
\end{array}
\right.
\]
We then obtain from loc. cit. that if $m\neq2$, then $\mathcal{C}_{25,m}^{1}(t,d)$ is $5$-isogenous to $\mathcal{C}_{25,2}^{1}(t,d)$ over $K(t,d)$, which concludes the proof.
\end{proof}

Proposition~\ref{prop:expclassg0} shows that the families $\mathcal{C}_{n,m}^1$ extend the strutural properties of the families $\mathcal{C}_{n,j}$ associated to the modular curves $X_0(n)$. In particular, for each fixed $n$, they form an isogenous family of pairwise non-isomorphic elliptic curves whose isogeny class realizes the prescribed partial isogeny graph. We now show that these families likewise parameterize the $K$-rational non-cuspidal points on the remaining genus $0$ modular curves.

\begin{theorem}\label{thm:explicitclassK}
Let $n\in\left\{  9,16,18,25\right\}$, and let $H_{n}$ be as defined in (\ref{eq:Hndef}). Let $E$ be an elliptic curve over a field $K$ of characteristic $0$ or relatively prime to $6n$, and suppose that $\rho_{E,n}(G_{K})$ is conjugate to a subgroup of $H_{n}\leq\operatorname*{GL}\nolimits_{2}(\mathbb{Z}/n\mathbb{Z})$. Then there exist $t\in K$ and $d\in K^{\times}/(K^{\times})^{2}$ such that the following hold:

\begin{enumerate}
\item $\{  [  \mathcal{C}_{n,m}^{1}(t,d)]  _{K}\}_{m}\subseteq V(\mathcal{G}(E/K))$, where the elliptic curves $\mathcal{C}_{n,m}^{1}(t,d)$ are given in Table~\ref{ta:curves};

\item The elliptic curve $E$ is $K$-isomorphic to $\mathcal{C}_{n,1}^{1}(t,d)$;

\item The field $K$ contains the set $S_{n}$, where $S_n$ is as given in \eqref{eq:explicitKWnHn}.

\item if $\operatorname*{End}_{K}E\cong\mathbb{Z}$, then the partial isogeny graph associated to $\{  [  \mathcal{C}_{n,m}^{1}(t,d)]_{K}\}  _{m}$, given in Table \ref{isographs}, is a subgraph of $\mathcal{G}(E/K)$.
\end{enumerate}
\end{theorem}

\begin{proof}
Since $H_{n}\leq H_{n}^{0}$, we have from Theorem~\ref{thm:barisothm1} and the proof of Proposition \ref{prop:expclassg0} that there exists $t\in K$ and $d\in K^{\times}/(K^{\times})^{2}$ such that $E$ is $K$-isomorphic to $\mathcal{C}_{n,1}^{1}(t,d)$, which establishes (2). By Corollary~\ref{cor:mainthmfieldofdef}, we have that $S_{n}\subseteq K$ if $n\neq25$ and $\sqrt{5}\in K$ if $n=25$. Since the family $\left\{  \mathcal{C}_{n,m}^{1}(t,d)\right\}  _{m}$ is defined over $K$ if $n\neq25$, we deduce from Proposition \ref{prop:expclassg0} that (1), (3), and (4) hold if $n\neq25$. 

For $n=25$, we observe that the proof of the converse of Corollary~\ref{cor:subgroup-subgraph-corr} shows that since $\rho_{E,25}(G_{K})$ is conjugate to a subgroup of $H_{25}=H_{25}^{1}$, it is the case that $E$ is $25$-isogenous over $K$ to an elliptic curve $E_{2}$ such that the $5$-isogenous elliptic curve $E_{1}$ through which the isogeny factors admits all of its six $5$-isogenies. Indeed, the techniques involving Galois representations used in the proof of loc. cit. continue to hold so long as the characteristic of $K$ is not $5$. Further, by Corollary \ref{cor:bloompotprop}, the $\operatorname{mod}25$ blooming invariant is an isogeny class invariant, and this also holds if the characteristic of $K$ is not $5$. Consequently, any elliptic curve that is $5$-isogenous to $E$ with the property that a $25$-isogeny from $E$ factors through it has the property that it admits all of its six $5$-isogenies over $K$. Since $\mathcal{C}_{25,2}^{1}(t,d)$ satisfies these assumptions, we have that of its six isogenies are defined over $K$. It follows from Proposition \ref{prop:expclassg0} that each curve in the family $\left\{  \mathcal{C}_{25,m}^{1}(t,d)\right\}  _{m}$ is defined over $K$. Necessarily, $S_{25}\subseteq K$ and so (1), (3), and (4) hold, which completes the proof.
\end{proof}

\begin{remark}\label{rmk:Qsqrt5flddef}
As an immediate consequence, we recover an explicit proof of Corollary~\ref{cor:mainthmfieldofdef} in the special case of $X_{0}^{1}(25)$. Namely, the collection of $\mathbb{Q}(\sqrt{5})$-rational non-cuspidal points on the modular curve $X_{0}^{1}(25)$ is empty.
\end{remark}

We conclude with the following application, which explicitly determines, for $n\in\{3,4,5\}$, the $n$-division field of elliptic curves whose $\operatorname{mod}n$ Galois image is diagonalizable.

\begin{corollary}\label{cor:ndivisionclass} 
Let $n\in\left\{  3,4,5\right\}  $ and suppose that $E$ is an elliptic curve over a field $K$ of characteristic $0$ or relatively prime to $6n$ such that $\rho_{E,n}(G_{K})$ is conjugate to a subgroup of the split Cartan subgroup of $\operatorname*{GL}\nolimits_{2}(\mathbb{Z}/n\mathbb{Z})$. Then, there exists a $t\in K$ and $d\in K^{\times}/(K^{\times})^{2}$ such that
\[
K(E[n])\cong\left\{
\begin{array}
[c]{cl}
K(\zeta_{3},\sqrt{d}) & \text{if }n=3,\\
K(i,\sqrt{d}) & \text{if }n=4,\\
K\!\left(  \zeta_{5},s,w\right)   & \text{if }n=5,
\end{array}
\right.
\]
where $s$ and $w$ are as given below:
\begin{align}
s  & =\sqrt{-2(5+\sqrt{5})(t^{2}+4)},\label{eq:valsn25}\\
w  & =\sqrt{d(((1-\zeta_{5})(11t^{5}+55t^{3}+55t+4))s+50\zeta_{5}^{3}(t^{2}+4)(t^{4}+3t^{2}+1))}.\label{eq:valcn25}
\end{align}

Furthermore, the following statements hold:

\begin{enumerate}
\item The elliptic curve $E$ is $K$-isomorphic to $\mathcal{C}_{n^{2},m}^{1}(t,d)$ where $m=2$ (resp. $4$) if $n=3,5$ (resp. $4$);

\item $\{[\mathcal{C}_{n^{2},j}(t,d)]_{K}\}_{j}\subseteq V(\mathcal{G}(E/K))$ where $\mathcal{C}_{n^{2},j}(t,d)$ is as defined in \cite[Table 5]{Bariso};

\item If $F/K$ is a field extension such that $S_{n}\subset F$ where $S_{n}$ is as given in (\ref{eq:explicitKWnHn}), then $\{[\mathcal{C}_{n^{2},m}^{1}(t,d)]_{K}\}_{m}\subseteq V(\mathcal{G}(E/F))$ where $\mathcal{C}_{n^{2},m}^{1}(t,d)$ is as defined in Table~\ref{ta:curves}.

\item The group $E[n]=\left\langle P_{n},Q_{n}\right\rangle $ where
{\footnotesize \begin{align*}
P_{n} &  =\left\{
\begin{array}
[c]{cl}
\left(  3t^{2}d,4(t^{3}-27)\sqrt{d^{3}}\right)  \hspace{3.8em} & \text{if }n=3,\\
\left(  -3d\left(  t^{4}-12t^{3}+24t^{2}-48t\right)  ,108it(t^{2}+4)(t-2)^{2}\sqrt{d^{3}})\right)  \hspace{3.8em}  & \text{if }n=4,\\
\left(  \frac{9ds}{5}\zeta_{5}^{2}(1-\zeta_{5})\left(  t^{5}+5t^{3}+5t-11\right)  +\left(  15t^{6}+105t^{4}+195t^{2}+60\right)  d,Y_{1}\right) \hspace{3.8em} & \text{if }n=5,
\end{array}
\right.  \\
Q_{n} &  =\left\{
\begin{array}
[c]{cl}
\left(  -d(t-6\zeta_{3}-6)^{2},4(2\zeta_{3}+1)^{3}\zeta_{3}^{2}(t-3)(t-3\zeta_{3})\sqrt{d^{3}}\right)   & \text{if }n=3,\\
\left(  -3d(t^{4}-12it^{3}-24t^{2}+48it+16),108t(t^{2}-4)(t-2i)^{2}\sqrt{d^{3}}\right)   & \text{if }n=4,\\
\left(  -\frac{9ds}{\zeta_{5}-\zeta_{5}^{4}}\left(  t^{4}+t^{3}+6t^{2}+6t+11\right)  -3d(t^{2}+4)\left(  t^{4}+6t^{3}+21t^{2}+36t+61\right),Y_{2}\right)   & \text{if }n=5,
\end{array}
\right.
\end{align*}}
where
{\footnotesize
\begin{align*}
Y_{1}  & =\frac{27dsw\zeta_{5}}{5(\zeta_{5}-\zeta_{5}^{4})}(t-1)\left(t^{4}+t^{3}+6t^{2}+6t+11\right)  ,\\
Y_{2}  & =\frac{27\zeta_{5}dw}{5(t-1)}\Biggl(\frac{s}{2}(1-3\zeta_{5}^{2}-3\zeta_{5}^{3})\left(  t^{3}+\left(  1+3(\zeta_{5}+\zeta_{5}^{4})\right) t^{2}+\left(  5-(\zeta_{5}+\zeta_{5}^{4})\right)  t+\left(  -2+8(\zeta_{5}+\zeta_{5}^{4})\right)  \right)  \cdot\\
& \left(  t^{2}+\left(  2(\zeta_{5}+\zeta_{5}^{4})-1\right)  t+5-2(\zeta_{5}+\zeta_{5}^{4})\right)  +  
(\zeta_5^3+7\zeta_5^2+8\zeta_5+4)t^6+(-4\zeta_5^3+22\zeta_5^2+18\zeta_5+9)t^5+\\
& (25\zeta_5^3+75\zeta_5^2+100\zeta_5+50)t^4 +(-50\zeta_5^3+150\zeta_5^2+100\zeta_5+50)t^3
+(125\zeta_5^3+225\zeta_5^2+350\zeta_5+175)t^2+\\
& (-136\zeta_5^3+248\zeta_5^2+112\zeta_5+56)t+164\zeta_5^3+148\zeta_5^2+312\zeta_5+156\Biggr).
\end{align*}}

\end{enumerate}
\end{corollary}

\begin{proof}
Since $\rho_{E,n}(G_{K})$ is conjugate to a subgroup of the split Cartan
subgroup of $\operatorname*{GL}\nolimits_{2}(\mathbb{Z}/n\mathbb{Z})$, it follows that $E$ admits at least two distinct $n$-isogenies. Consequently, there are elliptic curves $E_{1}$ and $E_{2}$ such that the following isogenies hold over $K$:
\[
\begin{tikzcd}
E_1 \arrow[r, "n", no head] & E \arrow[r, "n", no head] & E_2
\end{tikzcd}
\]
By Theorem \ref{thm:barisothm1} and the proof of Proposition~\ref{prop:expclassg0}, there exists $t\in K$ and $d\in K^{\times}/(K^{\times})^{2}$ such that the following isomorphisms hold over $K$:
\[
E_{1}\cong\mathcal{C}_{n^{2},m_{1}}^{1}(t,d)\qquad E\cong\mathcal{C}_{n^{2},m_{2}}^{1}(t,d)\qquad E_{2}\cong\mathcal{C}_{n^{2},m_{3}}^{1}(t,d).
\]
where
\[
\left(  m_{1},m_{2},m_{3}\right)  =\left\{
\begin{array}
[c]{cl}
\left(  1,2,3\right)   & \text{if }n=9,25,\\
\left(  1,4,8\right)   & \text{if }n=16.
\end{array}
\right.
\]
Consequently, (1) holds. We then obtain (2) and (3) from Theorem~\ref{thm:barisothm1} and Proposition~\ref{prop:expclassg0}. We note that proving (4) automatically gives that $K(E[n])$ is as claimed. The proof now follows, since in \cite{GitHubIsogenies}, we verify that the generators of $E[n]$ are as claimed. Specifically, see the three Jupyter Notebooks \texttt{n=3division.ipynb}, \texttt{n=4division.ipynb}, and \texttt{n=5division.ipynb}.
\end{proof}

\begin{example}\label{ex:nonunique}
Consider the elliptic curve $E$ with LMFDB label \href{https://www.lmfdb.org/EllipticCurve/Q/11/a/2}{11.a2}. The $\operatorname{mod}5$ Galois image is split Cartan, and the elliptic curve is $\mathbb{Q}$-isomorphic to
\[
\mathcal{C}_{5,2}^{1}(0,-2):y^{2}=x^{3}-214272x-69147648.
\]
Next, let $s$ and $w$ be as given in Corollary \ref{cor:ndivisionclass}. It is then computed that
\[
s=4\left(  \zeta_{5}^{3}+\zeta_{5}^{2}+2\zeta_{5}+1\right)  \qquad \text{and}\qquad w=\frac{4\left(  \zeta_{5}^{4}+\zeta_{5}+2\right)  \left(\zeta_{5}^{3}+\zeta_{5}-1\right)  \left(  \zeta_{5}^{3}+\zeta_{5}+2\right)}{\left(  \zeta_{5}^{3}+\zeta_{5}^{2}+2\zeta_{5}+1\right)  ^{3}}.
\]
In particular, $s,w\in\mathbb{Q}(\zeta_{5})$, and so loc. cit. implies that $\mathbb{Q}(E[5])\cong\mathbb{Q}(\zeta_{5})$. In particular, the base change of $E$ to $\mathbb{Q}(\zeta_{5})$ is $5$-bloomed, and its isogeny graph is, in fact, $\mathcal{G}(E/\mathbb{Q}(\zeta_{5}))\cong\mathcal{H}_{25}^{1}$.

Now consider the elliptic curve $E^{\prime}$ with LMFDB label \href{https://www.lmfdb.org/EllipticCurve/Q/18176/g/2}{18176.g2}. It is then checked that $E^{\prime}$ also has a split Cartan as its $\operatorname{mod}5$ Galois image, and it is $\mathbb{Q}$-isomorphic to
\[
\mathcal{C}_{5,2}^{1}(2,-2):y^{2}=x^{3}-22391424x-36949782528.
\]
We then compute $s=4\left(  \zeta_{5}^{3}+\zeta_{5}^{2}+2\zeta_{5}+1\right)
\sqrt{2}$. In particular, $s\not \in\mathbb{Q}(\zeta_{5})$, and so by Theorem~\ref{thm:explicitclassK}, $\mathcal{H}_{25}^{1}\not \hookrightarrow \mathcal{G}(E^{\prime}/\mathbb{Q}(\zeta_{5}))$. In fact, $\mathcal{G}(E^{\prime}/\mathbb{Q}(\zeta_{5}))\cong\mathcal{H}_{25}^{0}$. Consequently, over $\mathbb{Q}(\zeta_{5})$, the $5$-primary graphs $\mathcal{H}_{25}^{0}$ and $\mathcal{H}_{25}^{1}$ occur.
\end{example}

\section{Isogeny graphs over certain number fields} \label{sec:nffields}

Let $K$ be a number field, and let $\operatorname{IsogDeg}(K)$ denote
the set of degrees of cyclic $K$-rational isogenies of elliptic curves
defined over $K$. As explained in the introduction, the determination
of $\operatorname{IsogDeg}(\Q)$ was achieved through the work of several
authors. Among the fundamental contributions to this classification is
Mazur's landmark work \cite{MR482230}, while the determination was
completed by Kenku \cite{MR675184}. More precisely,
\begin{equation*}
\operatorname{IsogDeg}(\Q)
=
\{2,3,4,5,6,7,8,9,10,12,13,14,15,16,17,18,19,21,25,27,
37,43,67,163\}.
\end{equation*}

For number fields other than $\Q$, no analogous determination is known
in general. To the best of our knowledge, the only further instances for which
$\operatorname{IsogDeg}(K)$ has been completely determined are, conditional
on GRH, those obtained by Banwait, Najman, and Padurariu \cite{Banwait}. More precisely, they determined
$\operatorname{IsogDeg}(\Q(\sqrt d))$ for every $d\in\mathcal D$, where
\begin{equation*}
\mathcal D
=\left\{
\begin{array}{c}
-6846,-2289,213,834,1545,1885,1923,2517,2847,4569,\\
6537,7131,7302,7319,7635,7890,8383,9563,9903
\end{array}
\right\}.
\end{equation*}

The following theorem determines, conditional on GRH, all isogeny graphs
over $\Q(\sqrt d)$ for $d\in\mathcal D$. It is the analogue over these
quadratic fields of the characterization over $\Q$ given in
Theorem~\ref{ThmQ}.

\begin{theorem}
Assume GRH, and let $d\in\mathcal D$. The isogeny graph of an elliptic
curve defined over $\Q(\sqrt d)$ is uniquely determined by its isogeny
class degree. The possible isogeny graphs are precisely those appearing
in Table~\ref{isographsQ}, together with the additional graphs listed in
Table~\ref{isographsQd}.

In Table~\ref{isographsQd}, a checkmark in the column indexed by $d$
indicates that the graph in the corresponding row occurs over
$\Q(\sqrt d)$, whereas a cross indicates that it does not occur.
\end{theorem}

\begin{table}[htp]
\centering
\renewcommand{\arraystretch}{1.2}

\begin{tabular}{|c|cccccccccc|}
\hline
Graph $\backslash$ $d$
& $-6846$ & $-2289$ & $213$ & $834$ & $1545$
& $1885$ & $1923$ & $2517$ & $2847$ & $4569$ \\
\hline
$\mathcal H^1_4 \mathbin{\square} \mathcal H^0_5$
& $\times$ & $\times$ & \checkmark & \checkmark & $\times$
& \checkmark & \checkmark & \checkmark & $\times$ & $\times$ \\
\hline
$\mathcal H^1_8 \mathbin{\square} \mathcal H^0_3$
& $\times$ & \checkmark & \checkmark & \checkmark & \checkmark
& \checkmark & $\times$ & \checkmark & $\times$ & $\times$ \\
\hline
$\mathcal H^1_{32}$
& \checkmark & $\times$ & \checkmark & $\times$ & $\times$
& \checkmark & $\times$ & \checkmark & \checkmark & $\times$ \\
\hline
$\mathcal H^1_4 \mathbin{\square} \mathcal H^0_9$
& \checkmark & \checkmark & \checkmark & $\times$ & \checkmark
& \checkmark & $\times$ & \checkmark & $\times$ & \checkmark \\
\hline
\end{tabular}

\vspace{0.5cm}

\begin{tabular}{|c|ccccccccc|}
\hline
Graph $\backslash$ $d$
& $6537$ & $7131$ & $7302$ & $7319$ & $7635$
& $7890$ & $8383$ & $9563$ & $9903$ \\
\hline
$\mathcal H^1_4 \mathbin{\square} \mathcal H^0_5$
& \checkmark & \checkmark & $\times$ & \checkmark & \checkmark
& \checkmark & $\times$ & $\times$ & \checkmark \\
\hline
$\mathcal H^1_8 \mathbin{\square} \mathcal H^0_3$
& $\times$ & $\times$ & \checkmark & \checkmark & $\times$
& \checkmark & $\times$ & \checkmark & $\times$ \\
\hline
$\mathcal H^1_{32}$
& $\times$ & $\times$ & \checkmark & \checkmark & $\times$
& $\times$ & \checkmark & $\times$ & \checkmark \\
\hline
$\mathcal H^1_4 \mathbin{\square} \mathcal H^0_9$
& \checkmark & $\times$ & $\times$ & \checkmark & $\times$
& \checkmark & \checkmark & \checkmark & $\times$ \\
\hline
\end{tabular}

\caption{Additional isogeny graphs occurring over the quadratic fields
$\Q(\sqrt d)$ with $d\in\mathcal D$. A checkmark indicates that the
corresponding graph occurs over $\Q(\sqrt d)$, while a cross indicates
that it does not.}
\label{isographsQd}
\end{table}

\begin{proof}
Fix $d\in\mathcal D$, and consider an elliptic curve $E$ defined over $K=\Q(\sqrt d)$. We first verify that $E$ has no nontrivial endomorphisms defined over $K$, that is, that $\operatorname{End}_K \! E\cong \mathbb Z$. For $d$ positive, this is automatic by Proposition~\ref{Prop:BloomInvReal}, as $K$ admits a real embedding. The two remaining fields, corresponding to $d=-6846$ and $d=-2289$, are also excluded, since neither is the CM field of an order of class number one.

Moreover, $i\notin K$ and $\sqrt{p^{*}}\notin K$ for every prime $p$,
where $p^{*}=(-1)^{(p-1)/2}p$. Corollary~\ref{maincor1} therefore implies
that the isogeny class degree of $E$ uniquely determines its isogeny
graph $\mathcal G(E/K)$.

Finally, \cite[Table~7.1]{Banwait} gives
$\operatorname{IsogDeg}(\Q(\sqrt d))$ for every $d\in\mathcal D$,
conditional on GRH. In particular,
\begin{equation*}
\operatorname{IsogDeg}(\Q(\sqrt d))
\subseteq
\operatorname{IsogDeg}(\Q)\cup\{20,24,32,36\}.
\end{equation*}
The graphs arising from degrees in $\operatorname{IsogDeg}(\Q)$ are
listed in Table~\ref{isographsQ}, while the four remaining possibilities
are
$\mathcal H^1_4\mathbin{\square}\mathcal H^0_5$,
$\mathcal H^1_8\mathbin{\square}\mathcal H^0_3$,
$\mathcal H^1_{32}$, and
$\mathcal H^1_4\mathbin{\square}\mathcal H^0_9$,
corresponding respectively to degrees $20$, $24$, $32$, and $36$.
Their occurrence over each field $\Q(\sqrt d)$ is recorded in
Table~\ref{isographsQd}, completing the proof.
\end{proof}

\appendix
\section{An LMFDB zoo of isogeny graphs}\label{LMFDBzoo}
The \emph{L-functions and Modular Forms Database} (LMFDB) provides a database
for browsing elliptic curves over number fields, their isogeny classes, and the
arithmetic data attached to them. In particular, it offers a concrete way to
examine examples of isogenies arising in practice, through explicit elliptic
curve labels.

\underline{\it Labeling.} We now introduce a labelling convention for the isogeny graphs described in this
paper. Let \(E\) be an elliptic curve defined over a number field \(K\) such that \(\operatorname{End}_K(E)\cong \mathbb{Z}\). We have shown that the
isogeny graph of \(E\), denoted by \(\mathcal{G}(E/K)\), is the direct
product of its \(p\)-primary components \(\mathcal{G}_p(E/K)\). Moreover, each
of these \(p\)-primary graphs is isomorphic to one of the graphs
\(\mathcal{H}^{r}_{p^k}\). We label such a \(p\)-primary graph by
\[
p^k.r.
\]
Thus, to label \(\mathcal{G}(E/K)\), we concatenate the labels of the
non-trivial factors \(\mathcal{G}_p(E/K)\), ordered by increasing primes \(p\).

The purpose of Table~\ref{tab:lmfdb-zoo} is to exhibit one elliptic curve for
each isogeny graph appearing in the LMFDB. In this sense, the table may be
viewed as a small LMFDB zoo: a curated list of examples illustrating the
different isogeny graphs that occur among elliptic curves recorded in the
database.

For each row, the first column gives the label of the isogeny graph according
to the convention introduced above. The second column records the size of the
corresponding isogeny class, and the third column records the relevant isogeny
degree. The final column gives the LMFDB label of a representative elliptic
curve.

\vspace{1em}

\begin{longtable}{lccl}
\caption{A small zoo of isogeny configurations}\label{tab:lmfdb-zoo}\\

\hline
\textbf{isogeny label} & \textbf{size} & \textbf{degree} & \textbf{LMFDB label} \\
\hline
\endfirsthead

\hline
\textbf{isogeny label} & \textbf{size} & \textbf{degree} & \textbf{LMFDB label} \\
\hline
\endhead

\hline
\endfoot

\hline
\endlastfoot

\texttt{1.0} & $1$ & $1$ & \href{https://www.lmfdb.org/EllipticCurve/2.0.1016.1/9.2/c}{2.0.1016.1-9.2-c1} \\
\texttt{2.0} & $2$ & $2$ & \href{https://www.lmfdb.org/EllipticCurve/2.0.1003.1/9.1/a}{2.0.1003.1-9.1-a1} \\
\texttt{3.0} & $2$ & $3$ & \href{https://www.lmfdb.org/EllipticCurve/2.0.103.1/56.3/a}{2.0.103.1-56.3-a1} \\
\texttt{5.0} & $2$ & $5$ & \href{https://www.lmfdb.org/EllipticCurve/2.0.1003.1/9.1/c}{2.0.1003.1-9.1-c1} \\
\texttt{7.0} & $2$ & $7$ & \href{https://www.lmfdb.org/EllipticCurve/2.0.103.1/416.10/a}{2.0.103.1-416.10-a1} \\
\texttt{11.0} & $2$ & $11$ & \href{https://www.lmfdb.org/EllipticCurve/2.0.1023.1/9.1/a}{2.0.1023.1-9.1-a1} \\
\texttt{13.0} & $2$ & $13$ & \href{https://www.lmfdb.org/EllipticCurve/2.0.107.1/432.2/b}{2.0.107.1-432.2-b1} \\
\texttt{17.0} & $2$ & $17$ & \href{https://www.lmfdb.org/EllipticCurve/2.0.127.1/32.2/a}{2.0.127.1-32.2-a1} \\
\texttt{19.0} & $2$ & $19$ & \href{https://www.lmfdb.org/EllipticCurve/2.0.1007.1/9.1/a}{2.0.1007.1-9.1-a1} \\
\texttt{23.0} & $2$ & $23$ & \href{https://www.lmfdb.org/EllipticCurve/2.0.11.1/1587.1/b}{2.0.11.1-1587.1-b1} \\
\texttt{29.0} & $2$ & $29$ & \href{https://www.lmfdb.org/EllipticCurve/2.0.4.1/841.1/a}{2.0.4.1-841.1-a1} \\
\texttt{31.0} & $2$ & $31$ & \href{https://www.lmfdb.org/EllipticCurve/2.0.3.1/47089.3/a}{2.0.3.1-47089.3-a1} \\
\texttt{37.0} & $2$ & $37$ & \href{https://www.lmfdb.org/EllipticCurve/2.0.104.1/441.1/a}{2.0.104.1-441.1-a1} \\
\texttt{41.0} & $2$ & $41$ & \href{https://www.lmfdb.org/EllipticCurve/2.0.359.1/32.2/a}{2.0.359.1-32.2-a1} \\
\texttt{43.0} & $2$ & $43$ & \href{https://www.lmfdb.org/EllipticCurve/2.0.1247.1/9.1/a}{2.0.1247.1-9.1-a1} \\
\texttt{67.0} & $2$ & $67$ & \href{https://www.lmfdb.org/EllipticCurve/2.0.335.1/9.1/a}{2.0.335.1-9.1-a1} \\
\texttt{73.0} & $2$ & $73$ & \href{https://www.lmfdb.org/EllipticCurve/2.0.31.1/2450.16/g}{2.0.31.1-2450.16-g1} \\
\texttt{163.0} & $2$ & $163$ & \href{https://www.lmfdb.org/EllipticCurve/2.0.1956.1/9.1/a}{2.0.1956.1-9.1-a1} \\
\hline
\texttt{9.0} & $3$ & $9$ & \href{https://www.lmfdb.org/EllipticCurve/2.0.103.1/98.2/a}{2.0.103.1-98.2-a1} \\
\texttt{25.0} & $3$ & $25$ & \href{https://www.lmfdb.org/EllipticCurve/2.0.1012.1/11.1/a}{2.0.1012.1-11.1-a1} \\
\texttt{49.0} & $3$ & $49$ & \href{https://www.lmfdb.org/EllipticCurve/2.0.11.1/9153.2/a}{2.0.11.1-9153.2-a1} \\
\hline
\texttt{4.1} & $4$ & $4$ & \href{https://www.lmfdb.org/EllipticCurve/2.0.103.1/224.4/b}{2.0.103.1-224.4-b1} \\
\texttt{2.0-3.0} & $4$ & $6$ & \href{https://www.lmfdb.org/EllipticCurve/2.0.103.1/400.3/a}{2.0.103.1-400.3-a1} \\
\texttt{2.0-5.0} & $4$ & $10$ & \href{https://www.lmfdb.org/EllipticCurve/2.0.1016.1/9.2/a}{2.0.1016.1-9.2-a1} \\
\texttt{2.0-7.0} & $4$ & $14$ & \href{https://www.lmfdb.org/EllipticCurve/2.0.11.1/240.1/a}{2.0.11.1-240.1-a1} \\
\texttt{3.0-5.0} & $4$ & $15$ & \href{https://www.lmfdb.org/EllipticCurve/2.0.11.1/225.3/a}{2.0.11.1-225.3-a1} \\
\texttt{3.0-7.0} & $4$ & $21$ & \href{https://www.lmfdb.org/EllipticCurve/2.0.11.1/26244.5/b}{2.0.11.1-26244.5-b1} \\
\texttt{2.0-11.0} & $4$ & $22$ & \href{https://www.lmfdb.org/EllipticCurve/2.0.47.1/6.1/a}{2.0.47.1-6.1-a1} \\
\texttt{2.0-13.0} & $4$ & $26$ & \href{https://www.lmfdb.org/EllipticCurve/2.0.11.1/3600.1/a}{2.0.11.1-3600.1-a1} \\
\texttt{27.0} & $4$ & $27$ & \href{https://www.lmfdb.org/EllipticCurve/2.0.103.1/729.1/a}{2.0.103.1-729.1-a1} \\
\texttt{3.0-11.0} & $4$ & $33$ & \href{https://www.lmfdb.org/EllipticCurve/2.0.7.1/242.3/a}{2.0.7.1-242.3-a1} \\
\texttt{2.0-17.0} & $4$ & $34$ & \href{https://www.lmfdb.org/EllipticCurve/2.0.15.1/800.2/a}{2.0.15.1-800.2-a1} \\
\texttt{5.0-7.0} & $4$ & $35$ & \href{https://www.lmfdb.org/EllipticCurve/2.2.5.1/2401.1/b}{2.2.5.1-2401.1-b1} \\
\texttt{3.0-13.0} & $4$ & $39$ & \href{https://www.lmfdb.org/EllipticCurve/2.0.7.1/15876.2/e}{2.0.7.1-15876.2-e1} \\
\texttt{3.0-17.0} & $4$ & $51$ & \href{https://www.lmfdb.org/EllipticCurve/2.2.17.1/81.1/b}{2.2.17.1-81.1-b1} \\
\texttt{3.0-19.0} & $4$ & $57$ & \href{https://www.lmfdb.org/EllipticCurve/2.0.23.1/162.1/a}{2.0.23.1-162.1-a1} \\
\texttt{3.0-41.0} & $4$ & $123$ & \href{https://www.lmfdb.org/EllipticCurve/2.2.41.1/81.1/c}{2.2.41.1-81.1-c1} \\
\texttt{3.0-89.0} & $4$ & $267$ & \href{https://www.lmfdb.org/EllipticCurve/2.2.89.1/81.1/a}{2.2.89.1-81.1-a1} \\
\hline
\texttt{9.1} & $5$ & $9$ & \href{https://www.lmfdb.org/EllipticCurve/2.0.3.1/324.1/a}{2.0.3.1-324.1-a1} \\
\hline
\texttt{8.1} & $6$ & $8$ & \href{https://www.lmfdb.org/EllipticCurve/2.0.103.1/238.3/a}{2.0.103.1-238.3-a1} \\
\texttt{2.0-9.0} & $6$ & $18$ & \href{https://www.lmfdb.org/EllipticCurve/2.0.103.1/196.5/a}{2.0.103.1-196.5-a1} \\
\texttt{5.0-9.0} & $6$ & $45$ & \href{https://www.lmfdb.org/EllipticCurve/2.0.39.1/4.2/a}{2.0.39.1-4.2-a1} \\
\texttt{2.0-25.0} & $6$ & $50$ & \href{https://www.lmfdb.org/EllipticCurve/2.0.7.1/5632.18/a}{2.0.7.1-5632.18-a1} \\
\texttt{3.0-25.0} & $6$ & $75$ & \href{https://www.lmfdb.org/EllipticCurve/2.2.5.1/2025.1/c}{2.2.5.1-2025.1-c1} \\
\texttt{9.0-11.0} & $6$ & $99$ & \href{https://www.lmfdb.org/EllipticCurve/2.2.33.1/1.1/a}{2.2.33.1-1.1-a1} \\
\texttt{3.0-49.0} & $6$ & $147$ & \href{https://www.lmfdb.org/EllipticCurve/2.2.21.1/49.1/a}{2.2.21.1-49.1-a1} \\
\hline
\texttt{4.1-3.0} & $8$ & $12$ & \href{https://www.lmfdb.org/EllipticCurve/2.0.103.1/900.2/a}{2.0.103.1-900.2-a1} \\
\texttt{16.1} & $8$ & $16$ & \href{https://www.lmfdb.org/EllipticCurve/2.0.103.1/225.1/a}{2.0.103.1-225.1-a1} \\
\texttt{4.1-5.0} & $8$ & $20$ & \href{https://www.lmfdb.org/EllipticCurve/2.0.11.1/27.2/a}{2.0.11.1-27.2-a1} \\
\texttt{4.1-7.0} & $8$ & $28$ & \href{https://www.lmfdb.org/EllipticCurve/2.0.23.1/6.1/a}{2.0.23.1-6.1-a1} \\
\texttt{2.0-3.0-5.0} & $8$ & $30$ & \href{https://www.lmfdb.org/EllipticCurve/2.0.7.1/324.2/a}{2.0.7.1-324.2-a1} \\
\texttt{2.0-3.0-7.0} & $8$ & $42$ & \href{https://www.lmfdb.org/EllipticCurve/4.4.1600.1/1.1/a}{4.4.1600.1-1.1-a1} \\
\texttt{2.0-3.0-11.0} & $8$ & $66$ & \href{https://www.lmfdb.org/EllipticCurve/4.4.17424.1/1.1/a}{4.4.17424.1-1.1-a1} \\
\hline
\texttt{16.2} & $10$ & $16$ & \href{https://www.lmfdb.org/EllipticCurve/2.0.4.1/200.2/a}{2.0.4.1-200.2-a1} \\
\texttt{2.0-9.1} & $10$ & $18$ & \href{https://www.lmfdb.org/EllipticCurve/2.0.3.1/196.2/a}{2.0.3.1-196.2-a1} \\
\texttt{32.0} & $10$ & $32$ & \href{https://www.lmfdb.org/EllipticCurve/2.0.15.1/15.1/a}{2.0.15.1-15.1-a1} \\
\hline
\texttt{8.1-3.0} & $12$ & $24$ & \href{https://www.lmfdb.org/EllipticCurve/2.0.15.1/60.2/a}{2.0.15.1-60.2-a1} \\
\texttt{4.1-9.0} & $12$ & $36$ & \href{https://www.lmfdb.org/EllipticCurve/2.0.7.1/28.2/a}{2.0.7.1-28.2-a1} \\
\texttt{8.1-5.0} & $12$ & $40$ & \href{https://www.lmfdb.org/EllipticCurve/4.4.2048.1/1.1/a}{4.4.2048.1-1.1-a1} \\
\texttt{64.1} & $12$ & $64$ & \href{https://www.lmfdb.org/EllipticCurve/4.4.3600.1/225.1/c}{4.4.3600.1-225.1-c1} \\
\texttt{2.0-9.0-5.0} & $12$ & $90$ & \href{https://www.lmfdb.org/EllipticCurve/4.4.3600.1/1.1/b}{4.4.3600.1-1.1-b1} \\
\texttt{9.0-5.0-7.0} & $12$ & $315$ & \href{https://www.lmfdb.org/EllipticCurve/4.4.11025.1/1.1/a}{4.4.11025.1-1.1-a1} \\
\hline
\texttt{16.1-3.0} & $16$ & $48$ & \href{https://www.lmfdb.org/EllipticCurve/4.4.11025.1/25.3/b}{4.4.11025.1-25.3-b1} \\
\texttt{4.1-3.0-5.0} & $16$ & $60$ & \href{https://www.lmfdb.org/EllipticCurve/4.4.3600.1/1.1/a}{4.4.3600.1-1.1-a1} \\
\texttt{16.1-7.0} & $16$ & $112$ & \href{https://www.lmfdb.org/EllipticCurve/4.4.12544.1/1.1/a}{4.4.12544.1-1.1-a1} \\
\hline
\texttt{8.1-9.0} & $18$ & $72$ & \lmfdb{https://www.lmfdb.org/EllipticCurve/4.4.2304.1/1.1/a}{4.4.2304.1-1.1-a1} \\
\hline
\end{longtable}

\section{Tables for explicit classification of genus \texorpdfstring{$0$}{0} isogeny graphs}\label{appendix_para}

Below, we include the tables referenced in Section~\ref{sec:expcalssgenus0}. In particular, for $n\in \{9,16,18,25\}$, Tables~\ref{ta:jinv} and~\ref{ta:curves} give the Fricke parameterizations $(j_{n,1}(t),j_{n,2}(t)$ and Weierstrass models for the parameterized isogenous families of elliptic curves $\mathcal{C}_{n,m}^1(t,d)$, respectively. For a complete list of the Fricke parameterizations, see~\cite[Table~1]{Bariso},\cite[Table 3]{MR3084348},\cite[Tables 4 and 5]
{MR2514149}. We note that these tables are obtained from results in~\cite{MR0376533,MR1486831,MR3221641,MR2341166,MR675184,MR3838339}.

{\renewcommand*{\arraystretch}{1.55} \begin{longtable}{cC{3.5in}C{2in}}
\caption{The Fricke Parameterizations for $n\in\{9,16,18,25\}$}\label{ta:jinv}\\
\hline
$n$ & $j_{n,1}(  t)  $ & $j_{n,2}(  t)  $\\
\hline
\endfirsthead
\hline
$n$ & $j_{n,1}(  t)  $ & $j_{n,2}(  t)  $ \\
\hline
\endhead
\hline
\multicolumn{3}{r}{\emph{continued on next page}}
\endfoot
\hline
\endlastfoot
$9$ & $\frac{(t+6)^{3}(t^{3}+234t^{2}+756t+2160)^{3}}{(t-3)^{8}(t^{3}-27)}$ &
$\frac{t^{3}(t^{3}-24)^{3}}{t^{3}-27}$\\\hline
$16$ & $\frac{1}{t(t-2)^{16}(t+2)^{4}(t^{2}+4)}(t^{8}+240t^{7}+2160t^{6}+6720t^{5}+17504t^{4}+26880t^{3}+34560t^{2}+15360t+256)^{3}$ & $\frac
{(t^{8}-16t^{4}+16)^{3}}{t^{4}(t^{4}-16)}$\\\hline
$18$ & $\frac{(t^{3}+6t^{2}+4)^{3}}{t^{2}(t-2)^{18}(t+1)^{9}
(t^{2}-t+1)(t^{2}+2t+4)^{2}}(t^{9}+234t^{8}+756t^{7}+2172t^{6}
+1872t^{5}+3024t^{4}+48t^{3}+3744t^{2}+64)^{3}$ & $\frac{(t^{3}-2)^{3}(t^{9}-6t^{6}-12t^{3}-8)^{3}}{t^{9}(t^{3}-8)(t^{3}+1)^{2}}$\\\hline
$25$ & $\frac{1}{(t-1)^{25}
(t^{4}+t^{3}+6t^{2}+6t+11)}(t^{10}+240t^{9}+2170t^{8}+8880t^{7}+34835t^{6}+83748t^{5}
+206210t^{4}+313380t^{3}+503545t^{2}+424740t+375376)^{3}$ & $\frac{1}{t^{5}+5t^{3}+5t-11}(t^{10}+10t^{8}+35t^{6}-12t^{5}
+50t^{4}-60t^{3}+25t^{2}-60t+16)^{3}$
\end{longtable}}

{\renewcommand*{\arraystretch}{1.14} \footnotesize  \begin{longtable}{ccC{2.0in}C{3.45in}}
\caption{The elliptic curve $\mathcal{C}^1_{n,m}(t,d):y^2=x^3+d^2\mathcal{A}_{n,m}x+d^3\mathcal{B}_{n,m}$ \\ Note that for $n=25$, $a=\sqrt{5}$ and $s=\sqrt{-2(5+a)(t^2+4)}$}\\
\hline
$n$ & $m$ & $\mathcal{A}_{n,m}$ & $\mathcal{B}_{n,m}$\\
\hline
\endfirsthead
\caption[]{\emph{continued}}\\
\hline
$n$ & $m$ & $\mathcal{A}_{n,m}$ & $\mathcal{B}_{n,m}$ \\
\hline
\endhead
\hline
\multicolumn{4}{r}{\emph{continued on next page}}
\endfoot
\hline
\endlastfoot

$9$ & $1$ & $-3(t + 6) (t^{3} + 234 t^{2} + 756 t + 2160)$ & $-2 (t^{6} - 504 t^{5} - 16632 t^{4} - 123012 t^{3} - 517104 t^{2} - 1143072 t - 1475496)$\\\cmidrule{2-4}
& $2$ & $ -3 t^{4} - 648 t $  & $-2 t^{6} + 1080 t^{3} + 11664$ \\\cmidrule{2-4}
& $3$ & $-3t (t^{3} - 24) $ & $ -2(t^{6} - 36 t^{3} + 216)$  \\\cmidrule{2-4}
& $4$ & $-3(t - 6 \zeta_{3} - 6) (t^{3} + (-234 \zeta_{3} - 234) t^{2} + 756 \zeta_{3} t + 2160)$ & $-2(t^{6} + (504 \zeta_{3} + 504) t^{5} - 16632 \zeta_{3} t^{4} - 123012 t^{3} + (517104 \zeta_{3} + 517104) t^{2} - 1143072 \zeta_{3} t - 1475496)$\\\cmidrule{2-4}
& $5$ & $-3 t^{4} - 720 \zeta_{3} t^{3} + (6480 \zeta_{3} + 6480) t^{2} - 20088 t - 38880 \zeta_{3}$ & $-2(t^{6} - 504 \zeta_{3} t^{5} + (16632 \zeta_{3} + 16632) t^{4} - 123012 t^{3} - 517104 \zeta_{3} t^{2} + (1143072 \zeta_{3} + 1143072) t - 1475496)$ \\\hline

$16$ & $1$ & $-27  (t^{8} + 240 t^{7} + 2160 t^{6} + 6720 t^{5} + 17504 t^{4} + 26880 t^{3} + 34560 t^{2} + 15360 t + 256) $  & $-54  (t^{4} + 24 t^{3} + 24 t^{2} + 96 t + 16)  (t^{8} - 528 t^{7} - 3984 t^{6} - 14784 t^{5} - 31648 t^{4} - 59136 t^{3} - 63744 t^{2} - 33792 t + 256) $ \\\cmidrule{2-4}
& $2$ & $ -27  (t^{8} + 240 t^{6} + 2144 t^{4} + 3840 t^{2} + 256)$ & $ -54  (t^{4} - 24 t^{3} + 24 t^{2} - 96 t + 16)  (t^{4} + 24 t^{2} + 16)  (t^{4} + 24 t^{3} + 24 t^{2} + 96 t + 16)$ \\\cmidrule{2-4}
& $3$ & $-27  (t^{8} - 240 t^{7} + 2160 t^{6} - 6720 t^{5} + 17504 t^{4} - 26880 t^{3} + 34560 t^{2} - 15360 t + 256) $ & $-54  (t^{4} - 24 t^{3} + 24 t^{2} - 96 t + 16)  (t^{8} + 528 t^{7} - 3984 t^{6} + 14784 t^{5} - 31648 t^{4} + 59136 t^{3} - 63744 t^{2} + 33792 t + 256) $ \\\cmidrule{2-4}
& $4$ & $-27 t^{8} - 6048 t^{4} - 6912 $ & $-54 t^{12} + 28512 t^{8} + 456192 t^{4} - 221184 $ \\\cmidrule{2-4}
& $5$ & $-27 t^{8} + 6480 t^{6} - 57888 t^{4} + 103680 t^{2} - 6912 $ & $-54 t^{12} - 27216 t^{10} + 899424 t^{8} - 6096384 t^{6} + 14390784 t^{4} - 6967296 t^{2} - 221184 $ \\\cmidrule{2-4}
& $6$ & $-27 t^{8} + 432 t^{4} - 6912$ & $-54 t^{12} + 1296 t^{8} + 20736 t^{4} - 221184 $\\\cmidrule{2-4}
& $7$ & $-27  (t^{8} - 256 t^{4} + 4096) $  & $ -54  (t^{4} - 32)  (t^{8} + 512 t^{4} - 8192)$ \\\cmidrule{2-4}
& $8$ & $-27  (t^{8} - 16 t^{4} + 16) $  & $-54  (t^{4} - 8)  (t^{8} - 16 t^{4} - 8) $ \\\cmidrule{2-4}
& $9$ & $-27(t^{8} + 240 i t^{7} - 2160 t^{6} - 6720 i t^{5} + 17504 t^{4} + 26880 i t^{3} - 34560 t^{2} - 15360 i t + 256)$ & $-54(t^{4} + 24 i t^{3} - 24 t^{2} - 96 i t + 16) (t^{8} - 528 i t^{7} + 3984 t^{6} + 14784 i t^{5} - 31648 t^{4} - 59136 i t^{3} + 63744 t^{2} + 33792 i t + 256)$\\\cmidrule{2-4}
& $10$ & $-27(t^{8} - 240 i t^{7} - 2160 t^{6} + 6720 i t^{5} + 17504 t^{4} - 26880 i t^{3} - 34560 t^{2} + 15360 i t + 256)$ & $-54(t^{4} - 24 i t^{3} - 24 t^{2} + 96 i t + 16) (t^{8} + 528 i t^{7} + 3984 t^{6} - 14784 i t^{5} - 31648 t^{4} + 59136 i t^{3} + 63744 t^{2} - 33792 i t + 256)$
\\\hline

$18$ & $1$ & $-3  (t^{3} + 6 t^{2} + 4)  (t^{9} + 234 t^{8} + 756 t^{7} + 2172 t^{6} + 1872 t^{5} + 3024 t^{4} + 48 t^{3} + 3744 t^{2} + 64) $ & $-2  (t^{6} + 24 t^{5} + 24 t^{4} + 92 t^{3} - 48 t^{2} + 96 t - 8)  (t^{12} - 528 t^{11} - 3984 t^{10} - 14792 t^{9} - 27936 t^{8} - 42624 t^{7} - 37632 t^{6} - 52992 t^{5} - 25344 t^{4} - 43520 t^{3} - 6144 t^{2} - 6144 t - 512) $ \\\cmidrule{2-4}
& $2$ & $-3  (t^{3} + 6 t - 2)  (t^{9} + 234 t^{7} - 6 t^{6} + 756 t^{5} - 936 t^{4} + 2172 t^{3} - 1512 t^{2} + 936 t - 8) $ & $-2  (t^{6} + 24 t^{5} + 24 t^{4} + 92 t^{3} - 48 t^{2} + 96 t - 8)  (t^{12} - 24 t^{11} + 48 t^{10} - 680 t^{9} + 792 t^{8} - 3312 t^{7} + 4704 t^{6} - 10656 t^{5} + 13968 t^{4} - 14792 t^{3} + 7968 t^{2} - 2112 t - 8) $ \\\cmidrule{2-4}
& $3$ & $-3  (t^{12} + 232 t^{9} + 960 t^{6} + 256 t^{3} + 256) $ & $-2  (t^{18} - 516 t^{15} - 12072 t^{12} - 24640 t^{9} - 30720 t^{6} + 6144 t^{3} + 4096) $  \\\cmidrule{2-4}
& $4$ & $-3 (t^{3} - 2)  (t^{3} + 6 t - 2)  (t^{6} - 6 t^{4} - 4 t^{3} + 36 t^{2} + 12 t + 4) $ & $-2  (t^{2} + 2 t - 2)  (t^{4} - 2 t^{3} - 8 t - 2)  (t^{4} - 2 t^{3} + 6 t^{2} + 4 t + 4)  (t^{8} + 2 t^{7} + 4 t^{6} - 16 t^{5} - 14 t^{4} + 8 t^{3} + 64 t^{2} - 16 t + 4) $ \\\cmidrule{2-4}
& $5$ & $-3  (t^{3} + 4)  (t^{9} - 12 t^{6} + 48 t^{3} + 64) $  & $-2  (t^{6} - 4 t^{3} - 8)  (t^{12} - 8 t^{9} - 512 t^{3} - 512) $ \\\cmidrule{2-4}
& $6$ & $-3  (t^{3} - 2)  (t^{9} - 6 t^{6} - 12 t^{3} - 8) $ & $-2  (t^{6} - 4 t^{3} - 8)  (t^{12} - 8 t^{9} - 8 t^{3} - 8) $\\\cmidrule{2-4}
& $7$ & $-3(t^{3} + (-6 \zeta_{3} - 6) t^{2} + 4) (t^{9} + (-234 \zeta_{3} - 234) t^{8} + 756 \zeta_{3} t^{7} + 2172 t^{6} + (-1872 \zeta_{3} - 1872) t^{5} + 3024 \zeta_{3} t^{4} + 48 t^{3} + (-3744 \zeta_{3} - 3744) t^{2} + 64)$ & $-2(t^{6} + (-24 \zeta_{3} - 24) t^{5} + 24 \zeta_{3} t^{4} + 92 t^{3} + (48 \zeta_{3} + 48) t^{2} + 96 \zeta_{3} t - 8) (t^{12} + (528 \zeta_{3} + 528) t^{11} - 3984 \zeta_{3} t^{10} - 14792 t^{9} + (27936 \zeta_{3} + 27936) t^{8} - 42624 \zeta_{3} t^{7} - 37632 t^{6} + (52992 \zeta_{3} + 52992) t^{5} - 25344 \zeta_{3} t^{4} - 43520 t^{3} + (6144 \zeta_{3} + 6144) t^{2} - 6144 \zeta_{3} t - 512)$\\\cmidrule{2-4}
& $8$ & $-3(t^{3} + 6 \zeta_{3} t - 2) (t^{9} + 234 \zeta_{3} t^{7} - 6 t^{6} + (-756 \zeta_{3} - 756) t^{5} - 936 \zeta_{3} t^{4} + 2172 t^{3} + (1512 \zeta_{3} + 1512) t^{2} + 936 \zeta_{3} t - 8)$ & $-2(t^{6} + (-24 \zeta_{3} - 24) t^{5} + 24 \zeta_{3} t^{4} + 92 t^{3} + (48 \zeta_{3} + 48) t^{2} + 96 \zeta_{3} t - 8) (t^{12} + (24 \zeta_{3} + 24) t^{11} + 48 \zeta_{3} t^{10} - 680 t^{9} + (-792 \zeta_{3} - 792) t^{8} - 3312 \zeta_{3} t^{7} + 4704 t^{6} + (10656 \zeta_{3} + 10656) t^{5} + 13968 \zeta_{3} t^{4} - 14792 t^{3} + (-7968 \zeta_{3} - 7968) t^{2} - 2112 \zeta_{3} t - 8)$
\\\cmidrule{2-4}
& $9$ & $-3(t^{3} + 6 \zeta_{3} t^{2} + 4) (t^{9} + 234 \zeta_{3} t^{8} + (-756 \zeta_{3} - 756) t^{7} + 2172 t^{6} + 1872 \zeta_{3} t^{5} + (-3024 \zeta_{3} - 3024) t^{4} + 48 t^{3} + 3744 \zeta_{3} t^{2} + 64)$ & $-2(t^{6} + 24 \zeta_{3} t^{5} + (-24 \zeta_{3} - 24) t^{4} + 92 t^{3} - 48 \zeta_{3} t^{2} + (-96 \zeta_{3} - 96) t - 8) (t^{12} - 528 \zeta_{3} t^{11} + (3984 \zeta_{3} + 3984) t^{10} - 14792 t^{9} - 27936 \zeta_{3} t^{8} + (42624 \zeta_{3} + 42624) t^{7} - 37632 t^{6} - 52992 \zeta_{3} t^{5} + (25344 \zeta_{3} + 25344) t^{4} - 43520 t^{3} - 6144 \zeta_{3} t^{2} + (6144 \zeta_{3} + 6144) t - 512)$\\\cmidrule{2-4}
& $10$ & $-3(t^{3} + (-6 \zeta_{3} - 6) t - 2) (t^{9} + (-234 \zeta_{3} - 234) t^{7} - 6 t^{6} + 756 \zeta_{3} t^{5} + (936 \zeta_{3} + 936) t^{4} + 2172 t^{3} - 1512 \zeta_{3} t^{2} + (-936 \zeta_{3} - 936) t - 8)$ & $-2(t^{6} + 24 \zeta_{3} t^{5} + (-24 \zeta_{3} - 24) t^{4} + 92 t^{3} - 48 \zeta_{3} t^{2} + (-96 \zeta_{3} - 96) t - 8) (t^{12} - 24 \zeta_{3} t^{11} + (-48 \zeta_{3} - 48) t^{10} - 680 t^{9} + 792 \zeta_{3} t^{8} + (3312 \zeta_{3} + 3312) t^{7} + 4704 t^{6} - 10656 \zeta_{3} t^{5} + (-13968 \zeta_{3} - 13968) t^{4} - 14792 t^{3} + 7968 \zeta_{3} t^{2} + (2112 \zeta_{3} + 2112) t - 8)$\\\hline

$25$ & $1$ & $-27  (t^{2} + 4)  (t^{10} + 240 t^{9} + 2170 t^{8} + 8880 t^{7} + 34835 t^{6} + 83748 t^{5} + 206210 t^{4} + 313380 t^{3} + 503545 t^{2} + 424740 t + 375376) $ & $-54  (t^{2} + 4)^{2}  (t^{4} + 6 t^{3} + 21 t^{2} + 36 t + 61)  (t^{10} - 510 t^{9} - 13580 t^{8} - 36870 t^{7} - 190915 t^{6} - 393252 t^{5} - 1068040 t^{4} - 1508370 t^{3} - 2581955 t^{2} - 2087010 t - 1885124) $ \\\cmidrule{2-4}
& $2$ & $-27  (t^{2} + 4)  (t^{2} + 3 t + 1)  (t^{4} - 4 t^{3} + 11 t^{2} - 14 t + 31)  (t^{4} + t^{3} + 11 t^{2} - 4 t + 16) $ & $-54  (t^{2} - 2 t - 4)  (t^{2} + 4)^{2}  (t^{4} - 4 t^{3} + 21 t^{2} - 34 t + 41)  (t^{4} + 3 t^{2} + 1)  (t^{4} + 6 t^{3} + 21 t^{2} + 36 t + 61) $ \\\cmidrule{2-4}
& $3$ &  $-27  (t^{2} + 4)  (t^{10} + 10 t^{8} + 35 t^{6} - 12 t^{5} + 50 t^{4} - 60 t^{3} + 25 t^{2} - 60 t + 16) $ & $ -54  (t^{2} + 4)^{2}  (t^{4} + 3 t^{2} + 1)  (t^{10} + 10 t^{8} + 35 t^{6} - 18 t^{5} + 50 t^{4} - 90 t^{3} + 25 t^{2} - 90 t + 76)$\\\cmidrule{2-4}
& $4$ & $(-1620 t^{10} + (-7290a + 7290) t^{9} + (22680a - 40500) t^{8} + (-72900a + 191160) t^{7} + (204120a - 273780) t^{6} + (-247860a + 803520) t^{5} + (383940a - 497340) t^{4} + (-51030a - 44550) t^{3} + (-461700a + 675540) t^{2} + (962280a - 2666520) t - 732240a + 1846800) s - 27 t^{12} + (-1620a + 1620) t^{11} + (14580a + 14202) t^{10} + (24300a + 66420) t^{9} + (56700a + 291195) t^{8} + (396900a - 215136) t^{7} + (-814860a + 2326050) t^{6} + (2093040a - 4805244) t^{5} + (-5086800a + 8959005) t^{4} + (5248800a - 14176620) t^{3} + (-7322400a + 14084388) t^{2} + (5028480a - 8935920) t + 362880a + 1968192$ & $(6804 t^{16} + (112266a - 112266) t^{15} + (-830088a + 952560) t^{14} + (1908522a - 9100350) t^{13} + (-13281408a + 14186340) t^{12} + (13135122a - 66342402) t^{11} + (-29284416a + 32808888) t^{10} + (-55918674a + 133749630) t^{9} + (392345856a - 687374100) t^{8} + (-1038348234a + 2663034570) t^{7} + (2420958456a - 5000273208) t^{6} + (-3964122666a + 8897968422) t^{5} + (4793608512a - 11268172440) t^{4} + (-4644628128a + 9684133200) t^{3} + (2665453392a - 6270158160) t^{2} + (180768672a + 64393056) t - 721877184a + 1810299456) s - 54 t^{18} + (6804a - 6804) t^{17} + (-224532a - 225666) t^{16} + (-1524096a - 1796256) t^{15} + (2925720a - 11467980) t^{14} + (-28746900a + 54058752) t^{13} + (143475948a - 210435138) t^{12} + (-324795744a + 1026222048) t^{11} + (1309851648a - 1998896292) t^{10} + (-2071375740a + 5889279960) t^{9} + (4658188500a - 8949729042) t^{8} + (-5863020408a + 12837645288) t^{7} + (4908378384a - 13522320198) t^{6} + (-2345692608a + 589169052) t^{5} + (-8539836480a + 15626432520) t^{4} + (16095542400a - 34654762656) t^{3} + (-19834585344a + 47889239904) t^{2} + (16691899392a - 33649349184) t - 4800466944a + 9192410496$\\\cmidrule{2-4}
& $5$ & $(1620 t^{10} + (7290a - 7290) t^{9} + (-22680a + 40500) t^{8} + (72900a - 191160) t^{7} + (-204120a + 273780) t^{6} + (247860a - 803520) t^{5} + (-383940a + 497340) t^{4} + (51030a + 44550) t^{3} + (461700a - 675540) t^{2} + (-962280a + 2666520) t + 732240a - 1846800) s - 27 t^{12} + (-1620a + 1620) t^{11} + (14580a + 14202) t^{10} + (24300a + 66420) t^{9} + (56700a + 291195) t^{8} + (396900a - 215136) t^{7} + (-814860a + 2326050) t^{6} + (2093040a - 4805244) t^{5} + (-5086800a + 8959005) t^{4} + (5248800a - 14176620) t^{3} + (-7322400a + 14084388) t^{2} + (5028480a - 8935920) t + 362880a + 1968192$ & $(-6804 t^{16} + (-112266a + 112266) t^{15} + (830088a - 952560) t^{14} + (-1908522a + 9100350) t^{13} + (13281408a - 14186340) t^{12} + (-13135122a + 66342402) t^{11} + (29284416a - 32808888) t^{10} + (55918674a - 133749630) t^{9} + (-392345856a + 687374100) t^{8} + (1038348234a - 2663034570) t^{7} + (-2420958456a + 5000273208) t^{6} + (3964122666a - 8897968422) t^{5} + (-4793608512a + 11268172440) t^{4} + (4644628128a - 9684133200) t^{3} + (-2665453392a + 6270158160) t^{2} + (-180768672a - 64393056) t + 721877184a - 1810299456) s - 54 t^{18} + (6804a - 6804) t^{17} + (-224532a - 225666) t^{16} + (-1524096a - 1796256) t^{15} + (2925720a - 11467980) t^{14} + (-28746900a + 54058752) t^{13} + (143475948a - 210435138) t^{12} + (-324795744a + 1026222048) t^{11} + (1309851648a - 1998896292) t^{10} + (-2071375740a + 5889279960) t^{9} + (4658188500a - 8949729042) t^{8} + (-5863020408a + 12837645288) t^{7} + (4908378384a - 13522320198) t^{6} + (-2345692608a + 589169052) t^{5} + (-8539836480a + 15626432520) t^{4} + (16095542400a - 34654762656) t^{3} + (-19834585344a + 47889239904) t^{2} + (16691899392a - 33649349184) t - 4800466944a + 9192410496$\\\cmidrule{2-4}
& $6$ & $((810a - 810) t^{10} - 14580 t^{9} + (8910a + 36450) t^{8} + (-59130a - 86670) t^{7} + (34830a + 373410) t^{6} + (-277830a - 217890) t^{5} + (56700a + 711180) t^{4} + (47790a - 149850) t^{3} + (-106920a - 816480) t^{2} + (852120a + 1072440) t - 557280a - 907200) s - 27 t^{12} + (1620a + 1620) t^{11} + (-14580a + 14202) t^{10} + (-24300a + 66420) t^{9} + (-56700a + 291195) t^{8} + (-396900a - 215136) t^{7} + (814860a + 2326050) t^{6} + (-2093040a - 4805244) t^{5} + (5086800a + 8959005) t^{4} + (-5248800a - 14176620) t^{3} + (7322400a + 14084388) t^{2} + (-5028480a - 8935920) t - 362880a + 1968192$ & $((-3402a + 3402) t^{16} + 224532 t^{15} + (-61236a - 1598940) t^{14} + (3595914a + 221130) t^{13} + (-452466a - 26110350) t^{12} + (26603640a - 333396) t^{11} + (-1762236a - 56806596) t^{10} + (-38915478a - 72921870) t^{9} + (147514122a + 637177590) t^{8} + (-812343168a - 1264353300) t^{7} + (1289657376a + 3552259536) t^{6} + (-2466922878a - 5461322454) t^{5} + (3237281964a + 6349935060) t^{4} + (-2519752536a - 6769503720) t^{3} + (1802352384a + 3528554400) t^{2} + (-122580864a + 484118208) t - 544211136a - 899543232) s - 54 t^{18} + (-6804a - 6804) t^{17} + (224532a - 225666) t^{16} + (1524096a - 1796256) t^{15} + (-2925720a - 11467980) t^{14} + (28746900a + 54058752) t^{13} + (-143475948a - 210435138) t^{12} + (324795744a + 1026222048) t^{11} + (-1309851648a - 1998896292) t^{10} + (2071375740a + 5889279960) t^{9} + (-4658188500a - 8949729042) t^{8} + (5863020408a + 12837645288) t^{7} + (-4908378384a - 13522320198) t^{6} + (2345692608a + 589169052) t^{5} + (8539836480a + 15626432520) t^{4} + (-16095542400a - 34654762656) t^{3} + (19834585344a + 47889239904) t^{2} + (-16691899392a - 33649349184) t + 4800466944a + 9192410496$\\\cmidrule{2-4}
& $7$ & $((-810a + 810) t^{10} + 14580 t^{9} + (-8910a - 36450) t^{8} + (59130a + 86670) t^{7} + (-34830a - 373410) t^{6} + (277830a + 217890) t^{5} + (-56700a - 711180) t^{4} + (-47790a + 149850) t^{3} + (106920a + 816480) t^{2} + (-852120a - 1072440) t + 557280a + 907200) s - 27 t^{12} + (1620a + 1620) t^{11} + (-14580a + 14202) t^{10} + (-24300a + 66420) t^{9} + (-56700a + 291195) t^{8} + (-396900a - 215136) t^{7} + (814860a + 2326050) t^{6} + (-2093040a - 4805244) t^{5} + (5086800a + 8959005) t^{4} + (-5248800a - 14176620) t^{3} + (7322400a + 14084388) t^{2} + (-5028480a - 8935920) t - 362880a + 1968192$ & $((3402a - 3402) t^{16} - 224532 t^{15} + (61236a + 1598940) t^{14} + (-3595914a - 221130) t^{13} + (452466a + 26110350) t^{12} + (-26603640a + 333396) t^{11} + (1762236a + 56806596) t^{10} + (38915478a + 72921870) t^{9} + (-147514122a - 637177590) t^{8} + (812343168a + 1264353300) t^{7} + (-1289657376a - 3552259536) t^{6} + (2466922878a + 5461322454) t^{5} + (-3237281964a - 6349935060) t^{4} + (2519752536a + 6769503720) t^{3} + (-1802352384a - 3528554400) t^{2} + (122580864a - 484118208) t + 544211136a + 899543232) s - 54 t^{18} + (-6804a - 6804) t^{17} + (224532a - 225666) t^{16} + (1524096a - 1796256) t^{15} + (-2925720a - 11467980) t^{14} + (28746900a + 54058752) t^{13} + (-143475948a - 210435138) t^{12} + (324795744a + 1026222048) t^{11} + (-1309851648a - 1998896292) t^{10} + (2071375740a + 5889279960) t^{9} + (-4658188500a - 8949729042) t^{8} + (5863020408a + 12837645288) t^{7} + (-4908378384a - 13522320198) t^{6} + (2345692608a + 589169052) t^{5} + (8539836480a + 15626432520) t^{4} + (-16095542400a - 34654762656) t^{3} + (19834585344a + 47889239904) t^{2} + (-16691899392a - 33649349184) t + 4800466944a + 9192410496$

\label{ta:curves}	
\end{longtable}}

\bibliographystyle{plain}

\end{document}